\documentclass[b5,11pt,leqno]{article}

\usepackage{amsmath} 
\usepackage{amssymb} 
\usepackage{amscd} 
\usepackage{array} 

\usepackage{mathrsfs}  

\usepackage{amsthm} 

\usepackage[OT2,T1]{fontenc}

\numberwithin{equation}{section}

\theoremstyle{definition} 

\numberwithin{equation}{section}  

\allowdisplaybreaks               

\title{Projective representations and spin characters of complex reflection groups $G(m,p,n)$ and $G(m,p,\infty)$, III}

\author{Takeshi HIRAI and Akihito HORA}

\date*{}

\begin{document}

\maketitle

\setcounter{section}{0}
\setcounter{page}{1}
~\\[-10ex]
\begin{quotation}
{\small
{\bf Abstract.} Based on the hereditary property, studied in [18] and [21], from \lq mother groups' $G(m,1,n),\,4\le n\le\infty,$ to \lq child groups' $G(m,p,n),\,p|m,\,p>1,$ so-called complex reflection groups, we study in this paper classification and construction of irreducible projective representations (= spin representations) and their characters of the generalized symmetric groups $G(m,1,n)$, and spin characters of the inductive limit groups $G(m,1,\infty)$ in detail. By the heredity studied further, this gives the main kernel of the results for the child groups with $p|m, p>1$.}\footnote{
\\ 
symmetric groups, complex reflection groups, limit of characters
}

\end{quotation}

\begin{quotation}
\noindent
{\small {\bf Table of Contents}}

\noindent
{\small 
Part I\; Spin irreducible representations and their characters of generalized symmetric groups $G(m,1,n),\,4\le n<\infty$ 

1\;\; Representation groups of generalized symmetric groups\hfill $\cdots\;$ \;\;4\hspace{4ex}

2\;\; Preliminaries to Part I\hfill$\cdots\;$ \;\;7\hspace{4ex}

3\;\; Linear IRs of $G(m,1,n)$ and their characters\hfill$\cdots\;$ \;\;9\hspace{4ex}

4\;\; Spin IRs of Schur-Young subgroups of $\widetilde{{\mathfrak S}}_n$\hfill$\cdots\;$ 12\hspace{4ex}

5\;\; Spin IRs of $G(m,1,n)$ of spin type $z_1\to -1$\hfill$\cdots\;$ 17\hspace{4ex}

6\;\; Spin IRs and their characters of spin type $\chi^{\rm I}$\hfill$\cdots\;$ 19\hspace{4ex}

7\;\; Spin IRs and their characters of spin type $\chi^{\rm II}$\hfill$\cdots\;$ 25\hspace{4ex}

8\;\; Spin IRs of spin type $\chi^{\rm III}=(-1,\,1,-1)$\hfill$\cdots\;$ 32\hspace{4ex}

9\;\; Spin IRs of spin type $\chi^{\rm V}=(1,-1,-1)$\hfill$\cdots\;$ 38\hspace{4ex}

10\; Spin IRs and their characters of spin type $\chi^{\rm VI}$\hfill$\cdots\;$ 40\hspace{4ex}

11\; Spin IRs of spin type $\chi^{\rm VII}=(1,\,1,-1)$\hfill$\cdots\;$ 43\hspace{4ex}

\noindent
Part II\; Spin characters of infinite generalized symmetric groups $G(m,1,\infty)$ 

12\; Preliminaries to Part II\hfill$\cdots\;$ 47\hspace{4ex}

13\; Infinite symmetric groups ${\mathfrak S}_\infty$ and its covering group $\widetilde{{\mathfrak S}}_\infty$\hfill$\cdots\;$ 51\hspace{4ex}

14\; 
Characters of infinite generalized symmetric group $G(m,1,\infty)$\hfill$\cdots\;$ 56\hspace{4ex}

15\; Heredity from $R\big(G(m,1,\infty)\big)$ and the case of $R\big(G(m,p,\infty)\big)$\hfill$\cdots\;$ 59\hspace{4ex}


16\; Spin characters of $R\big(G(m,1,\infty)\big)$ of spin type $z_1\to -1$\hfill$\cdots\;$ 68\hspace{4ex}

17\; Spin characters of spin type $\chi^{\rm VII}=(1,\,1,-1)$\hfill$\cdots\;$ 71\hspace{4ex}

18\; Spin characters of spin type $\chi^{\rm I}=(-1,-1,-1)$\hfill$\cdots\;$ 73\hspace{4ex}

19\; Spin characters of spin type $\chi^{\rm II}=(-1,-1,1)$\hfill$\cdots\;$ 75\hspace{4ex}

20\; Spin characters of spin type $\chi^{\rm III}\!=(-1,\,1,-1)$ and $\chi^{\rm IV}\!=(-1,\,1,\,1)$\hfill$\cdots\;$ 78\hspace{4ex}

21\; Spin characters of spin type $\chi^{\rm V}=(1,-1,-1)$ and $\chi^{\rm VI}=(1,-1,\,1)$\hfill$\cdots\;$ 80\hspace{4ex}

22\, Parameter spaces and summary for the cases of $n=\infty$\hfill$\cdots\;$ 82\hspace{4ex}
 
\hspace*{2.7ex} References\hfill$\cdots\;$ 86\hspace{4ex}

}
\end{quotation}

\markboth{\hfil{\small\sc \;T.\ Hirai and A.\ Hora}\qquad\qquad\hfil}{{\thesection}\quad\mbox{\it\small Introduction}}

{\bf Introduction.} This paper is a continuation of the preceding papers [18] and [21] (referred as [I] and [II] in the following) in which we studied complex reflection groups $G(m,p,n)$ and $G(m,p,\infty),\,p|m,$ in general. We call the generalized symmetric groups $G(m,1,n)$ and $G(m,1,\infty)$ as \lq mother groups' and $G(m,p,n)$ and $G(m,p,\infty)$ for $p|m,\,p>1$, as \lq child groups', and studied hereditary property from a mother to her children. On the basis of this heredity, 
here in this paper we study the case of mother groups in total, and in addition to this we study the hereditary property furthermore, mainly in \S 15, for not only $G(m,p,n)$ but also for other normal subgroups of different types.

The paper consists of two Parts. In Part I, irreducible spin representations (= projective representations) of generalized symmetric groups $G(m,1,n),\,n\ge 4,$ are classified and constructed as induced representations, and their characters, called spin characters of $G(m,1,n)$, are calculated. 
In Part II, we studied spin characters of the inductive limits $G(m,1,\infty):=\lim_{n\to\infty}G(m,1,n)$ as pointwise limits of normalized irreducible spin characters of $G(m,1,n)$ as $n\to\infty$. Note that the relationship between characters and factor representations of finite type is clarified in [11].

To be more precise, we have  
$G(m,1,n)=D_n(T)\rtimes {\mathfrak S}_n$ with $D_n(T):=T^n,\,T={\boldsymbol Z}_m$, and 
the Schur multiplier $Z:={\mathfrak M}\big(G(m,1,n)\big)$ is equal to $Z_1$ for $m$ odd, and to $Z_1\times Z_2\times Z_3$ for $m$ even, where $Z_i=\langle z_i\rangle\cong{\boldsymbol Z}_2 ,\,z_i^{\;2}=e$, with the identity element $e$. 
The representation groups used here are special central extensions of $G(m,1,n)$ with the central subgroup $Z$ (which may be called a \lq universal covering group' of $G(m,1,n)$) given as 
\begin{eqnarray}
 R\big(G(m,1,n)\big)=\left\{
 \begin{array}{ll}
 \;D_n(T)\rtimes \widetilde{{\mathfrak S}}_n\quad&\mbox{\rm if $m$ is odd}, 
 \\[.3ex]
 \big(\widetilde{D}_n(T)\times Z_3\big)\rtimes\widetilde{{\mathfrak S}}_n&\mbox{\rm if $m$ is even,}
 \end{array}
\right.
\end{eqnarray}
where $\widetilde{{\mathfrak S}}_n$ is a representation group of the symmetric group ${\mathfrak S}_n$ and a central extension with the central subgroup $Z_1$, and $\widetilde{D}_n(T)$ is a double covering of $D_n(T)$ with central subgroup $Z_2$ (Theorem 1.1). Let $R\big(G(m,1,\infty)\big):=\lim_{n\to\infty}R\big(G(m,1,n)\big)$. 
A spin representation of $G(m,1,n)$ is essentially a  linear representation $\pi$ of $R\big(G(m,1,n)\big)$. If $\pi$ is irreducible, it has its own {\it spin type} $\chi^\pi_Z\in \widehat{Z}$ given by $\pi(z)=\chi^\pi_Z(z)I$, where $I$ denotes the identity operator. Suppose $m$ be even, then there are 8 different spin types containing the trivial one ${\bf 1}_Z\in\widehat{Z}$ which corresponds to linear irreducible representation (= IR) of the base group $G(m,1,n)$.  
 Take a non-trivial $\chi\in\widehat{Z}$ and put $Z^\chi:=\mathrm{Ker}(\chi)$, then any spin IR of spin type $\chi$ is essentially a linear representation of $R\big(G(m,1,n)\big)/Z^\chi$, a double covering of $G(m,1,n)$. Thus we can work, if we wish, mainly in the frame work of double coverings as was done in Stembidge [36],  
 Morris-Jones [29], [I] and [II]. However in the present paper, we work on the 8 times covering $R\big(G(m,1,n)\big)$ in principle and it enables us to simplify the methods used and clarify the relationships among different spin types.
  
In Part I, our important points are the following. 

(1) To classify and to construct spin IRs, we use \lq classical method for constructing IRs of a semidirect product group' in [7]. To apply this method to the semidirect product (0.1) above, {\it Step 1} $\sim$ {\it Step 4} are given in Theorem 2.2. 

(2) We apply the theory of the twisted central product of double covering groups, the twisted central product of their spin representations and its character formula, given in [19].

(3) To calculate spin characters, we apply a key lemma, Lemma 2.3, which is rather elementary but very effective to the concrete calculations of characters, and actually this enables us to simplify calculations in [II] very much. 

(4) Parametrization of spin dual of each spin type are respectively given by Young diagrams, shifted Young diagrams and certain extensions of them, such as Notations 3.2, 4.3, 4.6 and 7.6, with slight deviations in cases of spin types $\chi^{\rm III}$ and $\chi^{\rm VII}$. Irreducible spin characters are explicitly written by means of these parameters. 

In Part II, our important points are the following. Put 
$G'_\infty:=R\big(G(m,1,\infty)\big)$. 

(5) Evaluation in Table 12.1 of supports of $f \in K(G'_\infty; \chi)$, for each $\chi\in\widehat{Z}$, where $K(G'_\infty; \chi)$ denotes the set of all positive-definite invariant (or class) functions on $G'_\infty$ of spin type $\chi$, i.e., $f(zg')=\chi(z)f(g')\;(z\in Z, g'\in G'_\infty)$ (cf. [I, \S 10]). 

(6) Limit theorem, Theorem 12.5, asserts that every character of $G'_\infty:=R\big(G(m,1,\infty)\big)$ is a pointwise limit of a series of normalized irreducible characters of $R\big(G(m,1,n)\big)$ as $n\to\infty$. 

(7) Criteria (EF) and (EF$^{\chi^{\rm Y}}$) respectively in Theorems 12.3 and 12.4 for an $f\in K_1(G'_\infty;\chi)$ to be extremal or to be a character, where $K_1(G'_\infty;\chi)$ denotes the set of all $f\in K(G'_\infty;\chi)$ normalized as $f(e)=1$. 

(8) The data in (5) and (7) enable us to answer if the product of two characters is also a character or not (Theorem 12.8, see also Theorem 22.4). 

(9) Parameters for characters of spin type $\chi$ and for their sets $E(G'_\infty;\chi)$ are given by $A\in{\cal A}({\cal K})$ in Notation 14.3 and $C\in{\cal C}({\cal K})$ in Notation 16.5, depending on the spin type $\chi$, where ${\cal K}=\widehat{T}$ or ${\cal K}=\widehat{T}^0$ the half of $\widehat{T}$.

(10) 
The hereditary property is studied further, mainly in \S 15, for $G=G(m,1,\infty)$ and its normal subgroups $N$ such as her child groups $G(m,p,\infty),\,p|m,\,p>1,$ and $G^{\mathfrak A}(m,1,\infty):=(D_\infty(T)\times Z_3)\rtimes\widetilde{{\mathfrak A}}_\infty$ and $G^{\mathfrak A}(m,p,\infty):=G(m,p,\infty)\cap G^{\mathfrak A}(m,1,\infty)$. By the resriction map ${\rm Res}^G_N:\, E(G)\ni F \mapsto f=F|_N$, we can get the set $E(N)$ of all characters of $N$. Furthermore, going up to the spin level, we study the heredity from mother group $R\big(G(m,1,n)\big)$ to child groups $R\big(G(m,p,n)\big),\,p|m,\,p>1,$ particularly in detail for the case $n=\infty$, and clarified that the main kernel of spin characters of child groups can be obtained. In addition, in Case EE, where $p$ and $q=m/p$ are even, we propose Conjecture 15.17.

(11) 
The heredity at the non-spin level can be expressed by means of certain actions of finite groups (or symmetries) in the corresponding parameter spaces such as ${\cal A}({\cal K})$ etc. 
For example, in the case of child groups 
 $N=G(m,p,\infty)$, all the spin characters of $N$ are obtained by restriction map ${\rm Res}^G_N$ from $G=G(m,1,\infty)$ to $N$, and this can be expressed in a beautiful way by means of symmeries on the parameter spaces ${\cal A}({\cal K})$ with ${\cal K}=\widehat{T}$ or $\widehat{T}^0$ (see Theorems 15.1, 15.18, 15.19 and 15.20).

(12) Mutual relations between the sets of characters $E(G'_\infty;\chi)$ of different spin types $\chi$, for certain pair of $\chi$'s, are displayed by diagrams in Theorems 17.6, 19.5, 20.2 and 21.2.  
\vskip1em

This paper is organized as shown in the table of contents above, and the final results on characters of $G'_\infty$ are summarized in \S 22, in particular in Tables 22.2 and 22.5.


\part{\Large Spin irreducible representations and  their characters of generalized symmetric groups $G(m,1,n),\,4\le n<\infty$}

\section{\large Representation groups of generalized symmetric groups}

{\bf 1.1. Representation groups of $G(m,1,n)$ and their inductive limits.} 
\vskip.5em

We introduce our main objects, quoting our previous papers \cite{[HHH2013]} and \cite{[HHoH2013b]} (which we quote as [I] and [II] respectively in the following). 
Generalized symmetric groups $G(m,1,n),\,4\le n<\infty$, and their representation groups $R\big(G(m,1,n)\big)$, and also their inductive limits: 
\begin{eqnarray}
\label{2015-05-12-2}
&&\quad
\begin{array}{l}
G(m,1,\infty)=\lim_{n\to\infty}G(m,1,n)={\bigcup}_{n\geqslant 4}G(m,1,n),
\\[.8ex]
R\big(G(m,1,\infty)\big)=\lim_{n\to\infty} R\big(G(m,1,n)\big)={\bigcup}_{n\geqslant 4}R\big(G(m,1,n)\big).  
\end{array}
\end{eqnarray}
  Put\;  
$
T:={\boldsymbol Z}_m,\,D_n(T):={\prod}_{j\in{\boldsymbol I}_n}\!T_j,\,T_j=T,\,{\boldsymbol I}_n:=\{1,2,\ldots,n\},
$ 
and let $\widetilde{D}_n=\widetilde{D}_n(T)$ be the central extension of $D_n=D_n(T)$, given below, with central subgroup $Z_2=\langle z_2\rangle\cong {\boldsymbol Z}_2,\;z_2^{\;2}=e$, and the canonical homomorphism $\Phi_D$\,: 
\begin{eqnarray}
\label{2016-05-02-2}
\{e\}\longrightarrow Z_2\longrightarrow \widetilde{D}_n(T)\stackrel{\Phi_D\,}{\longrightarrow} D_n(T)\longrightarrow \{e\}. 
\end{eqnarray}
Let $y$ be a fixed generator of $T$ and $y_j\in T_j$ its $j$-th copy, and $\eta_j\in \widetilde{D}_n$ be its preimage such that 
\begin{eqnarray}
\label{2016-05-02-3}
&&
\Phi_D(\eta_j)=y_j\;\big(j\in {\boldsymbol I}_n\big),\;\;\eta_j\eta_k=z_2\eta_k\eta_j\;(j\ne k).
\end{eqnarray}

Moreover $\widetilde{{\mathfrak S}}_n=R({\mathfrak S}_n)$ is the representation group ${\mathfrak T}'_n$ of ${\mathfrak S}_n$ given by Schur as 
\begin{eqnarray}
\label{2016-05-02-4}
&&
\begin{array}{l}
\;\{e\}\longrightarrow Z_1\longrightarrow \widetilde{{\mathfrak S}}_n\stackrel{\Phi_{\mathfrak S}\,}{\longrightarrow} {\mathfrak S}_n\longrightarrow \{e\}, \;\;Z_1=\langle z_1\rangle,\;z_1^{\;2}=e,
\\[.5ex]
\widetilde{{\mathfrak S}}_n=\langle z_1,r_i\;(i\in {\boldsymbol I}_{n-1})\rangle,\;\;
\Phi_{\mathfrak S}(r_i)=s_i=(i\;\,i{\hspace{-.1ex}+\hspace{-.1ex}}1)\;\;(i\in {\boldsymbol I}_{n-1}),\end{array}
\end{eqnarray}
with the fundamental relations as 
\begin{eqnarray}
\label{2016-05-02-5}
\left\{
\begin{array}{l}
r_i^{\;2}=e\;\;(i\in{\boldsymbol I}_{n-1}), \;\; (r_ir_{i+1})^3 =e\;\;(i\in{\boldsymbol I}_{n-2}),
\\[.5ex]
r_ir_j=z_1r_jr_i\quad(|i-j|\ge 2).
\end{array}
\right.
\end{eqnarray}

For $n\ge 1$, put $G(m,1,n)=D_n(T)\rtimes{\mathfrak S}_n$, with the action of $\sigma\in{\mathfrak S}_n$ on $d=(t_j)_{j\in{\boldsymbol I}_n}\in D_n(T),\,t_j\in T_j,$ as \,$\sigma(d):=(t'_j)_{j\in{\boldsymbol I}_n},\,t'_j=t_{\sigma^{-1}(j)}\;(j\in{\boldsymbol I}_n)$. For $n\ge 4$, 
the Schur multiplier ${\mathfrak M}(G_n)$ for the generalized symmetric group $G_n=G(m,1,n)$ is given as \cite{[DaMo1974]}\;
\begin{eqnarray}
\label{2016-05-04-1}
{\mathfrak M}(G_n)=\left\{
\begin{array}{ll}
Z_1=\langle z_1\rangle \cong {\boldsymbol Z}_2\quad &\quad\mbox{\rm in case $m$ is odd,}
\\[.5ex]
\prod_{1\leqslant i\leqslant 3}Z_i\cong {\boldsymbol Z}_2^{\;3} \quad &\quad\mbox{\rm in case $m$ is even,}
\end{array}
\right.
\end{eqnarray}
where $Z_i=\langle z_i\rangle,\;z_i^{\;2}=e\;(1\le i\le 3)$. 
A representation group  
\begin{eqnarray}
\label{2017-11-16-1}
&&
\{e\}\longrightarrow Z\longrightarrow R\big(G(m,1,n)\big)\stackrel{\Phi\,}{\longrightarrow} G(m,1,n)\longrightarrow \{e\}\quad({\rm exact}),
\end{eqnarray}
with $Z:={\mathfrak M}(G_n)$, is given as follows. 
\vskip1em 

{\bf Theorem 1.1.}  [II, \S 3.1]\; {\it 
Let $4\le n<\infty$.\; 

{\rm (i)} 
In case $m$ is odd,  
\begin{eqnarray}
\label{2016-05-04-2}
R\big(G(m,1,n)\big)=D_n(T)\rtimes \widetilde{{\mathfrak S}}_n,
\end{eqnarray}
where $\sigma'\in\widetilde{{\mathfrak S}}_n$ acts on $d=(t_j)_{j\in{\boldsymbol I}_n}\in D_n(T)$ through $\sigma=\Phi_{\mathfrak S}(\sigma')\in{\mathfrak S}_n$ as \;$\sigma'(d)=\sigma'd{\sigma'}^{\,-1}=(d'_j)_{j\in{\boldsymbol I}_n},\;d'_j=d_{\sigma^{-1}(j)}$. 

{\rm (ii)} 
In case $m$ is even, 
\begin{eqnarray}
\label{2016-05-04-3}
R\big(G(m,1,n)\big)=\widetilde{D}^\vee_n(T)\rtimes \widetilde{{\mathfrak S}}_n,\quad 
\widetilde{D}^\vee_n(T):=\widetilde{D}_n(T)\times Z_3,
\end{eqnarray}
where $r_i\in\widetilde{{\mathfrak S}}_n\;(i\in{\boldsymbol I}_{n-1})$ acts on the generators $\eta_j\in \widetilde{D}_n(T)$ as 
 \begin{eqnarray}
\label{2016-05-04-4}
r_i(\eta_j)=\left\{
\begin{array}{ll}
z_3\eta_{s_i(j)}\quad & (j=i,i+1), 
\\[.5ex]
\;\eta_j &(j\ne i,i+1).
\end{array}
\right.
\end{eqnarray}
Or, with $\widehat{\eta}_j:=z^{\;j-1}_3\eta_j\;(j\in{\boldsymbol I}_n)$\, and\, 
$D^\wedge_n(T):=\big\langle z_2, \widehat{\eta}_j\;(j\in{\boldsymbol I}_n)\big\rangle$, 
\begin{gather}
\label{2016-05-04-5}
\left\{
\begin{array}{l}
r_i\,\widehat{\eta}_j\,r_i^{\;-1}=z_3\,\widehat{\eta}_{s_i(j)}\quad(i\in{\boldsymbol I}_{n-1},\,j\in{\boldsymbol I}_n), 
\\[.5ex]
\widetilde{D}^\vee_n(T)\cong D^\wedge_n(T)\times Z_3.
\qquad
\end{array}
\right.
\end{gather}
}

We abbreviate the notations as 
 $D_n:=D_n(T),\,\widetilde{D}_n:=\widetilde{D}_n(T),\,D^\wedge_n:=D^\wedge_n(T),\,S={\mathfrak S}_n,\,S'=\widetilde{{\mathfrak S}}_n,\,G_n=G(m,1,n),$ and $G'_n=R\big(G(m,1,n)\big).$
Put $Z={\mathfrak M}\big(G(m,1,n)\big)$, or $Z=Z_1$ for $m$ odd and $Z=Z_1\times Z_2\times Z_3$ for $m$ even. 

The semidirect product structure at $n=\infty$ is given naturally as  
\begin{eqnarray}
\label{2015-05-12-2-9}
&&
G_\infty:=G(m,1,\infty)=D_\infty\rtimes {\mathfrak S}_\infty,\;D_\infty=D_\infty(T)=\lim_{n\to\infty}D_n(T),
 \\[.5ex]
\label{2017-07-25-1}
&&
G'_\infty:=R\big(G(m,1,\infty)\big)=
\left\{
\begin{array}{ll}
\; D_\infty\rtimes \widetilde{{\mathfrak S}}_\infty\;&
\mbox{\rm for $m$ odd},
\\[.3ex]
(\widetilde{D}_\infty\times Z_3)\rtimes \widetilde{{\mathfrak S}}_\infty\;&
\mbox{\rm for $m$ even},
\end{array}
\right.
\end{eqnarray}
where $\widetilde{D}_\infty=\widetilde{D}_\infty(T)=\lim_{n\to\infty}\widetilde{D}_n(T)$.  

\vskip1.2em 

{\bf 1.2. Projective representation.}\; Let $G$ be a topological group. A unitary {\it projective representation} (=\;{\it spin representation}) of $G$ is a continuous map $G\ni g \mapsto \pi(g)$, a unitary transformation on a complex Hilbert space $V(\pi)$, satisfying 
\begin{eqnarray}
\pi(g)\pi(h)=\alpha_{g,h}\pi(gh)\quad (\alpha_{g,h}\in{\boldsymbol C}^\times),
\end{eqnarray}
where the function $\alpha_{g,h}$ on $G\times G$ is called the {\it factor set} of $\pi$.  

A spin irreducible representation (=IR) $\pi$ of $G(m,1,n)$ has its spin type $\chi\in\widehat{Z}$ defined as follows. We can lift up $\pi$ to a linear representation $\pi'$ of $R\big(G(m,1,n)\big)$ by taking a section of $G(m,1,n)$ in $R\big(G(m,1,n)\big)$ [7, \S 2.3]. 
By abuse of terminology we say $\pi'$ a spin IR of $G(m,1,n)$ and also of $R\big(G(m,1,n)\big)$. 
 Each IR $\Pi$ of $R\big(G(m,1,n)\big)$ has its central character $\chi\in \widehat{Z}$ such as $\Pi(z)=\chi(z)I\;(z\in Z)$ with the identity operator $I$, and we call $\chi$ the {\it spin type} of $\Pi$. We call the spin type of $\pi'$ also that of $\pi$. In case $m$ is odd, the non-trivial spin type is denoted by $\chi^{\rm odd}$, and in case $m$ is even, the spin type $\chi\in\widehat{Z}$ is expressed as $\chi^{\rm Y}=(\beta_1,\beta_2,\beta_3)$ with $\beta_i=\chi^{\rm Y}(z_i)=\pm 1$, where the classifying number Y\,=\,I, II, \ldots, VIII are given as follows (see Table 9.1 in [II]): 
\vskip.5em
$
\begin{array}{llll}
\chi^{\rm I}=(-1,-1,-1),&\chi^{\rm II}=(-1,-1,\,1),&\chi^{\rm III}=(-1,\,1,-1),&\chi^{\rm IV}=(-1,\,1,\,1),
\\
\chi^{\rm V}=(1,-1,-1),&
\chi^{\rm VI}=(1,-1,\,1),&
\chi^{\rm VII}=(1,\,1,-1),&
\chi^{\rm VIII}=(1,\,1,\,1).
\end{array}
$

\vskip1.2em

{\bf 1.3. Spin characters of $G(m,1,n)$ and of $G(m,1,\infty)$.}\; 
Let $G$ be a topological group satisfying Hausdorff separation axiom, and denote by $K(G)$ the set of all continuous positive-definite central  functions on $G$ and by $K_1(G)$ the set of those $f\in K(G)$ normalized as $f(e)=1$ at the identity element $e$, and by $E(G)$ the set of all extremal elements of the convex set $K_1(G)$. Every $f\in E(G)$ is called a normalized {\it character} of $G$. 
Let $\pi$ be a finite dimensional representation,  and $\chi_\pi$ its trace character and $\widetilde{\chi}_\pi:=\chi_\pi/\dim \pi$ its normalized character. 
When $G$ is compact, $E(G)$ is equal to the set of all normalized irreducible characters $\widetilde{\chi}_\pi$ of IRs $\pi$ of $G$.  In general, every element $f\in E(G)$ corresponds canonically and bijectively to a quasi-equivalence class of factor representations of finite type [11, Theorem B].

Now take $G=G'_n=R\big(G(m,1,n)\big),\,4\le n\le \infty$. 
A function $f$ on 
$G'_n$ is said to have a {\it spin type} $\chi\in\widehat{Z}$, if 
\begin{eqnarray}
\label{2017-07-23-1}
 f(zg')=\chi(z)f(g')\quad(z\in Z,\,g'\in G'_n).
\end{eqnarray}
By the Fourier expansion with respect to $Z$, we see that any function $f$ on $G'_n$ is expressed as $f=\sum_{\chi\in\widehat{Z}}f^\chi$, with $f^\chi$ of spin  type $\chi$. On the other hand, we have, for $4\le n\le \infty$, a decomposition
\begin{eqnarray}
\label{2015-05-07-1}
E(G'_n)={\bigsqcup}_{\chi\in\widehat{Z}}\,E(G'_n;\chi), 
\end{eqnarray}
where $E(G'_n;\chi)$ denotes the set of all $f\in E(G'_n)$ with spin type $\chi$. 
\vskip1em 

{\bf 1.4.\; Limits of irreducible characters as $n\to\infty$.}\; 
Along with \,$G_n\to G_\infty=\bigcup_{n\geqslant 4}G_n$\, or with \,$G'_n\to G'_\infty=\bigcup_{n\geqslant 4}G'_n$\, as $n\to\infty$, we study pointwise limits of normalized irreducible characters $\widetilde{\chi}_{\pi_n}$ for a series of IRs $\pi_n$ of $G_n$ or of $G'_n$: $f_\infty=\lim_{n\to\infty}\widetilde{\chi}_{\pi_n}$. This kind of study started in the case of the symmetric groups ${\mathfrak S}_n\to{\mathfrak S}_\infty$ by A.\,Vershik\,-\,S.\,Kerov (\cite{[VK1981]}, \cite{[VK1988]}), and succeeded by M.\,Nazarov \cite{[Naz1992]} for the case of representation groups $\widetilde{{\mathfrak S}}_n\to\widetilde{{\mathfrak S}}_\infty$.  

An abstract asymptotic theory of characters is generalized in [23, \S\S 2--3] 
to cover the case of inductive limits of compact groups $H_n\to H_\infty=\bigcup_{n\geqslant 1}H_n$. It is  
explicitly applied to the case of wreath product of ${\mathfrak S}_n$ with any compact group $T$ as \,$H_n={\mathfrak S}_n(T)=D_n(T)\rtimes {\mathfrak S}_n \to H_\infty=D_\infty(T)\rtimes {\mathfrak S}_\infty$, and proved that \,{\it any character $f\in E(H_\infty)$ can be obtained as the pointwise limit $f_\infty=\lim_{n\to\infty}\widetilde{\chi}_{\pi_n}$ of a series of normalized characters $\widetilde{\chi}_{\pi_n}$ of IRs $\pi_n$ of $H_n$} [ibid., Theorem 4.3], that is, $E(H_\infty) \subset {\rm Lim}(H_\infty)$, where ${\rm Lim}(H_\infty)$ denotes the set of all such limits. 

 On the other hand, in case $T$ is not finite, we obtain, as weak limits $f_\infty=\lim_{n\to\infty}\widetilde{\chi}_{\pi_n}$, many everywhere discontinuous functions (called {\it bad}\, limits) [17, \S 5],  
 apart from the totality $E(H_\infty)$ of characters of $H_\infty$ given as {\it good}\, limits. We know that the set of these bad limits correspond to the difference of Martin boundary $\partial {\mathbb Y}(\widehat{T})$ and minimal Martin boundary $\partial_m {\mathbb Y}(\widehat{T})$ in the notation of Theorem 3.4 in \cite{[HoH2014]}.  Recently we proved (but not yet published) the bad limits are also pointwise limits of $\widetilde{\chi}_{\pi_n}$ for well-chosen series $\pi_n$'s, and just corresponds to the difference ${\rm Lim}(H_\infty)\setminus E(H_\infty)$ of two limit sets. 

In this paper, we calculate all the pointwise limits $f_\infty=\lim_{n\to\infty}\widetilde{\chi}_{\pi_n}$,   
and prove that their set ${\rm Lim}(G_\infty)$ or ${\rm Lim}(G'_\infty)$ is just equal to the set of characters $E(G_\infty)$ or to $E(G'_\infty)$ respectively, and so there is no bad limit in the present case. For any $\chi\in\widehat{Z}$, if $\pi_n$'s are of spin type $\chi$ or $\widetilde{\chi}_{\pi_n}\in E(G'_n;\chi)$, then the limit $f_\infty$ covers the totality of $E(G'_\infty;\chi)$. Thus we obtain explicit character formula for any $f\in E(G'_\infty)$. Note that, in the case of $\chi=\chi^{\rm Y}$ with Y=I, II, VII and VIII, this result has been already reported in [I] and [II].

\section{\large Preliminaries to Part I}

\markboth{\hfil{\small\sc \;T.\ Hirai and A.\ Hora}\qquad\qquad\hfil}{{\thesection}\quad\mbox{\it\small Introduction}} 

We study spin irreducible representations (=IRs) and their characters of generalized symmetric groups, using results in [I] and [II]. 
The results in this part will be applied in the second part to obtain characters of $R\big(G(m,1,\infty)\big)$ by taking limits as $n\to\infty$.
\vskip1.2em

{\bf 2.1. Construction of IRs of semidirect product groups.}

We apply the method given below for constructing IRs of semidirect product groups. 
Let $G$ be a finite group of semidirect product type\;  
$
G=U\rtimes S,
$ 
where $U$ is a normal subgroup of $G$, and $S$ acts on $U$.  The action of $s\in S$ on $u\in U$ is denoted by $s(u)$. We quote from \cite{[Hir2013]}  
 some results necessary in the following.

Take a unitary IR $\rho$ of $U$, on a Hilbert space $V(\rho)$, and consider its equivalence class $[\rho]\in\widehat{U}$. Every $s\in S$ acts on $\rho$ as \;$s(\rho)(u):=\rho\big(s^{-1}(u)\big)\;(u\in U)$, realised on the same space $V(\rho)$, and on equivalent classes as $[\rho]\to [s(\rho)]$. Denote by $\widehat{U}/S$ the set of $S$-orbits in the dual $\widehat{U}$ of $U$. 
Take the stationary subgroup $S([\rho])=\{s\in S\,;\,s(\rho)\cong \rho\}$ of a representative $[\rho]$ of an $S$-orbit in $\widehat{U}/S$. Put $H:=U\rtimes S([\rho])$. For $s \in S([\rho])$, we determine explicitly an intertwining operator $J_\rho(s): V(\rho)\to V(\rho)$ as   
\begin{eqnarray}
\label{2011-04-05-21}
\rho\big(s(u)\big)=J_\rho(s)\,\rho(u)\,J_\rho(s)^{-1}\quad(u\in U). 
\quad
\end{eqnarray}
Then it is determined up to a non-zero scalar factor, and   
we have naturally a spin (or projective) representation $S([\rho])\ni s \mapsto J_\rho(s)$. Let $\alpha_{s,t}\in {\boldsymbol T}^1=\{z\in{\boldsymbol C}\,;\,|z|=1\}$ be its factor set given as \;
\begin{eqnarray}
\label{2017-06-07-1}
 J_\rho(s)J_\rho(t)= \alpha_{s,t}\,J_\rho(st)\quad 
 \big(s,t \in S([\rho])\big).
\end{eqnarray}

Let $Z$ be the finite subgroup of ${\boldsymbol T}^1$ generated by $\{\alpha_{s,t}\}$ and put $K'=Z\times K$ for $K=S([\rho])$ with the product $(z,s)(z',t):=(zz'\alpha_{s,t},st)$. Then $K'=S([\rho])'$ is a central extension of $K=S([\rho])$ associated to the cocycle $\alpha_{s,t}$ 
[ibid., Lemma 2.2]:   
 \begin{eqnarray}
\label{2011-04-05-22}
\nonumber
1\longrightarrow Z \longrightarrow S([\rho])' \stackrel{\Phi_S}{\longrightarrow} S([\rho]) \longrightarrow 1\quad\mbox{\rm (exact)}, 
\quad
\end{eqnarray}
where $\Phi_S$ denotes the canonical homomorphism. 
Then $J_\rho$ can be lifted up to a linear representation $J'_\rho$ of $S([\rho])'$, acting on the same space $V(\rho)$. 
  Put $H':=U\rtimes S([\rho])'$ with the action $s'(u):=s(u),\,s'\in S([\rho])',\;s=\Phi_S(s')$, and   
\begin{eqnarray}
\label{2011-04-05-23}
\nonumber
\pi^0\big((u,s')\big):=\rho(u)\cdot J'_\rho(s')\quad\big(u\in U,\,s'\in S([\rho])'\big). 
\quad
\end{eqnarray}
Then $\pi^0=\rho\cdot J'_\rho$ is an IR of $H'$ on $V(\rho)$.
\vskip1em 

Now take an IR $\pi^1$ of $S([\rho])'$ and consider it as a representation of $H'$ through the quotient map $H' \to H'/U\cong S([\rho])' $, denoted as $\widetilde{\pi}^1$, and consider inner tensor product $\pi^0\otimes \widetilde{\pi}^1$ as representation of $H'$, which we denote as 
$\pi^0\boxdot \pi^1$ in short. Take a $\pi^1$ with the inverse factor set $\beta_{s,t}=\alpha_{s,t}^{\;-1}$. Then that of $\pi:=\pi^0\boxdot \pi^1$ is trivial and $\pi$ becomes a linear (or non-spin) representation of the base group $H=U\rtimes S([\rho])$ of $H'$.  
 Thus we obtain a representation of $G$ as 
\begin{eqnarray}
\label{2012-11-19-41}
\Pi(\pi^0,\pi^1) :=\mathrm{Ind}^G_H(\pi^0\boxdot \pi^1). 
\quad
\end{eqnarray}

{\bf Theorem 2.1.} ([7, Theorem 4.1])\;    
{\it 
For an IR $\rho$ of $U$ and a spin IR $\pi^1$ of $S([\rho])$ with the inverse factor set, 
the induced representation\; 
$\Pi(\pi^0,\pi^1)=\mathrm{Ind}^G_H (\pi^0\boxdot\pi^1)$\; 
of\, $G$ is irreducible. 
}
\bigskip

Put $\Pi=\mathrm{Ind}_H^G\pi$, $\pi=\pi^0\boxdot\pi^1$, then its  character is given by 
\begin{eqnarray}
\label{2011-11-19-21}
&&
\; 
\chi_\Pi(g) = \sum_{\dot{k}\in H\backslash G}\chi_\pi(kgk^{-1})= \frac{1}{|H|}\sum_{g\in G}\chi_\pi(kgk^{-1}). 
\end{eqnarray}
Here $\chi_\pi$ is extended from $H$ to $G$ by putting 0 outside $H$, and $\dot{k}=Hk$. 

Assume that \,$\big\{\rho_i\,;\,\mbox{\rm IR of $U$},\,i\in I_{U;\,S}\big\}$\, gives a complete set of representatives (= CSR) for $\widehat{U}/S$, and for each $i\in I_{U;\,S}$, let 
\;$\big\{\pi^1_{i,j}\;;\;j\in J_i\big\}$ be a CSR of equivalence classes of  
 IRs of $S([\rho_i])'$ with the factor set inverse to that of $J_{\rho_i}$. Put \,$H_i:=U\rtimes S([\rho_i]),\;H'_i:=U\rtimes S([\rho_i])'$,\, $\pi_{i,j}=\pi^0_i\boxdot \pi^1_{i,j}$, and \,$\Pi_{i,j}=\Pi(\pi^0_i,\pi^1_{i,j})=\mathrm{Ind}^G_{H_i}\pi_{i,j}$. 
 Define a set of IRs of $G$ as  
 \begin{eqnarray}
\label{2012-11-20-42}
\Omega(G):=\big\{\Pi_{i,j}:=\Pi(\pi^0_i,\pi^1_{i,j})\;;\;i\in I_{U;\,S},\, j\in J_i\big\}.
\end{eqnarray}

\vskip.2em

{\bf Theorem 2.2.} ([7, Theorem 5.1])\; 
{\it 
The set of characters \,$\big\{\chi_{\Pi_{i,j}}\,;\,\Pi_{i,j}\in\Omega(G)\big\}$\, gives a complete orthonormal system in the subspace of central elements of $L^2(G)$ with respect to the normalized Haar measure on $G$.  
 Consequently the set\, $\Omega(G)$ gives a complete set of representatives for the dual \,$\widehat{G}$ of $G$.  

}

\vskip1.2em

Thus, to construct a complete set of IRs of $G$, we follow the following steps: 
\vskip.7em

{\it Step 1.} Choose a complete set $\{\rho_i\,;\,i\in I_{U;S}\}$ of representatives for $\widehat{U}/S$. 

{\it Step 2.} Determine the stationary subgroup $S([\rho_i])$ in $S$. 

{\it Step 3.} Calculate explicitly the intertwining operators $J_{\rho_i}(s)\;\big(s\in S([\rho_i]\big).$

{\it Step 4.} Determine a complete set $\big\{\pi^1_{i,j}\,;\,j\in J_i\big\}$ of representatives of spin IRs of $S([\rho_i])'$ with the factor set inverse to that of   
 $J'_{\rho_i}$. Then take 
the induced representations $\Pi_{i,j}=\mathrm{Ind}^G_{H_i}(\pi^0_i\boxdot\pi^1_{i,j})$. 

\vskip1.2em

{\bf 2.2. A key lemma for induced representations.}\; 
Let $G$ be a finite group and $H$ its subgroup. 
The following property of induced representations will work as a key lemma for calculating characters of spin IRs of $G'_n=R\big(G(m,1,n)\big)$.  
\vskip1em

{\bf Lemma 2.3.}\,  
{\it 
Let \,${\cal T}$ and $\pi$ be linear representations of $G$ and $H$ respectively. Then 
\begin{eqnarray}
\label{2016-05-27-1}
&&
\mathrm{Ind}^G_H\big(({\cal T}|_H)\otimes\pi\big)\,\cong\, {\cal T}\otimes \mathrm{Ind}^G_H\,\pi.
\end{eqnarray}
}

{\it Proof.} 
The character of $\Pi=\mathrm{Ind}^G_H\big(({\cal T}|_H)\otimes\pi\big)$ is calculated as, for $h\in H$, 
\begin{eqnarray}
\label{2016-05-27-3}
&&
\chi_\Pi(h)=\frac{1}{|H|}\,\sum_{g\in G:\,ghg^{-1}\in H}\chi_{{\cal T}}(ghg^{-1})\,\chi_{\pi}(ghg^{-1}),
\end{eqnarray}
and $\chi_\Pi(g_0)=0$ for $g_0\in G$ not conjugate to any element in $H$. 
Since $\chi_{{\cal T}}$ is $G$-invariant, it turns out to be 
\begin{eqnarray*}
\label{2016-05-27-4}
&&
\chi_\Pi(h)
=
\chi_{{\cal T}}(h)\cdot\frac{1}{|H|}\sum_{g\in G:\,ghg^{-1}\in H}\chi_{\pi}(ghg^{-1})
=
\chi_{{\cal T}}(h)\times \chi_{\mathrm{Ind}^G_H \pi}(h)\,.
\end{eqnarray*}
Hence we get\; $\chi_\Pi=\chi_{{\cal T}}\cdot \chi_{\mathrm{Ind}^G_H \pi}$. Consequently \;$\Pi\cong {\cal T}\otimes \big(\mathrm{Ind}^G_H \pi\big)$. 
\hfill 
$\Box$\;
\vskip1em

{\bf Example 2.4.}\; For any linear representation ${\cal T}$ of $G$,  
$$
\mathrm{Ind}^G_H\big({\cal T}|_H\big)=\mathrm{Ind}^G_H\big(\mathrm{Res}^G_H {\cal T}\big)\cong {\cal T}\otimes \mathrm{Ind}^G_H{\boldsymbol 1}_H, 
$$
where $\mathrm{Ind}^G_H{\boldsymbol 1}_H$ is the quasi-regular representation on $\ell^2(H\backslash G)$. 
More specifically, take the case of $H=\{e\}$ the trivial group, then the 2nd component in the most right hand side is the regular representation ${\cal R}$ of $G$ on $\ell^2(G)$, and we have   
\begin{eqnarray}
\label{2016-06-09-1}
{\cal T}\otimes {\cal R} \cong [\dim {\cal T}]\cdot {\cal R}, 
\end{eqnarray}
where $[\dim {\cal T}]\cdot$ denotes the multiplicity.

\section{\large Linear IRs of $G(m,1,n)$ and their characters}

Let $\zeta^{(a)},\,0\le a<m,$ be a  character of $T={\boldsymbol Z}_m$ defined as 
\begin{eqnarray}
\label{2015-02-08-11-8}
\zeta^{(a)}(y):=\omega^a,\quad\omega=e^{2\pi i/m}.
\end{eqnarray}
The dual of $T$ is expressed as $\widehat{T}=\{\zeta^{(a)}\,;\,0\le a<m\}$, and we give a linear order in it according to $a$. Put, for even $m=2m^0$, the half of $\widehat{T}$ as  
\begin{eqnarray}
\label{2016-04-11-41-8}
\widehat{T}^0:=\{\zeta^{(a)}\,;\,0\le a<m^0:=m/2\} \,\subset\,\widehat{T}.
\end{eqnarray}

{\bf 3.1. Linear irreducible representations.}\; 
 Let us review shortly this case since it will become   indispensable later, and on the way we introduce necessary notations. 
We apply the general method of constructing IRs of semidirect product groups, given in {\bf 2.1}, to $G_n=D_n\rtimes {\mathfrak S}_n$. Denote this group as $G=U\rtimes S$. The spin type of these linear (non-spin) IRs is 
$\chi^{\rm VIII}=(1,1,1)$. 
\vskip.6em

{\it Step 1.}\;    
The dual $\widehat{U}$ consists of $[\rho]=\{\rho\}$ of one-dimensional characters $\rho$ given by a partition ${\cal I}_n=(I_{n,\zeta})_{\zeta\in\widehat{T}}\,:$ ${\boldsymbol I}_n={\bigsqcup}_{\zeta\in\widehat{T}}\,I_{n,\zeta}\;(I_{n,\zeta}=\emptyset$ is admitted), in such a way that, for  
$d=(t_i)_{i\in{\boldsymbol I}_n}\in D_n,\,t_i\in T_i=T,$ 
\begin{eqnarray}
\label{2014-03-18-1-8}
\rho(d):={\prod}_{\zeta\in\widehat{T}}{\prod}_{\;i\in I_{n,\zeta}}\zeta(t_i).
\end{eqnarray}
We denote $\rho$ by $\rho_{{\cal I}_n}$. The action of $\sigma \in S={\mathfrak S}_n$ on $\rho$ is defined as 
$(\sigma\rho)(d):=\rho\big(\sigma^{-1}d\big)\;(d\in D_n)$, for which a partition $(I'_\zeta)_{\zeta\in\widehat{T}},\;I'_\zeta=\sigma^{-1}(I_{n,\zeta})$, is associated. 
Therefore a CSR of $\widehat{U}/S$ is given as follows. 
According to the order of $\zeta\in\widehat{T}$, take subintervals $I_{n,\zeta}$ of ${\boldsymbol I}_n$ successively, then we get a so-called {\it standard partition} ${\cal I}_n:=\big(I_{n,\zeta}\big)_{\zeta\in\widehat{T}}$.     
Put $n_\zeta:=|I_{n,\zeta}|$, then \,${\boldsymbol \nu}:=(n_\zeta)_{\zeta\in\widehat{T}}$ is a partition of $n$ indexed by $\zeta$. 
For ${\cal K}=\widehat{T}$ or ${\cal K}=\widehat{T}^0$, put 
\begin{eqnarray}
\label{2016-05-05-1}
P_n({\cal K}):=\big\{ {\boldsymbol \nu}=(n_\zeta)_{\zeta\in{\cal K}}\;;\;{\sum}_{\zeta\in{\cal K}}\,n_\zeta=n\;(n_\zeta=0\; \mbox{\rm admitted})
\big\}.
\end{eqnarray}
A standard partition ${\cal I}_n$ of $n$ and a ${\boldsymbol \nu}\in P_n(\widehat{T})$ correspond each other bijectively.   

\vskip.5em

{\it Step 2.}\;  Let ${\mathfrak S}_{I_{n,\zeta}}$ be the subgroup of ${\mathfrak S}_n$ acting on $I_{n,\zeta}\subset{\boldsymbol I}_n$, and call\, 
${\displaystyle 
{\prod}_{\zeta\in\widehat{T}}\,{\mathfrak S}_{I_{n,\zeta}}
}$\, {\it Frobenius-Young subgroup} of ${\mathfrak S}_n$ corresponding to the partition ${\boldsymbol \nu}\in P_n(\widehat{T})$, which will be identified with ${\mathfrak S}_{\boldsymbol \nu}:={\displaystyle {\prod}_{\zeta\in\widehat{T}}\,{\mathfrak S}_{n_\zeta}}$. We admit $P_n(\widehat{T}^0)\subset P_n(\widehat{T})$ naturally. 

\vskip1em

{\bf Lemma 3.1.}\;\; 
{\it 
The stationary subgroup $S([\rho])$ of $[\rho]$ is given as   
\begin{eqnarray}
\label{2014-03-18-2}
S([\rho])={\prod}_{\zeta\in\widehat{T}}{\mathfrak S}_{I_{n,\zeta}}
\cong{\mathfrak S}_{\boldsymbol \nu},\quad {\boldsymbol \nu}=(n_\zeta)_{\zeta\in\widehat{T}},\,n_\zeta:=|I_{n,\zeta}|.
\end{eqnarray}
}

{\it Step 3.}\; An intertwining operator $J_\rho(\sigma)$ satisfying \;$\rho\big(\sigma(d)\big)=J_\rho(\sigma)\rho(d)J_\rho(\sigma)^{-1}$\; is trivial: $J_\rho(\sigma)={\bf 1}$ (1-dimensional identity operator). So,  
\begin{eqnarray}
\label{2016-02-04-1-8}
\pi^0:\,U\rtimes S([\rho])\ni(d,\sigma)\mapsto \rho(d)\cdot J'_\rho(\sigma)=\rho(d).
\end{eqnarray}

{\bf Notation 3.2.}\, Let ${\boldsymbol Y}_{\!n}$ be the set of all Young diagram of size $n$, and put
${\boldsymbol Y}=\bigsqcup_{n\geqslant 0}{\boldsymbol Y}_{\!n}$ with ${\boldsymbol Y}_{\!0}=\{\varnothing\}$. A Young diagram $D_{\boldsymbol \lambda}$ with $i$-th row of length $\lambda_i$ and a partition ${\boldsymbol \lambda}=(\lambda_1,\lambda_2,\ldots),\,\lambda_1\ge \lambda_2\ge \ldots\,$ of $n$ will be identified if convenient. 
For ${\boldsymbol \lambda}\in{\boldsymbol Y}_{\!n}$ and ${\boldsymbol \mu}=(\mu_1,\mu_2,\ldots)\in{\boldsymbol Y}_{\!n+1}$, we denote by ${\boldsymbol \lambda} \nearrow {\boldsymbol \mu}$ if  
\begin{eqnarray}
\label{2016-02-04-11}
\mu_{i_0}= \lambda_{i_0}+1=\mu_{i_0}\;(\exists i_0),\;\;\lambda_i=\mu_i\;(\forall i\ne i_0),
\end{eqnarray}
where we supply $\lambda_i=0$ if necessary.
For ${\boldsymbol \nu}=(n_\zeta)_{\zeta\in\widehat{T}}\in P_n(\widehat{T})$ and ${\cal K}=\widehat{T}$ or $\widehat{T}^0$, put 
\begin{gather}
\nonumber
{\boldsymbol Y}_{\!n}({\boldsymbol \nu}):=\big\{\Lambda^n=({\boldsymbol \lambda}^\zeta)_{\zeta\in\widehat{T}}\;;\;{\boldsymbol \lambda}^\zeta\in{\boldsymbol Y},\;|{\boldsymbol \lambda}^\zeta|=n_\zeta\;(\forall\zeta\in\widehat{T})\big\},
\\
\label{2014-03-27-1}
{\boldsymbol Y}_{\!n}({\cal K}):=\bigsqcup_{{\boldsymbol \nu}\in P_n({\cal K})}
 {\boldsymbol Y}_{\!n}({\boldsymbol \nu}),\quad \;\;
 {\boldsymbol Y}({\cal K}):=\bigsqcup_{n\geqslant 0}{\boldsymbol Y}_{\!n}({\cal K}).
\end{gather}

Any $\Lambda^n=({\boldsymbol \lambda}^\zeta)_{\zeta\in\widehat{T}}\in{\boldsymbol Y}_{\!n}(\widehat{T})$\, determines naturally 
a standard partition ${\cal I}_n=(I_{n,\zeta})_{\zeta\in\widehat{T}}:{\boldsymbol I}_n=\bigsqcup_{\zeta\in\widehat{T}}I_{n,\zeta}$\,, and ${\boldsymbol \nu}=(n_\zeta)_{\zeta\in\widehat{T}}\in P_n(\widehat{T}),\, |I_{n,\zeta}|=n_\zeta$. This is expressed as\; \lq\lq\hspace{.3ex}$\Lambda^n\rightsquigarrow {\cal I}_n,\,{\boldsymbol \nu}$\hspace{.5ex}''\, in the following.    

For $\Lambda^n\in{\boldsymbol Y}_{\!n}({\cal K})$ and 
$M^{n+1}=({\boldsymbol \mu}^\zeta)_{\zeta\in{\cal K}}\in{\boldsymbol Y}_{\!n+1}({\cal K})$, we denote $\Lambda^n\nearrow M^{n+1}$ if \,${\boldsymbol \lambda}^{\zeta_0}\nearrow{\boldsymbol \mu}^{\zeta_0}\;(\exists \zeta_0\in{\cal K})$\, and\,${\boldsymbol \lambda}^\zeta={\boldsymbol \mu}^\zeta\;(\forall \zeta\ne \zeta_0)$. 
By this definition, 
${\boldsymbol Y}({\cal K})$ becomes a branching graph in Definition 2.1 in \cite{[HoHH2008]}.

\vskip.5em

{\it Step 4.}\; For $S([\rho])=\prod_{\zeta\in\widehat{T}}{\mathfrak S}_{I_{n,\zeta}}\cong{\mathfrak S}_{\boldsymbol \nu}$, we choose a CSR of its dual as follows. For each $\zeta\in\widehat{T}$, take an IR $\pi_{{\boldsymbol \lambda}^\zeta}$ of ${\mathfrak S}_{I_{n,\zeta}}\cong {\mathfrak S}_{n_\zeta}$ corresponding to Young diagram ${\boldsymbol \lambda}^\zeta$ of size $n_\zeta=|I_{n,\zeta}|$, and then take their outer tensor product as 
$$
\pi_{\Lambda^n}:=\boxtimes_{\zeta\in\widehat{T}}\,\pi_{{\boldsymbol \lambda}^\zeta}\quad{\rm with}\quad    
\Lambda^n=({\boldsymbol \lambda}^\zeta)_{\zeta\in\widehat{T}}\in{\boldsymbol Y}_{\!n}(\widehat{T}). 
$$
Put $\pi^1:=\pi_{\Lambda^n}$, then 
 the IR $\pi=\pi^0\boxdot \pi^1$ of the subgroup $H=U\rtimes S([\rho])$ is determined by \,$
\Lambda^n:=\big({\boldsymbol \lambda}^\zeta\big)_{\zeta\in\widehat{T}}\in{\boldsymbol Y}_{\!n}(\widehat{T})$\, since $\Lambda^n\rightsquigarrow {\boldsymbol \nu}\in P_n(\widehat{T})$. Put 
\begin{eqnarray}
\label{2016-02-04-21}
\breve{\Pi}_{\Lambda^n}:=\mathrm{Ind}^G_{H}(\pi^0\boxdot \pi^1).
\end{eqnarray}
\vskip.3em

{\bf Theorem 3.3.}\; 
{\it 
For generalized symmetric group $G_n=G(m,1,n),\, 4\le n<\infty,$ with $T={\boldsymbol Z}_m$,\, the set of IRs 
\begin{eqnarray}
\label{2016-05-06-1}
\mbox{\rm IR}\big(G(m,1,n)\big)
:=
\big\{\breve{\Pi}_{\Lambda^n}\,;\,\Lambda^n=\big({\boldsymbol \lambda}^\zeta\big)_{\zeta\in\widehat{T}}\in{\boldsymbol Y}_{\!n}(\widehat{T})\big\}
\end{eqnarray}
gives a complete set of representatives of the dual \,$\widehat{G_n}$ of $G_n$. 
}
\vskip1.2em

{\bf 3.2. Irreducible characters.} 
We give the character of $\breve{\Pi}_{\Lambda^n}$. 
For a $g=(d,\sigma)\in G_n=G(m,1,n)=D_n(T)\rtimes {\mathfrak S}_n$, we put
\begin{eqnarray}
\label{2015-05-03-11}
&&
\left\{
\begin{array}{l}
 {\rm supp}(d) := \{i \in {\boldsymbol I}_n\;;\; t_i \ne e_T\},\;\;  d=(t_i)_{i\in {\boldsymbol I}_n},\; t_i \in T_i=T,
 \\[.5ex] 
{\rm supp}(\sigma):= \{j \in {\boldsymbol I}_n\;;\; \sigma(j)\ne j\},\;\;
 {\rm supp}(g) := {\rm supp}(d) \bigcup {\rm supp}(\sigma), 
 \end{array}
 \right.
 \end{eqnarray}
where $e_T$ denotes the identity element of $T$.  
A standard decomposition of $g=(d,\sigma)\in G_n$ is 
\begin{eqnarray}
\label{2014-03-18-21-8}
&&
g=(d, \sigma) = \xi_{q_1}\xi_{q_2} \cdots \xi_{q_r}g_1g_2 \cdots g_s,\;\;\xi_q=\big(t_q,(q)\big), \;\;g_j=(d_j,\sigma_j), 
\end{eqnarray}
where for each $j \in J:={\boldsymbol I}_s$, $\sigma_j$ is a cycle of length $\ell_j:=\ell(\sigma_j)$, $\mathrm{supp}(d_j)\subset \mathrm{supp}(\sigma_j)$, and $K_j:=\mathrm{supp}(\sigma_j)\;(j\in J)$ are mutually disjoint. Then \,$g \in H_n=H$ if and only if \,$K_j\subset I_{n,\zeta}\;(\exists \zeta\in\widehat{T})$ for each $j\in J$.  
We put \,$Q:=\{q_1,q_2,\ldots,q_r\}$.  

Denote by $\chi(\pi_{\boldsymbol \lambda}|\,\cdot)$ the character of IR $\pi_{\boldsymbol \lambda},\,{\boldsymbol \lambda} \in {\boldsymbol Y}_{\!k}(\widehat{T})$, of the symmetric group ${\mathfrak S}_k$. A $\sigma\in{\mathfrak S}_k$ is decomposed into a product of mutually disjoint cycles $\sigma_i$, and the character value for $\sigma$ is determined by the set of lengths $(\ell_i),\,\ell_i=\ell(\sigma_i),\,\large\sum_i\ell_i\le k,$ which will be denoted as $\chi\big(\pi_{\boldsymbol \lambda}\big|(\ell_i)\big)$. Moreover we prepare notations for partitions, with parameter $\zeta\in\widehat{T}$, of $Q$ and $J$ respectively as  
\begin{eqnarray}
\label{2014-12-16-1-8}
{\cal Q}:=(Q_\zeta)_{\zeta\in\widehat{T}},\quad
{\cal J}:=(J_\zeta)_{\zeta\in\widehat{T}}.
\end{eqnarray} 
\vskip.2em 

{\bf Theorem 3.4.}\;([9, Theorem 4, p.312])\, 
{\it 
For IR $\breve{\Pi}_{\Lambda^n},\,\Lambda^n=({\boldsymbol \lambda}^\zeta)_{\zeta\in\widehat{T}}\in{\boldsymbol Y}_{\!n}(\widehat{T}),$ of  
$G_n=D_n(T)\rtimes{\mathfrak S}_n$ its character $\chi\big(\breve{\Pi}_{\Lambda^n}|\cdot\big)$ is given as follows. 
If $g=(d,\sigma)\in G_n$ is not conjugate to an element in $H_n:=D_n(T)\rtimes {\large\prod}_{\zeta\in\widehat{T}}\,{\mathfrak S}_{I_{n,\zeta}}$, then $\chi\big(\breve{\Pi}_{\Lambda^n}|g\big)=0$. 
  For $g=(d,\sigma)\in H_n$, let its standard decomposition be as in (\ref{2014-03-18-21-8}).   Then  
\begin{eqnarray}
\label{2014-03-17-11}
\qquad\;
\chi\big(\breve{\Pi}_{\Lambda^n}|g\big)
\!\!&=&\!\!
\sum_{{\cal Q},{\cal J}} 
b(\Lambda^n;{\cal Q},{\cal J};g) X(\Lambda^n;{\cal Q},{\cal J};g), 
\\
\nonumber
b(\Lambda^n;{\cal Q},{\cal J};g)
\!\!&=&\!\! 
\frac{(n-|\mathrm{supp}(g)|)!}{%
{\displaystyle 
{\prod}_{\zeta\in\widehat{T}}
\Big(n_\zeta - |Q_\zeta| - 
{\sum}_{j \in J_\zeta}\ell_j\Big)!
}
}\,.
\\[.3ex]
\label{2010-06-11-1-2}
\nonumber
X(\Lambda^n;{\cal Q},{\cal J};g)
\!\!&=&\!\!
\prod_{\zeta\in\widehat{T}} 
\Big(\prod_{q\in Q_\zeta}\zeta(t_q)\cdot\prod_{j\in J_\zeta}\zeta\big(P(d_j)\big)\cdot \chi\big(\pi_{{\boldsymbol \lambda}^\zeta}|(\ell_j)_{j\in J_\zeta}\big)\Big), 
\end{eqnarray}
where \,$\ell_j:=\ell(\sigma_j)$,\, $P(d_j)=\large\prod_{i\in K_j}\!t_i$ for $d_j=(t_i)_{i\in K_j}$, and \,${\cal Q}=(Q_\zeta)_{\zeta\in\widehat{T}},\,{\cal J}=(J_\zeta)_{\zeta\in\widehat{T}}$, are partitions of $Q$ and $J$ respectively and the sum runs over all pairs $({\cal Q},{\cal J})$ which satisfy the condition 
\begin{eqnarray}
\label{2014-03-17-11-3}
\nonumber
\mbox{\bf (QJ)}\qquad\qquad\qquad
|Q_\zeta|+{\sum}_{j\in J_\zeta}\ell_j\le n_\zeta=|{\boldsymbol \lambda}^\zeta|\;\;\;
(\zeta\in\widehat{T}). \qquad\qquad\qquad\qquad
\end{eqnarray}
}

\section{\large Spin IRs of Schur-Young subgroups of $\widetilde{{\mathfrak S}}_n$}

For ${\boldsymbol \nu}=(n_\zeta)_{\zeta\in{\cal K}}\in P_n({\cal K})$, ${\cal K}=\widehat{T}$ or $\widehat{T}^0$, we define {\it Schur-Young subgroup}\, $\widetilde{{\mathfrak S}}_{\boldsymbol \nu}:=\Phi^{\;\;-1}_{\mathfrak S}\big({\mathfrak S}_{\boldsymbol \nu})\subset \widetilde{{\mathfrak S}}_n$ under the identification ${\mathfrak S}_{\boldsymbol \nu}$ with $\large\prod_{\zeta\in {\cal K}} {\mathfrak S}_{I_{n,\zeta}}\subset{\mathfrak S}_n,\,|I_{n,\zeta}|=n_\zeta$ (cf. {\bf 3.1}). To describe spin IRs of this group, we introduce notations of Schur's type. 
\vskip1.2em

{\bf 4.1. Hauptdarstellung of Schur for ${\mathfrak S}_n$.}  
Schur defined Hauptdarstellung $\Delta_n$ of ${\mathfrak T}_n$ in \cite{[Sch1911]}, and we rewrite it to $\Delta'_n$ of $\widetilde{{\mathfrak S}}_n={\mathfrak T}'_n$. 
Put $2\times 2$ matrices 
$\varepsilon,a,b,c$ as 
\begin{eqnarray}
\label{2014-04-13-1-8}
&&\quad\;
\varepsilon = 
{\small 
\begin{pmatrix}
1&0 \\
0&1
\end{pmatrix}
}, \;\;
a= 
{\small 
\begin{pmatrix}
0 & 1 \\
1 & 0
\end{pmatrix}
}, \;\;
b=
{\small 
\begin{pmatrix}
0 & -i \\
i & 0
\end{pmatrix}
}, \;\;
c=
{\small 
\begin{pmatrix}
1 & 0 \\
0 & -1
\end{pmatrix}
}
. 
\end{eqnarray}
For $N:=[(n-1)/2]$, we define matrices 
 $X_j\;(j\in{\boldsymbol I}_{2N+1})$ of degree $2^{N}$\,: for  
$i\in{\boldsymbol I}_{N}$,
\begin{eqnarray}
\label{2017-06-19-1}
&&
\left\{
\begin{array}{l}
X_{2i-1} =  c^{\otimes(i-1)}\otimes a\otimes \varepsilon^{\otimes(N-i)}, \, 
\\[.5ex]
X_{2i} \;\;\;=\;  c^{\otimes(i-1)}\otimes b\otimes \varepsilon^{\otimes(N-i)}, \;\;
\\[.5ex] 
X_{2N+1} = c^{\otimes N}\,,
\end{array}
\right.
\end{eqnarray}
and define $\Delta'_n(r_j)$ with coefficients $a_0=0,\,a_j,\,b_j$ as, for $j\in{\boldsymbol I}_{n-1}$,  
\begin{gather}
\label{2017-06-19-2}
\Delta'_n(r_j):=a_{j-1}X_{j-1}+b_jX_j, 
\\
a_j=-\sqrt{\frac{j}{2(j\!+\!1)}}\,,\;\;
b_j=\sqrt{\frac{j\!+\!1}{2j}}\,.
\end{gather}
Then $\Delta'_n$ gives a spin IR of $\widetilde{{\mathfrak S}}_n$ (cf. [19, Theorem 8.3]). 

 Let $P_n$ be the set of all partitions ${\boldsymbol \nu}=(\nu_p)_{p\geqslant 1},\,\nu_1\ge \nu_s\ge \ldots \ge\nu_t\ge 1,\,\sum_{p\in{\boldsymbol I}_t}\nu_p=n,$ and $O\!P_n$ that of ${\boldsymbol \nu}=(\nu_p)_{p\geqslant 1}$ with all $\nu_p$ odd. Here $\ell({\boldsymbol \nu}):=t$\, is called the {\it length}\, of the partition ${\boldsymbol \nu}$. 
For an interval $K=[a,b]\subset{\boldsymbol I}_n$, put $\sigma'_K:=r_ar_{a+1}\cdots r_{b-1}$, and for  ${\boldsymbol \nu}=(\nu_p)_{p\in{\boldsymbol I}_t}\in P_n$, put  
\begin{eqnarray}
\label{2015-05-22-21}
&&
\qquad
\sigma'_{\boldsymbol \nu}=\sigma'_1\sigma'_2\cdots\sigma'_t,\;\;\sigma'_j= \sigma'_{K_j},\;K_j=[\nu_1\!+\!\cdots\!+\!\nu_{j-1},\nu_1\!+\!\cdots\!+\!\nu_j],  
\end{eqnarray} 
and call this a {\it standard representative} of conjugacy classes of $\widetilde{{\mathfrak S}}_n$ modulo $Z_1$. 
Thus we have a complete system of representatives $\{\sigma'_{\boldsymbol \nu}\,;\,{\boldsymbol \nu}\in P_n\}$.   
As symbols, we put $\sigma_{\boldsymbol \nu}=\Phi_\mathfrak{S}(\sigma'_{\boldsymbol \nu})=\sigma_1\sigma_2\cdots\sigma_t,\,\sigma_i=\sigma_{K_i}=\Phi_\mathfrak{S}(\sigma'_i)\in\mathfrak{S}_{K_i}$. Then, 
 $\mathrm{sgn}(\sigma'_{\boldsymbol \nu}):=\mathrm{sgn}(\sigma_{\boldsymbol \nu})=(-1)^{d({\boldsymbol \nu})}=(-1)^{n-t}$\,
 with 
$d({\boldsymbol \nu}):=\sum_{p\in{\boldsymbol I}_t}(\nu_p-1)=n-l({\boldsymbol \nu})$.
\bigskip

{\bf Theorem 4.1.} ([35, Abs.VI], 
[19, Theorem 8.4])\; 

{\it 
Let $n\ge 4$. Then the character of Schur's Hauptdarstellung $\Delta'_n$ is given as follows. 

{\rm (i)}\;
Suppose $\sigma'_{\boldsymbol \nu}$ is even. Then\, $\chi_{\Delta'_n}(\sigma'_{\boldsymbol \nu})\ne 0$\,  if and only if ${\boldsymbol \nu}\in O\!P_n$, and in that case 
\begin{eqnarray}
\label{2013-04-29-1}
\chi_{\Delta'_n}(\sigma'_{\boldsymbol \nu})=(-1)^{(n-t)/2}\,2^{[(t-1)/2]}
=
(-1)^{d({\boldsymbol \nu})/2}\,2^{[(\ell({\boldsymbol \nu})-1)/2]}. 
\end{eqnarray}

{\rm (ii)}\; Suppose $\sigma'_{\boldsymbol \nu}$ is odd, and $n$ is even. Then,  
 $\chi_{\Delta'_n}(\sigma'_{\boldsymbol \nu})\ne 0$\, if and only if \,$\ell({\boldsymbol \nu})=1$, or equivalently $\sigma_{\boldsymbol \nu}$ is a cycle of length $n$. For ${\boldsymbol \nu}=(n)$,\, $\sigma'_{(n)}=r_1r_2\cdots r_{n-1},\,\sigma_{(n)}=(1\;2\;3\;\ldots\;n),$ and 
\begin{eqnarray}
\label{2013-04-29-2}
\chi_{\Delta'_n}(\sigma'_{(n)})=i^{n/2-1}\hspace{.2ex}\sqrt{n/2}\;\;\;(i=\sqrt{-1})\,.
\end{eqnarray}

{\rm (iii)}\; Suppose $\sigma'_{\boldsymbol \nu}$ is odd, and $n$ is odd. Then, $\chi_{\Delta'_n}(\sigma'_{\boldsymbol \nu})= 0$. 

{\rm (vi)}\; $\Delta'_n$ is self-associate or not (i.e., $\mathrm{sgn}_{\mathfrak S}\cdot \Delta'_n\cong \Delta'_n$\, or\;  
$\mathrm{sgn}_{\mathfrak S}\cdot \Delta'_n\not\cong \Delta'_n$), according as $n$ is even or odd. 

}

\vskip1em

{\bf 4.2. Spin IRs of ${\mathfrak S}_n$.} 
  Denote by $S\!P_n$ the set of all strict partitions ${\boldsymbol \lambda} =(\lambda_j)_{1\leqslant j\leqslant l}$ of $n$\,:  $\lambda_1>\lambda_2>\ldots>\lambda_l>0,\;\large\sum_{1\leqslant j\leqslant l}\lambda_j=n$. To each ${\boldsymbol \lambda}\in S\!P_n$, there corresponds bijectively a shifted Young diagram $F_{\boldsymbol \lambda}$ of size $n$ as follows: in  
 $(x,y)$-plane, take $y$-axis downwards and represents by $(i,j)\in{\boldsymbol N}^2$ a unit box with center $(i,j)$, then 
\begin{eqnarray}
\label{2013-06-07-1-2}
F_{\boldsymbol \lambda}:=\{(i,j)\,;\;j\le i\le \lambda_j+(j-1)\}.
\end{eqnarray}
$F_{\boldsymbol \lambda}$ will be identified with ${\boldsymbol \lambda}$ if there is no fear of confusion. Let ${\boldsymbol Y}^{\rm sh}_{\!n}$ be the set of all shifted Young diagrams of size $n$ (${\boldsymbol Y}^{\rm sh}_{\!n}=\varnothing$ for $n=0$), and put 
\begin{eqnarray}
\label{2017-06-12-1}
{\boldsymbol Y}^{\rm sh}={\bigsqcup}_{n\geqslant 0}{\boldsymbol Y}^{\rm sh}_{\!n}.
\end{eqnarray}
 
For ${\boldsymbol \lambda}=(\lambda_j)_{1\leqslant j\leqslant l} \in {\boldsymbol Y}^{\rm sh}_{\!n},$ put
\begin{eqnarray}
\label{2016-04-10-23-8}
\left\{
\begin{array}{l}
l({\boldsymbol \lambda}):=l\;(\mbox{\rm the length of ${\boldsymbol \lambda}$}),\quad  d({\boldsymbol \lambda}):=n-l({\boldsymbol \lambda}), 
\\
\varepsilon({\boldsymbol \lambda}):=0,1,\;\;\mbox{\rm according as $d({\boldsymbol \lambda})$ is even or odd.}
\end{array}
\right.
\end{eqnarray}

The Schur-Young subgroup $\widetilde{{\mathfrak S}}_{\boldsymbol \lambda}$ is canonically isomorphic to the twisted central product [19, Definition 1.1] 
of $\Phi^{\;-1}_{\mathfrak S}\big({\mathfrak S}_{I_{n,\zeta}}\big)\cong\widetilde{{\mathfrak S}}_{\lambda_j}$ as 
\begin{eqnarray}
\label{2016-05-09-1}
\widetilde{{\mathfrak S}}_{\boldsymbol \lambda}\cong \widetilde{{\mathfrak S}}_{\lambda_1}\,\widehat{*}\,\widetilde{{\mathfrak S}}_{\lambda_2}\,\widehat{*}\,\cdots\,\widehat{*}\,\widetilde{{\mathfrak S}}_{\lambda_l},
\end{eqnarray}
where $\widetilde{{\mathfrak S}}_n:={\mathfrak S}_n\times Z_1\;(n=1,2,3)$. Take for each $j$, Schur's Hauptdarstellung $\Delta'_{\lambda_j}$ of $\widetilde{{\mathfrak S}}_{\lambda_j}$ and consider their twisted central product $\Delta'_{\boldsymbol \lambda}:=\Delta'_{\lambda_1}\,\widehat{*}\,\Delta'_{\lambda_2}\,\widehat{*}\,\cdots\,\widehat{*}\,\Delta'_{\lambda_l}$ (cf. [ibid., \S 8.3]). Denote by $\tau_{\boldsymbol \lambda}$ the unique top irreducible component of the standard induced representation 
$\Pi_{\boldsymbol \lambda}:=\mathrm{Ind}^{\widetilde{{\mathfrak S}}_n}_{\widetilde{{\mathfrak S}}_{\boldsymbol \lambda}}\Delta'_{\boldsymbol \lambda}$. Then we know the following due to Schur\;\cite{[Sch1911]} (cf. [19., \S 8]):  
\vskip1em

{\bf (4-i)}\;\; 
{\it $\tau_{\boldsymbol \lambda}$ is self-associate, or\, $\tau_{\boldsymbol \lambda}\cong \mathrm{sgn}_{\mathfrak S}\cdot\tau_{\boldsymbol \lambda}$, if and only if $\varepsilon({\boldsymbol \lambda})=0$. 
}

{\bf (4-ii)}\; 
{\it The complete set of representatives of spin IRs of $\widetilde{{\mathfrak S}}_n$ is given by 
\begin{eqnarray}
\label{2016-05-09-11}
&&
\big\{\tau_{\boldsymbol \lambda}\,;\,{\boldsymbol \lambda}\in {\boldsymbol Y}^{\rm sh}_{\!n},\;\varepsilon({\boldsymbol \lambda})=0\big\}\bigsqcup\big\{\tau_{\boldsymbol \lambda},\;\mathrm{sgn}_{\mathfrak S}\cdot\tau_{\boldsymbol \lambda}\,;\,{\boldsymbol \lambda}\in {\boldsymbol Y}^{\rm sh}_{\!n},\;\varepsilon({\boldsymbol \lambda})=1\big\}.
\end{eqnarray}
}

\vskip.2em

{\bf Lemma 4.2.} (cf. [19, \S 10.3])\;  
{\it 
For the support of character $\chi_{\tau_{\boldsymbol \lambda}}$ of spin IR $\tau_{\boldsymbol \lambda},\,{\boldsymbol \lambda}=(\lambda_j)_{1\leqslant j\leqslant l}\in {\boldsymbol Y}^{\rm sh}_{\!n}$, of $\widetilde{{\mathfrak S}}_n, 4\le n <\infty,$ there holds the following.   
\vskip.2em

{\rm (i)}\; Denote by $\widetilde{\Sigma}_k$ the set of  elements $\sigma'\!\in\widetilde{{\mathfrak S}}_k$ for which the cycle decompositions of $\sigma=\Phi_{\mathfrak S}(\sigma')$ consist of only even cycles,  and put 
\begin{eqnarray}
\label{2014-12-19-11-8}
\nonumber
\widetilde{\Sigma}_{\boldsymbol \lambda}:=\widetilde{\Sigma}_{\lambda_1}\widetilde{\Sigma}_{\lambda_2}\cdots \widetilde{\Sigma}_{\lambda_l}\subset \widetilde{{\mathfrak S}}_{\lambda_1}\widehat{*}\,\widetilde{{\mathfrak S}}_{\lambda_2}\widehat{*}\,\cdots\,\widehat{*}\,\widetilde{{\mathfrak S}}_{\lambda_l}=\underset{1\leqslant j\leqslant l}{\widehat{*}}\widetilde{{\mathfrak S}}_{\lambda_j}\hookrightarrow\widetilde{{\mathfrak S}}_n\,,
\end{eqnarray}
\\[-5ex]
then
\\[-4ex]
\begin{eqnarray}
\label{2016-05-10-1}
\mathrm{supp}(\chi_{\tau_{\boldsymbol \lambda}})\bigcap\widetilde{{\mathfrak A}}_n\,\subset\,\widetilde{\Sigma}_{\boldsymbol \lambda}^{\;\widetilde{{\mathfrak S}}_n}:={\bigcup}_{\sigma'_0\in\widetilde{{\mathfrak S}}_n}\sigma'_0\,\widetilde{\Sigma}_{\boldsymbol \lambda}{\sigma'_0}^{\,-1}.
\end{eqnarray}
Moreover, for $\sigma'\in\widetilde{{\mathfrak A}}_n$, let $\sigma=\Phi_{\mathfrak S}(\sigma')=\sigma_1\sigma_2\cdots\sigma_p$ be a decompositions into disjoint cycles, and $\sigma'=\sigma'_1\sigma'_2\cdots\sigma'_p$ with $\Phi_{\mathfrak S}(\sigma'_j)=\sigma_j\;(\forall j)$, then for any permutation $\kappa$ of $\{1,2,\ldots,p\}$, 
 $
\chi_{\tau_{\boldsymbol \lambda}}(\sigma')=\chi\big(\sigma'_{\kappa(1)}\sigma'_{\kappa(2)}\cdots\sigma'_{\kappa(p)}\big).
$

{\rm (ii)}\;\; If $\varepsilon({\boldsymbol \lambda})=0$, then $\tau_{\boldsymbol \lambda}$ is self-associate, and \,$\chi_{\tau_{\boldsymbol \lambda}}=0$\, on $\widetilde{{\mathfrak C}}_n=\widetilde{{\mathfrak S}}_n\setminus \widetilde{{\mathfrak A}}_n$. 

{\rm (iii)}\; If $\varepsilon({\boldsymbol \lambda})=1$, then $\tau_{\boldsymbol \lambda}$ is non-self-associate, and 
\begin{eqnarray}
\label{2016-05-10-11}
&&
 \qquad 
 \mathrm{supp}(\chi_{\tau_{\boldsymbol \lambda}})\bigcap\widetilde{{\mathfrak C}}_n=\widetilde{\Xi}_{\boldsymbol \lambda}^{\;\widetilde{{\mathfrak S}}_n}:={\bigcup}_{\sigma'_0\in\widetilde{{\mathfrak S}}_n}\,\sigma'_0\,\widetilde{\Xi}_{\boldsymbol \lambda}{\sigma'_0}^{\,-1},\qquad {\rm with}
 \\
 \nonumber
 &&
\left\{
\begin{array}{l}
\widetilde{\Xi}_k:=\{\sigma'\in\widetilde{{\mathfrak S}}_k\;;\;\sigma=\Phi_{\mathfrak S}(\sigma')\;\,\mbox{is a cycle of length $k$}\}\subset\widetilde{{\mathfrak S}}_k,
\\[1.2ex]
\widetilde{\Xi}_{\boldsymbol \lambda}:=\widetilde{\Xi}_{\lambda_1}\widetilde{\Xi}_{\lambda_2}\cdots \widetilde{\Xi}_{\lambda_l}\subset \widetilde{{\mathfrak S}}_{\lambda_1}\widehat{*}\,\widetilde{{\mathfrak S}}_{\lambda_2}\widehat{*}\,\cdots\,\widehat{*}\,\widetilde{{\mathfrak S}}_{\lambda_l}=\underset{1\leqslant j\leqslant l}{\widehat{*}}\widetilde{{\mathfrak S}}_{\lambda_j}\hookrightarrow\widetilde{{\mathfrak S}}_n.
\end{array}
\right.
\end{eqnarray}
}

\vskip.2em 

{\bf 4.3. Spin IRs of Schur-Young subgroup $\widetilde{{\mathfrak S}}_{\boldsymbol \nu}$\,.}\; 
\vskip.3em

For the later use, we put for ${\boldsymbol \nu}=(n_\zeta)_{\zeta\in\widehat{T}}\in P_n(\widehat{T})$ and ${\cal K}=\widehat{T}$ or $\widehat{T}^0$,   
\begin{gather}
\nonumber
{\boldsymbol Y}^{\rm sh}_{\!n}({\boldsymbol \nu}):=\big\{\Lambda^n=({\boldsymbol \lambda}^\zeta)_{\zeta\in\widehat{T}}\;;\;{\boldsymbol \lambda}^\zeta\in{\boldsymbol Y}^{\rm sh}_{n_\zeta}\big\},
\\
\label{2016-05-09-21}
{\boldsymbol Y}^{\rm sh}_{\!n}({\cal K}):={\bigsqcup}_{{\boldsymbol \nu}\in P_n({\cal K})}{\boldsymbol Y}^{\rm sh}_{\!n}({\boldsymbol \nu}).
\end{gather}
and for $\Lambda^n=({\boldsymbol \lambda}^\zeta)_{\zeta\in{\cal K}}\in {\boldsymbol Y}^{\rm sh}_{\!n}({\cal K})$, put
\begin{gather}
\label{2014-11-21-21-8}
\tau_{\Lambda^n}:=\underset{\zeta\in{\cal K}}{\widehat{*}}\,\tau_{{\boldsymbol \lambda}^\zeta}, 
\qquad
\\[.5ex]
\label{2016-04-11-31-8}
\left\{
\begin{array}{l}
s(\Lambda^n):=\sharp\{\zeta\in{\cal K}\;;\;\varepsilon({\boldsymbol \lambda}^\zeta)=1\},
\\[1ex]
d(\Lambda^n):=n-{\sum}_{\zeta\in{\cal K}}\,l({\boldsymbol \lambda}^\zeta),
\end{array}
\right.
\end{gather}
where the twisted central product \,$\widehat{*}_{\zeta\in{\cal K}}\,\tau_{{\boldsymbol \lambda}^\zeta}$ of $\tau_{{\boldsymbol \lambda}^\zeta}$'s is taken according to the order in ${\cal K}\subset\widehat{T}$. 
Then \,$s(\Lambda^n)\equiv d(\Lambda^n)\;(\mathrm{mod}\;2)$, and 
$\tau_{\Lambda^n}$ is self-associate or not according as $s(\Lambda^n)$ is even or not. 
\vskip1em

{\bf Notation 4.3.}  To describe the spin dual of\, $\widetilde{{\mathfrak S}}_{\boldsymbol \nu}=\Phi^{\;\;-1}_{\mathfrak S}\big({\mathfrak S}_{\boldsymbol \nu}\big)$ etc., we introduce 
\begin{gather}
\label{2014-05-05-1-8}
\left\{
\begin{array}{ll}
{\boldsymbol Y}^{\rm sh}_{\!n}({\boldsymbol \nu})^{\rm ev}\,:=\,
\big\{\Lambda^n=({\boldsymbol \lambda}^\zeta)_{\zeta\in{\cal K}}\in{\boldsymbol Y}^{\rm sh}_{\!n}({\boldsymbol \nu})\;;\;s(\Lambda^n)\equiv d(\Lambda^n)\equiv 0\big\},  
\\[1ex]
{\boldsymbol Y}^{\rm sh}_{\!n}({\boldsymbol \nu})^{\rm odd}:=
\big\{\Lambda^n=({\boldsymbol \lambda}^\zeta)_{\zeta\in{\cal K}}\in{\boldsymbol Y}^{\rm sh}_{\!n}({\boldsymbol \nu})\;;\;s(\Lambda^n)\equiv d(\Lambda^n)\equiv 1\big\},  
 \end{array}
 \right.
 \\
 \label{2014-12-16-11-8}
{\mathscr Y}_n({\boldsymbol \nu})={\mathscr Y}_n({\boldsymbol \nu})^{\rm ev}\bigsqcup{\mathscr Y}_n({\boldsymbol \nu})^{\rm odd},\quad{\rm with}
 \\
 \nonumber
{\mathscr Y}_n({\boldsymbol \nu})^{\rm ev}:={\boldsymbol Y}^{\rm sh}_{\!n}({\boldsymbol \nu})^{\rm ev},\;\,
{\mathscr Y}_n({\boldsymbol \nu})^{\rm odd}:={\boldsymbol Y}^{\rm sh}_{\!n}({\boldsymbol \nu})^{\rm odd}\times\{\pm 1\}. 
 \\
 \label{2016-05-09-22}
{\mathscr Y}_n({\cal K})={\mathscr Y}_n({\cal K})^{\rm ev}\bigsqcup{\mathscr Y}_n({\cal K})^{\rm odd}, \quad{\rm with}
 \\
 \nonumber
{\mathscr Y}_n({\cal K})^{\rm ev}:=\bigsqcup_{{\boldsymbol \nu}\in P_n({\cal K})}{\mathscr Y}_n({\boldsymbol \nu})^{\rm ev},\;\,
{\mathscr Y}_n({\cal K})^{\rm odd}:=\bigsqcup_{{\boldsymbol \nu}\in P_n({\cal K})}{\mathscr Y}_n({\boldsymbol \nu})^{\rm odd}.
\end{gather}

\vskip.2em

{\bf Theorem 4.4.}\;[19, Theorem 9.1]\; 
{\it 
For ${\boldsymbol \nu}=(n_\zeta)_{\zeta\in{\cal K}}\in P_n({\cal K})$, the twisted central product \,${\widehat{*}}_{\,\zeta\in{\cal K}}\,\widetilde{{\mathfrak S}}_{n_\zeta}$ is naturally isomorphic to $\widetilde{{\mathfrak S}}_{\boldsymbol \nu}=\Phi^{\;\;-1}_{\mathfrak S}\big({\mathfrak S}_{\boldsymbol \nu}\big)$, and a complete set of representatives for spin dual of $\widetilde{{\mathfrak S}}_{\boldsymbol \nu}$ is given by 
\begin{gather}
\label{2016-05-09-31}
\big\{\tau_{M^n}\,;\,M^n\in {\mathscr Y}_n({\boldsymbol \nu})\big\},\quad {\rm with} 
\\
\label{2016-05-09-3}
\tau_{M^n}:=\left\{
\begin{array}{ll}
\quad\;
\tau_{\Lambda^n},\qquad & M^n=\Lambda^n\in {\mathscr Y}_n({\boldsymbol \nu})^{\rm ev}, 
\\[.5ex]
\mathrm{sgn}_{\mathfrak S}^{(1-\mu)/2}\cdot\tau_{\Lambda^n},\quad 
& M^n=(\Lambda^n,\mu)\in {\mathscr Y}_n({\boldsymbol \nu})^{\rm odd}, 
\,\mu=\pm 1.
\end{array}
\right.
\end{gather}
Thus ${\mathscr Y}_n({\boldsymbol \nu})$ gives a parameter space for the spin dual of Schur-Young subgroup $\widetilde{{\mathfrak S}}_{\boldsymbol \nu}$ of $\widetilde{{\mathfrak S}}_n$.

}

\vskip1em

{\bf Lemma 4.5.} (cf.\;Lemma 4.2 and [19, \S 10.4]) 
{\it 
For $\Lambda^n=({\boldsymbol \lambda}^\zeta)_{\zeta\in{\cal K}}\in {\boldsymbol Y}_{\!n}({\cal K})$ with ${\cal K}=\widehat{T}$, 
the character of spin IR $\tau_{\Lambda^n}=\widehat{*}_{\zeta\in{\cal K}}\hspace{.3ex}\tau_{{\boldsymbol \lambda}^\zeta}$ of\, $\widetilde{{\mathfrak S}}_{\boldsymbol \nu}\cong\widehat{*}_{\zeta\in{\cal K}}\hspace{.3ex}\widetilde{{\mathfrak S}}_{n_\zeta},\, {\boldsymbol \nu}=(n_\zeta)_{\zeta\in{\cal K}},\,n_\zeta=|{\boldsymbol \lambda}^\zeta|,$ has the following property. For an element $\sigma'\in \widetilde{{\mathfrak S}}_{\boldsymbol \nu}$, let 
$\sigma'=\sigma'_1\sigma'_2\cdots \sigma'_s$ be a decomposition such that $\sigma_j=\Phi_{\mathfrak S}(\sigma'_j)$ are disjoint cycles.
 
In case $\sigma'$ is even, let\, ${\textstyle \large\prod}_{j\in J_\zeta}\sigma_j\in{\mathfrak S}_{I_{n,\zeta}}\cong\widetilde{{\mathfrak S}}_{n_\zeta},\;\bigsqcup_{\zeta\in{\cal K}}J_\zeta={\boldsymbol I}_s$, then  
\begin{eqnarray*}
\chi(\tau_{\Lambda^n}|\,\sigma')=
{\prod}_{\zeta\in{\cal K}}\,\chi\big(\tau_{{\boldsymbol \lambda}^\zeta}\,\big|\,{\textstyle \large\prod}_{j\in J_\zeta}\sigma'_j\big), 
\end{eqnarray*}
where the order in each product ${\prod}_{j\in J_\zeta}\sigma'_j$ does not change the character value. 

In case $\sigma'$ is odd, the character value $\chi(\tau_{\Lambda^n}|\,\sigma')$ is zero if\, $s(\Lambda^n)=0$ or\, $|\mathrm{supp}(\sigma')|<n$. Otherwise it can be calculated by means of the character formula in {\rm [ibid., Proposition 4.3]} with reference to {\rm [ibid., Theorem 10.1]}.
}

\vskip1em

{\bf 4.4. Spin IRs of the subgroup $\widetilde{{\mathfrak A}}_n\bigcap\widetilde{{\mathfrak S}}_{\boldsymbol \nu}$ of Schur-Young subgroup $\widetilde{{\mathfrak S}}_{\boldsymbol \nu}$\,.}

To get a CSR of spin IRs of the normal subgroup $\widetilde{{\mathfrak A}}_{\boldsymbol \nu}:=\widetilde{{\mathfrak A}}_n\bigcap\widetilde{{\mathfrak S}}_{\boldsymbol \nu}$ of $\widetilde{{\mathfrak S}}_{\boldsymbol \nu}$, generally of index 2, we restrict spin IRs from $\widetilde{{\mathfrak S}}_{\boldsymbol \nu}$ to $\widetilde{{\mathfrak A}}_{\boldsymbol \nu}$ and pick up its irreducible components. 
Note that, for $\Lambda^n=({\boldsymbol \lambda}^\zeta)_{\zeta\in{\cal K}} \in {\boldsymbol Y}^{\rm sh}_{\!n}({\boldsymbol \nu})$, spin IR $\tau_{\Lambda^n}$ is self-associate or not according as $\Lambda^n\in {\boldsymbol Y}^{\rm sh}_{\!n}({\boldsymbol \nu})^{\rm ev}$ or $\Lambda^n\in {\boldsymbol Y}^{\rm sh}_{\!n}({\boldsymbol \nu})^{\rm odd},$ and in turn, this corresponds to whether the restriction $\tau_{\Lambda^n}|_{\widetilde{{\mathfrak A}}_{\boldsymbol \nu}}$ splits or not. 
We call 
$\Lambda^n=({\boldsymbol \lambda}^\zeta)_{\zeta\in\widehat{T}}$ {\it quasi-trivial} \,if $n_\zeta=|{\boldsymbol \lambda}^\zeta|\le 1\;(\forall \zeta)$. In that case, $\prod_{\zeta\in\widehat{T}}{\mathfrak S}_{I_{n,\zeta}}\cong{\mathfrak S}_{\boldsymbol \nu}$ is trivial, and we understand $\widetilde{{\mathfrak S}}_{\boldsymbol \nu}=\widetilde{{\mathfrak A}}_{\boldsymbol \nu}\cong Z_1$ and 
$\tau_{\Lambda^n}=\mathrm{sgn}_{Z_1}$. 
\vskip1em

{\bf Notation 4.6.}  To describe the spin dual of \,$\widetilde{{\mathfrak A}}_n\cap\widetilde{{\mathfrak S}}_{\boldsymbol \nu}$ etc., we introduce 
\begin{eqnarray}
\nonumber
&&
{\boldsymbol Y}^{{\mathfrak A},{\rm sh}}_n({\boldsymbol \nu})^{\rm ev,non}
:=\big\{ (\Lambda^n,\kappa)\;;\;\Lambda^n\in 
{\boldsymbol Y}^{\rm sh}_{\!n}({\boldsymbol \nu})^{\rm ev}\;\mbox{\rm non quasi-trivial},\,\kappa=0, 1\big\},
\\[.5ex]
\nonumber
&&
{\boldsymbol Y}^{{\mathfrak A},{\rm sh}}_n({\boldsymbol \nu})^{\rm ev,tri}
:=
\big\{\Lambda^n\;;\;\Lambda^n\in {\boldsymbol Y}^{\rm sh}_{\!n}({\boldsymbol \nu})^{\rm ev}\;\mbox{\rm quasi-trivial}\big\}, 
\\[.3ex]
\label{2015-02-11-22-2}
&&
\quad
{\mathscr Y}^{\mathfrak A}_n({\boldsymbol \nu}):={\boldsymbol Y}^{{\mathfrak A},{\rm sh}}_{\!n}({\boldsymbol \nu})^{\rm ev,non}\bigsqcup{\boldsymbol Y}^{{\mathfrak A},{\rm sh}}_{\!n}({\boldsymbol \nu})^{\rm ev,tri}\bigsqcup{\boldsymbol Y}^{\rm sh}_{\!n}({\boldsymbol \nu})^{\rm odd},
\\[.3ex]
&&
\label{2017-07-04-1}
\nonumber
\quad
{\displaystyle 
{\mathscr Y}^{\mathfrak A}_n({\cal K}):={\bigsqcup}_{{\boldsymbol \nu}\in P_n({\cal K})}\,{\mathscr Y}^{\mathfrak A}_n({\boldsymbol \nu})\quad{\rm for}\;\;{\cal K}=\widehat{T}\;\,{\rm or}\;\,\widehat{T}^0.
}
\end{eqnarray}

Note that ${\boldsymbol Y}^{{\mathfrak A},{\rm sh}}_n({\boldsymbol \nu})^{\rm ev,tri}=\emptyset$ if $n>m$ for ${\boldsymbol \nu}\in P_n(\widehat{T})$ (resp. if $n>m^0=m/2$ for ${\boldsymbol \nu}\in P_n(\widehat{T}^0)$).
By using Theorem 4.4, we obtain the following. 
\vskip1em

{\bf Theorem 4.7.}\;
{\it 
Put ${\cal K}=\widehat{T}$ or $\widehat{T}^0$. 
Let 
${\boldsymbol \nu}=(n_\zeta)_{\zeta\in {\cal K}}\in P_n({\cal K})$. 
A complete set of representatives of spin dual of   
\,$\widetilde{{\mathfrak A}}_n\cap\widetilde{{\mathfrak S}}_{\boldsymbol \nu}$ is given by 
\begin{eqnarray}
\label{2015-02-13-1-4} 
&& 
\left\{
\begin{array}{ll}
\widetilde{\rho}_{M^n}:= \widetilde{\rho}^{\,(\kappa)}_{\Lambda^n}
\quad 
&
\big(M^n=(\Lambda^n,\kappa)\in{\boldsymbol Y}^{{\mathfrak A},{\rm sh}}_n({\boldsymbol \nu})^{\rm ev,non} \big), 
\\[.5ex]
\widetilde{\rho}_{M^n}:=\mathrm{sgn}_{Z_1} 
\quad 
&
\big(M^n=\Lambda^n\in{\boldsymbol Y}^{{\mathfrak A},{\rm sh}}_n({\boldsymbol \nu})^{\rm ev,tri}\big), 
\\[.5ex]
\widetilde{\rho}_{M^n}\!:=\!\tau_{\Lambda^n}\big|_{\widetilde{{\mathfrak A}}_n\cap\widetilde{{\mathfrak S}}_{\boldsymbol \nu}} \quad
&
\big(M^n\!=\!\Lambda^n\in{\boldsymbol Y}^{\rm sh}_{\!n}({\boldsymbol \nu})^{\rm odd} \big),
\end{array}
\right.
\end{eqnarray}
Here $\widetilde{\rho}_{M^n}=\widetilde{\rho}^{\,(\kappa)}_{\Lambda^n}\;(\kappa=0,1)$ are irreducible components of $\tau_{\Lambda^n}\big|_{\widetilde{{\mathfrak A}}_n\cap\widetilde{{\mathfrak S}}_{\boldsymbol \nu}}$. 
}
\hfill
$\Box$\;

\section{\large Spin IRs of $G(m,1,n)$ of spin type $z_1\to -1$}

For $2\le n <\infty$,  
consider a central extension $\widetilde{G}(m,1,n)$ of $G(m,1,n)=D_n(T)\rtimes{\mathfrak S}_n=U\rtimes S$ by the central subgroup $Z_1=\langle z_1\rangle,\,z_1^{\;2}=e,$ given as  
\begin{eqnarray}
\label{2016-12-15-1}
&&\;
\widetilde{G}(m,1,n):=D_n(T)\rtimes \widetilde{{\mathfrak S}}_n=U\rtimes S',\;\;\sigma'(d):=\sigma(d)\;(\sigma'\in S'),\,\end{eqnarray}
where $\sigma=\Phi_{\mathfrak S}(\sigma')\in S={\mathfrak S}_n$. 
Then $\widetilde{G}(m,1,n)=R\big(G(m,1,n)\big)$ if $m$ is odd, and $\widetilde{G}(m,1,n)\cong R\big(G(m,1,n)\big)/Z_{23},\,Z_{23}:=\langle z_2,z_3\rangle,$ if $m$ is even, which are also denoted by $\widetilde{G}^{\rm Y}(m,1,n)$ for Y\,=\,odd and IV respectively, where $\chi^{\rm IV}=(-1,1,1)$. A spin IR of $G(m,1,n)$ of spin type $z_1\to -1$ is nothing but an IR $\Pi$ of $\widetilde{G}(m,1,n)$ such that $\Pi(z_1)=-I$ ($I=$\,the identity operator).
Here we study this kind of spin IRs of the generalized symmetric group $G(m,1,n)$. 
\vskip.8em

{\bf 5.1. Spin IRs of $\widetilde{G}(m,1,n)$.}\; 
  To construct this kind of spin IRs, we follow the steps given in {\bf 2.1}. 
\vskip.5em

{\it Step 1.}\; For a CSR of\, $\widehat{U}/S$, take standard partitions ${\cal I}_n=\big(I_{n,\zeta}\big)_{\zeta\in\widehat{T}}$ of the interval ${\boldsymbol I}_n=[1,n]$ into subintervals $I_{n,\zeta}$ as in {\bf 3.1},   
 and characters $\rho=\rho_{{\cal I}_n}\in \widehat{U}$ in (\ref{2014-03-18-1-8}). 
\vskip.5em

{\it Step 2.}\; 
The stationary subgroup $S'([\rho])$ in $S'$ of $[\rho]=\{\rho\}$ is the full inverse image under $\Phi_{\mathfrak S}$ of the stationary subgroup $S([\rho])$ of $[\rho]$ in $S=\Phi_{\mathfrak S}(S')$\,:
\begin{eqnarray}
\label{2014-04-05-2-8}
&&
\quad 
S'([\rho])=\widetilde{{\mathfrak S}}_{\boldsymbol \nu}\cong \underset{\zeta\in\widehat{T}}{\widehat{*}}\widetilde{{\mathfrak S}}_{n_\zeta},\;\;{\boldsymbol \nu}=(n_\zeta)_{\zeta\in\widehat{T}},\;\;n_\zeta=|I_{n,\zeta}|.  
\end{eqnarray}

{\it Step 3.}\;   
For $\sigma'\in S'([\rho])$, $J_\rho(\sigma')={\bf 1}$ and \,$\pi^0:\,U\rtimes S'([\rho])\ni(d,\sigma')\mapsto \rho(d)J_\rho(\sigma')=\rho(d)$\, is a one-dimensional character. 

\vskip.5em

{\it Step 4.}\;   
By Theorem 4.4, we know that 
 
{\bf (*)} 
{\it for ${\boldsymbol \nu}=(n_\zeta)_{\zeta\in\widehat{T}}\in P_n(\widehat{T})$, a CSR of spin IRs of $S'([\rho])=\widetilde{{\mathfrak S}}_{\boldsymbol \nu}$ is given by 
\begin{eqnarray}
\label{2016-05-09-51}
\big\{\tau_{M^n}\,;\,M^n\in {\mathscr Y}_n({\boldsymbol \nu})\big\}.
\end{eqnarray}
}

Such is the classification of spin IRs of the spin type $z_1\to -1$. However, to give a good parametrization of such IRs, we reverse the steps. 
First start from $M^n\in{\mathscr Y}_n(\widehat{T})$ which contains   $\Lambda^n\in{\boldsymbol Y}^{\rm sh}_{\!n}(\widehat{T})$. Let ${\boldsymbol \nu}=(n_\zeta)_{\zeta\in\widehat{T}}$ and $\rho=\rho_{{\cal I}_n}\in\widehat{U}$ be associated with $\Lambda^n$ canonically. Then we have an IR $\pi=\pi^0\boxdot\pi^1$ of $H'_n=U\rtimes S'([\rho])$, with $\pi^0=\rho\cdot {\bf 1}$ of $H'_n$ and $\pi^1=\tau_{M^n}$ of $S'([\rho])$. Thus we obtain a spin IR of $\widetilde{G}_n=\widetilde{G}(m,1,n)$ as 
\begin{eqnarray}
\label{2015-02-22-12-8}
&&
\Pi^{\rm Y}_{M^n}:=\Pi(\pi^0,\pi^1)=\mathrm{Ind}^{\widetilde{G}_n}_{H'_n}(\rho\boxdot \tau_{M^n})\quad \mbox{\rm for\; Y\,=\,odd or IV}. 
\end{eqnarray}

\vskip.2em

{\bf Theorem 5.1.}\;  
{\it 
For the double covering group\, $\widetilde{G}_n=\widetilde{G}(m,1,n)$ of generalized symmetric group $G(m,1,n),\, 4\le n<\infty$, a complete set of representatives of equivalence classes of spin IRs of $G(m,1,n)$ of spin type $\chi^{\rm odd}$ (for $m$ odd) or $\chi^{\rm IV}=(-1,1,1)$ (for $m$ even) 
is given by the set, with\, {\rm Y\,=\,odd} or\, {\rm IV},  
\begin{eqnarray}
\label{2014-04-10-2}
&&
\mbox{\rm {\footnotesize spin}IR}^{\rm Y}\big(G(m,1,n)\big)
\;=\;
\big\{\Pi^{\rm Y}_{M^n}\;;\;M^n\in{\mathscr Y}_n(\widehat{T})\big\}. 
\end{eqnarray}
So ${\mathscr Y}_n(\widehat{T})$ gives a parameter space for the spin dual with this spin type. 
More in detail, 
put \,$\Pi^{\rm Y}_{\Lambda^n}:=\mathrm{Ind}^{\widetilde{G}_n}_{H'_n}(\rho\boxdot \tau_{\Lambda^n})$, then the above set (\ref{2014-04-10-2}) is equal to 
\begin{eqnarray*}
\big\{\Pi^{\rm Y}_{\Lambda^n}\;;\;\Lambda^n\in{\boldsymbol Y}^{\rm sh}_{\!n}(\widehat{T})^{\rm ev}\big\}\bigsqcup  
\big\{\Pi^{\rm Y}_{\Lambda^n},\;\,\mathrm{sgn}\cdot\Pi^{\rm Y}_{\Lambda^n}\;;\,\Lambda^n\in{\boldsymbol Y}^{\rm sh}_{\!n}(\widehat{T})^{\rm odd}\big\}.
\end{eqnarray*}
}

{\bf 5.2. Spin irreducible characters of $\widetilde{G}(m,1,n)$.} 
We denote the character of a representation $\tau$ by $\chi_\tau$ or $\chi(\tau|\cdot)$.  
Put $G_n=G(m,1,n),\,\widetilde{G}_n=\widetilde{G}(m,1,n).$ 
The calculation of the character of induced representation 
 $\Pi^{\rm Y}_{\Lambda^n}=\mathrm{Ind}^{\widetilde{G}_n}_{H'_n} (\rho\boxdot\tau_{\Lambda^n})$ proceeds as usual. 
 A $g'=(d,\sigma')\in \widetilde{G}_n$ has a standard decomposition as  
\begin{eqnarray}
\label{2014-04-07-1-8}
&&\;
g'=(d, \sigma') = \xi_{q_1}\xi_{q_2} \cdots \xi_{q_r}g'_1g'_2 \cdots g'_s,\;\;\xi_q=\big(t_q,(q)\big), \;\;g'_j=(d_j,\sigma'_j), 
\end{eqnarray}
where $d\in D_n(T),\,\sigma',\sigma'_j\in\widetilde{{\mathfrak S}}_n$, $\mathrm{supp}(d_j)\subset K_j:=\mathrm{supp}(\sigma_j),\,\sigma_j=\Phi_{\mathfrak S}(\sigma'_j)$ a cycle. Then, such a decomposition of $g=\Phi(g')\in G_n$ is given by $\xi_q$ and $g_j:=(d_j,\sigma_j)$. 
Put 
\begin{eqnarray}
\label{2016-05-09-63}
Q:=\{q_1,q_2,\ldots,q_r\},\;J:={\boldsymbol I}_s,\;\ell_j:=\ell(\sigma'_j):=\ell(\sigma_j)\;(j\in J).
\end{eqnarray}
Decompositions of $Q$ and $J$ indexed by $\zeta\in\widehat{T}$ are denoted by ${\cal Q}:=(Q_\zeta)_{\zeta\in\widehat{T}}$ and ${\cal J}:=(J_\zeta)_{\zeta\in\widehat{T}}$\, respectively.

\vskip1em 

{\bf Theorem 5.2.} 
{\it 
Let $\widetilde{G}_n=\widetilde{G}(m,1,n)=D_n(T)\rtimes\widetilde{{\mathfrak S}}_n,\,T={\boldsymbol Z}_m,\,4\le n <\infty$. The character $\chi\big(\Pi^{\rm Y}_{\Lambda^n}|\cdot\big)$ of spin IR\, $\Pi^{\rm Y}_{\Lambda^n},\,\Lambda^n=({\boldsymbol \lambda}^\zeta)_{\zeta\in\widehat{T}}\in{\boldsymbol Y}^{\rm sh}_{\!n}(\widehat{T})$, with\, {\rm Y\,=\,odd} or\, {\rm IV}, 
is given as follows: 
let ${\boldsymbol \nu}=(n_\zeta)_{\zeta\in\widehat{T}},\,n_\zeta=|{\boldsymbol \lambda}^\zeta|,$ and ${\boldsymbol I}_n=\bigsqcup_{\zeta\in\widehat{T}}I_{n,\zeta}$ the associated standard partition, and put $H'_n=D_n(T)\rtimes \widetilde{{\mathfrak S}}_{\boldsymbol \nu},\,\widetilde{{\mathfrak S}}_{\boldsymbol \nu}={\widehat{*}}_{\zeta\in\widehat{T}}\,\widetilde{{\mathfrak S}}_{I_{n,\zeta}}$.

{\rm (i)}\;   
If $g'\in \widetilde{G}_n$ is not conjugate to an element in $H'_n$, then \,$\chi\big(\Pi^{\rm Y}_{\Lambda^n}|g'\big)=0$.

{\rm (ii)}\;  
If $g'=(d,\sigma')\in H'_n$ has a standard decomposition given in (\ref{2014-04-07-1-8}),  
then 
\begin{eqnarray}
\label{2014-03-17-11-2}
\chi\big(\Pi^{\rm Y}_{\Lambda^n}|g'\big)
=
\sum_{{\cal Q},{\cal J}} 
b(\Lambda^n;{\cal Q},{\cal J};g) X(\Lambda^n;{\cal Q},{\cal J};g'), 
\end{eqnarray}
~\\[-7ex]
\begin{gather}
\nonumber
b(\Lambda^n;{\cal Q},{\cal J};g)
= 
\frac{(n-|\mathrm{supp}(g)|)!}{%
\prod_{\zeta\in\widehat{T}}
\Big(n_\zeta - |Q_\zeta| - 
\sum_{j \in J_\zeta}\ell_j\Big)!
}\,, 
\\[.5ex]
\nonumber 
X(\Lambda^n;{\cal Q},{\cal J};g')
=
\prod_{\zeta\in\widehat{T}} 
\Big(\prod_{q\in Q_\zeta}\zeta(t_q)\cdot\prod_{j\in J_\zeta}\zeta\big(P(d_j)\big)\Big)
 \cdot \chi\big(\tau_{\Lambda^n}|\,\sigma'_{{\cal Q}{\cal J}}\sigma'(\sigma'_{{\cal Q}{\cal J}})^{-1}\big),
\end{gather}
where a pair $\big({\cal Q},{\cal J}\big)$ runs over such ones that satisfy the condition 
\begin{eqnarray}
\label{2014-04-07-11-8}
\nonumber
\mbox{\bf (QJ)}\qquad\qquad\quad
|Q_\zeta|+{\sum}_{j\in J_\zeta}\ell_j\le n_\zeta\;\;\;
(\forall \zeta\in\widehat{T}), \qquad\qquad\qquad\qquad
\end{eqnarray}
and $\chi\big(\tau_{\Lambda^n}|\,\cdot)$ denotes spin character of\, $\widetilde{{\mathfrak S}}_{\boldsymbol \nu}$ of spin IR $\tau_{\Lambda^n}$, and further 
$\sigma'_{{\cal Q}{\cal J}}\in\widetilde{{\mathfrak S}}_n$ is an element whose image 
$\sigma_{{\cal Q}{\cal J}}=\Phi_{\mathfrak S}(\sigma'_{{\cal Q}{\cal J}})\in{\mathfrak S}_n$ satisfies, for each $\zeta\in\widehat{T}$, 
\begin{eqnarray}
\label{2014-04-07-21-8}
&&
\sigma_{{\cal Q}{\cal J}}Q_\zeta\subset I_{n,\zeta},\;\; \sigma_{{\cal Q}{\cal J}}K_j\subset I_{n,\zeta}\;\;(j\in J_\zeta). 
\end{eqnarray}
}

\vskip.2em

Taking into account Lemmas 4.2 and 4.3, we know the following

\vskip1em

{\bf Proposition 5.3.} 
{\it 
In the character formula {\rm (\ref{2014-03-17-11-2})} of spin IR \,$\Pi^{\rm Y}_{\Lambda^n}$, the factor $X(\Lambda^n;{\cal Q},{\cal J};g')$ can be rewritten, in case $\sigma'$ is even, as 
in the character formula {\rm (\ref{2014-03-17-11})} as 
$$
X(\Lambda^n;{\cal Q},{\cal J};g')
=\!
\prod_{\zeta\in\widehat{T}} 
\Big\{\prod_{q\in Q_\zeta}\zeta(t_q)\,\prod_{j\in J_\zeta}\zeta\big(P(d_j)\big)\Big\} 
\chi\big(\tau_{{\boldsymbol \lambda}^\zeta}\,\big|\,\sigma'_{{\cal Q}{\cal J}}\big({\textstyle \large\prod}_{j\in J_\zeta}\sigma'_j\big)(\sigma'_{{\cal Q}{\cal J}})^{-1}\big). 
$$
 
}

\section{\large Spin IRs and their characters of spin type $\chi^{\rm I}$} 

Let $4\le n<\infty$, $m$ even, and spin type $\chi^{\rm I}=(-1,-1,-1)$. This case was treated in [II, Part II], and here we rewrite the results there in a different,  short way for the later use. Schur multiplier for $G_n=G(m,1,n)$ is $Z=\langle z_1,z_2,z_3\rangle$ and the representation group $G'_n:=R\big(G(m,1,n)\big)$ is given in Theorem 1.1\,(ii).  
For $g'=(d',\sigma')\in G'_n=\widetilde{D}^\vee_n\rtimes \widetilde{{\mathfrak S}}_n$, put $g=\Phi(g')=(d,\sigma)\in G_n=D_n\rtimes{\mathfrak S}_n$ and 
$$
\mathrm{supp}(d'):=\mathrm{supp}(d),\,\mathrm{supp}(\sigma'):=\mathrm{supp}(\sigma),\,\mathrm{supp}(g'):=\mathrm{supp}(g),\,L(\sigma'):=L(\sigma).
$$
\vskip.2em

{\bf 6.1. Spin IRs of $G'_n$ of spin type $\chi^{\rm I}=(-1,-1,-1)$.}\; 
  Write the semidirect product group $G'_n=\widetilde{D}^\vee_n\rtimes \widetilde{{\mathfrak S}}_n$ by 
 $G'=U'\rtimes S'$, and $G_n=D_n\rtimes {\mathfrak S}_n$ by $G=U\rtimes S$. 
\vskip.5em

{\it Step 1.}\; Let $\varepsilon,a,b,c$ be $2\times 2$ matrices 
given in (\ref{2014-04-13-1-8}). 
For $k:=[n/2]$, we define matrices 
 $Y_j\;(j\in{\boldsymbol I}_{2k+1})$ of degree $2^{k}$\,: for  
$i\in{\boldsymbol I}_{k}$,
\begin{eqnarray}
\label{2014-04-13-2-8}
&&
\left\{
\begin{array}{l}
Y_{2i-1} =  c^{\otimes(i-1)}\otimes a\otimes \varepsilon^{\otimes(k-i)}, \, 
\\[.5ex]
Y_{2i} \;\;\;=\;  c^{\otimes(i-1)}\otimes b\otimes \varepsilon^{\otimes(k-i)}, 
\;\;
\\[.5ex] 
Y_{2k+1} = c^{\otimes k}\,.
\end{array}
\right.
\end{eqnarray}

As parameter space, put  
\begin{eqnarray}
\label{2014-04-13-3-8}
&&\quad
\left\{
\begin{array}{l}
\Gamma_n:=\big\{
\gamma=(\gamma_1,\gamma_2,\ldots,\gamma_n)\,;\; 
0\le \gamma_j\le m-1\;\;(j\in{\boldsymbol I}_n)\big\},
\\[.5ex]
\Gamma^0_n:=
\big\{\gamma=(\gamma_j)_{j\in{\boldsymbol I}_n}\in\Gamma_n\,;\;0\le \gamma_j <m^0\!:=\!m/2\,(j\in{\boldsymbol I}_n) \big\}, 
\\[.5ex] 
\Gamma^1_n:=
\big\{\gamma=(\gamma_j)_{j\in{\boldsymbol I}_n}\in\Gamma_n\,;\;0\le \gamma_1\le \gamma_2\le \ldots\le \gamma_n\big\}.
\end{array}
\right.
\end{eqnarray}
We define two kinds of actions on $\gamma\in \Gamma_n$ as follows: for $p\in{\boldsymbol I}_n,\;\sigma\in{\mathfrak S}_n,$
\begin{eqnarray}
\label{2014-12-20-11-8}
&&
\quad
\left\{
\begin{array}{l}
\tau_p\gamma:=(\gamma'_j)_{j\in{\boldsymbol I}_n},\;\gamma'_p\equiv\gamma_p+m^0\;(\mathrm{mod}\;m), \;\gamma'_j=\gamma_j\;(j\ne p);
\\[.5ex]
 \sigma\gamma:=\big(\gamma_{\sigma^{-1}(j)}\big)_{j\in{\boldsymbol I}_n}.
\end{array}
\right. 
\end{eqnarray}

For $\gamma\in\Gamma_n$, define a one-dimensional character of $U=D_n=\langle y_j\;(j\in{\boldsymbol I}_n)\rangle$ as 
\begin{eqnarray}
\label{2014-04-13-11-8}
&&
\zeta_\gamma=(\zeta_j)_{j\in{\boldsymbol I}_n},\;\;
\zeta_j(y_j):=\omega^{\gamma_j}\;(j\in{\boldsymbol I}_n),\;\;\omega:=e^{2\pi i/m}, 
\end{eqnarray}
and induce it up to $U'=D^\wedge_n\times Z_3$ through the homomorphism $U'\to U\cong U'/Z_{2,3}$ with $Z_{2,3}:=\langle z_2,z_3\rangle$, and denote it again by $\zeta_\gamma$ for brevity:   
\begin{eqnarray}
\label{2017-11-20-1}
\zeta_\gamma(\widehat{\eta}_j):=\zeta_\gamma(y_j),\;\;
\zeta_\gamma(z_2):=1,\;\;\zeta_\gamma(z_3):=1.
\end{eqnarray}
For IRs of\, $U'$ of spin type $(z_2,z_3)\to (-1,-1)$, we first put 
\begin{eqnarray}
\label{2014-04-17-91-8}
P_\gamma(z_2):=-E,\quad P_\gamma(\widehat{\eta}_j):=\zeta_\gamma(\widehat{\eta}_j)\,\widehat{Y}_j\;\;(j\in{\boldsymbol I}_n),
\end{eqnarray}
where $\widehat{Y}_j:=(-1)^{j-1}Y_j\;(j\in{\boldsymbol I}_n)$. 
Then $P_\gamma$ is a spin IR of $D^\wedge_n$ of dimension $2^k$. Their characters are given in [II, Theorem 6.3].   By [ibid., Lemma 6.4], we know  that $P_\gamma\cong P_{\gamma'}$ if and only if\,  

(1) in case $n$ is even, $\gamma, \gamma'$ are congruent under $\langle \tau_p\;(p\in{\boldsymbol I}_n)\rangle$,\, and 
 
(2) in case $n$ is odd, $\gamma, \gamma'$ are congruent under $\langle \tau_p\tau_q\;(p,q\in{\boldsymbol I}_n)\rangle$. 
\\
In case $n$ is odd, we put 
\begin{eqnarray*}
P^+_\gamma=P_\gamma,\;P^-_\gamma=P_{\tau_n\gamma}\quad (\gamma\in\Gamma^0_n). 
\end{eqnarray*}
Denote by $\widehat{U'}^{(-1,-1)}$ the set of equivalence classes of IRs of $U'$ with spin type $(z_2,z_3)\to (-1,-1)$.  
A complete set of representatives (= CSR) of the spin dual $D^\wedge_n$ of spin type $z_2\to -1$ was given in [ibid., Theorem 6.5], from which we get the following.  
Let $\mathrm{sgn}_{Z_3}$ be the non-trivial character of $Z_3=\langle z_3\rangle$. 

\vskip1em 

{\bf Lemma 6.1.} {\it 
A CSR of the spin dual \;$\widehat{U'}^{(-1,-1)}$ of spin type $(z_2,z_3)\to (-1,-1)$ 
is given for $U'=D^\wedge_n\times Z_3$ as  
\vskip.5em

\quad 
for $n\ge 2$ even, 
\; 
$\big\{P_\gamma\cdot\mathrm{sgn}_{Z_3} \,;\,\gamma\in\Gamma^0_n\big\},$

\quad 
for $n\ge 3$ odd, 
\;\;\,
$\big\{P^\epsilon_\gamma\cdot\mathrm{sgn}_{Z_3}\,;\,\gamma\in\Gamma^0_n,\;\epsilon=+,-\big\}.
$
}

\vskip.8em

 Put 
\begin{eqnarray}
\label{2016-05-13-1}
\left\{
\begin{array}{ll}
P^0:=P_{\bf 0}\;(\gamma={\bf 0})\quad\;&
\mbox{\rm in case $n$ is even,} 
\\[.5ex]
P^+=P_{\bf 0},\;P^-:=P_{\tau_n{\bf 0}}\quad\;&
\mbox{\rm 
 in case $n$ is odd.}
\end{array}
\right.
\end{eqnarray}
Then \,$P_\gamma=\zeta_\gamma  P^0,\, 
P^\epsilon_\gamma:=\zeta_\gamma P^\pm\;(\epsilon=+,-)$.

The action of $\sigma'\in S'=\widetilde{{\mathfrak S}}_n$ on $\widehat{\eta}_j \in D^\wedge_n$ is through $\sigma=\Phi(\sigma')\in{\mathfrak S}_n$ as 
\begin{eqnarray}
\label{2016-06-24-1}
{\sigma'}(\widehat{\eta}_j)=z_3^{\;L(\sigma')}\widehat{\eta}_{\sigma(j)}\quad(j\in{\boldsymbol I}_n).
\end{eqnarray}
Calculating the character of $(r_i\rho)(u')=\rho(r_i^{-1}u')\;(u'\in U', i\in{\boldsymbol I}_{n-1})$ for $\rho=P_\gamma\cdot \mathrm{sgn}_{Z_3}$ etc., we obtain the following. 

\vskip1em

{\bf Lemma 6.2.}\; (cf. [II, Theorem 7.1])\;  
{\it 
$\sigma'\in \widetilde{{\mathfrak S}}_n$ acts on IRs listed in Lemma 6.1 as follows: with\, $\sigma=\Phi_{\mathfrak S}(\sigma')$,\; 
$\sigma'(P_\gamma\cdot\mathrm{sgn}_{Z_3})\cong P_{\sigma\gamma}\cdot\mathrm{sgn}_{Z_3}\;(\sigma'\in\widetilde{{\mathfrak S}}_n).$\,  
In case $n\ge 3$ odd, this gives\; 
$
\sigma'(P^\epsilon_\gamma\cdot\mathrm{sgn}_{Z_3})\cong P^\epsilon_{\sigma\gamma}\cdot\mathrm{sgn}_{Z_3}\;\;(\epsilon=+,-)$. 
}

\vskip1em
 
{\bf Lemma 6.3.}\; 
{\it 
A CSR of the spin dual \;$\widehat{U'}^{(-1,-1)}$ of spin type $(z_2,z_3)\to (-1,-1)$ under the action of $S'=\widetilde{{\mathfrak S}}_n$ is given by 
\begin{eqnarray}
\label{2014-04-16-13}
&&
\rho=
\left\{
\begin{array}{lll}
 P_\gamma\cdot\mathrm{sgn}_{Z_3}\;\;\quad &(\gamma\in \Gamma^0_n\cap\Gamma^1_n), \quad&
  \mbox{\it for $n$ even,}
\\[.5ex]
  P^\pm_\gamma\cdot\mathrm{sgn}_{Z_3}\;&(\gamma\in \Gamma^0_n\cap\Gamma^1_n),
  \quad &
\mbox{\it for $n$ odd.}  
\end{array}
\right.
\end{eqnarray}
}

{\it Step 2.} 
$\gamma=(\gamma_j)_{j\in{\boldsymbol I}_n}\in \Gamma_n$ determines 
$\zeta_\gamma=(\zeta_j)_{j\in{\boldsymbol I}_n},\,\zeta_j=\zeta^{(\gamma_j)}\in\widehat{T}$, and then a partition ${\cal I}_n$ of ${\boldsymbol I}_n=[1,n]$ as     
\begin{eqnarray}
\label{2016-05-11-21}
&&
{\cal I}_n:=(I_{n,\zeta})_{\zeta\in\widehat{T}},\;\;
I_{n,\zeta}:=\{j\in{\boldsymbol I}_n\;;\;\zeta_j=\zeta\},\;\;
{\boldsymbol I}_n={\bigsqcup}_{\zeta\in\widehat{T}}I_{n,\zeta}. 
\end{eqnarray}

{\bf Lemma 6.4.} (cf. [II, Theorem 7.2])\; 
{\it 
A parameter $\gamma\in \Gamma_n^0\cap\Gamma^1_n$ corresponds bijectively to a standard partition $(I_{n,\zeta})_{\zeta\in\widehat{T}^0}$ of ${\boldsymbol I}_n$, or to ${\boldsymbol \nu}=(n_\zeta)_{\zeta\in{\cal K}}\in P_n({\cal K})$ with $n_\zeta=|I_{n,\zeta}|$ and ${\cal K}=\widehat{T}^0$ (i.e., $\gamma\leftrightarrow{\boldsymbol \nu}$). In this situation, the stationary subgroup $S'([\rho])$ of $\rho=P_\gamma\cdot\mathrm{sgn}_{Z_3},\,P^+_\gamma\cdot\mathrm{sgn}_{Z_3}$ and $P^-_\gamma\cdot\mathrm{sgn}_{Z_3}$ in $S'=\widetilde{{\mathfrak S}}_n$ is given as   
\begin{eqnarray}
\label{2016-05-11-22}
&&
S'([\rho])=\widetilde{{\mathfrak S}}_{\boldsymbol \nu}:=\Phi_{{\mathfrak S}}^{\;-1}\Big({\prod}_{\zeta\in{\cal K}}{\mathfrak S}_{I_{n,\zeta}}\Big)\;\cong\;\underset{\zeta\in{\cal K}}{\widehat{*}}\,\widetilde{{\mathfrak S}}_{n_\zeta}\,.
\end{eqnarray}

}

\vskip.2em

{\it Step 3.}\;   
Let $\rho$ be a spin IR of $U'$ listed in Lemma 6.1. For $\sigma'\in S'([\rho])$, the equivalence $\sigma'(\rho)\cong \rho$ is realised by an intertwining operator $J_\rho(\sigma')$ satisfying 
\begin{eqnarray}
\label{2016-05-12-1}
\rho\big({\sigma'}^{\,-1}(d')\big)=J_\rho(\sigma')^{-1}\rho(d')J_\rho(\sigma')\quad(d'\in U').
\end{eqnarray}

 Put\, $\widehat{Y}_j=(-1)^{j-1}Y_j\;(j\in{\boldsymbol I}_n)$\, and
\begin{eqnarray}
\label{2014-04-16-23-8}
&&
\nabla_n(r_j):=\frac{1}{\sqrt{2}}(\widehat{Y}_j-\widehat{Y}_{j+1})\;\hspace{8.7ex} (j\in{\boldsymbol I}_{n-1}),
\\
&&
\label{2014-04-20-11}
\nabla^-_n(r_j):=-(i\widehat{Y}_n)\,\nabla_n(r_j)\,(i\widehat{Y}_n)^{-1}\quad\big(j\in{\boldsymbol I}_{n-1},\;i=\sqrt{-1}\big). 
\end{eqnarray}
Then, \;
 $\nabla^-_n(r_j)=\nabla_n(r_j)\;(j\in{\boldsymbol I}_{n-2}),\;\nabla^-_n(r_{n-1})={\textstyle  
\frac{1}{\sqrt{2}}}
\big(\widehat{Y}_{n-1}+\widehat{Y}_n\big)$, and 
$\nabla_n$ and $\nabla^-_n$ both give spin representations of $\widetilde{{\mathfrak S}}_n$.

\bigskip

{\bf Lemma 6.5.} [II, Theorem 9.2, p.166]\; 

{\it   
{\rm (i)}\; 
Let $\rho$ be one of $P_\gamma\cdot\mathrm{sgn}_{Z_3},\,P^+_\gamma\cdot\mathrm{sgn}_{Z_3}\,(=P_\gamma\cdot\mathrm{sgn}_{Z_3})$ in Lemma 6.1. An intertwining operator in (\ref{2016-05-12-1}) is given as 
$J_\rho(\sigma')=c_{\sigma'}\nabla_n(\sigma')\;(\exists c_{\sigma'}\in{\boldsymbol C}^\times)$. Put $c_{\sigma'}=1\;(\forall \sigma')$, then $J_\rho(\sigma'):=\nabla_n(\sigma')$ gives a spin representation of $S'([\rho])$. 

{\rm (ii)}\; When $n$ is odd, let $\rho=P^-_\gamma\cdot\mathrm{sgn}_{Z_3}$ in Lemma 6.1. An intertwining operator in (\ref{2016-05-12-1}) is given as 
$J_\rho(\sigma')=c_{\sigma'}\nabla^-_n(\sigma')\;(\exists c_{\sigma'}\in{\boldsymbol C}^\times)$. Put $c_{\sigma'}=1\;(\forall \sigma')$, then $J_\rho(\sigma'):=\nabla^-_n(\sigma')$ gives a spin representation of $S'([\rho])$. 

}

 \vskip1em

{\it Step 4.}\;    
As the counter part $\pi^1$ of $\pi^0=\rho\cdot J_\rho$, we take non-spin IRs (with respect to $Z_1$) of the stationary subgroup $S'([\rho])$, or linear IRs of the quotient group 
$$
S'([\rho])/Z_1=\Phi_{{\mathfrak S}}^{\;-1}\Big({\prod}_{\zeta\in{\cal K}}{\mathfrak S}_{I_{n,\zeta}}\Big)\big/Z_1\cong{\prod}_{\zeta\in{\cal K}}\,{\mathfrak S}_{n_\zeta},\;n_\zeta=|I_{n,\zeta}|. 
$$

A CSR of the dual of $S'([\rho])/Z_1$ is given similarly as in Step 4 in \S 3.1, but in case of spin type $\chi^{\rm I}$ here, we take ${\cal K}=\widehat{T}^0$, instead of ${\cal K}=\widehat{T}$ there\,: 
\begin{eqnarray*}
\label{2014-04-19-1-8}
&&
{\boldsymbol Y}_{\!n}({\cal K})=\Big\{\Lambda^n=({\boldsymbol \lambda}^\zeta)_{\zeta\in{\cal K}}\;;\;{\boldsymbol \lambda}^\zeta\in{\boldsymbol Y},\,{\boldsymbol \nu}=(n_\zeta)_{\zeta\in{\cal K}}\in P_n({\cal K}),\,n_\zeta=|{\boldsymbol \lambda}^\zeta|\Big\}. 
\end{eqnarray*}
For $\Lambda^n=({\boldsymbol \lambda}^\zeta)_{\zeta\in{\cal K}}\in{\boldsymbol Y}_{\!n}({\cal K})$, put $\pi_{\Lambda^n}:=\boxtimes_{\zeta\in{\cal K}}\,\pi_{{\boldsymbol \lambda}^\zeta}$, an outer direct tensor product, and take $\pi^1:=\pi_{\Lambda^n}$. 
 From $\Lambda^n$, the standard partition ${\cal I}_n=(I_{n,\zeta})_{\zeta\in{\cal K}}$ of ${\boldsymbol I}_n$\, and a partition ${\boldsymbol \nu}=(n_\zeta)_{\zeta\in{\cal K}}\in P_n({\cal K})$ are determined. With ${\boldsymbol \nu}\leftrightarrow\gamma$, this situation is expressed as $\Lambda^n\rightsquigarrow \gamma,{\cal I}_n, {\boldsymbol \nu}$. We put, for $\pi^0, \pi^1$ above, 
\begin{eqnarray}
\label{2017-07-08-1}
&&
 \Pi(\pi^0,\pi^1)=\mathrm{Ind}^{G'_n}_{H'_n}(\pi^0\boxdot\pi^1)\quad  {\rm with}\quad H'_n:=U'\rtimes S'([\rho]).
\end{eqnarray}
 
\vskip.2em 

{\bf Theorem 6.6.} (cf. [II, Theorem 11.5])\;  
{\it 
Suppose $4\le n<\infty$. A complete system of representatives 
$\mbox{\rm {\footnotesize spin}IR}^{\rm I}\big(G(m,1,n)\big)$ of spin 
IRs of $R\big(G(m,1,n)\big)$ of spin type $\chi^{\rm I}=(-1,-1,-1)$ is given, with ${\cal K}=\widehat{T}^0$, as follows.  

{\rm (i)}\; In case $n$ is even, let $\Lambda^n \in{\boldsymbol Y}_{\!n}({\cal K}),\,\Lambda^n\rightsquigarrow \gamma$, and put
\begin{gather*}
\Pi^{\rm I}_{\Lambda^n}:=\Pi(\pi^0,\pi^1)\;\; {\rm with}\; \pi^0=\big(P_\gamma\cdot\mathrm{sgn}_{Z_3}\big)\cdot \nabla_n,\;\pi^1=\pi_{\Lambda^n},
\\
\mbox{\rm {\footnotesize spin}IR}^{\rm I}\big(G(m,1,n)\big)\;=\;\big\{\Pi^{\rm I}_{\Lambda^n}\;;\; \Lambda^n \in{\boldsymbol Y}_{\!n}({\cal K})\big\}.
\end{gather*}

{\rm (ii)}\; In case $n$ is odd,  
\begin{gather*}
\Pi^{\rm \,I+}_{\Lambda^n}:=\Pi(\pi^0,\pi^1)\;\;{\rm with}\;
\pi^0=\big(P^+_\gamma\cdot\mathrm{sgn}_{Z_3}\big)\cdot \nabla_n,\;\pi^1=\pi_{\Lambda^n}, 
\\
\Pi^{\rm \,I-}_{\Lambda^n}:=\Pi(\pi^0,\pi^1)\;\; {\rm with}\;
\pi^0=\big(P^-_\gamma\cdot\mathrm{sgn}_{Z_3}\big)\cdot \nabla^-_n,\;\pi^1=\pi_{\Lambda^n}, 
\\
\mbox{\rm {\footnotesize spin}IR}^{\rm I}\big(G(m,1,n)\big)
=
\big\{\Pi^{\rm I\epsilon}_{\Lambda^n}\,;\, 
(\epsilon,\Lambda^n) \in \{+,-\}\times{\boldsymbol Y}_{\!n}({\cal K})\big\}.
\end{gather*}

}

\vskip.2em

{\bf 6.2. Special spin IRs and tensor products with them.}\; 
As an important consequence of Lemma 6.1, we have the following structural information about spin IRs listed in Theorem 6.6. In (\ref{2016-05-13-1}), we took special spin IRs of $\widetilde{D}_n$  as\;  
$P^0:=P_{{\bf 0}}\;(\gamma={\bf 0})$ for $n$ even, and  
$P^+:=P_{{\bf 0}},\;P^-:=P_{\tau_n{\bf 0}}$\, for $n$ odd. Here     
 take also, for ${\cal K}=\widehat{T}^0$,
\begin{eqnarray}
\label{2014-12-24-21-8}
&&
\quad
\Lambda^n_0:=({\boldsymbol \lambda}^\zeta)_{\zeta\in{\cal K}}\in{\boldsymbol Y}_{\!n}({\cal K})\;\;{\rm with}\;\; {\boldsymbol \lambda}^{\zeta_0}=(n),\;{\boldsymbol \lambda}^\zeta=\varnothing\;(\zeta\ne\zeta_0:=\zeta^{(0)}).  
\end{eqnarray}
For this special $\Lambda^n_0$, the subgroup $H'_n$ is equal to $G'_n$ itself and the corresponding spin IRs of $G'_n$ is nothing else but  
\begin{eqnarray}
\label{2014-04-24-52}
\left\{
 \begin{array}{ll}
 \Pi^{\rm I0}_n
\,:=\Pi^{\rm I}_{\Lambda^n_0}= \big(P^0\cdot\mathrm{sgn}_{Z_3}\big)\cdot \nabla_n\,,\;\;\;& \mbox{\rm for $n\ge 4$ even,}
\\[.5ex]
\Pi^{\rm I\pm}_n
:=\Pi^{\rm I\pm}_{\Lambda^n_0}= \big(P^\pm\cdot\mathrm{sgn}_{Z_3}\big)\cdot \nabla^\pm_n\,,\qquad & \mbox{\rm for $n\ge 5$ odd.} 
\end{array}
\right.
\end{eqnarray}

\vskip.2em 

{\bf Proposition 6.7.}\; 
{\it Let $4\le n<\infty$. 
Spin IRs of spin type $\chi^{\rm I}$ have the following equivalence relation containing (non-spin) linear IRs\,: for $\Lambda^n\in{\boldsymbol Y}_{\!n}({\cal K}),\,{\cal K}=\widehat{T}^0,$ 
\begin{eqnarray*}
\label{2011-06-15-1}
\Pi^{\rm I}_{\Lambda^n}\,\cong\, \Pi^{\rm I0}_n\otimes \breve{\Pi}_{\Lambda^n}\,,\quad  
\Pi^{\rm I\pm}_{\Lambda^n}\,\cong\, \Pi^{\rm I\pm}_n\otimes \breve{\Pi}_{\Lambda^n}\,.
\end{eqnarray*}
}

{\it Proof.}  For the first equivalence, 
take \,${\cal T}=\Pi^{\rm I0}_n$ and $\pi=(\zeta_\gamma\,{\bf 1}_{Z_3})\cdot \pi_{\Lambda^n}$, where $\Lambda^n\rightsquigarrow\gamma\in \Gamma^0_n\cap\Gamma^1_n$.  
Then, by Lemma 2.3, we obtain the desired equivalence. Other two equivalences can be obtained similarly.
\hfill 
$\Box$\; 

\vskip1.2em

{\bf 6.3. Spin irreducible characters of spin type $\chi^{\rm I}=(-1,-1,-1)$.} 
\vskip.5em 

{\bf 6.3.1. Conjugacy classes modulo $Z$ in $G'_n=R\big(G(m,1,n)\big)=\widetilde{D}^\vee_n\rtimes\widetilde{{\mathfrak S}}_n$.}  

For $g'_1,g'_2\in G'_n$, $g'_1$ is said to be  
{\it conjugate to $g_2$ modulo $Z$} if 
$$
g'_1=zg'_0g'_2{g'_0}^{\,-1}\quad(\exists z\in Z, \exists g'_0\in G'). 
$$
\vskip.2em 

{\bf Lemma 6.8.} (cf. [II, Lemma 16.1])\; 
{\it 
Any element of $G'_n=R\big(G(m,1,n)\big)=\widetilde{D}^\vee_n\rtimes\widetilde{{\mathfrak S}}_n$ is conjugate modulo $Z$ to an element $g'$ satisfying  
the normalization condition:  

{\bf (NC)}\; $g'=(d',\sigma')\in G'_n$ has a decomposition\, $g'=\xi'_{q_1}\cdots \xi'_{q_r}\,g'_1\cdots g'_s,$ such that 
\begin{eqnarray}
\label{2016-05-12-11}
&&
\left\{
\begin{array}{l}
\xi'_q = \big(t'_q, (q)\big)={\widehat{\eta}_q}^{\;a_q},\;\;
g'_j=(d'_j, \sigma'_j),\;K_j:=\mathrm{supp}(\sigma'_j) \supset \mathrm{supp}(d'_j),
\\[.7ex]
\mbox{\rm $K_j$ is a subinterval $[a_j, b_j]$ of ${\boldsymbol I}_n=[1,n]$},
\\[.7ex]  
\sigma'_j=\sigma'_{K_j}:=r_{a_j}r_{a_j+1}\cdots r_{b_j-1},\;
d'_j =\widehat{\eta}_{a_j}^{\;\,\mathrm{ord}(d'_j)},\;\;\mbox{\rm for $j\in J:={\boldsymbol I}_s$},
\\[.7ex]
 \mbox{\rm $Q:=\{q_1,q_2,\ldots,q_r\}$ and $K_j\;(j\in J)$ are mutually disjoint,}
\end{array}
\right.
\end{eqnarray}
where \,$\mathrm{ord}\big(z_2^{\;a}z_3^{\;b}\widehat{\eta}_1^{\;a_1}\widehat{\eta}_2^{\;a_n}\cdots\widehat{\eta}_n^{\;a_n}\big):=\sum_{j\in{\boldsymbol I}_n}a_j\;(\mathrm{mod}\;m)$.
}
 \vskip1em

Such an element $g'\in G'_n$ is called a {\it standard representative} of conjugacy classes modulo $Z$ of $G'_n$.

Take $\Lambda^n=({\boldsymbol \lambda}^\zeta)_{\zeta\in{\cal K}}\in {\boldsymbol Y}_{\!n}({\cal K}),\,{\cal K}=\widehat{T}^0$, and let $\Lambda^n\rightsquigarrow\gamma, {\boldsymbol \nu}, {\cal I}_n$.   

\vskip1em

{\bf Lemma 6.9.}\; 
{\it 
Let $\rho\in \widehat{U'}^{(-1,-1)}$ be as in Lemma 6.3, and put $H'_n:=U'\rtimes S'([\rho])$. Then, a standard representative $g'\in G'_n$ in (\ref{2016-05-12-11}) is contained in $H'_n$ \,if and only if\;   
\lq\lq\,$
K_j \subset I_{n,\zeta}\;(\forall j\in J,\,\exists \zeta \in{\cal K}=\widehat{T}^0)$.''
 }

\vskip1em

{\bf 6.3.2. Characters of spin IRs of $G'_n$ of spin type $\chi^{\rm I}$.}
 
 We introduce two conditions concerning the supports of characters as follows. Express $g'=(d',\sigma')\in G'_n$ in the form of $g'=\xi'_{q_1}\cdots \xi'_{q_r}\,g'_1\cdots g'_s,$ where $\xi'_q$ and $g'_j=(d'_j,\sigma'_j)$ are basic components with disjoint supports. Then    
\\[1ex] 
\mbox{\bf Condition (I-00)}\quad 
$
\left\{
\begin{array}{l}
\mathrm{ord}(d')+L(\sigma')\equiv 0\;(\mathrm{mod}\;2)\,,
\\[.5ex]
\mathrm{ord}(\xi'_{q_i})\equiv 0\;(\forall i),\;\mathrm{ord}(d'_j)+L(\sigma'_j)\equiv 0\;(\mathrm{mod}\;2) \;\;(\forall j)\,; 
\end{array}
\right. 
$
\\[1ex] 
\mbox{\bf Condition (I-11)}\quad
$
\left\{
\begin{array}{l}
|\mathrm{supp}(g')|=n\;{\rm odd},\;\;\mathrm{ord}(d')+L(\sigma')\equiv 1\;(\mathrm{mod}\;2)\,,
\\[.5ex]
\mathrm{ord}(\xi'_{q_i})\equiv 1\;(\forall i),\;\;\mathrm{ord}(d'_j)\equiv 1\;(\mathrm{mod}\;2) \;\;(\forall j). 
\end{array}
\right. 
$
\vskip1em

The characters of special spin IRs are given as follows.
\vskip1em

{\bf Theorem 6.10.}\; 
{\it 
Let $g'=(d',\sigma')\in G'_n=\widetilde{D}^\vee_n\rtimes \widetilde{{\mathfrak S}}_n$.

{\rm (i)}\; {\bf Case of even $n=2n'\ge 4$.}\;
 
If $g'=(d',\sigma')\in G'$ satisfies the condition\,: $\mathrm{ord}(d')+L(\sigma')\equiv 0\;(\mathrm{mod}\;2)$, then
\vskip.6em

\qquad 
$\chi\big(\Pi^{\rm I0}_n\big|g'\big)=\mathrm{tr} \big(\Pi^{\rm I0}_n(g')\big)\ne 0\;\Longleftrightarrow\;g'$ satisfies Condition (I-00). 
\vskip.6em 
\noindent
Under Condition (I-00), for a standard representative $g'\in G'_n$ of conjugacy classes modulo $Z$ in (\ref{2016-05-12-11}), with $\ell_j=\ell(\sigma'_j):= \ell(\sigma_j),\,\sigma_j=\Phi_{\mathfrak S}(\sigma'_j)$, 
\begin{eqnarray}
\label{2011-06-16-1-8}
\quad
\chi\big(\Pi^{\rm I0}_n\big|g'\big)=
2^{n'}\cdot{\prod}_{j\in J}(-1)^{[(\ell_j-1)/2]}\, 2^{-(\ell_j-1)/2}. 
\end{eqnarray}

If $g'$ satisfies the condition\,:\, $\mathrm{ord}(d')+L(\sigma')\equiv 1\;(\mathrm{mod}\;2)$, then\; $\chi\big(\Pi^{\rm I0}_n\big|g'\big)= 0$.\; 
\vskip.5em

{\rm (ii)}\; {\bf Case of odd $n=2n'+1\ge 5$.}\;

If $g'=(d',\sigma')\in G'_n$ satisfies the condition\,:  $\mathrm{ord}(d')+L(\sigma')\equiv 0\;(\mathrm{mod}\;2)$, then
\vskip.6em

\quad 
$\chi\big(\Pi^{\rm I-}_n\big|g'\big)=\chi\big(\Pi^{\rm I+}_n\big|g'\big)=\chi\big(\Pi^{\rm I0}_n\big|g'\big)$,\; and 
 is given by the formula (\ref{2011-06-16-1-8}).
\vskip.6em

If $g'$ satisfies\,: $\mathrm{ord}(d')+L(\sigma')\equiv 1\;(\mathrm{mod}\;2)$, then\, $\chi\big(\Pi^{\rm I-}_n\big|g'\big)=-\chi\big(\Pi^{\rm I+}_n\big|g'\big),$ and
$$
\chi\big(\Pi^{\rm I+}_n\big|g'\big)\ne 0\;\Longleftrightarrow\;g'\;\;\mbox{\rm satisfies Condition (I-11)}.
$$
Under Condition (I-11), for $g'\in G'_n$ a standard representative of conjugacy classes modulo $Z$ given in (\ref{2016-05-12-11}),    
\begin{eqnarray} 
\label{2011-06-16-2}
\nonumber
\chi\big(\Pi^{\rm I+}_n\big|g'\big)
=
\varepsilon^{\rm I}(g')\, (2i)^{n'}\cdot {\prod}_{j\in J}(-1)^{[(\ell_j-1)/2]}\, 2^{-(\ell_j-1)/2},
\end{eqnarray}
where the sign $\varepsilon^{\rm I}(g')=\pm 1$ is determined explicitly by a formula in 
{\rm \,[II, Theorem 19.3, p.230]}. 
}

\vskip1em

Now we give spin irreducible characters of $G'_n$ of spin type $\chi^{\rm I}$ in general.

\vskip.8em

{\bf Theorem 6.11.} [II, Theorem 19.9]  
{\it 
Let 
$\Lambda^n=\big({\boldsymbol \lambda}^\zeta)_{\zeta\in\widehat{T}^0} \in {\boldsymbol Y}_{\!n}(\widehat{T}^0)$. For  
$g'=(d',\sigma')\in G'_n=\widetilde{D}^\vee_n\rtimes\widetilde{{\mathfrak S}}_n$, put\,   
$g=\Phi(g')\in G_n=G(m,1,n)$.  

{\rm (i)}\; {\bf Case of even $n\ge 4$.}\; 

In general,\, 
$\chi\big(\Pi^{\rm \;I}_{\Lambda^n}\big|g'\big)
=\chi\big(\Pi^{\rm I0}_n\big|g'\big)\times
\chi\big(\breve{\Pi}_{\Lambda^n}|g\big).$ 

\quad 
If $g'$ satisfies\; $\mathrm{ord}(d')+L(\sigma')\equiv 0\;(\mathrm{mod}\;2)$, then 
\\
\qquad\qquad
$\chi\big(\Pi^{\rm \,I}_{\Lambda^n}\big|g'\big)\ne 0\;\Longrightarrow$ $g'$ satisfies Condition (I-00).

\quad 
If $g'$ satisfies\; $\mathrm{ord}(d')+L(\sigma')\equiv 1\;(\mathrm{mod}\;2)$, then\, $\chi\big(\Pi^{\rm \,I}_{\Lambda^n}\big|g'\big)= 0$. 
 \vskip.3em

{\rm (ii)}\; {\bf Case of odd $n\ge 5$.}\; 

In general,\,  
$\chi\big(\Pi^{\rm I\epsilon}_{\Lambda^n}\big|g'\big)
=
\chi\big(\Pi^{\rm I\epsilon}_n\big|g'\big)
\times  
 \chi\big(\breve{\Pi}_{\Lambda^n}|g\big)\;\;(\epsilon=+,-).$

\quad 
If $g'$ satisfies\, $\mathrm{ord}(d')+L(\sigma')\equiv 0\;(\mathrm{mod}\;2)$, then 

\qquad \quad
$\chi\big(\Pi^{\rm I\pm}_{\Lambda^n}\big|g'\big)\ne 0\;\Longrightarrow$ $g'$ satisfies Condition (I-00). 

\quad 
If $g'$ satisfies\, $\mathrm{ord}(d')+L(\sigma')\equiv 1\;(\mathrm{mod}\;2)$, then 

\qquad \quad
$\chi\big(\Pi^{\rm I-}_{\Lambda^n}\big|g'\big)=-\,\chi\big(\Pi^{\rm I+}_{\Lambda^n}\big|g'\big)$,\, and 

\qquad \quad
$\chi\big(\Pi^{\rm I\pm}_{\Lambda^n}\big|g'\big)\ne 0\;\Longrightarrow$ $g'$ satisfies Condition (I-11).

}

\section{\large Spin IRs and their characters of spin type $\chi^{\rm II}$}

\quad
{\bf 7.1. Spin IRs of spin type type $\chi^{\rm II}$.}\; 
Let $4\le n<\infty$, $m$ even, and spin type $\chi^{\rm II}=(-1,-1,\,1)$. 
We use the notations in Theorem 1.1\,(ii), and put\; $G'_n=R\big(G(m,1,n)\big)=\widetilde{D}^{\vee}_n\rtimes \widetilde{{\mathfrak S}}_n=U'\rtimes S',\; \widetilde{D}^{\vee}_n=\widetilde{D}_n\times Z_3\;(=D^\wedge_n\times Z_3)$. 
The details of construction of IRs of spin type $\chi^{\rm II}$ are given in [II, \S 12], and here we summarize the results. We put ${\cal K}=\widehat{T}^0$ for the case of spin type $\chi^{\rm II}$. 
\vskip.5em

{\it Step 1.}\; 
 For IRs $\rho$ of $U'=\widetilde{D}^{\vee}_n=D^\wedge_n\times Z_3$ of spin type $(z_2,z_3)\to (-1,1)$, we refer Lemma 6.1\,: with the trivial character ${\bf 1}_{Z_3}$ of $Z_3$, 

{\it 
A CSR of the spin dual\, $\widehat{U'}^{(-1,1)}$ of $U'$ of spin type $(z_2,z_3)\to (-1,1)$ is given by 
\begin{eqnarray}
\label{2014-04-16-13-9}
&&
\rho=
\left\{
\begin{array}{lll}
 P_\gamma\cdot{\bf 1}_{Z_3}\;\;\quad &(\gamma\in \Gamma^0_n), \quad&
  \mbox{\it for $n$ even,}
\\[.5ex]
  P^\epsilon_\gamma\cdot{\bf 1}_{Z_3}\;&(\gamma\in \Gamma^0_n,\;\epsilon=+,-),
  \quad &
\mbox{\it for $n$ odd.}  
\end{array}
\right.
\end{eqnarray}
}

\vskip.2em

{\bf Lemma 7.1.}\,   
{\it 
The action of  
$S'$ on this spin dual is through $\Phi_{\mathfrak S}: \widetilde{{\mathfrak S}}_n\to {\mathfrak S}_n$ as 
\vskip.5em
\quad 
for $n\ge 2$ even, \hspace{3.2ex}$\sigma'(P_\gamma{\bf 1}_{Z_3})\cong P_{\sigma\gamma}{\bf 1}_{Z_3}\quad\;\big(\sigma=\Phi_{\mathfrak S}(\sigma')\in{\mathfrak S}_n\big)$, 
\vskip.2em
\quad 
for $n\ge 3$ odd,
$
\left\{
\begin{array}{lll}
\sigma'(P^+_\gamma{\bf 1}_{Z_3})\cong P^+_{\sigma\gamma}{\bf 1}_{Z_3},&
\sigma'(P^-_\gamma{\bf 1}_{Z_3})\cong 
P^-_{\sigma\gamma}{\bf 1}_{Z_3}\; & (\sigma'\in\widetilde{{\mathfrak A}}_n), 
\\[.5ex]
\sigma'(P^+_\gamma{\bf 1}_{Z_3})\cong P^-_{\sigma\tau_n\gamma}{\bf 1}_{Z_3},&
\sigma'(P^-_\gamma{\bf 1}_{Z_3})\cong 
P^+_{\sigma\tau_n\gamma}{\bf 1}_{Z_3}\; & (\sigma'\not\in\widetilde{{\mathfrak A}}_n).
\end{array}
\right.
$
}

\vskip1em

{\bf Lemma 7.2.}  
{\it 
A CSR for\, $\widehat{U'}^{(-1,1)}/S'$ is given as   
\vskip.3em
  
for $n\ge 2$ even,\; 
$\big\{P_\gamma{\bf 1}_{Z_3}\,;\,\gamma\in\Gamma^0_n\cap\Gamma^1_n\big\}$,

for $n\ge 3$ odd,\;  
$
\big\{P^+_\gamma{\bf 1}_{Z_3}\,;\,\gamma\in\Gamma^0_n\cap\Gamma^1_n\big\}\bigsqcup 
\big\{P^-_\gamma{\bf 1}_{Z_3}\,;\,\gamma\in \Gamma^0_n\cap\Gamma^1_n,\;\,\gamma_j\;\mbox{\rm all different}\big\}. $
}

\vskip1em

{\it Step 2.}\; 
We determine stationary subgroup $S'([\rho])$ as follows. 

\vskip.7em 

{\bf Lemma 7.3.}\; 
{\it 
The stationary subgroup $S'([\rho])$ of spin IRs $\rho= P_\gamma{\bf 1}_{Z_3},\;P^+_\gamma{\bf 1}_{Z_3}$ and $P^-_\gamma{\bf 1}_{Z_3}$ $(\gamma\in\Gamma^0_n\cap\Gamma^1_n)$ of\,  
$U'=\widetilde{D}^\vee_n$ is given as follows: 
\vskip.3em
For $n\ge 2$ even,
\hspace{3ex}
$S'\big([P_\gamma{\bf 1}_{Z_3}]\big)=\Phi_{\mathfrak S}^{\;-1}\Big({\prod}_{\zeta\in\widehat{T}^0}{\mathfrak S}_{I_{n,\zeta}}\Big)= \widetilde{{\mathfrak S}}_{\boldsymbol \nu}.$

 For $n\ge 3$ odd,
\hspace{4.2ex}
$
S'\big([P^\pm_\gamma{\bf 1}_{Z_3}]\big)=\widetilde{{\mathfrak A}}_n\bigcap\Phi_{\mathfrak S}^{\;-1}\Big({\prod}_{\zeta\in\widehat{T}^0}{\mathfrak S}_{I_{n,\zeta}}\Big)= \widetilde{{\mathfrak A}}_n\cap\widetilde{{\mathfrak S}}_{\boldsymbol \nu}\,.
$
}

\vskip.8em

{\it Step 3.}\; 
Two  
spin IRs $\nabla'_n,\,\nabla''_n$
 of $\widetilde{{\mathfrak S}}_n$ are defined as (cf. [II, \S 8.1, p.158])   
\begin{eqnarray}
\label{2014-12-25-31}
&&
\left\{
\begin{array}{l}
{\displaystyle 
\nabla'_n(r_i):=\mbox{$\frac{1}{\sqrt{2}}$}
(Y_i-Y_{i+1})\;\; (i\in{\boldsymbol I}_{n-1}),\quad
\nabla'_n(z_1):=-E, 
}
\\[.5ex]
\nabla''_n(r'_i):=-Y_{2n'+1} \nabla'_n(r'_i)Y_{2n'+1}\;\; (i\in{\boldsymbol I}_{n-1}),\;\;
\nabla''_n(z_1):=-E.  
\end{array}
\right.
\end{eqnarray}

 Also, for $n=2n'$ even, a spin representation of $\widetilde{{\mathfrak S}}_n=\widetilde{{\mathfrak S}}_{2n'}$ is defined by 
\begin{eqnarray}
\label{2014-04-29-12}
&&
\nabla^{\rm II}_n(r_j)=(iY_{2n'+1})\nabla'_n(r_j)=(iY_{2n'+1})\cdot{\textstyle \frac{1}{\sqrt{2}}}(Y_j-Y_{j+1})\;\;(j\in{\boldsymbol I}_{2n'-1}). 
\end{eqnarray}

For $n=2n'\!+\!1$ odd, two spin representations of $\widetilde{{\mathfrak A}}_n=\widetilde{{\mathfrak A}}_{2n'+1}$ are defined by restricting $\nabla'_n$ and $\nabla''_n$\,: for $\sigma'\in\widetilde{{\mathfrak A}}_n$,  
\begin{eqnarray}
\label{2014-04-29-11}
&&
\mho^+_n:=\nabla'_n|_{\widetilde{{\mathfrak A}}_n},\quad 
\mho^-_n:=\nabla''_n|_{\widetilde{{\mathfrak A}}_n}. 
\end{eqnarray}
\vskip.2em

{\bf Proposition 7.4.}\; [II, Theorem 10.2, p.168]
{\it 

{\bf (i)}\; Case of $n\ge 2$ even.\; 
For $\rho=P_\gamma{\bf 1}_{Z_3}\;\big(\gamma\in\Gamma^0_n\big)$, an intertwining operator in (\ref{2016-05-12-1}) for 
$\sigma'\in S'([\rho])$ is (up to a scalar multiple) $J_\rho(\sigma')=\nabla^{\rm II}_n(\sigma')$, and this gives a spin representation of the stationary subgroup  $S'([\rho])$. 
\\[.5ex]
\indent
 {\bf (ii)}\; Case of $n\ge 3$ odd.\;  
For $\rho=P^\varepsilon_\gamma\,{\bf 1}_{Z_3}\,\big(\gamma\in\Gamma^0_n,\,\varepsilon=+,-\big)$, an intertwining operator for 
$\sigma'\in S'([\rho])\subset\widetilde{{\mathfrak A}}_n$ is (up to a scalar multiple) $J_\rho(\sigma')=\mho^\varepsilon_n(\sigma')$. In a special case where $S'([\rho])=Z_1$, we understand as $J_\rho=\mho^\varepsilon_n|_{Z_1}=\mathrm{sgn}_{Z_1}\!\cdot I\;(I=$ the identity operator). 
}

\vskip1em 

{\it Step 4.}\; 
A representation $\pi^0=\rho\cdot J_\rho$ of the subgroup $H'_n=U'\rtimes S'([\rho])$ has spin type $(z_1,z_2,z_3)\to (-1,-1,\,1)$, and as its counter part, IR $\pi^1$ of $S'([\rho])$ should have the trivial spin type $z_1\to 1$ and is a linear IR of $S'([\rho])/Z_1$. 
Recall that $\Lambda^n=({\boldsymbol \lambda}^\zeta)_{\zeta\in\widehat{T}^0}\in{\boldsymbol Y}_{\!n}(\widehat{T}^0)$\, determines naturally 
a standard partition ${\cal I}_n=(I_{n,\zeta})_{\zeta\in\widehat{T}^0}$ of ${\boldsymbol I}_n$, ${\boldsymbol \nu}=(n_\zeta)_{\zeta\in\widehat{T}^0}\in P_n(\widehat{T}^0),\, |I_{n,\zeta}|=n_\zeta$, and $\gamma\in\Gamma^0_n\cap\Gamma^1_n$, i.e., $\Lambda^n\rightsquigarrow {\cal I}_n,\,{\boldsymbol \nu},\,\gamma$.

\vskip.3em
{\it Step 4-1.}\; 
Case of $n$ even.\, 
For 
$\rho=P_\gamma$, we have 
$S'([\rho])/Z_1\cong{\large\prod}_{\zeta\in\widehat{T}^0}{\mathfrak S}_{n_\zeta}={\mathfrak S}_{\boldsymbol \nu}$\,. 
A CSR of the dual of ${\mathfrak S}_{\boldsymbol \nu}$ is given by the set  $\{\pi_{\Lambda^n}\,;\,\Lambda^n\in {\boldsymbol Y}_{\!n}(\widehat{T}^0)\}$.

\vskip.3em
{\it Step 4-2.}\; Case of $n$ odd.\, 
   
{\sc Case odd-1.} Case where $n_\zeta\le 1\;(\forall \zeta\in\widehat{T}^0)$, and so\; $n\le m^0=m/2$. For $\rho=P^+_\gamma{\bf 1}_{Z_3}$ (or $P^-_\gamma{\bf 1}_{Z_3}$), $\prod_{\zeta\in\widehat{T}^0}{\mathfrak S}_{I_{n,\zeta}}=\{e\}$ and so $S'([\rho])=Z_1$. We can take $J_\rho=\mho^+_n|_{S'([\rho])}=\mathrm{sgn}_{Z_1}\!\cdot I$. Then representation of $S'([\rho])/Z_1\cong \{e\}$ is $\pi^1={\bf 1}$ (trivial representation).

{\sc Case odd-2.} Case where\, $n_\zeta\ge 2\;(\exists \zeta\in\widehat{T}^0)$.\;

For $\rho=P^+_\gamma{\bf 1}_{Z_3}$, we have\,    
$
S'([\rho])/Z_1\cong {\displaystyle {\mathfrak A}_n\bigcap{\mathfrak S}_{\boldsymbol \nu}.
}
$ 
A CSR of IRs of the direct product group ${\mathfrak S}_{\boldsymbol \nu}$ is given by the set $\{\pi_{\Lambda^n}\,;\,\Lambda^n=({\boldsymbol \lambda}^\zeta)_{\zeta\in\widehat{T}^0}\in {\boldsymbol Y}_{\!n}(\widehat{T}^0)\}$. From this, we obtain necessary IRs by restriction from ${\mathfrak S}_{\boldsymbol \nu}$ to ${\mathfrak A}_n\bigcap{\mathfrak S}_{\boldsymbol \nu}$ of index 2 as follows.

\vskip1em
  
{\bf Lemma 7.5.} 
{\it 
{\rm (i)}\;  For an IR $\pi_{\boldsymbol \lambda}\;({\boldsymbol \lambda}\in {\boldsymbol Y}_{\!q})$ of ${\mathfrak S}_q$, there holds \;$\mathrm{sgn}\cdot \pi_{\boldsymbol \lambda}\cong \pi_{({}^t{\boldsymbol \lambda})}$. 

{\rm (ii)}\; When\, 
$\mathrm{sgn}\cdot \pi_{\boldsymbol \lambda}\cong \pi_{{\boldsymbol \lambda}}$ (\,$\Leftrightarrow\,{}^t{\boldsymbol \lambda} ={\boldsymbol \lambda}$), 
 the restriction 
$\pi_{\boldsymbol \lambda}|_{{\mathfrak A}_n}$ splits into two mutually non-equivalent IRs as\;   
 $\rho^{(0)}_{\boldsymbol \lambda}\oplus \rho^{(1)}_{\boldsymbol \lambda}$.
 
{\rm (iii)}\; 
When $\mathrm{sgn}\cdot \pi_{\boldsymbol \lambda}\not\cong \pi_{{\boldsymbol \lambda}}$ (\,$\Leftrightarrow\,{}^t{\boldsymbol \lambda} \ne{\boldsymbol \lambda}$), 
the restriction \,$\rho_{\boldsymbol \lambda}:=\pi_{\boldsymbol \lambda}|_{{\mathfrak A}_n}\cong (\pi_{\,{}^t{\boldsymbol \lambda}})|_{{\mathfrak A}_n}$\, is irreducible.
}
 
\vskip1.2em 
{\bf Notation 7.6.}\; For the dual of  
${\mathfrak A}_n\bigcap{\mathfrak S}_{\boldsymbol \nu}$, let ${\cal K}=\widehat{T}^0$ or $\widehat{T}$ and    
\begin{gather*}
\label{2010-10-14-1-8}
{\boldsymbol Y}_{\!n}^{\mathfrak{A}}({\cal K})\;:=\;{\boldsymbol Y}_{\!n}^{\mathfrak{A}}({\cal K})^{\rm asym}\bigsqcup {\boldsymbol Y}_{\!n}^{\mathfrak{A}}({\cal K})^{\rm sym}\bigsqcup {\boldsymbol Y}_{\!n}({\cal K})^{\rm tri},
\qquad\qquad
\\[.7ex]
\quad 
\begin{array}{l}
{\boldsymbol Y}_{\!n}^{\mathfrak{A}}({\cal K})^{\rm asym}:=
\big\{\{\Lambda^n,{}^t\Lambda^n\}\;;\;\Lambda^n=({\boldsymbol \lambda}^\zeta)_{\zeta\in{\cal K}}\in{\boldsymbol Y}_{\!n}({\cal K}),\,{}^t\Lambda^n\ne\Lambda^n\big\},
\\[1ex]
{\boldsymbol Y}_{\!n}^{\mathfrak{A}}({\cal K})^{\rm sym}:=
\big\{(\Lambda^n,\kappa)\;;\;\Lambda^n\in{\boldsymbol Y}_{\!n}({\cal K}),\,{}^t\Lambda^n=\Lambda^n,\,
|{\boldsymbol \lambda}^\zeta|\ge 2\;(\exists \zeta),\,\kappa=0,1\big\},
\\[1ex]
{\boldsymbol Y}_{\!n}({\cal K})^{\rm tri}:=
\big\{\Lambda^n\;;\;\Lambda^n=({\boldsymbol \lambda}^\zeta)_{\zeta\in{\cal K}},\;\mbox{\rm quasi-trivial, i.e,}\;|{\boldsymbol \lambda}^\zeta|\le 1\;(\forall \zeta)\big\}.
\end{array}
\end{gather*}
Here, ${\boldsymbol Y}_{\!n}^{\mathfrak{A}}({\cal K})^{\rm tri}\ne \emptyset$ only when $n\le m^0=m/2$ (resp. $n\le m$) according as ${\cal K}=\widehat{T}^0$ (resp. ${\cal K}=\widehat{T}^0$), and  
$\Lambda^n$ is non quasi-trivial for ${\boldsymbol Y}_{\!n}^{\mathfrak{A}}({\cal K})^{\rm asym}$ and ${\boldsymbol Y}_{\!n}^{\mathfrak{A}}({\cal K})^{\rm sym}$, in the terminology of \S 4.4. 
 \vskip1.2em 
 
 {\bf Lemma 7.7.}\; 
 {\it 
 Put $G^0={\mathfrak S}_{\boldsymbol \nu},\,H^0={\mathfrak A}_n\cap G^0={\mathfrak A}_n\cap{\mathfrak S}_{\boldsymbol \nu}$.  
 
 {\rm (i)}\, For $\Lambda^n=({\boldsymbol \lambda}^\zeta)_{\zeta\in{\cal K}}\in {\boldsymbol Y}_{\!n}({\cal K})$, put ${\boldsymbol \nu}=(n_\zeta)_{\zeta\in{\cal K}}\in P_n({\cal K}),\,n_\zeta=|{\boldsymbol \lambda}^\zeta|$.  
For an IR $\pi_{\Lambda^n}$ of $G^0=\large\prod_{\zeta\in{\cal K}}{\mathfrak S}_{n_\zeta}$, \;$\mathrm{sgn}\cdot \pi_{\Lambda^n}\cong \pi_{({}^t\Lambda^n)}$ with ${}^t(\Lambda^n)={}^t\Lambda^n:=({}^t{\boldsymbol \lambda}^\zeta)_{\zeta\in{\cal K}}$. 

\vskip.3em
 
{\rm (ii)}\, Assume $n_\zeta\ge 2\;(\exists \zeta\in{\cal K})$.  
Then $[G^0:H^0]=2$. 
 \vskip.3em 
{\rm (a)} When ${}^t\Lambda^n=\Lambda^n$, i.e., ${}^t{\boldsymbol \lambda}^\zeta ={\boldsymbol \lambda}^\zeta\;(\forall \zeta\in{\cal K})$,\; 
$\pi_{\Lambda^n}|_{H^0}\cong \rho^{(0)}_{\Lambda^n}\oplus \rho^{(1)}_{\Lambda^n},\; 
\rho^{(0)}_{\Lambda^n}\not\cong \rho^{(1)}_{\Lambda^n}$\,. 
\vskip.4em 

{\rm (b)} When ${}^t\Lambda^n\ne\Lambda^n$, i.e., ${}^t{\boldsymbol \lambda}^\zeta \ne{\boldsymbol \lambda}^\zeta\;(\exists \zeta\in{\cal K})$,\; 
 $\rho_{\Lambda^n}:=\pi_{\Lambda^n}|_{H^0} \cong (\pi_{\,{}^t\Lambda^n})|_{H^0}$\, is irreducible. 
\vskip.4em 

{\rm (iii)}\, For a fixed ${\boldsymbol \nu}=(n_\zeta)_{\zeta\in{\cal K}}\in P_n({\cal K}),\, \exists n_\zeta\ge 2,$ the set of IRs listed in (a) and (b) gives a complete set of representatives of IRs of $H^0$. 
}
\vskip1em

 Taking induced representations $\Pi(\pi^0,\pi^1)=\mathrm{Ind}^{G'_n}_{H'_n}(\pi^0\boxdot\pi^1)$ from $H'_n=U'\rtimes S'([\rho])$, we get spin IRs of $G'_n$ of spin type $\chi^{\rm II}=(-1,-1,\,1)$\, as follows, with ${\cal K}=\widehat{T}^0$.  

\vskip1em

{\bf (II.even)}\; For $n\ge 4$ even, 
\vskip.3em 

\hspace*{5ex} 
$\Pi^{\rm II}_{\Lambda^n}=\Pi(\pi^0, \pi^1)$,\quad 
$\pi^0=\big(P_\gamma{\bf 1}_{Z_3}\cdot\nabla^{\rm II}_n\big)|_{H'_n},\;\;
\pi^1 =  \pi_{\Lambda^n},\;\; 
\Lambda^n\in{\boldsymbol Y}_{\!n}(\widehat{T}^0).$  
\vskip.5em

{\bf (II.odd)}\;\, For $n\ge 5$ odd, with $\mho^+_n,\, \mho^-_n$ in (\ref{2014-04-29-11}), and $H'_n\subset 
G^{\prime\hspace{.15ex}{\mathfrak A}}_n:=\widetilde{D}^\vee_n\rtimes\widetilde{{\mathfrak A}}_n$, 
\begin{eqnarray*}
\begin{array}{lll}
\mbox{\rm 1.}& 
\{\Lambda^n,{}^t\Lambda^n\}\in {\boldsymbol Y}_{\!n}^{\mathfrak{A}}(\widehat{T}^0)^{\rm asym},&
\Pi^{\mho+}_{\Lambda^n}=\Pi(\pi^0, \pi^1),\, 
\pi^0=\big(P^+_\gamma\hspace{.2ex}{\bf 1}_{Z_3}\cdot\mho^+_n\big)|_{H'_n},\,\pi^1 = \rho_{\Lambda^n}, 
\\[.3ex]
\mbox{\rm 2.}& 
(\Lambda^n,\kappa)\in {\boldsymbol Y}_{\!n}^{\mathfrak{A}}(\widehat{T}^0)^{\rm sym},&
\Pi^{\mho+}_{\Lambda^n,\,\kappa}=\Pi(\pi^0, \pi^1),\, 
\pi^0=\big(P^+_\gamma\hspace{.2ex}{\bf 1}_{Z_3}\cdot\mho^+_n\big)|_{H'_n},\,\pi^1\! = \rho^{(\kappa)}_{\Lambda^n},\, 
 \\[.3ex]
\mbox{\rm 3.}& 
(\epsilon,\Lambda^n)\in \{\pm\}\times{\boldsymbol Y}_{\!n}(\widehat{T}^0)^{\rm tri}& 
\Pi^{\mho\hspace{.15ex}\epsilon}_{\Lambda^n}=\Pi(\pi^0, \pi^1),\,  
\pi^0=\big(P^\epsilon_\gamma\hspace{.2ex}{\bf 1}_{Z_3}\cdot\mho^\epsilon_n\big)|_{H'_n},\, 
\pi^1 = \rho_{\Lambda^n}\,,
 \end{array}
\end{eqnarray*}
where, for $(\epsilon,\Lambda^n)$ above, $S'([\rho])=Z_1,\,H'_n=U'\times Z_1,\,\rho_{\Lambda^n}={\bf 1}_{Z_1}$\,. 
\vskip1em

{\bf Special spin IRs of spin type $(z_1,z_2,z_3)\to(-1,-1,\,1)$.}\; 

We define special spin IRs of spin type $\chi^{\rm II}$ as
\begin{eqnarray}
\label{2017-07-16-21}
&&
\left\{
\begin{array}{lll}
\mbox{\rm for $n\ge 4$ even}, \quad 
&\Pi^{\rm II\hspace{.15ex}0}_n:=(P^0\hspace{.15ex}{\boldsymbol 1}_{Z_3})\cdot \nabla^{\rm II}_n,\quad 
&\mbox{\rm spin IR of $G'_n$},
\\[.3ex]
\mbox{\rm for $n\ge 5$ odd}, \quad 
& 
\Pi^{\mho^\epsilon,\hspace{.15ex}{\mathfrak A}}_n:=(P^\epsilon\hspace{.15ex}{\boldsymbol 1}_{Z_3})\cdot \mho^\epsilon_n\;\;(\epsilon =\pm), \quad\;
&\mbox{\rm spin IR of $G^{\prime\hspace{.15ex}{\mathfrak A}}_n$}. 
\end{array}
\right.
\end{eqnarray}

\vskip.1em

{\bf Proposition 7.8.}\;  
{\it 

{\rm (i)}\; For $n\ge 4$ even,\; 
$\Pi^{\rm II}_{\Lambda^n}\cong \Pi^{\rm II\hspace{.15ex}0}_n\otimes\breve{\Pi}_{\Lambda^n},\;\Lambda^n\in{\boldsymbol Y}_{\!n}(\widehat{T}^0).$
\vskip.3em

{\rm (ii)}\; For $n\ge 5$ odd, let\, $G^{\mathfrak A}_n:=\Phi(G^{\prime\hspace{.15ex}{\mathfrak A}}_n),\,H_n:=\Phi(H'_n)$\, 
with 
\begin{eqnarray}
\label{2017-07-29-1}
G^{\prime\hspace{.15ex}{\mathfrak A}}_n:=\widetilde{D}^\vee_n\rtimes\widetilde{{\mathfrak A}}_n=U'\rtimes\widetilde{{\mathfrak A}}_n, \quad 
H'_n:=U'\rtimes S'([\rho])\subset G^{\prime\hspace{.15ex}{\mathfrak A}}_n, 
\end{eqnarray}
and\, $\Pi^{\rm \mho^\epsilon\!,{\mathfrak A}}_n:=P^\epsilon{\bf 1}_{Z_3}\!\cdot\mho^\epsilon_n\;(\epsilon=+,-)$ a spin IR of $G^{\prime\hspace{.15ex}{\mathfrak A}}_n$. Then, $H_n=D_n\rtimes\big({\mathfrak A}_n\cap{\mathfrak S}_{\boldsymbol \nu}\big),\,\Lambda^n\rightsquigarrow{\boldsymbol \nu},\gamma$, and 
according to Cases 1, 2 and 3 in (II.odd),  
\begin{eqnarray*}
\begin{array}{ll}
\;\;
\mbox{\rm Case 1}.\qquad 
&
\Pi^{\mho^+}_{\Lambda^n}\,\cong\,\mathrm{Ind}^{G'_n}_{G^{\prime\hspace{.15ex}{\mathfrak A}}_n}\Big(\Pi^{\mho^+\!,{\mathfrak A}}_n\otimes\mathrm{Ind}^{G^{\mathfrak A}_n}_{H_n}\!\big( \zeta_\gamma\hspace{.2ex}{\bf 1}_{Z_3}\cdot\rho_{\Lambda^n}\big)\Big),  
\\[1ex]
\;\;
\mbox{\rm Case 2}.\;\;
&
\Pi^{\mho^+}_{\Lambda^n,\kappa}\cong \mathrm{Ind}^{G'_n}_{G^{\prime\hspace{.15ex}{\mathfrak A}}_n}\Big(\Pi^{\mho^+\!,{\mathfrak A}}_n\otimes
\mathrm{Ind}^{G^{\mathfrak A}_n}_{H_n}\!\big(\zeta_\gamma\hspace{.2ex}{\bf 1}_{Z_3}\cdot\rho^{(\kappa)}_{\Lambda^n}\big)\Big),  
\\[1ex]
\;\;
\mbox{\rm Case 3}.\;
&
\Pi^{\mho^\epsilon}_{\Lambda^n}\,\cong\, \mathrm{Ind}^{G'_n}_{G^{\prime\hspace{.15ex}{\mathfrak A}}_n}\Big(\Pi^{\mho^\epsilon\!,{\mathfrak A}}_n\otimes\mathrm{Ind}^{G^{\mathfrak A}_n}_{H_n}\!\big(\zeta_\gamma\hspace{.2ex}{\bf 1}_{Z_3}\cdot\rho_{\Lambda^n}\big)\Big).  
\end{array}
\end{eqnarray*}

}

\vskip.1em

{\it Proof.}\,  
(i) Apply Lemma 2.3 for the pair $(G, H)$ with $G=G'_n=U'\rtimes S'$,  $H=U'\rtimes S'([\rho])$, and ${\cal T}=(P^0{\bf 1}_{Z_3})\cdot\nabla^{\rm II}_n=\Pi^{\rm II\hspace{.15ex}0}_n$,\, $\pi=\zeta_\gamma\boxdot \pi_{\Lambda^n}.$ 

(ii)   
Apply Lemma 2.3 for $G=G^{\prime\,{\mathfrak A}}_n$,  $H=U'\rtimes S'([\rho])\subset G$, and ${\cal T}=(P^\epsilon{\bf 1}_{Z_3})\cdot\mho^\epsilon_n\;(\epsilon=\pm)$ and $\pi=\zeta_\gamma\boxdot\rho_{\Lambda^n},\,\zeta_\gamma\boxdot\rho^{(\kappa)}_{\Lambda^n},\,\zeta_\gamma\boxdot\rho_{\Lambda^n}$ respectively for three cases. Then we obtain\, 
$
\mathrm{Ind}^{G^{\prime\hspace{.15ex}{\mathfrak A}}_n}_{H'_n}(\pi^0\boxdot\pi^1)\cong \Pi^{\mho^\epsilon\!,{\mathfrak A}}_n\otimes\mathrm{Ind}^{G^{\prime\hspace{.15ex}{\mathfrak A}}_n}_{H'_n}\!\big(\zeta_\gamma\boxdot\rho_{\Lambda^n}\big)
$\, etc. We induce up these equivalences from $G^{\prime\hspace{.15ex}{\mathfrak A}}_n$ to the whole $G'_n$, getting asserted equivalences respectively. 
\hfill 
$\Box$ 
\vskip1.2em

{\bf Theorem 7.9.}\; 
{\it  
Let \,$n\ge 4$.  
A complete set of representatives $\mbox{\rm {\footnotesize spin}IR}^{\rm II}\big(G(m,1,n)\big)$ of spin IRs of $R\big(G(m,1,n)\big)$ of spin type $\chi^{\rm II}=(-1,-1,\,1)$ is given as follows: 
\\[.5ex]
\indent
{\rm (i)}\;\; {\sc Case of $n$ even.} \;\; 
$\mbox{\rm {\footnotesize spin}IR}^{\rm II}\big(G(m,1,n)\big)\;=\;\big\{\Pi^{\rm II}_{\Lambda^n}\;;\; \Lambda^n \in{\boldsymbol Y}_{\!n}(\widehat{T}^0)\big\}$.
\\[.5ex]
\indent
{\rm (ii)}\; {\sc Case of $n$ odd}, $n>m^0=m/2$.\; 

\qquad\quad
$\mbox{\rm {\footnotesize spin}IR}^{\rm II}\big(G(m,1,n)\big)=\big\{\Pi^{\mho+}_{\Lambda^n}\;;\; 
\{\Lambda^n,{}^t\Lambda^n\} \in {\boldsymbol Y}_{\!n}^\mathfrak{A}(\widehat{T}^0)^{\rm asym}\big\}\bigsqcup 
\\[.3ex]
\hspace*{46ex}
\big\{
\Pi^{\mho+}_{\Lambda^n,\kappa}\;;\; 
(\Lambda^n,\kappa) \in {\boldsymbol Y}_{\!n}^\mathfrak{A}(\widehat{T}^0)^{\rm sym}\big\}$.

{\rm (iii)}\, {\sc Case of $n$ odd}, $n\le m^0=m/2$. \;

\qquad\quad 
$\mbox{\rm {\footnotesize spin}IR}^{\rm II}\big(G(m,1,n)\big)= 
\big\{
\Pi^{\mho+}_{\Lambda^n}\;;\;\{\Lambda^n,{}^t\Lambda^n\}\in {\boldsymbol Y}_{\!n}^{\mathfrak{A}}(\widehat{T}^0)^{\rm asym}\big\}\bigsqcup 
\\[.3ex]
\hspace*{12ex} 
\big\{
\Pi^{\mho+}_{\Lambda^n,\kappa}\;;\;(\Lambda^n,\kappa)\in {\boldsymbol Y}_{\!n}^{\mathfrak{A}}(\widehat{T}^0)^{\rm sym}\big\}
\bigsqcup
\big\{ 
\Pi^{\mho\hspace{.15ex}\epsilon}_{\Lambda^n}\;;\; (\epsilon,\Lambda^n)\in\{+,-\}\times{\boldsymbol Y}_{\!n}(\widehat{T}^0)^{\rm tri}\big\}$.

}

\vskip1.2em

{\bf
 7.2. Spin irreducible characters of spin type $\chi^{\rm II}$ for $n=2n'\ge 4$ even.}

We calculate the character of 
$\Pi=\mathrm{Ind}^{G'_n}_{H'_n}\pi,\;\pi=\pi^0\boxdot\pi^1$. Since the results are quite different, we divide the cases into two, depending on $n$ even or odd. 
 
Among IRs $\Pi^{\rm II}_{\Lambda^n}$, pick up a special one  
$\Pi^{\rm II\hspace{.15ex}0}_n:=(P^0{\bf 1}_{Z_3})\cdot \nabla^{\rm II}_n$ for $\Lambda^n=({\boldsymbol \lambda}^\zeta)_{\zeta\in\widehat{T}^0}$ with ${\boldsymbol \lambda}^{\zeta^{(0)}}=(n)$ and ${\boldsymbol \lambda}^\zeta=\varnothing$ for $\zeta\ne \zeta^{(0)}$. 
\vskip1em 

{\bf Lemma 7.10.}\; 
{\it 
The character of spin IR\,  
$\Pi^{\rm II}_{\Lambda^n}$ for $\Lambda^n\in{\boldsymbol Y}_{\!n}(\widehat{T}^0)$ is given as 
\begin{eqnarray}
\label{2016-06-25-11}
&&\;
\chi\big(\Pi^{\rm II}_{\Lambda^n}\big|g'\big)=\chi\big(\Pi^{\rm II\hspace{.15ex}0}_n\big|g'\big)\times
\chi\big(\breve{\Pi}_{\Lambda^n}|g\big)\quad\big(g'\in \widetilde{G}^{\rm \,II}_n,\,g=\Phi(g')\in G_n\big). 
\end{eqnarray}

}

{\it Proof.}\; Apply Proposition 7.8\,(i).
\hfill 
$\Box$\;

\vskip1em

The non-spin character $\chi\big(\breve{\Pi}_{\Lambda^n}|g\big)$ for $\Lambda^n\in{\boldsymbol Y}_{\!n}(\widehat{T}^0)\subset {\boldsymbol Y}_{\!n}(\widehat{T})$ is given by the formula (\ref{2014-03-17-11}). 
On the other side, we have $\chi\big(\Pi^{\rm II\hspace{.15ex}0}_n\big|g'\big)=\chi_{P^0}(d')\cdot\chi_{\nabla^{\rm II}_n}(\sigma')$ for $g'=(d',\sigma')\in G'_n=U'\rtimes S'$.  
Take $g'$ with standard decomposition $g'=\xi'_{q_1}\cdots \xi'_{q_r}\,g'_1\cdots g'_s,\,g'_j=(d'_j,\sigma'_j),$  
and put 
 $Q:=\{q_1,q_2,\ldots,q_r\}$, $J:={\boldsymbol I}_s$. Consider two conditions on $G'$ as 
 
\vskip1em
\noindent
\mbox{\bf Condition (II-00)}
\quad 
$
\left\{
\begin{array}{l}
\mathrm{ord}(d')\equiv 0,\quad L(\sigma')\equiv 0\;(\mathrm{mod}\;2)\;;
\\[.5ex]
\mathrm{ord}(\xi'_{q_i})\equiv 0\;(\forall i),\;\mathrm{ord}(d'_j)+L(\sigma'_j)\equiv 0\;(\mathrm{mod}\;2) \;\;(\forall j),
\end{array}
\right. 
$

\vskip.5em
\noindent
\mbox{\bf Condition (II-11)}
\;\; 
$
\left\{
\begin{array}{l}
|\mathrm{supp}(g')|=n=2n',\; \mathrm{ord}(d')\equiv L(\sigma')\equiv 1,\; r+s\equiv 1,
\\[.5ex]
\mathrm{ord}(\xi'_{q_i})\equiv 1\;(i\in {\boldsymbol I}_r),\;\; \mathrm{ord}(d'_j)\equiv 1\;(j\in{\boldsymbol I}_s)\;\;(\mathrm{mod}\;2).
\end{array}
\right.
$
\vskip1em
\noindent
Put $\mathrm{sgn}(\sigma'):=(-1)^{L(\sigma')}$, and 
divide the set of suffices 
$J=\{1,2,\ldots,s\}$ as  
\begin{eqnarray}
\label{2009-06-26-52}
&&
J=J_+\bigsqcup J_-,\quad J_{\pm}:=\{j \in J\;;\; \mathrm{sgn}(\sigma'_j)=\pm 1\}. 
\qquad\qquad
\end{eqnarray}

 {\bf Lemma 7.11.}\; [II, Theorem 20.2, p.237] 
{\it 
Let $n\ge 4$ be even. 

{\rm (i)}  
Case of $L(\sigma')\equiv 0\;(\mathrm{mod}\;2):$   
\vskip.5em
\qquad\quad 
$\chi\big(\Pi^{\rm II\hspace{.15ex}0}_n\big|g'\big)\ne 0\;\Longleftrightarrow\;g'$ satisfies Condition (II-00). 
\vskip.5em

 For such a $g'$, but normalised as in (\ref{2016-05-12-11}), 
 \begin{eqnarray}
\label{2011-06-17-1}
\chi\big(\Pi^{\rm II\hspace{.15ex}0}_n\big|g'\big)
=2^{n'}\cdot{\prod}_{j\in J}(-1)^{(\ell_j-1)/2}\, 2^{-(\ell_j-1)/2},
\end{eqnarray}
where, for $j\in J_-$, $\ell_j$ is even and $(-1)^{1/2}(-1)^{1/2}:=-1$.

{\rm (ii)}  
Case of $L(\sigma')\equiv 1\;(\mathrm{mod}\;2):$ 
\vskip.5em

\qquad\quad  
$\chi\big(\Pi^{\rm II\hspace{.15ex}0}_n\big|g'\big)\ne 0\;\Longleftrightarrow\;g'$ satisfies Condition (II-11). 
\vskip.5em

For such a $g'$, but normalised as in (\ref{2016-05-12-11}), 
\begin{eqnarray}
\label{2011-06-17-2}
\nonumber
\chi\big(\Pi^{\rm II\hspace{.15ex}0}_n\big|g'\big)
&=&
\varepsilon^{\rm II}(g')\cdot 2^{n'}i^{n'-1}\cdot{\prod}_{j\in J} (-1)^{\ell_j-1}\,2^{-(\ell_j-1)/2}
\\
\nonumber 
&=& 
-\varepsilon^{\rm II}(g')\cdot 2^{n'}i^{n'-1}\cdot{\prod}_{j\in J}2^{-(\ell_j-1)/2},
\end{eqnarray}
where the sign $\varepsilon^{\rm II}(g')=\pm 1$ is determined by a formula in [ibid., Theorem 17.7]. 
}
\vskip1em 

{\bf Note.} For the formula (\ref{2011-06-17-1}), 
the sign factor $\varepsilon\big({\textstyle \prod_{j\in J_-}g'_j}\big)$ in front of the right hand side in the general formula in [II, Theorem 20.2] disappeared here thanks to the normalization (\ref{2016-05-12-11}).

\vskip1em

{\bf Theorem 7.12.} [II, Theorem 20.2]\, 
{\it 
 Let $n\ge 4$ be even,   
 $\Lambda^n=({\boldsymbol \lambda}^\zeta)_{\zeta\in\widehat{T}^0}\in {\boldsymbol Y}_{\!n}(\widehat{T}^0)$, $\Lambda^n\rightsquigarrow\gamma, {\boldsymbol \nu}, {\cal I}_n$, and $\rho=P_\gamma{\bf 1}_{Z_3}$.  
 For $g'=(d',\sigma')\in G'_n$, put $g=\Phi(g')\in G_n=G(m,1,n),\,\sigma=\Phi_{\mathfrak S}(\sigma')\in{\mathfrak S}_n,\;\mathrm{sgn}_{\mathfrak S}(\sigma'):=\mathrm{sgn}_{\mathfrak S}(\sigma)$. Then,\;    
 \begin{gather*}
 \chi(\Pi^{\rm \,II}_{\Lambda^n}| g')=
\chi(\Pi^{\rm II\hspace{.15ex}0}_n|g')\times  
\chi(\breve{\Pi}_{\Lambda^n}|g),\quad
\chi(\Pi^{\rm \,II}_{{}^t\Lambda^n}| g')=\mathrm{sgn}_{\mathfrak S}(\sigma')\chi(\Pi^{\rm \,II}_{\Lambda^n}| g').
\end{gather*}

For $g'=(d',\sigma')\in H'_n=\widetilde{D}_n\stackrel{\rm II}{\rtimes}S'([\rho]),\;S'([\rho])=\widetilde{{\mathfrak S}}_{\boldsymbol \nu}$,   
\\
If $L(\sigma')\equiv 0\;(\mathrm{mod}\;2)$, then\;\; 
$\chi(\Pi^{\rm \;II}_{\Lambda^n}|g')\ne 0\Longrightarrow g'$ satisfies Condition (II-00).
\\
If $L(\sigma')\equiv 1\;(\mathrm{mod}\;2)$, then\;\; 
$\chi(\Pi^{\rm \;II}_{\Lambda^n}|g')\ne 0\Longrightarrow g'$ satisfies Condition (II-11). 

If $g'$ is not conjugate to an element in $H'_n$, then  \;$\chi\big(\Pi^{\rm \,II}_{\Lambda^n}\big|g'\big)=0,\;\chi\big(\breve{\Pi}_{\Lambda^n}|g)=0$. 
}

\vskip1.2em

{\bf 7.3. Spin irreducible characters of spin type $\chi^{\rm II}$ for $n=2n'+1\ge 5$ odd.}

Directly from Proposition 7.8, we obtain  

\vskip.8em 
{\bf Lemma 7.13.}\; 
{\it 
For spin IRs of $G'_n=R\big(G(m,1,n)\big),\,n$ odd, of spin type $\chi^{\rm II}=(-1,-1,1)$, there hold the following relations: according to the cases in (II.odd), 
\begin{eqnarray*}
\label{2016-08-23-31-2}
\;\,
\begin{array}{ll}
\mbox{\rm Case 1}.\qquad\quad &
\chi\big(\Pi^{\mho^+}_{\Lambda^n}\big)=\mathrm{Ind}^{G'_n}_{G^{\prime\hspace{.15ex}{\mathfrak A}}_n}\Big(\chi\big(\Pi^{\mho^+\!,{\mathfrak A}}_n\big)\cdot\chi\big(\mathrm{Ind}^{G^{\mathfrak A}_n}_{H_n}\!( \zeta_\gamma\hspace{.2ex}{\bf 1}_{Z_3}\cdot\rho_{\Lambda^n})\big)\Big),  
\\[.5ex]
\mbox{\rm Case 2}.\qquad &
\chi\big(\Pi^{\mho^+}_{\Lambda^n,\kappa}\big)\!=\! \mathrm{Ind}^{G'_n}_{G^{\prime\hspace{.15ex}{\mathfrak A}}_n}\Big(\chi\big(\Pi^{\mho^+\!,{\mathfrak A}}_n\big)\cdot\chi\big(
\mathrm{Ind}^{G^{\mathfrak A}_n}_{H_n}\!(\zeta_\gamma\hspace{.2ex}{\bf 1}_{Z_3}\cdot\rho^{(\kappa)}_{\Lambda^n})\big)\Big),  
\\[.5ex]
\mbox{\rm Case 3}.\qquad &
\chi\big(\Pi^{\mho^\epsilon}_{\Lambda^n}\big)\!=\! \mathrm{Ind}^{G'_n}_{G^{\prime\hspace{.15ex}{\mathfrak A}}_n}\Big(\chi\big(\Pi^{\mho^\epsilon\!,{\mathfrak A}}_n\big)\cdot\chi\big(\mathrm{Ind}^{G^{\mathfrak A}_n}_{H_n}\!(\zeta_\gamma\hspace{.2ex}{\bf 1}_{Z_3}\cdot\rho_{\Lambda^n})\big)\Big),   
\end{array}
\end{eqnarray*}
where $H_n=\Phi(H'_n)=D_n\rtimes\big({\mathfrak A}_n\cap{\mathfrak S}_{\boldsymbol \nu}\big),\,\Lambda^n\rightsquigarrow\gamma,{\boldsymbol \nu}$, and for a class function $F$ on $G^{\prime\hspace{.15ex}{\mathfrak A}}_n$,\, $\big(\mathrm{Ind}^{G'_n}_{G^{\prime\hspace{.15ex}{\mathfrak A}}_n}F\big)(g')=F(g')+F(\sigma'g'{\sigma'}^{\,-1})\;(g'\in G'_n)$ for a $\sigma'\in\widetilde{{\mathfrak S}}_n\setminus\widetilde{{\mathfrak A}}_n$\,. 
}

\vskip.8em

We give characters of the special spin IRs $\Pi^{\rm \mho^\epsilon\!,{\mathfrak A}}_n$ of $G^{\prime\hspace{.15ex}{\mathfrak A}}_n,\,n=2n'+1$ [II, \S 17.5].
\vskip.8em

{\bf Lemma 7.14.}\; 
{\it 
{\rm (i)}\; 
For spin IRs $\Pi^{\mho^\epsilon\!,{\mathfrak A}}_n=P^\epsilon{\bf 1}_{Z_3}\cdot\mho^\epsilon_n\;(\epsilon=+,-)$ of $G^{\prime\hspace{.15ex}{\mathfrak A}}_n$,
the representation matrices of $\Pi^{\rm \mho^-\!,{\mathfrak A}}_n(g')\;(g'\in G'_n)$ are given by replacement $Y_{2n'+1}\to -Y_{2n'+1}$ from those of $\Pi^{\rm \mho^+\!,{\mathfrak A}}_n(g')$. 

{\rm (ii)}\; 
For an $s'_0\in \widetilde{{\mathfrak S}}_n\setminus\widetilde{{\mathfrak A}}_n$, \;$s'_0\big(\Pi^{\rm \mho^+\!,{\mathfrak A}}_n\big)\cong \Pi^{\rm \mho^-\!,{\mathfrak A}}_n$\,. Put\; 
$
\Pi^{\rm IIodd}_n:=\mathrm{Ind}^{G'_n}_{G^{\prime\hspace{.15ex}{\mathfrak A}}_n}\Pi^{\rm \mho^+\!,{\mathfrak A}}_n$\, a spin IR of\, $G'_n$, then\; $\mathrm{Ind}^{G'_n}_{G^{\prime\hspace{.15ex}{\mathfrak A}}_n}\Pi^{\rm \mho^-\!,{\mathfrak A}}_n\cong \Pi^{\rm IIodd}_n$\, and \;$
\Pi^{\rm IIodd}_n|_{G^{\prime\hspace{.15ex}{\mathfrak A}}_n}\cong \Pi^{\rm \mho^+\!,{\mathfrak A}}_n\oplus \Pi^{\rm \mho^-\!,{\mathfrak A}}_n\,.
$

}

\vskip1em

For 
$g'=(d',\sigma')\in G^{\prime\hspace{.15ex}{\mathfrak A}}_n$, take a standard decomposition $g'=\xi'_{q_1}\cdots \xi'_{q_r}\,g'_1\cdots g'_s,\,g'_j=(d'_j,\sigma'_j),$ and divide $J={\boldsymbol I}_s$ as
$J=J_+\bigsqcup J_-,\;
J_\epsilon:=\{j\in J\,;\,\mathrm{sgn}_{\mathfrak S}(\sigma'_j)=\epsilon 1\},\;
\epsilon=+,-.$ For characters of the special spin IRs of $G^{\prime\hspace{.15ex}{\mathfrak A}}_n=\big(\widetilde{D}_n\times Z_3\big)\rtimes\widetilde{{\mathfrak A}}_n$, we have the following.

\vskip.8em

{\bf Theorem 7.15.} [II, \S 17.5]\; 

 {\it Let $n\ge 5$ be odd, and  
$g'=(d',\sigma')\in G^{\prime\hspace{.15ex}{\mathfrak A}}_n$ be as above.

{\rm (i)} Case of $\mathrm{ord}(d')\equiv 0\;(\mathrm{mod}\;2)$:
\\[.7ex]
\qquad\qquad
$
\chi\big(\Pi^{\mho^\epsilon\!,{\mathfrak A}}_n\big|g'\big)\ne 0\;
\Longleftrightarrow \;g'$ satisfies the condition (II-00). 
\vskip.8em
\noindent
In this case, $|J_-|$ is even, and if $g'$ is a standard representative modulo $Z$ in (\ref{2016-05-12-11}), 
\begin{eqnarray}
\label{2016-09-03-21}
\nonumber
\chi\big(\Pi^{\rm \mho^\epsilon\!,{\mathfrak A}}_n\big|g'\big)
\!\!&=&\!\!
2^{n'}\!\prod_{j\in J_+}\!2^{-(\ell_j-1)/2}(-1)^{(\ell_j-1)/2}\cdot
\prod_{j\in J_-}\!2^{-(\ell_j-1)/2}(-1)^{\ell_j/2-1}
\\
\label{2016-09-03-22}
\!\!&=&\!\!
 2^{n'}\prod_{j\in J}(-1)^{[(\ell_j-1)/2]}\,2^{-(\ell_j-1)/2}.
\end{eqnarray}

{\rm (ii)} Case of $\mathrm{ord}(d')\equiv 1\;(\mathrm{mod}\;2)$:
\\[.7ex]
\qquad\qquad 
$\chi\big(\Pi^{\mho^\pm\!,{\mathfrak A}}_n\big|g'\big)\ne 0\;\Longleftrightarrow\;g'$ 
 satisfies the condition ($\mho$-11) below: 
\\[.8ex]
{\bf ($\mho$-11)}
\quad 
$
\left\{
\begin{array}{l}
|\mathrm{supp}(g')|\!=\!n\!=\!2n'\!+\!1,\, \mathrm{ord}(d')\equiv 1,\, L(\sigma')\equiv 0\;\;(\mathrm{mod}\;2),
\\[.5ex]
\mathrm{ord}(\xi'_q)\equiv 1\;(q\in Q),\, \mathrm{ord}(d'_j)\equiv 1\;(j\in J)
\quad
(\therefore\; r\!+\!s\!\equiv\! 1). 
\end{array}
\right.
$
\vskip.8em
\noindent
In this case, $|J_-|$ is even, and if $g'$ is a standard representative modulo $Z$ in (\ref{2016-05-12-11}), 
\begin{eqnarray}
\label{2010-06-26-22}
\nonumber
\left\{
\begin{array}{l}
\chi\big(\Pi^{\mho^-\!,{\mathfrak A}}_n\big|g'\big)=-\chi\big(\Pi^{\mho^+\!,{\mathfrak A}}_n\big|g'\big),
\\[.5ex]
\chi\big(\Pi^{\mho^+\!,{\mathfrak A}}_n\big|g'\big)
= 
\varepsilon^\mho(g')\cdot(2i)^{n'}\cdot \prod_{j\in J}2^{-(\ell_j-1)/2}\,,
\end{array}
\right.
\end{eqnarray}
where the sign $\varepsilon^\mho(g')=\pm 1$ is determined by putting $Y_{K_j}:=\prod_{p\in K_j}\!Y_p$ (product according to the order of $p$), 
\begin{eqnarray}
\label{2010-06-26-11}
{\prod}_{q\in Q}Y_q\times {\prod}_{j\in J}Y_{K_j}\;=\;\varepsilon^\mho(g')\cdot Y_1Y_2\cdots Y_n\,.
\end{eqnarray}
}

\vskip.2em 

Now we give the character of spin IR (of spin type $\chi^{\rm II}$) $\Pi^{\rm IIodd}_n=\mathrm{Ind}^{G'_n}_{G^{\prime\hspace{.15ex}{\mathfrak A}}_n}\Pi^{\mho^+\!,{\mathfrak A}}_n\cong \mathrm{Ind}^{G'_n}_{G^{\prime\hspace{.15ex}{\mathfrak A}}_n}\Pi^{\mho^-\!,{\mathfrak A}}_n$ of $G'_n=R\big(G(m,1,n)\big)\supset G^{\prime\hspace{.15ex}{\mathfrak A}}_n,\, n=2n'+1$. 

\vskip1em 

{\bf Theorem 7.16.}  
{\it 
Let $n\ge 5$ be odd. Then, for $g'=(d',\sigma')\in G'_n=R\big(G(m,1,n)\big)$,  
\begin{eqnarray*}
\label{2016-09-04-11}
&&
\chi\big(\Pi^{\rm IIodd}_n\big|g'\big)=
\left\{
\begin{array}{ll}
2\,\chi\big(\Pi^{\rm \mho^+\!,{\mathfrak A}}_n\big|g'\big),\qquad
&{\rm if}\; g'\in G^{\prime\hspace{.15ex}{\mathfrak A}}_n\;(\mbox{\rm i.e. $L(\sigma')\equiv 0\;(\mathrm{mod}\;2)$}),
\\[.5ex]
\;\;0\,, 
& {\rm if}\;g'\not\in G^{\prime\hspace{.15ex}{\mathfrak A}}_n\;(\mbox{\rm i.e. $L(\sigma')\equiv 1\;(\mathrm{mod}\;2)$}).
\end{array}
\right.
\end{eqnarray*}
}

\vskip.2em

Further, for irreducible characters of spin type $\chi^{\rm II}$, we have (cf. [II, \S 20.2]), 
\vskip1em

{\bf Theorem 7.17.}\,  
{\it 
Let $n\ge 5$ be odd, and\, 
$\Pi_n$ be one of\, $\Pi^{\mho+}_{\Lambda^n},\,\Pi^{\mho+}_{\Lambda^n,\,\kappa}$ and\, $\Pi^{\mho\hspace{.15ex}\epsilon}_{\Lambda^n}$ in Cases 1, 2 and 3 in (II.odd) just before  Proposition 7.8. Then its character satisfies the following. 
Let\, $s'_0\in\widetilde{{\mathfrak S}}_n\setminus\widetilde{{\mathfrak A}}_n$ and put\, $s_0=\Phi_{\mathfrak S}(s'_0)$. 

{\rm (i)}\,  For $g'\not\in G^{\prime\hspace{.15ex}{\mathfrak A}}_n$,\, 
$\chi(\Pi_n|g')=0$. 

{\rm (ii)}\,  For $g'=(d',\sigma')\in G^{\prime\hspace{.15ex}{\mathfrak A}}_n=U'\rtimes\widetilde{{\mathfrak A}}_n$, let\, $g=(d,\sigma)\in G^{\mathfrak A}_n=U\rtimes {\mathfrak A}_n$ be its canonical image. Recall\; 
$H_n=\Phi(H'_n)=D_n\rtimes\big({\mathfrak A}_n\cap{\mathfrak S}_{\boldsymbol \nu}\big),\,\Lambda^n\rightsquigarrow\gamma,{\boldsymbol \nu}$. 
 
In case\, $\mathrm{ord}(d')=\mathrm{ord}(d)\equiv 0\;(\mathrm{mod}\;2)$,\; 
\begin{eqnarray*}
\begin{array}{lll}
\mbox{\rm 1.}\quad& 
\{\Lambda^n,{}^t\Lambda^n\}\in {\boldsymbol Y}_{\!n}^{\mathfrak{A}}(\widehat{T}^0)^{\rm asym},&
\chi(\Pi^{\mho^+}_{\Lambda^n}|g')=\chi\big(\Pi^{\mho^+\!,{\mathfrak A}}_n\big|g'\big)\cdot\chi\big(\mathrm{Ind}^{G_n}_{H_n}\!( \zeta_\gamma\cdot\rho_{\Lambda^n})\big|g\big), 
\\[.3ex]
\mbox{\rm 2.}& 
(\Lambda^n,\kappa)\in {\boldsymbol Y}_{\!n}^{\mathfrak{A}}(\widehat{T}^0)^{\rm sym},&
\chi(\Pi^{\mho^+}_{\Lambda^n,\kappa}|g')=\chi\big(\Pi^{\mho^+\!,{\mathfrak A}}_n\big|g'\big)\cdot\chi\big(\mathrm{Ind}^{G_n}_{H_n}\!( \zeta_\gamma\cdot\rho^{(\kappa)}_{\Lambda^n})\big|g\big), 
\\[.3ex]
\mbox{\rm 3.}& 
(\epsilon,\Lambda^n)\in \{\pm\}\times{\boldsymbol Y}_{\!n}(\widehat{T}^0)^{\rm tri},& 
\chi(\Pi^{\mho^\epsilon}_{\Lambda^n}|g')=\chi\big(\Pi^{\mho^\epsilon\!,{\mathfrak A}}_n\big|g'\big)\cdot\chi\big(\mathrm{Ind}^{G_n}_{H_n}\!( \zeta_\gamma\cdot\rho_{\Lambda^n})\big|g\big). 
\end{array}
\end{eqnarray*}
In case $\mathrm{ord}(d')\equiv 1\;(\mathrm{mod}\;2)$, with an $s_0=\Phi_{\mathfrak S}(s'_0)\in{\mathfrak S}_n\setminus{\mathfrak A}_n$, 
\begin{eqnarray*}
\begin{array}{ll}
\mbox{\rm 1.}& 
\chi(\Pi^{\mho^+}_{\Lambda^n}|g')=\chi\big(\Pi^{\mho^+\!,{\mathfrak A}}_n\big|g'\big)\Big(\chi\big(\mathrm{Ind}^{G^{\mathfrak A}_n}_{H_n}\!( \zeta_\gamma\cdot\rho_{\Lambda^n})\big|g\big)-\chi\big(\mathrm{Ind}^{G^{\mathfrak A}_n}_{H_n}\!( \zeta_\gamma\cdot\rho_{\Lambda^n})\big|s_0gs_0^{\,-1}\big)\Big), 
\\[.9ex]
\mbox{\rm 2.}& 
\chi(\Pi^{\mho^+}_{\Lambda^n,\kappa}|g')
\!=\!\chi\big(\Pi^{\mho^+\!,{\mathfrak A}}_n\big|g'\big)\Big(\chi\big(\mathrm{Ind}^{G^{\mathfrak A}_n}_{H_n}\!( \zeta_\gamma\cdot\rho^{(\kappa)}_{\Lambda^n})\big|g\big)\!-\!\chi\big(\mathrm{Ind}^{G^{\mathfrak A}_n}_{H_n}\!( \zeta_\gamma\cdot\rho^{(\kappa)}_{\Lambda^n})\big|s_0gs_0^{\,-1}\big)\Big),
\\[.9ex]
\mbox{\rm 3.}& 
\chi(\Pi^{\mho^\epsilon}_{\Lambda^n}|g')=\chi\big(\Pi^{\mho^\epsilon\!,{\mathfrak A}}_n\big|g'\big)\Big(\chi\big(\mathrm{Ind}^{G^{\mathfrak A}_n}_{H_n}\!( \zeta_\gamma\cdot\rho_{\Lambda^n})\big|g\big)-\chi\big(\mathrm{Ind}^{G^{\mathfrak A}_n}_{H_n}\!( \zeta_\gamma\cdot\rho_{\Lambda^n})\big|s_0gs_0^{\,-1}\big)\Big).
\end{array}
\end{eqnarray*}
}
\vskip.2em

{\it Proof.}\; For the 1st part of (ii), it suffices to note  
\;$\chi\big(\Pi^{\mho^+\!,{\mathfrak A}}_n\big|g'\big)=\chi\big(\Pi^{\mho^+\!,{\mathfrak A}}_n\big|s'_0g'{s'_0}^{\,-1}\big)$ 
and 
$\chi\big(\mathrm{Ind}^{G_n}_{H_n}\!( \zeta_\gamma\cdot\rho_{\Lambda^n})\big|g\big)=\chi\big(\mathrm{Ind}^{G^{\mathfrak A}_n}_{H_n}\!( \zeta_\gamma\cdot\rho_{\Lambda^n})\big|g\big)+\chi\big(\mathrm{Ind}^{G^{\mathfrak A}_n}_{H_n}\!( \zeta_\gamma\cdot\rho_{\Lambda^n})\big|s_0gs_0^{\,-1}\big)$. 
For the 2nd part of (ii), it suffices to note\; \;$\chi\big(\Pi^{\mho^+\!,{\mathfrak A}}_n\big|g'\big)=-\chi\big(\Pi^{\mho^+\!,{\mathfrak A}}_n\big|s'_0g'{s'_0}^{\,-1}\big)$.
\hfill 
$\Box$\;

\vskip1.2em

 Recall linear IRs 
$\breve{\Pi}_{\Lambda^n}$ of $G_n=G(m,1,n)$ 
in (\ref{2016-05-06-1}).  
\vskip1.2em

{\bf Lemma 7.18.} 
{\it 
In Cases 1, 2 and 3 in Theorem 7.17 respectively, we have
\begin{eqnarray*}
\begin{array}{lll}
\mbox{\rm 1.}& 
\{\Lambda^n,{}^t\Lambda^n\}\in {\boldsymbol Y}_{\!n}^{\mathfrak{A}}(\widehat{T}^0)^{\rm asym},&
\mathrm{Ind}^{G_n}_{H_n}\!( \zeta_\gamma\cdot\rho_{\Lambda^n})
\cong \breve{\Pi}_{\Lambda^n}\oplus \breve{\Pi}_{\hspace{.15ex}{}^t\hspace{-.15ex}\Lambda^n},\;\;{\rm where}
\\[.2ex]
&&
\breve{\Pi}_{\hspace{.15ex}{}^t\hspace{-.15ex}\Lambda^n}(g)= \mathrm{sgn}_{\mathfrak S}(\sigma)\breve{\Pi}_{\Lambda^n}(g)\;\;\mbox{\rm for $g=(d,\sigma)$}, 
\\[.3ex]
\mbox{\rm 2.}& 
(\Lambda^n,\kappa)\in {\boldsymbol Y}_{\!n}^{\mathfrak{A}}(\widehat{T}^0)^{\rm sym},&
\mathrm{Ind}^{G_n}_{H_n}\!( \zeta_\gamma\cdot\rho^{(\kappa)}_{\Lambda^n})\cong \breve{\Pi}_{\Lambda^n}\quad(\kappa=0,1),  
\\[.3ex]
\mbox{\rm 3.}& 
(\epsilon,\Lambda^n)\in \{\pm\}\times{\boldsymbol Y}_{\!n}(\widehat{T}^0)^{\rm tri},& 
\mathrm{Ind}^{G_n}_{H_n}\!( \zeta_\gamma\cdot\rho_{\Lambda^n})
\cong \breve{\Pi}_{\Lambda^n}. 
\end{array}
\end{eqnarray*}
}

\vskip.2em

{\bf Proposition 7.19.}\; 
{\it 
Let $n\ge 5$ be odd. For characters of spin IRs $\Pi^{\mho+}_{\Lambda^n},\,\Pi^{\mho+}_{\Lambda^n,\,\kappa}$ and\, $\Pi^{\mho\hspace{.15ex}\epsilon}_{\Lambda^n}$ of spin type $\chi^{\rm II}$, we have the following expression. 
 
On the subset defined by\, $\mathrm{ord}(d')=\mathrm{ord}(d)\equiv 0\;(\mathrm{mod}\;2)$,\; 
\begin{eqnarray*}
\begin{array}{lll}
\mbox{\rm 1.}\quad& 
\{\Lambda^n,{}^t\Lambda^n\}\in {\boldsymbol Y}_{\!n}^{\mathfrak{A}}(\widehat{T}^0)^{\rm asym},&
\chi(\Pi^{\mho^+}_{\Lambda^n})=
\tfrac{1}{2}\chi(\Pi^{\rm IIodd}_n)\cdot\big(\chi(\breve{\Pi}_{\Lambda^n})+\chi(\breve{\Pi}_{\hspace{.15ex}{}^t\hspace{-.15ex}\Lambda^n})\big), 
\\[.3ex]
\mbox{\rm 2.}& 
(\Lambda^n,\kappa)\in {\boldsymbol Y}_{\!n}^{\mathfrak{A}}(\widehat{T}^0)^{\rm sym},&
\chi(\Pi^{\mho^+}_{\Lambda^n,\kappa})=
\tfrac{1}{2}\chi(\Pi^{\rm IIodd}_n)\cdot\chi\big(\breve{\Pi}_{\Lambda^n}), 
\\[.3ex]
\mbox{\rm 3.}& 
(\epsilon,\Lambda^n)\in \{\pm\}\times{\boldsymbol Y}_{\!n}(\widehat{T}^0)^{\rm tri},& 
\chi(\Pi^{\mho^\epsilon}_{\Lambda^n})=
\tfrac{1}{2}\chi(\Pi^{\rm IIodd}_n)\cdot\chi\big(\breve{\Pi}_{\Lambda^n}\big). 
\end{array}
\end{eqnarray*}

On the subset defined by\, $\mathrm{ord}(d')\equiv 1\;(\mathrm{mod}\;2)$, 
\begin{eqnarray*}
\begin{array}{ll}
\mbox{\rm 1.}& 
\chi(\Pi^{\mho^+}_{\Lambda^n}|g')=\tfrac{1}{2}\chi\big(\Pi^{\rm IIodd}_n\big|g'\big)\Big(\chi\big(\mathrm{Ind}^{G^{\mathfrak A}_n}_{H_n}\!( \zeta_\gamma\cdot\rho_{\Lambda^n})\big|g\big)-\chi\big(\mathrm{Ind}^{G^{\mathfrak A}_n}_{H_n}\!( \zeta_\gamma\cdot\rho_{\,{}^t\hspace{-.15ex}\Lambda^n})\big|g\big)\Big), 
\\[.9ex]
\mbox{\rm 2.}& 
\chi(\Pi^{\mho^+}_{\Lambda^n,\kappa}\big|g'\big)
\!=\!\tfrac{1}{2}\chi\big(\Pi^{\rm IIodd}_n\big|g'\big)\Big(\chi\big(\mathrm{Ind}^{G^{\mathfrak A}_n}_{H_n}\!( \zeta_\gamma\cdot\rho^{(\kappa)}_{\Lambda^n})\big|g\big)\!-\!\chi\big(\mathrm{Ind}^{G^{\mathfrak A}_n}_{H_n}\!( \zeta_\gamma\cdot\rho^{(\kappa+1)}_{\Lambda^n})\big|g\big)\Big),
\\[.9ex]
\mbox{\rm 3.}& 
\chi(\Pi^{\mho^\epsilon}_{\Lambda^n}\big|g'\big)=
\tfrac{1}{2}\chi\big(\Pi^{\rm IIodd}_n\big|g'\big)\Big(\,\underset{\sigma\in{\mathfrak A}_n}{\large\sum}\!\zeta_\gamma\big(\sigma(d)\big)-\underset{\sigma\in{\mathfrak S}_n\setminus{\mathfrak A}_n}{\large\sum}\!\zeta_\gamma\big(\sigma(d)\big)\Big)\cdot\delta^{{\mathfrak S}_n}_e(\sigma),
\end{array}
\end{eqnarray*}
where $\kappa+1$ is understood in $\mathrm{mod}\;2$, and $\delta^{{\mathfrak S}_n}_e$ is the delta function on ${\mathfrak S}_n$ supported on the identity element $e\in{\mathfrak S}_n$.
}
\vskip1em

{\it Proof.} 
For the second part, in Cases 1 and 2, take $s'_0$ such that $s_0\in {\mathfrak S}_{\boldsymbol \nu}\setminus{\mathfrak A}_{\boldsymbol \nu}$, then $s_0\gamma=\gamma$.    Note that, in Case 1 and Case 2 respectively, for $\sigma\in{\mathfrak S}_{\boldsymbol \nu}$,  
$$
\chi\big(\rho_{\Lambda^n}|s_0\sigma s_0^{\,-1}\big)=\chi\big(\rho_{\,{}^t\hspace{-.15ex}\Lambda^n}|\sigma\big),\quad \chi\big(\rho^{(\kappa)}_{\Lambda^n})|s_0\sigma s_0^{\,-1}\big)
=\chi\big(\rho^{(\kappa+1)}_{\Lambda^n})|\sigma\big).
$$
~\\[-7ex]
\hfill 
$\Box$\quad

\vskip1.2em

{\bf Corollary 7.20.}\; 
{\it 
In general, spin IRs\, $\Pi^{\mho+}_{\Lambda^n},\,\Pi^{\mho+}_{\Lambda^n,\,\kappa}$ and\, $\Pi^{\mho\hspace{.15ex}\epsilon}_{\Lambda^n}$ of spin type $\chi^{\rm II}$ 
are not equivalent to tensor products \,$\Pi^{(1)}\otimes\Pi^{(2)}$ of two IRs\, $\Pi^{(i)}$'s of\, $G'_n$.
}
\vskip1em

\section{\large Spin IRs of spin type $\chi^{\rm III}=(-1,\,1,-1)$}

Here we classify and construct spin IRs of 
$G(m,1,n)$ of spin type $\chi^{\rm III}=(-1,\,1,-1)$. 
We keep to the notations in \S 6.1, and  
$G'_n=R\big(G(m,1,n)\big)=\widetilde{D}^\vee_n\rtimes \widetilde{{\mathfrak S}}_n$ is abbreviated as $G'=U'\rtimes S'$, where $U'=\widetilde{D}^\vee_n=D^\wedge_n\times Z_3,\,D^\wedge_n=\langle z_2,\widehat{\eta}_j\;(j\in{\boldsymbol I}_n)\rangle$.  

\vskip.5em

{\it Step 1.} IRs of the group $U'$ of spin type $(z_2,z_3)\to (1,-1)$ are one-dimensional characters. A character of $U'$ is called {\it spin} or {\it non-spin} depending on whether 
$(z_2,z_3)\to (1,-1)$ or $(z_2,z_3)\to (1,1)$. 
The sets of such characters are denoted by $\widehat{U'}^{(1,-1)}$ and $\widehat{U'}^{(1,1)}$ respectively. 
Put 
$\omega:=e^{2\pi i/m}$ and $m^0=m/2$. 
For $\gamma\in \Gamma_n$, define characters $\zeta^\wedge_\gamma$ of $D^\wedge_n$, and $Z_\gamma$ of $U'=\widetilde{D}^\vee_n=D^\wedge_n\times Z_3$ as 
\begin{eqnarray}
\label{2017-07-15-1}
&&
\zeta^\wedge_\gamma(z_2^{\,b}\widehat{\eta}_1^{\,a_1}\widehat{\eta}_2^{\,a_2}\cdots\widehat{\eta}_n^{\,a_n}):=\omega^{\gamma_1a_1}\omega^{\gamma_2a_2}\cdots \omega^{\gamma_na_n},
\quad
Z_\gamma:=\zeta^\wedge_\gamma\cdot\mathrm{sgn}_{Z_3}\,.
\end{eqnarray}
Then\; $\widehat{U'}^{(1,-1)}=\big\{Z_\gamma\,;\,\gamma\in\Gamma_n\big\}$, and for ${\boldsymbol m}^0:=(m^0,m^0,\ldots,m^0)\in\Gamma_n$, $Z_{\gamma+{\boldsymbol m}^0}(d')=(-1)^{\mathrm{ord}(d')}Z_\gamma(d')\;(d'\in \widetilde{D}^\vee_n)$. 
 Moreover $\widehat{U'}^{(1,1)}=\big\{\zeta_\gamma\,;\,\gamma\in\Gamma_n\big\}$, where $\zeta_\gamma$ is naturally induced from $U=D_n$ up  to $U'$ as in (\ref{2017-11-20-1}).  Note that $\zeta_\gamma$ can be expressed by  
$(\zeta_j)_{j\in{\boldsymbol I}_n}$ 
with $\zeta_j(\widehat{\eta}_j):=\zeta_j(y_j)=\omega^{\gamma_j}\,,\; \zeta_j\in \widehat{T}\cong\langle \widehat{\eta}_j\rangle$, i.e., $\zeta_\gamma=(\zeta_j)_{j\in{\boldsymbol I}_n}$.

\vskip1.2em

{\bf Lemma 8.1.}\; 
{\it 
 The action of \,$S'=\widetilde{{\mathfrak S}}_n$ on spin character $Z_\gamma$ is given as follows: 
for $\gamma=(\gamma_j)_{j\in{\boldsymbol I}_n}\in\Gamma_n$, as equality in $\mathrm{mod}\;m$,
\begin{eqnarray}
\label{2011-09-09-2}
&&\quad\;
\left\{
\begin{array}{l}
r_iZ_\gamma = Z_{(r_i^{\,{\rm III}}\gamma)}\,,\;\quad 
r_i^{\,{\rm III}}\gamma:=s_i\gamma+{\boldsymbol m}^0,
\\
{\sigma'}^{\,{\rm III}}\gamma =\sigma\gamma +L(\sigma){\boldsymbol m}^0,\,\sigma\gamma=(\gamma_{\sigma^{-1}(j)})_{j\in{\boldsymbol I}_n}\;\big(\sigma'\in{\mathfrak S}_n,\,\sigma=\Phi_{\mathfrak S}(\sigma')\big).
\end{array}
\right. 
\end{eqnarray}
}

{\it Proof.} 
For $i\in{\boldsymbol I}_{n-1}$,\; 
$
r_i(Z_\gamma)(\widehat{\eta}_j)
=Z_\gamma\big(r_i^{\;-1}(\widehat{\eta}_j)\big)=
Z_\gamma\big(z_3\widehat{\eta}_{s_i(j)}\big)=-\omega^{\gamma_{s_i(j)}}
=-Z_{s_i\gamma}(\widehat{\eta}_j)=Z_{s_i\gamma+{\boldsymbol m}^0}(\widehat{\eta}_j)\;\;(j\in{\boldsymbol I}_n).
$
~\hfill
$\Box$\;\;
\vskip1.2em

Under this action of $S'$, a CSR of spin characters for $\widehat{U'}^{(1,-1)}/S'$ is given as follows. For $\gamma\in \Gamma^1_n$, choose $0\le p\le n$ as 
\begin{eqnarray}
\label{2012-04-19-02}
\quad
\left\{
\begin{array}{l}
\gamma_0:=0\le\gamma_1\le \gamma_2\le \ldots \le \gamma_p<m^0=m/2,
\\[.5ex]
m^0\le \gamma_{p+1}\le \ldots \le \gamma_n<m, 
\end{array}
\right.
\end{eqnarray}
and put \,$q=n\!-\!p$. Put $\gamma'_{p+j}:=\gamma_{p+j}-m^0\;(j\in{\boldsymbol I}_q)$, and  
\begin{eqnarray}
\label{2012-04-19-3}
\quad
\left\{
\begin{array}{l}
\Psi_0(\gamma):=(\gamma'_{p+1},\ldots,\gamma'_n,\gamma_1+m^0,\ldots,\gamma_p+m^0)\in\Gamma^1_n,
\\[1ex]
\Psi_1(\gamma):=(\gamma_{p+1},\,\ldots\,,\gamma_n,\gamma_1,\,\ldots\,,\gamma_p)\in \Gamma_n.
\end{array}
\right.
\end{eqnarray}
Then 
\;$\Psi_1(\gamma)\equiv\Psi_0(\gamma)+{\boldsymbol m}^0\;\;(\mbox{\rm elementwise }\mathrm{mod}\;m)$.
 For $\gamma,\gamma'\in\Gamma_n$, denote by $\gamma\sim\gamma'$ if  they are mutually conjugate under $\widetilde{{\mathfrak S}}_n$.

\vskip1.2em

{\bf Lemma 8.2.}\; 
{\it 
{\rm (i)}\; 
For any $\gamma\in\Gamma^1_n$\,,
\begin{eqnarray}
\label{2012-04-19-62}
\gamma\sim\Psi_1(\gamma)+pq{\boldsymbol m}^0\equiv \Psi_0(\gamma)+(pq+1){\boldsymbol m}^0\;\;(\mathrm{mod}\;m).
\end{eqnarray}
If $pq$ is even, then\; $\Psi_1(\gamma)\sim\gamma,$\, and  
if $pq$ is odd ($\therefore n$ even), then\; 
$\Psi_0(\gamma)\sim\gamma.$

{\rm (ii)}\; For $\gamma\in\Gamma^1_n$\,,
if $\gamma_j=\gamma_{j+1}\;(\exists j)$, then\;   
$r_j^{\,{\rm III}}\gamma=\gamma+{\boldsymbol m}^0$,\; 
$
\Psi_0(\gamma)\sim\gamma,\; 
\Psi_1(\gamma)\sim\gamma. 
$

{\rm (iii)}\; Let $\gamma\in\Gamma^1_n$. 
If all $\gamma_j$ are distinct, then $n\le m$.   
 If $\Psi_0(\gamma)=\gamma$, then\, $n$ is even.

}

\vskip1.2em

{\bf Proposition 8.3.}\; 
{\it 
The set \,$\widehat{U'}^{ (1,-1)}\!\big/S'$ of conjugacy classes of characters of spin type $(z_2,z_3)\to(1,-1)$ of\, $U'$ under the action of\, $S'=\widetilde{{\mathfrak S}}_n$ corresponds 1-1 way to $\Gamma_n/\!\sim$, and has a complete set of representatives $\big\{Z_\gamma\,;\,\gamma\in \Gamma^{{\rm III}}_n\big\}$, where the subset \,$\Gamma^{{\rm III}}_n\subset \Gamma_n$ is given as follows. Denote by $\preccurlyeq$ the lexicographic order in 
 $\Gamma_n$.

{\rm (i)}\; {\bf Case of $n$ odd.}  
\;$\Gamma^{{\rm III}}_n=\Gamma^{{\rm III},0}_n\bigsqcup\Gamma^{{\rm III},1}_n$\;\;\mbox{\it (if $n>m$, then\, $\Gamma^{{\rm III},1}_n=\emptyset)$}, with  
\begin{eqnarray}
\nonumber
&&
\Gamma^{{\rm III},0}_n:=\big\{\gamma\in\Gamma^1_n\;;\;\gamma\preccurlyeq\Psi_0(\gamma),\;\gamma_j=\gamma_{j+1}\;(\exists j)\big\},
\\
\nonumber
&&
\Gamma^{{\rm III},1}_n:=\big\{\gamma\in\Gamma^1_n\;;\;\mbox{\rm 
$\gamma_j$'s are all distinct}\,
\big\}. 
\end{eqnarray}

{\rm (i)}\; {\bf Case of $n$ even.}\; $
\Gamma^{{\rm III}}_n=\Gamma^{{\rm III},0}_n\bigsqcup\Gamma^{{\rm III},1}_n\bigsqcup\Gamma^{{\rm III},2}_n,\;\;
\Gamma^{{\rm III},0}_n=\Gamma^{{\rm III},0,\prec}_n\bigsqcup\Gamma^{{\rm III},0,=}_n$ \;
\\
\hspace*{21.1ex}
(if $n>m$, then \,$\Gamma^{{\rm III},1}_n=\Gamma^{{\rm III},2}_n=\emptyset$), 
\\[-4ex]
\begin{eqnarray*}
\nonumber
&&
\qquad
\Gamma^{{\rm III},0,\prec}_n:=
\big\{\gamma\in\Gamma^1_n\;;\;\gamma\prec\Psi_0(\gamma),\;\gamma_j=\gamma_{j+1}\;(\exists j)\big\},
\\
&&
\qquad
\Gamma^{{\rm III},0,=}_n:=
\big\{\gamma\in\Gamma^1_n\;;\;\gamma=\Psi_0(\gamma),\;\gamma_j=\gamma_{j+1}\;(\exists j)\big\},
\\
\nonumber
&&
\qquad
\Gamma^{{\rm III},1}_n:=\big\{\gamma\in\Gamma^1_n\;;\;\mbox{\rm 
$\gamma_j$ all distinct, $pq$ even}\,
\big\},
\\
\nonumber
&&
\qquad
\Gamma^{{\rm III},2}_n:=\big\{\gamma',\Psi_1(\gamma')\;;\;\gamma'\in\Gamma^1_n,\;\gamma'\preccurlyeq\Psi_0(\gamma'),\; 
\mbox{\rm $\gamma'_j$ all distinct, $pq$ odd}\,
\big\}.
\end{eqnarray*}
}

{\it Proof.}\; {\sc (a)\; Case of $\gamma_j=\gamma_{j+1}\;(\exists j)$.} By Lemmas 8.1 and 8.2. 

{\sc (b)\; 
Case of $\gamma_j$ all different.}\, For $\gamma\in\Gamma_n$, there exists  
$\sigma\in{\mathfrak S}_n$ such that $\gamma':=\sigma\gamma\in\Gamma^1_n$. So, depending on whether 
$\sigma=\Phi_{\mathfrak S}(\sigma')$ is even or odd, ${\sigma'}^{\,{\rm III}}\gamma\equiv\sigma\gamma$ or $\equiv\sigma\gamma+{\boldsymbol m}^0\;(\mathrm{mod}\;m)$. 
Accordingly, $\exists \gamma'\in\Gamma^1_n$ such that 
\;$\gamma\sim\gamma'$ \;or\; $\gamma\sim\gamma'+{\boldsymbol m}^0$. 
In the latter case, for $\gamma'+{\boldsymbol m}^0$, since $\gamma'=\Psi_0(\gamma^0)\;(\exists \gamma^0\in\Gamma^1_n)$, we have $\gamma'+{\boldsymbol m}^0=\Psi_0(\gamma^0)+{\boldsymbol m}^0=\Psi_1(\gamma^0)$. So, 
$\gamma$ is connected by $\sim$ to an element in $B=\{\gamma', \Psi_1(\gamma')\,;\,\gamma'\in\Gamma^1_n\}$. 
When $pq$ is even, inside $B$, $\Psi_1(\gamma^0)\sim \gamma^0$ and no other relation $\sim$ is possible. 
When $pq$ is odd, let $\tau'_{pq}\in\widetilde{{\mathfrak S}}_n$ be such element that $\tau_{pq}=\Phi_{\mathfrak S}(\tau'_{pq})$ is the exchange of two subintervals $[1,\ldots,p],\,[p+1, \ldots, n]$ of $[1,2,\ldots,n]$. Then only possible action of $\widetilde{{\mathfrak S}}_n$ inside $B$ is by $\tau'_{pq}$: ${\tau'_{pq}}^{\rm \!III}\gamma'=\Psi_0(\gamma')$, and ${\tau'_{pq}}^{\rm \!III}\big(\Psi_1(\gamma')\big)=\gamma'+{\boldsymbol m}^0$ which is in $B$ only when $\gamma'=\Psi_0(\gamma')$ and equals to $\Psi_1(\gamma')$ itself. 
\\~
\hfill
$\Box$\;

\vskip1em

{\it Step 2.}\; We determine the stationary subgroup in $S'=\widetilde{{\mathfrak S}}_n$ of spin character $\rho=Z_\gamma\in\widehat{U'}^{(1,-1)}\;(\gamma\in \Gamma^{{\rm III}}_n)$. 
As notations necessary in the following, we put 
\vskip.5em
$
\begin{array}{ll}
\widetilde{{\mathfrak S}}_n(\gamma):=\big\{\sigma'\in\widetilde{{\mathfrak S}}_n\;;\;{\sigma'}^{\,{\rm III}}\gamma=\gamma\big\},\quad &
{\mathfrak S}_n(\gamma):=\big\{\sigma\in{\mathfrak S}_n\,;\,\sigma\gamma=\gamma\big\},
\\[.5ex]
\widetilde{{\mathfrak A}}_n(\gamma):=\big\{\sigma'\in\widetilde{{\mathfrak A}}_n\,;\,{\sigma'}^{\,{\rm III}}\gamma\;(=\sigma\gamma)=\gamma\big\},\quad &
{\mathfrak A}_n(\gamma):=\big\{\sigma\in{\mathfrak A}_n\,;\,\sigma\gamma=\gamma\big\}.
\end{array}
$
\vskip1em

{\bf Proposition 8.4.}\; 
{\it  
For $\rho=Z_\gamma, \;\gamma=(\gamma_j) \in \Gamma^{{\rm III}}_n$, the stationary subgroup $S'([\rho])$ in $S'=\widetilde{{\mathfrak S}}_n$ of equivalence class $[\rho]$ is given as follows.  

{\rm (i)}\;\; {\bf Case of $n$ odd:}\hspace{20ex} 
$S'([\rho])=\widetilde{{\mathfrak A}}_n(\gamma)$, 
\\
\hspace*{15ex}
in particular, for $\gamma\in \Gamma^{\rm III,1}$, 
\hspace{1ex} 
$S'([\rho])=\widetilde{{\mathfrak A}}_n(\gamma)=Z_1$.

{\rm (ii)}\; {\bf Case of $n$ even:\;\; $n=2n'$.}

\quad
{\rm (ii-1)}\; For $\gamma\in\Gamma^{{\rm III},0,\prec}_n$ $\big(here\; \gamma\ne \Psi_0(\gamma)\big)$,\,  
$S'([\rho])=\widetilde{{\mathfrak A}}_n(\gamma)$. 

\quad
{\rm (ii-2)}\; For $\gamma\in \Gamma^{{\rm III},0,=}_n$ $\big(here\; \gamma=\Psi_0(\gamma)\big)$, 
\\
put 
 $\tau_0:=(1\;\,n'\!+\!1)(2\;\,n'\!+\!2)\cdots(n'\;n)\in{\mathfrak S}_n$\,, and take its preimage $\tau'_0\in\widetilde{{\mathfrak S}}_n$. Then
 
\quad 
 $(\tau'_0)^2=z_1^{\,n'(n'-1)/2}$,\, $\tau_0\gamma=\Psi_1(\gamma)\equiv\Psi_0(\gamma)+{\boldsymbol m}^0=\gamma+{\boldsymbol m}^0\,(\mathrm{mod}\;m)$. 
\vskip.3em

\qquad
{\rm (1)}\; If $n'$ is odd, then\,    
$\tau'_0$ is odd, and \hspace{.8ex} $S'([\rho])=\widetilde{{\mathfrak A}}_n(\gamma)\bigsqcup \tau'_0\,\widetilde{{\mathfrak A}}_n(\gamma)$.
\vskip.3em

\qquad
{\rm (2)}\; If $n'$ is even, then\,  
$\tau'_0$ is even, $\not\in S'([\rho])$, 
 and 
$\gamma_j=\gamma_{j+1}\;(\exists j)$, 
\\[.3ex]
\hspace{46ex}
$S'([\rho])=\widetilde{{\mathfrak A}}_n(\gamma)\bigsqcup (r_j\tau'_0)\,\widetilde{{\mathfrak A}}_n(\gamma)$.  
\vskip.2em

\quad{\rm (ii-3)}\; 
If $\gamma_j$'s are all distinct, then  
\hspace{5.5ex}  $S'([\rho])=\widetilde{{\mathfrak A}}_n(\gamma)= Z_1$.
}
\vskip1em

{\it Proof.}\; Assume 
$\sigma'\in\widetilde{{\mathfrak S}}_n$ be odd. Since ${\sigma'}^{\,{\rm III}}\gamma\equiv\sigma\gamma+{\boldsymbol m}^0$, to be 
$\sigma'\in S'([\rho])$ or to be ${\sigma'}^{\,{\rm III}}\gamma\equiv\gamma$, it is necessary to be $\sigma\gamma\equiv\gamma+{\boldsymbol m}^0$. Since $\sigma\gamma$ and $\gamma$ have the same components, the right hand side should be equal to $\Psi_1(\gamma)=\Psi_0(\gamma)+{\boldsymbol m}^0$. Hence $\Psi_0(\gamma)=\gamma$. Accordingly, if $\Psi_0(\gamma)\ne \gamma$, we have $S'([\rho])\subset\widetilde{{\mathfrak A}}_n$. This proves (ii-1).

If $n$ is odd, then\, $\Psi_0(\gamma)\ne \gamma$ by Lemma 8.2. So we get (i). 

(ii-2)\;\, $(\tau'_0)^2=z_1^{\,n'(n'-1)/2}$ follows from the next lemma.  Corresponding to $n'$ even or odd, $\tau'_0$ is even or odd, and then ${\tau'_0}^{\,{\rm III}}\gamma\equiv\tau_0\gamma\equiv\Psi_1(\gamma),$ or \;${\tau'_0}^{\,{\rm III}}\gamma\equiv \tau_0\gamma+{\boldsymbol m}^0\equiv \Psi_0(\gamma)$. 
Succeeding detailed discussions are omitted here. 
\hfill 
$\Box$\;
\vskip1.2em

{\bf Lemma 8.5.}\;(cf. [II, Lemma 7.2])\; 
{\it  
 For $j<k$, put
 
\quad 
$r_{jk} := (r_jr_{j+1}\cdots r_{k-2})r_{k-1}(r_{k-2}\cdots r_{j+1}r_j),\quad r_{kj} := r_{jk}$, 
\\[.5ex]
then, \,$(r_{jk})^2=e, \; 
 \Phi_{\mathfrak S}(r_{jk}) = s_{jk}:= (j\;\,k)$, transposition of $j,k$,\, and 
\begin{eqnarray}
 \label{2009-02-13-2-2}
&&
\qquad\quad 
 \left\{
 \begin{array}{ll}
  r_ir_{jk}r_i^{\;-1} = z_1r_{jk}&  (j,k \ne i, i+1),
 \qquad \qquad\quad 
 \\[.5ex]
 r_ir_{ji}r_i^{\;-1} = r_{j,i+1}, \quad 
 r_ir_{j,i+1}r_i^{\;-1} = r_{ji} & (j \ne i,i+1). 
 \hspace{11ex}\Box
   \end{array}
 \right.
\end{eqnarray}
 }

\vskip.5em 

 {\it Step 3.}\; 
Since $\rho=Z_\gamma$ is one-dimensional, we can put $J_\rho={\bf 1}$ (trivial). 

\vskip.7em

 {\it Step 4.}\; Let us choose the counter part $\pi^1$ for $\pi^0=\rho\cdot J_\rho=Z_\gamma\cdot{\bf 1}\;(\gamma\in\Gamma^1_n)$.
 \vskip.3em

  {\it Step 4-a.}\; {\bf Case of $S'([\rho])=\widetilde{{\mathfrak A}}_n(\gamma)$,} cf. Proposition 8.4 (i), (ii-1), (ii-3){\bf :} 
 \vskip.3em
In the case of distinct $\gamma_j$'s,\quad 
\;$S'([\rho])=\widetilde{{\mathfrak A}}_n(\gamma)=Z_1$, \;\;$\pi^1=\mathrm{sgn}_{Z_1}$.   
\vskip.2em

In the contrary case, 
with $\zeta^\vee_\gamma=(\zeta^\vee_j)_{j\in{\boldsymbol I}_n}$ in 
(\ref{2017-07-15-1}), put 
\begin{eqnarray}
\label{2015-02-06-1}
&&
\,
{\boldsymbol I}_n={\bigsqcup}_{\zeta\in\widehat{T}}\,I_{n,\zeta},\;I_{n,\zeta}=\{j\in{\boldsymbol I}_n\,;\,\zeta^\vee_j=\zeta\}, \;
{\boldsymbol \nu}=(n_\zeta)_{\zeta\in\widehat{T}},\;n_\zeta=|I_{n,\zeta}|,
\end{eqnarray}
and $\widetilde{{\mathfrak S}}_{\boldsymbol \nu}:=\Phi_{\mathfrak S}^{\;-1}\big({\prod}_{\zeta\in\widehat{T}}{\mathfrak S}_{I_{n,\zeta}}\big)$. Then $\prod_{\zeta\in\widehat{T}}{\mathfrak S}_{n_\zeta}\cong {\mathfrak S}_{\boldsymbol \nu}$ and 
\begin{eqnarray}
\label{2015-02-06-2} 
S'([\rho])=
\widetilde{{\mathfrak A}}_n(\gamma)=\widetilde{{\mathfrak S}}_{\boldsymbol \nu}\bigcap\widetilde{{\mathfrak A}}_n,
\quad 
\widetilde{{\mathfrak S}}_{\boldsymbol \nu}\cong\underset{\zeta\in\widehat{T}}{\widehat{*}}\widetilde{{\mathfrak S}}_{n_\zeta}.
\end{eqnarray}

A CSR of spin IRs of $\widetilde{{\mathfrak S}}_{\boldsymbol \nu}$ for ${\boldsymbol \nu}\in P(\widehat{T})$ is given by Theorem 4.4, and that of spin IRs of its subgroup $S'([\rho])=\widetilde{{\mathfrak A}}_n(\gamma)=\widetilde{{\mathfrak S}}_{\boldsymbol \nu}\bigcap\widetilde{{\mathfrak A}}_n$ by Theorem 4.7. 

 \vskip1em 
 
{\it Step 4-b.}\; {\bf Case of $S'([\rho])=\widetilde{{\mathfrak A}}_n(\gamma)\bigsqcup s'_0\widetilde{{\mathfrak A}}_n(\gamma)
$,}  cf. Proposition 8.4 (ii-2){\bf :}
\vskip.3em

Put 
$s'_0:=\tau'_0$ or $r_j\tau'_0$ (resp. $s_0:=\Phi_{\mathfrak S}(s_0')=\tau_0$ or $s_j\tau_0$) if $n'$ is odd or even.  
Then 
$s_0$ is odd, and $s_0^{\,2}=e$\;(unit element). Since 
$\Psi_0(\gamma)=\tau_0\gamma+{\boldsymbol m}^0\;(\mathrm{mod}\;m)$, there exists  $\delta=(\delta_j)_{j\in{\boldsymbol I}_{n'}}\in \Gamma^1_{n'}$ of half size $n'=n/2$ such that 
\begin{eqnarray}
\label{2015-02-08-3}
&&
\gamma=(\delta,\delta\!+\!{\boldsymbol m}^0),\;\;\delta\!+\!{\boldsymbol m}^0:=(\delta_1\!+\!m^0, \delta_2\!+\!m^0, \ldots, \delta_{n'}\!+\!m^0),
\end{eqnarray}
and \,$s_0\gamma=\tau_0\gamma=(\delta\!+{\boldsymbol m}^0,\delta)$. With $\zeta^{(a)}$ in (\ref{2015-02-08-11-8}), put $\widehat{T}^0=\{\zeta^{(a)}\,;\,0\le a<m^0\}$, ${\boldsymbol \nu}^0:=({\boldsymbol \nu}_\zeta)_{\zeta\in\widehat{T}^0}$, then \,$n_{\zeta^{(m^0)}\zeta}=n_\zeta\;(\zeta\in\widehat{T}^0)$, and 
\begin{eqnarray}
\label{2015-02-08-5}
\quad
\widetilde{{\mathfrak S}}_{\boldsymbol \nu}\cong\widetilde{{\mathfrak S}}^{(l)}_{{\boldsymbol \nu}^0}\,\widehat{*}\,\widetilde{{\mathfrak S}}^{(r)}_{{\boldsymbol \nu}^0},\;\;\;
 \widetilde{{\mathfrak A}}_n(\gamma)=\widetilde{{\mathfrak A}}_n\bigcap\widetilde{{\mathfrak S}}_{\boldsymbol \nu}\,.
\end{eqnarray}
Here the left hand side $\widetilde{{\mathfrak S}}^{(l)}_{{\boldsymbol \nu}^0}\cong\widetilde{{\mathfrak S}}_{{\boldsymbol \nu}^0}$ and the right hand side $\widetilde{{\mathfrak S}}^{(r)}_{{\boldsymbol \nu}^0}\cong\widetilde{{\mathfrak S}}_{{\boldsymbol \nu}^0}$\, act respectively on the left half ${\boldsymbol I}^{(l)}_n:=\{1,2,\ldots,n'\}$ and on the right half ${\boldsymbol I}^{(r)}_n:=(1+n',2+n',\ldots,n)$ of ${\boldsymbol I}_n$. 
We can and do identify ${\boldsymbol Y}^{\rm sh}_{\!n}({\boldsymbol \nu})$ with ${\boldsymbol Y}^{\rm sh}_{\!n'}({\boldsymbol \nu}^0)\times{\boldsymbol Y}^{\rm sh}_{\!n'}({\boldsymbol \nu}^0)$.

Put ${\cal S}:=\widetilde{{\mathfrak S}}_{\boldsymbol \nu},\,{\cal A}:=\widetilde{{\mathfrak A}}_n(\gamma),\,{\cal C}:=S'([\rho])={\cal A}\sqcup s'_0{\cal A}$. 
To get a CSR of spin IRs of 
 ${\cal C}$, we first start from a spin IR of ${\cal S}$, and restrict it to the subgroup ${\cal A}$ of index 2, and then induce its irreducible components up to ${\cal C}$:\, ${\cal S}\searrow {\cal A}\nearrow {\cal C}$. 
\vskip1em

{\bf Lemma 8.6.} (Restriction ${\cal S}\downarrow{\cal A}$)\; 
{\it 
Let $\Lambda^n=(\Lambda^{n',1}, \Lambda^{n',2})$ be an element of  \,${\boldsymbol Y}^{\rm sh}_{\!n}({\boldsymbol \nu})={\boldsymbol Y}^{\rm sh}_{\!n'}({\boldsymbol \nu}^0)\times{\boldsymbol Y}^{\rm sh}_{\!n'}({\boldsymbol \nu}^0)$. 
 Take spin IR $\tau_{\Lambda^n}\cong \tau_{(\Lambda^{n',1})}\,\widehat{*}\,\tau_{(\Lambda^{n',2})}$ of\, 
${\cal S}=
\widetilde{{\mathfrak S}}^{(l)}_{{\boldsymbol \nu}^0}\,\widehat{*}\,\widetilde{{\mathfrak S}}^{(r)}_{{\boldsymbol \nu}^0}$. 
 If $s(\Lambda^n)$ is even, then $\tau_{\Lambda^n}$ is self-associate and $\tau_{\Lambda^n}|_{\cal A}\cong \widetilde{\rho}_{(\Lambda^n,0)}\,\bigoplus\,\widetilde{\rho}_{(\Lambda^n,1)}\,,\,\widetilde{\rho}_{(\Lambda^n,0)}\not\cong\widetilde{\rho}_{(\Lambda^n,1)}.$  
 If $s(\Lambda^n)$ is odd, then $\tau_{\Lambda^n}$ is non self-associate and $\tau_{\Lambda^n}|_{\cal A}=: \widetilde{\rho}_{\Lambda^n}$ is irreducible. 

The set of spin IRs obtained in this way gives a CSR (= complete system of representatives) 
 of spin IRs of ${\cal A}$. 
}
\vskip1em 

{\bf Lemma 8.7.} (Action of $s'_0$ on ${\cal A}$)\;   
{\it 

 Put\, $s'_0:=\tau'_0$\, for $n'=n/2$ odd,\, and\, $s'_0:=r_j\tau'_0$\, for $n'$ even. 

{\rm (i)}\;  
Conjugacy action of $s'_0=\tau'_0$ on 
 $\sigma'=\sigma'_1\,\widehat{*}\,\sigma'_2 \in {\cal A}=\widetilde{{\mathfrak A}}(\gamma)\subset\widetilde{{\mathfrak S}}^{(l)}_{{\boldsymbol \nu}^0}\,\widehat{*}\,\widetilde{{\mathfrak S}}^{(r)}_{{\boldsymbol \nu}^0}$ is \;$s'_0\sigma'{s'_0}^{\,-1}=\sigma'_2\,\widehat{*}\,\sigma'_1$ (transposition of components). That of $s'_0=r_j\tau'_0$ with $r_j\in \widetilde{{\mathfrak S}}^{(l)}_{{\boldsymbol \nu}^0}$ is \,$s'_0\sigma'{s'_0}^{\,-1}=\sigma'_2\,\widehat{*}\,(r_j\sigma'_1r_j^{\,-1}).$
\vskip.3em

{\rm (ii)}\; For $\Lambda^n=(\Lambda^{n',1}, \Lambda^{n',2})\in{\boldsymbol Y}^{\rm sh}_{\!n}({\boldsymbol \nu})={\boldsymbol Y}^{\rm sh}_{\!n'}({\boldsymbol \nu}^0)\times{\boldsymbol Y}^{\rm sh}_{\!n'}({\boldsymbol \nu}^0)$, put \,$\overline{\tau_0}(\Lambda^n):=(\Lambda^{n',2}, \Lambda^{n',1})$.\,   
If $s(\Lambda^n)$ is odd, then $s'_0\big(\widetilde{\rho}_{\Lambda^n}\big)=\widetilde{\rho}_{\,\overline{\tau_0}(\Lambda^n)}$.\,  
If $s(\Lambda^n)$ is even, then $s'_0\big(\widetilde{\rho}^{\,(\kappa)}_{\Lambda^n}\big)\cong\widetilde{\rho}^{\,(\kappa+1)}_{\overline{\tau_0}(\Lambda^n)}$, 
where $\kappa+1=0,1$ is under $\mathrm{mod}\;2$, and in case $\overline{\tau_0}(\Lambda^n)\ne \Lambda^n$, the suffices $\kappa'=0,1$ of\, $\widetilde{\rho}^{\,(\kappa')}_{\,\overline{\tau_0}(\Lambda^n)}$ will be appropriately adjusted. 

}

\vskip1em

{\it Proof.}\; 
(i)\; For generators $r_jr_k$ of ${\cal A}$, we calculate explicitly 
$\tau'_0(r_jr_k){\tau'_0}^{\,-1}$ by means of Lemma 8.5. \; 
(ii) follows from (i).
\hfill
$\Box$\;

\vskip1.2em

{\bf Lemma 8.8.} 
(Inducing up ${\cal A}\uparrow{\cal C}$)\; 
{\it  
Let 
$\Lambda^n=(\Lambda^{n',1}, \Lambda^{n',2})\in{\boldsymbol Y}^{\rm sh}_{\!n'}({\boldsymbol \nu}^0)\times{\boldsymbol Y}^{\rm sh}_{\!n'}({\boldsymbol \nu}^0)$. 

{\rm (i)}\; Case of \,$\overline{\tau_0}(\Lambda^n)\ne \Lambda^n$. \; 

For $M^n\in{\mathscr Y}^{\mathfrak A}_n$ containing $\Lambda^n$, 
$T^{\mathfrak A}_{M^n}:=\mathrm{Ind}^{\cal C}_{\cal A}\,\widetilde{\rho}_{M^n}$\; 
is irreducible and equivalent to $T^{\mathfrak A}_{M^{\prime\,n}}$ with $M^{\prime\, n}$ obtained from $M^n$ by replacing $\Lambda^n$ by $\overline{\tau_0}(\Lambda^n)$.

{\rm (ii)}\; Case of $\overline{\tau_0}(\Lambda^n)= \Lambda^n$. 
\\
\quad
For $M^n\in{\mathscr Y}^{\mathfrak A}_n$ containing $\Lambda^n$, 
$T^{\mathfrak A}_{M^n}:=\mathrm{Ind}^{\cal C}_{\cal A}\,\widetilde{\rho}_{M^n}$\; 
is reducible and splits into two non-equivalent irreducible components: \;  
$T^{\mathfrak A}_{M^n}\cong T^{\mathfrak A}_{M^n;0}\bigoplus T^{\mathfrak A}_{M^n;1}.$
 
{\rm (iii)}\; The set of IRs of\, ${\cal C}=S'([\rho])=\widetilde{{\mathfrak A}}_n(\gamma)\bigsqcup s'_0\widetilde{{\mathfrak A}}_n(\gamma)
$ listed above gives a CSR of spin IRs of ${\cal C}$. 
}

\vskip1em

Start with 
$\Lambda^n=({\boldsymbol \lambda}^\zeta)_{\zeta\in\widehat{T}}\in {\boldsymbol Y}^{\rm sh}_{\!n}(\widehat{T})$, then $\Lambda^n\rightsquigarrow \gamma\in\Gamma^1_n,\, {\boldsymbol \nu}=(n_\zeta)_{\zeta\in\widehat{T}}\in P_n(\widehat{T})$, $n_\zeta=|{\boldsymbol \lambda}^\zeta|$. So we have naturally 
 $\rho=Z_\gamma$ and $\pi^0=\rho\cdot J_\rho=Z_\gamma\cdot {\bf 1}$. Its counter part $\pi^1$, a spin IR of 
$S'([\rho])$, is given by $\Lambda^n$ as explained in Theorem 4.7, Lemmas 8.6, and 8.7. 
Thus a CSR of spin IRs of $G'_n$ is given by $
\Pi(\pi^0,\pi^1)=\mathrm{Ind}^{G'_n}_{H'_n}(\pi^0\boxdot \pi^1)$ with $H'_n=U'\rtimes S'([\rho])$. 
 We list up them giving appropriate names.
\vskip1em

{\bf List 8.9.\; Spin IRs of $G'_n=R\big(G(m,p,n)\big)$ of spin type $\chi^{\rm III}$.}

\vskip.3em

{\rm (i)}\; {\bf Case of $n$\,odd.}  
\,$\Gamma^{{\rm III}}_n=\Gamma^{{\rm III},0}_n\bigsqcup\Gamma^{{\rm III},1}_n.$ 

$\bullet$ 
Case of 
$\gamma \in \Gamma^{{\rm III},0}_n=\big\{\gamma\in\Gamma^1_n;\gamma\preccurlyeq\Psi_0(\gamma),\gamma_j=\gamma_{j+1}\;(\exists j)\big\}
$:\quad\,  
$S'([\rho])=\widetilde{{\mathfrak A}}_n(\gamma)$. 

 \vskip.2em

\quad
 For $M^n\in{\mathscr Y}^{\mathfrak A}_{\!n}({\boldsymbol \nu})$ containing $\Lambda^n$,\hspace{1ex} 
$\pi^1=\widetilde{\rho}_{M^n}$,\; 
$\Pi^{\rm III}_{M^n}:=\mathrm{Ind}^{G'_n}_{H'_n}\big(Z_\gamma\hspace{.2ex}{\bf 1}_{S'([\rho])}\boxdot \widetilde{\rho}_{M^n}\big)$. 
 \vskip.2em

$\bullet$ 
Case of $
\gamma\in \Gamma^{{\rm III},1}_n=\big\{\gamma\in\Gamma^1_n;\mbox{\rm $\gamma_j$ are all distinct}\big\}  
$: 
\quad\;\;\, 
$S'([\rho])=\widetilde{{\mathfrak A}}_n(\gamma)=Z_1$, \;
 \vskip.2em

\quad 
 For $M^n=\Lambda^n$ quasi-trivial, \hspace{2ex}\;$\pi^1=\mathrm{sgn}_{Z_1},$ \hspace{.1ex} $\Pi^{\rm III}_{M^n}:=\mathrm{Ind}^{G'_n}_{H'_n}\,\big(Z_\gamma\hspace{.2ex}{\bf 1}_{S'([\rho])}\boxdot \mathrm{sgn}_{Z_1}\big)$.

\vskip.8em

{\rm (i)}\; {\bf Case of $n$ even.}\; $
\Gamma^{{\rm III}}_n=\Gamma^{{\rm III},0}_n\bigsqcup\Gamma^{{\rm III},1}_n\bigsqcup\Gamma^{{\rm III},2}_n$.
 
 \vskip.2em

$\bullet$ 
Case of $\gamma\in \Gamma^{{\rm III},0,\prec}_n=\!
\big\{\gamma\in\Gamma^1_n\,;\gamma\prec\Psi_0(\gamma),\,\gamma_j=\gamma_{j+1}\,(\exists j)\big\}$: 
\;\;$S'([\rho])=\widetilde{{\mathfrak A}}_n(\gamma)$. 
 \vskip.2em
 
\quad
 For $M^n\in{\mathscr Y}^{\mathfrak A}_{\!n}({\boldsymbol \nu})$ containing $\Lambda^n$,\hspace{1ex} 
$\pi^1=\widetilde{\rho}_{M^n}$,\; 
$\Pi^{\rm III}_{M^n}:=\mathrm{Ind}^{G'_n}_{H'_n}\big(Z_\gamma\hspace{.2ex}{\bf 1}_{S'([\rho])}\boxdot \widetilde{\rho}_{M^n}\big)$. 

\vskip.2em

$\bullet$ 
Case of $\gamma\in \Gamma^{{\rm III},0,=}_n=\!
\big\{\gamma\in\Gamma^1_n\,;\gamma=\Psi_0(\gamma),\,\gamma_j=\gamma_{j+1}\;(\exists j)\big\}$: \;
 \vskip.2em
\hspace{50ex}
$S'([\rho])=\widetilde{{\mathfrak A}}_n(\gamma)\bigsqcup s'_0\widetilde{{\mathfrak A}}_n(\gamma)$. \;
\\
\qquad
For $M^n\in{\boldsymbol Y}^{{\mathfrak A},{\rm sh}}_n({\boldsymbol \nu})$ containing $\Lambda^n$, 

\quad 
(1)\;  
If $\overline{\tau_0}(\Lambda^n)\ne \Lambda^n$, then \hspace{7ex}  
$\pi^1=T^{\mathfrak A}_{M^n}, 
\hspace{.7ex}
\Pi^{\rm III}_{M^n}:=\mathrm{Ind}^{G'_n}_{H'_n}\big(Z_\gamma\hspace{.2ex}{\bf 1}_{S'([\rho])}\boxdot T^{\mathfrak A}_{M^n}\big). 
$
\\[.3ex]
\qquad\qquad 
Let $\overline{\tau_0}(M^n)$ be obtained from $M^n$ by replacing $\Lambda^n$ by $\overline{\tau_0}(\Lambda^n)$, then 
\\
\qquad\hspace{3ex} 
$\Pi^{\rm III}_{M^n}\cong \Pi^{\rm III}_{\overline{\tau_0}(M^n)}$.
\vskip.2em

\quad 
(2)   
If $\overline{\tau_0}(\Lambda^n)\!=\! \Lambda^n$, then\,   
$\pi^1\!=\!T^{\mathfrak A}_{M^n;\,\iota}\,(\iota=0,1),\, 
\Pi^{\rm III}_{M^n;\,\iota}\!:=\!\mathrm{Ind}^{G'_n}_{H'_n}\big(Z_\gamma\hspace{.2ex}{\bf 1}_{S'([\rho])}\boxdot T^{\mathfrak A}_{M^n;\,\iota}\big). 
$
\vskip.3em

$\bullet$ 
Case of $\gamma\in \Gamma^{{\rm III},1}_n\bigsqcup 
\Gamma^{{\rm III},2}_n$: 
\hspace{21.6ex} 
$S'([\rho])=Z_1$,\; 

\quad 
 For $M^n=\Lambda^n$ quasi-trivial,\;\; $\pi^1=\mathrm{sgn}_{Z_1},$\; 
$\Pi^{\rm III}_{M^n}:=\mathrm{Ind}^{G'_n}_{H'_n}\,\big(Z_\gamma\hspace{.2ex}{\bf 1}_{S'([\rho])}\boxdot \mathrm{sgn}_{Z_1}\big)$. 
\vskip1.2em 

{\bf Theorem 8.10.} $\big($Spin IRs of spin type $\chi^{\rm III}=(-1,\,1,-1)\big)$\;  
 {\it 
 A complete set of representatives of spin IRs of the representation group 
$R\big(G(m,1,n)\big),\; 
n\ge 4,$ of spin type $\chi^{\rm III}$ is given as follows. 
\\[.5ex]
\indent
{\rm (i)}\; {\bf Case of $n$ odd:} \;
 $\mbox{\rm {\footnotesize spin}IR}^{\rm III}\big(G(m,1,n)\big)
=\big\{\Pi^{\rm III}_{M^n}\,;\, M^n \in{\mathscr Y}^{\mathfrak A}_n(\widehat{T}),\;\gamma\in\Gamma^{\rm III}_n\big\}.
$
\vskip.5em

{\rm (ii)} {\bf Case of $n$ even:}\; 
\\[.5ex]
\hspace*{7ex} 
$\mbox{\rm {\footnotesize spin}IR}^{\rm III}\big(G(m,1,n)\big)=
\big\{\Pi^{\rm III}_{M^n}\;;\;M^n\in{\mathscr Y}^{\mathfrak A}_n(\widehat{T}),\;\gamma\in \Gamma^{\rm III,0,\prec}\bigsqcup\Gamma^{\rm III,1}\bigsqcup\Gamma^{\rm III,2}\big\}
\\[.5ex]
\hspace*{19.2ex} 
\bigsqcup 
\big\{\Pi^{\rm III}_{M^n}\cong \Pi^{\rm III}_{\overline{\tau_0}(M^n)}\;;\,M^n\in{\mathscr Y}^{\mathfrak A}_n(\widehat{T}),\;\gamma\in \Gamma^{\rm III,0,=},\,\overline{\tau_0}(\Lambda^n)\ne \Lambda^n\big\}
\\[.5ex]
\hspace*{19.2ex} 
\bigsqcup 
\big\{\Pi^{\rm III}_{M^n;\,\iota}\;;\,M^n\in{\mathscr Y}^{\mathfrak A}_n(\widehat{T}),\;\gamma\in \Gamma^{\rm III,0,=},\,\overline{\tau_0}(\Lambda^n)= \Lambda^n,\;\iota=0,1\big\}.
$
}

\vskip1.2em
{\bf Note 8.10.}\; 
 Characters of spin IRs of 
$R\big(G(m,1,n)\big)$ of spin type $\chi^{\rm III}$ can be given by explicit calculations from the above classification and construction.  
In particular, in the case of $n$ odd and $\gamma\in\Gamma^{\rm III}_n$, and also in the case of $n$ even and $\gamma\in \Gamma^{\rm III,0,\prec}_n$, take $\Lambda^n\in {\boldsymbol Y}^{\rm sh}_{\!n}(\widehat{T})$ such that $\Lambda^n\rightsquigarrow \gamma$. If $\Lambda^n \in {\boldsymbol Y}^{\rm sh}_n(\widehat{T})^{\rm asym}$, then the calculation of character is easy.

\section{\large Spin IRs of spin type $\chi^{\rm V}=(1,-1,-1)$}

We\footnote{Spin IRs of spin type $\chi^{\rm IV}=(-1,\,1,\,1)$ and their characters have been given in \S 5.
} keep to the notations in \S 6.1, and  
abbreviate $G'_n=R\big(G(m,1,n)\big)=\widetilde{D}^\vee_n\rtimes \widetilde{{\mathfrak S}}_n$ as $G'=U'\rtimes S'$. 
For the case of spin type $\chi^{\rm V}$, we put ${\cal K}=\widehat{T}^0$. 
 
\vskip.7em

{\it Step 1.}\; 
By Lemma 6.2, we take,  
as a CSR of spin dual of $U'=\widetilde{D}^\vee_n$ of spin type $(z_2,z_3)\to (-1,-1)$ modulo the action of $S'$, i.e., a CSR for $\widehat{U'}^{(-1,-1)}/S'$, the following sets, 
 with $P_\gamma=\zeta_\gamma\,P^0$ and $P^\epsilon_\gamma= \zeta_\gamma P^\epsilon\;(\epsilon=+,-)$: 
\vskip.5em

\qquad 
for $n$ even, \quad $\rho=P_\gamma\cdot\mathrm{sgn}_{Z_3}\;\;(\gamma\in\Gamma^0_n\cap\Gamma^1_n),$ 

\qquad 
for $n$ odd, \hspace{2.3ex} $\rho=P^\epsilon_\gamma\cdot\mathrm{sgn}_{Z_3}\;\;(\gamma\in\Gamma^0_n\cap\Gamma^1_n,\,\epsilon=+,-)$.

\vskip.8em

{\it Step 2.}\; 
For the stationary subgroup $S'([\rho])$ in $S'$ of $[\rho]$, we refer to Lemma 6.4.

\vskip.7em

{\bf Proposition 9.1.}  
{\it 
The stationary subgroup $S'([\rho])$ in $S'$ of $[\rho]$ is given as follows. 

{\rm (i)\; Case of $n\ge 2$ even:}\;   
for $\rho=P_\gamma\cdot\mathrm{sgn}_{Z_3}$,\, 
$
S'([\rho])=\widetilde{{\mathfrak S}}_n(\gamma)\cong \widetilde{{\mathfrak S}}_{\boldsymbol \nu}\cong{\widehat{*}}_{\,\zeta\in\widehat{T}^0}\,\widehat{{\mathfrak S}}_{n_\zeta}\,.
$

{\rm (ii) Case of $n\ge 3$ odd:} \,     
for $\rho=P^\epsilon_\gamma\cdot\mathrm{sgn}_{Z_3}$,\, 
 $
S'([\rho])=\widetilde{{\mathfrak S}}_n(\gamma).$ 
}
\vskip1em

{\it Step 3.}\; 
Take a spin IRs $\rho$ listed in Step 1 with $\gamma\in\Gamma^0_n\cap\Gamma^1_n$.  
In case $n$ is even, for $\sigma'\in S'([\rho])$, an intertwining operator $J_\rho(\sigma')$ from $\rho$ to $\sigma'\rho$\, is given by 
\begin{eqnarray}
\label{2010-01-19-1-3}
\hspace*{12ex} 
\nabla_n(r_i)
={\textstyle  
\frac{1}{\sqrt{2}}}
\big(\widehat{Y}_i-\widehat{Y}_{i+1}\big)
 \;\;
(i\in{\boldsymbol I}_{n-1}) 
\;\;\;
{\rm with}\;\;\, 
\widehat{Y}_j=(-1)^{j-1}Y_j\;\;(j\in{\boldsymbol I}_n). 
\end{eqnarray}
For $n$ odd, we put\; $\nabla^+_n=\nabla_n$, and 
\begin{eqnarray}
\label{2010-01-19-2-3}
&&
\nabla^-_n(r_j)
=
-(i\widehat{Y}_n)\,\nabla_n(r_j) \,(i\widehat{Y}_n)^{-1}\;\; \big(j\in{\boldsymbol I}_{n-1},\;i=\sqrt{-1}\big).
\end{eqnarray}
\vskip.2em

{\bf Proposition 9.2.}\; 

(i)\;\; Case of $n\ge 2$ even.\;   
{\it 
Let $\rho=P_\gamma\cdot\mathrm{sgn}_{Z_3}\;\big(\gamma\in\Gamma^0_n\cap\Gamma^1_n\big)$.\,  
Intertwining operators in (\ref{2016-05-12-1}) for $\sigma'\in S'([\rho])$ is $J_\rho(\sigma')=\nabla_n(\sigma')$\, up to a scalar  multiple. 
\\[.5ex]
\indent
{\rm (ii)\; Case of $n\ge 3$ odd.}\;  
Let $\rho=P^\epsilon_\gamma\cdot\mathrm{sgn}_{Z_3}\;(\gamma\in\Gamma^0_n\cap\Gamma^1_n,\,\epsilon=+,-)$. \; 
Intertwining operators for $\sigma'\in S'([\rho])$ is  $J_\rho(\sigma')=\nabla^\epsilon_n(\sigma')$\, up to a scalar multiple.

}

\vskip1em

{\it Proof.}\; 
(i)\; Denote $A^{-1}BA$ by  $\iota(A)B$. We have for $i\in{\boldsymbol I}_{n-1},\,j\in{\boldsymbol I}_n$, 
\vskip.4em 
\,$\iota\big(\nabla_n(r_i)\big)\widehat{Y}_j=-\widehat{Y}_{s_i(j)}$\;\; $\therefore\,\iota\big(\nabla_n(\sigma')\big)\widehat{Y}_j=\mathrm{sgn}(\sigma)\widehat{Y}_{\sigma^{-1}(j)}\;\big(\sigma'\in\widetilde{{\mathfrak S}}_n, \sigma=\Phi_{\mathfrak S}(\sigma')\big)$,\; 
\\[1ex]
and so\; $\iota\big(\nabla_n(\sigma')\big)P^0= \sigma^{\rm \prime\,V}(P^0).$
Moreover\; $\zeta_\gamma=\sigma^{\rm \prime\,V}(\zeta_{\sigma^{-1}\gamma})$\, because $\sigma^{\prime\,-1}(\widehat{\eta}_j)=\mathrm{sgn}(\sigma)\widehat{\eta}_{\sigma^{-1}(j)}$, whence $\sigma^{\rm \prime\,V}(\zeta_\gamma)=\zeta_{\sigma\gamma}$.

(ii)\, We prove for $\rho=P^-_\gamma$. For $r_i$ with \;$i<n-1$, the calculation is the same as above, and for $r_{n-1}$,
\vskip-2em   
\begin{gather*}
\iota\big(\nabla^-_n(r_{n-1})\big)\widehat{Y}_j=
\left\{ 
\begin{array}{ll}
-\widehat{Y}_j\,,\quad& j<n-1,\hspace*{6ex}
\\
\,\,\widehat{Y}_n\,,\quad &j=n-1,
\\
\widehat{Y}_{n-1}\,,\quad&j=n, 
\end{array}
\right.
\\
\therefore\;\; 
\iota\big(\nabla^-_n(r_{n-1})\big)P^-_\gamma(\widehat{\eta}_{n-1})=
\iota\big(\nabla^-_n(r_{n-1})\big)\big(\zeta_{\tau_n\gamma}(\widehat{\eta}_{n-1})\widehat{Y}_{n-1}\big)\hspace{8ex}\mbox{\rm (cf. (\ref{2014-12-20-11-8}))}
\\
\hspace*{36ex}
=\omega^{\gamma_{n-1}}\,\widehat{Y}_n
=P^-_{s_{n-1}\gamma}\big(r_{n-1}(\widehat{\eta}_{n-1})\big),
\\[.5ex]
\hspace*{6ex}  
\iota\big(\nabla^-_n(r_{n-1})\big)P^-_\gamma(\widehat{\eta}_n)=
\omega^{\gamma_n+m^0}\widehat{Y}_{n-1}=P^-_{s_{n-1}\gamma}\big(r_{n-1}(\widehat{\eta}_n)\big).\hspace{17.5ex} \Box
\end{gather*}

Now, $\pi^0=\rho\cdot J_\rho$ is a spin IR of 
$H'_n=U'\rtimes S'([\rho])$ of spin type $\chi^{\rm V}$ given as 
\begin{eqnarray}
\label{2014-04-21-20-3}
&&\;\;
\left\{
\begin{array}{lll}
{\rm for}\;\;\rho=P_\gamma\cdot\mathrm{sgn}_{Z_3}\;\;(n\;{\rm even}), \qquad&
\pi^0=(P_\gamma\cdot\mathrm{sgn}_{Z_3})\cdot \nabla_n|_{S'([\rho])}\,, 
\\[.5ex]
{\rm for}\;\;\rho=P^\epsilon_\gamma\cdot\mathrm{sgn}_{Z_3}\;\;(n\;{\rm odd}), \quad 
&
\pi^0=(P^\epsilon_\gamma\cdot\mathrm{sgn}_{Z_3})\cdot \nabla^\epsilon_n|_{S'([\rho])}\,.
\end{array}
\right.
\end{eqnarray}

\vskip.2em

 {\it Step 4.}\; As a counter part $\pi^1$ of $\pi^0$, it should be a spin representation of $S'([\rho])$, to get non-spin property $z_1\to 1$ of $\pi=\pi^0\boxdot \pi^1$. As in (\ref{2014-04-24-52}), we put 
 \begin{eqnarray}
\label{2014-04-24-52-99}
&&
\left\{
 \begin{array}{ll}
 \Pi^{\rm I0}_n
:= \big(P^0\cdot\mathrm{sgn}_{Z_3}\big)\cdot \nabla_n\,,\;\;\;& \mbox{\rm for $n\ge 4$ even,}
\\[.5ex]
\Pi^{\rm I\epsilon}_n
:= \big(P^\epsilon\cdot\mathrm{sgn}_{Z_3}\big)\cdot \nabla^\epsilon_n\,,\qquad & \mbox{\rm for $n\ge 5$ odd, $\epsilon\in\{+,-\}$.} 
\end{array}
\right.
\end{eqnarray}

\vskip.1em

{\bf Case of $n$ even.}\; 
Take 
$\rho=(P_\gamma\cdot\mathrm{sgn}_{Z_3})\cdot \nabla_n|_{S'([\rho])}\;(\gamma\in\Gamma^0_n\cap\Gamma^1_n)$. Then, from $\gamma$, we determine ${\boldsymbol \nu}=(n_\zeta)_{\zeta\in \widehat{T}^0}\in P_n(\widehat{T}^0)$ canonically and then take 
$\Lambda^n=({\boldsymbol \lambda}^\zeta)_{\zeta\in \widehat{T}^0}\in{\boldsymbol Y}^{\rm sh}_{\!n}({\boldsymbol \nu})$. Conversely, if $\Lambda^n\in {\boldsymbol Y}^{\rm sh}_{\!n}(\widehat{T}^0)$ is given, then \,$\Lambda^n\rightsquigarrow {\boldsymbol \nu}\in P_n(\widehat{T}^0),\,\gamma\in\Gamma^0_n\cap\Gamma^1_n$\,.  

Take an element $M^n=\Lambda^n$ or $M^n=(\Lambda^n,\mu),\,\mu=\pm 1$ of 
${\mathscr Y}_n(\widehat{T}^0)$ in \S 4.3. Then we have a spin IR of $S'([\rho])\cong \widetilde{{\mathfrak S}}_{\boldsymbol \nu}=\widehat{*}_{\zeta\in\widehat{T}^0}\widetilde{{\mathfrak S}}_{n_\zeta}$ as 
\;$
\pi^1:=\tau_{M^n}$\,. 
 Then, with  
$\pi^0=\rho\cdot J_\rho$\,, the IR $\pi=\pi^0\boxdot\pi^1$ of $H'_n=U'\rtimes S'([\rho])$ becomes of spin type $\chi^{\rm V}=(1,-1,-1)$, cancelling out spin properties of $J_\rho$ and $\pi^1$ with respect to $z_1$. Put 
\begin{eqnarray}
\label{2014-05-09-21}
&&\quad 
\Pi^{\rm V}_{M^n}:=\mathrm{Ind}^{G'_n}_{H'_n} \big(\pi^0\boxdot\pi^1\big)
=\mathrm{Ind}^{G'_n}_{H'_n}\Big(\big((P_\gamma\cdot\mathrm{sgn}_{Z_3})\cdot \nabla_n\big)\big|_{H'_n}\boxdot \tau_{M^n}\Big).
\end{eqnarray}
Since 
\;$\big((P_\gamma\cdot\mathrm{sgn}_{Z_3})\cdot \nabla_n\big)\big|_{H'_n}\boxdot \tau_{M^n}\cong \big((P^0\cdot\mathrm{sgn}_{Z_3})\cdot \nabla_n\big)\big|_{H'_n}\otimes (\zeta_\gamma\,{\bf 1}_{Z_3}\cdot\tau_{M^n}\big), 
$
we have a tensor product expression, with $\Pi^{\rm I0}_n$ in (\ref{2014-04-24-52}),  
\begin{eqnarray}
\label{2017-07-16-1}
&&
\Pi^{\rm V}_{M^n}\cong  \Pi^{\rm I0}_n\otimes\Pi^{\rm IV}_{M^n}\;\;\;{\rm with}\;\;\;\Pi^{\rm IV}_{M^n}:= \mathrm{Ind}^{G'_n}_{H'_n}(\zeta_\gamma\,{\bf 1}_{Z_3}\cdot\tau_{M^n}),
\end{eqnarray}
where $\Pi^{\rm IV}_{M^n}\;\big(M^n\in{\mathscr Y}(\widehat{T}^0)\big)$ is a spin IR of spin type $\chi^{\rm IV}=(-1,\,1,\,1)$. 

\vskip1em

{\bf Case of $n$ odd.}\; 
Start with $M^n\in{\mathscr Y}_n(\widehat{T}^0),$ containing $\Lambda^n=({\boldsymbol \lambda}^\zeta)_{\zeta\in\widehat{T}^0}\in {\boldsymbol Y}^{\rm sh}_{\!n}(\widehat{T}^0)$, then $\Lambda^n\rightsquigarrow {\boldsymbol \nu},\, \gamma$.  Here ${\boldsymbol \nu}=(n_\zeta)_{\zeta\in{\cal K}}\in P_n(\widehat{T}^0)$, $n_\zeta=|{\boldsymbol \lambda}^\zeta|\;(\zeta\in\widehat{T}^0),$ and $\gamma\in\Gamma^0_n\cap\Gamma^1_n$.  
Take $\rho=P^+_\gamma\cdot\mathrm{sgn}_{Z_3}$ or $P^-_\gamma\cdot\mathrm{sgn}_{Z_3}$, and put $\pi^0=\rho\cdot J_\rho$. 
Also take as $\pi^1$ a spin IR $\tau_{M^n}$ of  
$S'([\rho])\cong\widetilde{{\mathfrak S}}_{\boldsymbol \nu}$\,. Then we get IR $\pi=\pi^0\boxdot\pi^1$ of $H'_n=U'\rtimes S'([\rho])$ of spin type $\chi^{\rm V}$. Put, with $\epsilon=+,-$,  
\begin{eqnarray}
\label{2014-05-09-31}
&&\Pi^{\rm V}_{M^n}=\mathrm{Ind}^{G'_n}_{H'_n}\Big(\big((P^\epsilon_\gamma\cdot\mathrm{sgn}_{Z_3})\cdot \nabla^\epsilon_n|_{S'([\rho])}\big)\boxdot \tau_{M^n}\Big)\quad\big(M^n\in{\mathscr Y}_n(\widehat{T}^0)\big).
\end{eqnarray}
Since 
\;$\big((P^\epsilon_\gamma\cdot\mathrm{sgn}_{Z_3})\cdot \nabla^\epsilon_n\big)\big|_{H'_n}\boxdot \tau_{M^n}\cong \big((P^\epsilon\cdot\mathrm{sgn}_{Z_3})\cdot \nabla^\epsilon_n\big)\big|_{H'_n}\otimes (\zeta_\gamma\,{\bf 1}_{Z_3}\cdot\tau_{M^n}\big), 
$
we have, with reference to \S 4 for spin type $\chi^{\rm IV}=(-1,\,1,\,1)$ and with $\Pi^{\rm I\epsilon}_n\;(\epsilon=\pm)$ in (\ref{2014-04-24-52}),  
\begin{eqnarray}
\label{2017-07-16-1-9}
&&
\Pi^{\rm V\epsilon}_{M^n}\cong  \Pi^{\rm I\,\epsilon}_n\otimes\Pi^{\rm IV}_{M^n}\;\;\;{\rm with}\;\;\;\Pi^{\rm IV}_{M^n}:= \mathrm{Ind}^{G'_n}_{H'_n}(\zeta_\gamma\,{\bf 1}_{Z_3}\cdot\tau_{M^n}).
\end{eqnarray}

\vskip.2em

{\bf Theorem 9.3.}\; 
{\it 
Assume $n\ge 4,\,m$ even. A complete set of representatives of spin IRs of 
 $R\big(G(m,1,n)\big)$ of spin type $\chi^{\rm V}=(1,-1,-1)$ is given as follows.
\\[.5ex]
\indent
{\rm (i)}\;\; Case of $n$ even: \;
 $\mbox{\rm {\footnotesize spin}IR}^{\rm V}\big(G(m,1,n)\big)
=\Big\{\Pi^{\rm V}_{M^n}\,;\, M^n \in{\mathscr Y}^{\mathfrak A}_n(\widehat{T}^0)\Big\}.
$
\vskip.2em

{\rm (ii)}\; Case of $n$ odd:\; 
 $\mbox{\rm {\footnotesize spin}IR}^{\rm V}\big(G(m,1,n)\big)
=\Big\{\Pi^{\rm V\epsilon}_{M^n}\,;\,\epsilon\in\{+,-\},\, M^n \in{\mathscr Y}^{\mathfrak A}_n(\widehat{T}^0)\Big\}.
$
}

\section{\large Spin IRs and their characters of spin type $\chi^{\rm VI}$}

We abbreviate the notation $G'_n=R\big(G(m,1,n)\big)=\widetilde{D}^\vee_n\rtimes \widetilde{{\mathfrak S}}_n$ as $G'=U'\rtimes S'$. 

\vskip1.2em 
{\bf 10.1. Classification of IRs of spin type $\chi^{\rm VI}=(1,-1,\,1)$.}
\vskip.5em  

{\it Step 1.}\; 
 By Lemma 7.2 we know the following.  
 \vskip.8em 
 
 {\bf Lemma 10.1.} {\it A CSR of \,$\widehat{U'}^{(-1,1)}/S'$  is given by   
\vskip.3em

\quad  
for $n\ge 2$ even,\; 
$\big\{P_\gamma\hspace{.2ex}{\bf 1}_{Z_3}\,;\,\gamma\in\Gamma^0_n\cap\Gamma^1_n\big\}$,
\vskip.3em

\quad 
for $n\ge 3$ odd,\;  
$
\big\{P^+_\gamma{\bf 1}_{Z_3}\,;\gamma\in\Gamma^0_n\cap\Gamma^1_n\big\}\bigsqcup 
\big\{P^-_\gamma{\bf 1}_{Z_3}\,;\gamma\in \Gamma^0_n\cap\Gamma^1_n,\,\gamma_j\;\mbox{\rm all different}\big\}. $
}

\vskip1em

{\it Step 2.}\; 
By Proposition 7.4 we know also the following.
\vskip.7em

{\bf Lemma 10.2.}\; 
{\it 
The stationary subgroup $S'([\rho])$ of spin IRs $\rho= P_\gamma{\bf 1}_{Z_3},\;P^+_\gamma{\bf 1}_{Z_3}$ and $P^-_\gamma{\bf 1}_{Z_3}$ $(\gamma\in\Gamma^0_n\cap\Gamma^1_n)$ of\,  
$U'=\widetilde{D}^\vee_n$ is given as follows: 
with ${\cal K}=\widehat{T}^0$, 

For $n\ge 2$ even,
\hspace{3ex}
$S'\big([P_\gamma\hspace{.2ex}{\bf 1}_{Z_3}]\big)=\Phi_{\mathfrak S}^{\;-1}\Big({\prod}_{\zeta\in{\cal K}}{\mathfrak S}_{I_{n,\zeta}}\Big)= \widetilde{{\mathfrak S}}_{\boldsymbol \nu}.$

 For $n\ge 3$ odd,
\hspace{4.2ex}
$
S'\big([P^\pm_\gamma{\bf 1}_{Z_3}]\big)=\widetilde{{\mathfrak A}}_n\bigcap\Phi_{\mathfrak S}^{\;-1}\Big({\prod}_{\zeta\in{\cal K}}{\mathfrak S}_{I_{n,\zeta}}\Big)= \widetilde{{\mathfrak A}}_n\cap\widetilde{{\mathfrak S}}_{\boldsymbol \nu}\,.
$
}

\vskip1em

{\it Step 3.}\; By Lemma 6.5, we have the following. 
  
\vskip.7em

{\bf Proposition 10.3.}\;  
{\it 
{\rm (i)}\; {\bf Case of $n\ge 2$ even.}\; 
For 
$\rho=P_\gamma\hspace{.2ex}{\boldsymbol 1}_{Z_3}\;\big(\gamma\in\Gamma^0_n\cap\Gamma^1_n\big)$ and  $\sigma'\in S'([\rho])\cong\widetilde{{\mathfrak S}}_{\boldsymbol \nu}$, an intertwining operator (\ref{2016-05-12-1}) from $\rho$ to ${\sigma'}\rho$ is, up to a scalar multiple,  
\,$
J_\rho(\sigma')=\nabla^{\rm II}_n(\sigma')$.\,  
This gives a spin representation of\, $S'([\rho])$.
\\[.5ex]
\indent
 {\rm (ii)}\, {\bf Case of $n\ge 3$ odd.}\;  
Let $\rho=P^\epsilon_\gamma\hspace{.2ex}{\boldsymbol 1}_{Z_3}$ be as in Lemma 10.1. 
An intertwining operator $J_\rho(\sigma')$ for $\sigma'\in S'([\rho])\cong\widetilde{{\mathfrak A}}_n\cap\widetilde{{\mathfrak S}}_{\boldsymbol \nu}$ is, up to a scalar multiple, \,$J_\rho(\sigma')=\mho^\epsilon_n(\sigma')$ in (\ref{2014-04-29-11}). 
 This gives a spin representation of\, $S'([\rho])\cong\widetilde{{\mathfrak A}}_n\bigcap\widetilde{{\mathfrak S}}_{\boldsymbol \nu}$. 
In the special case where $S'([\rho])=Z_1$, we understand  \,$S'([\rho])=Z_1$ and $J_\rho=\mathrm{sgn}_{Z_1}\!\cdot I.$ 
}

\vskip1em

{\it Step 4.}\;  
The spin type of IR $\pi^0=\rho\cdot J_\rho$ of $H'_n=U'\rtimes S'([\rho])$ is $(z_1,z_2,z_3)\to (-1,-1,\,1)=\chi^{\rm II}$.  
To get the spin type $\chi^{\rm VI}=(1,-1,\,1)$ for $\pi=\pi^0\boxdot \pi^1$, we should take, as the counter part of $\pi^0$, a spin IR $\pi^1$ of $S'([\rho])$, that is, of spin type $z_1\to -1$.

\vskip.5em

{\bf Case of $n$ even.}\; 
We apply the results in {\bf 4.3}, in particular Theorem 4.4. Take an  
$M^n\in{\mathscr Y}_n(\widehat{T}^0)$.  Then it contains $\Lambda^n$, which gives $\gamma\in\Gamma^0_n\cap\Gamma^1_n$ and ${\boldsymbol \nu}\in P_n(\widehat{T}^0)$, i.e., $\Lambda^n\rightsquigarrow \gamma, {\boldsymbol \nu}$. So we have  $\rho=P_\gamma\hspace{.2ex}{\bf 1}_{Z_3}$, and $\pi^0=\rho\cdot J_\rho$ and $\pi^1=\tau_{M^n}$. 
Thus \,$\pi=\pi^0\boxdot\pi^1=
\big(P_\gamma\hspace{.2ex}{\bf 1}_{Z_3}\cdot\nabla^{\rm II}_n|_{\widetilde{{\mathfrak S}}_{\boldsymbol \nu}}\big)\boxdot \tau_{M^n}$\, is of spin type $\chi^{\rm VI}=(1,-1,\,1)$. Put  
\begin{eqnarray}
\label{2015-02-14-11}
\Pi^{\rm VI}_{M^n}:=\mathrm{Ind}^{G'_n}_{H'_n}(\pi^0\boxdot\pi^1).
\end{eqnarray}
Since\; $\big(P_\gamma\hspace{.2ex}{\bf 1}_{Z_3}\cdot\nabla^{\rm II}_n|_{\widetilde{{\mathfrak S}}_{\boldsymbol \nu}}\big)\boxdot \tau_{M^n}\cong\big(P^0\hspace{.2ex}{\bf 1}_{Z_3}\cdot\nabla^{\rm II}_n\big)\big|_{H'_n}\otimes\big((\zeta_\gamma\hspace{.2ex}{\bf 1}_{Z_3})\cdot\tau_{M^n}\big)
$, we have by Lemma 2.3, with $\Pi^{\rm II0}_n$ in  (\ref{2017-07-16-21}) and \,$\Pi^{\rm IV}_{M^n}=\mathrm{Ind}^{G'_n}_{H'_n}(\zeta_\gamma\hspace{.2ex}{\bf 1}_{Z_3}\cdot\tau_{M^n})$ in (\ref{2015-02-22-12-8}),    
\begin{eqnarray}
\label{2017-07-16-11}
\Pi^{\rm VI}_{M^n}\cong \Pi^{\rm II0}_n\otimes \Pi^{\rm IV}_{M^n}\,.
\end{eqnarray}

\vskip.2em

{\bf Case of $n$ odd.}\; 
We apply Theorem 4.7.  
Start with $\rho=P^+_\gamma\hspace{.2ex}{\bf 1}_{Z_3}$ or $\rho=P^-_\gamma\hspace{.2ex}{\bf 1}_{Z_3}$ in Lemma 10.1. Here $\gamma$ gives ${\boldsymbol \nu}=(n_\zeta)_{\zeta\in \widehat{T}^0}\in P_n(\widehat{T}^0)$. Take a spin IR $\pi^1$ of $S'([\rho])=\widetilde{{\mathfrak A}}_n\bigcap\widetilde{{\mathfrak S}}_{\boldsymbol \nu}$, then we get an IR $\pi^0\boxdot\pi^1$ of $H'_n=U'\rtimes S'([\rho])$. 
\vskip1em 

{\bf Lemma 10.4.}\; 
{\it 
IR $\pi^0\boxdot \pi^1$ is given explicitly according to three cases as follows:  
\begin{eqnarray}
\label{2015-02-15-1}
\nonumber
\hspace*{.5ex}
\quad 
\begin{array}{lll}
\mbox{\rm Case 1.}&
\big((P^+_\gamma\hspace{.2ex}{\bf 1}_{Z_3})\cdot\mho^+_n|_{S'([\rho])}\big)\boxdot \widetilde{\rho}_{M^n}\,, & 
M^n:=(\Lambda^n,\kappa) \in {\mathscr Y}^{\mathfrak A}_n(\widehat{T}^0)^{\rm ev,non}\,,
\\[1ex]
\mbox{\rm Case 2.}&
\big((P^\epsilon_\gamma\hspace{.2ex}{\bf 1}_{Z_3})\cdot\mho^\epsilon_n|_{S'([\rho])}\big)\boxdot\widetilde{\rho}_{\Lambda^n},& 
M^n:=(\epsilon,\Lambda^n)\in\{+,-\}\!\times\!{\mathscr Y}^{\mathfrak A}_n(\widehat{T}^0)^{\rm ev,tri}, 
\\[1ex]
\mbox{\rm Case 3.}&
\big((P^+_\gamma\hspace{.2ex}{\bf 1}_{Z_3})\cdot\mho^+_n|_{S'([\rho])}\big)\boxdot\widetilde{\rho}_{\Lambda^n}\,, & 
M^n:=\Lambda^n\in{\boldsymbol Y}^{\rm sh}_{\!n}(\widehat{T}^0)^{\rm odd}\,.
\end{array}
\end{eqnarray}
}

\vskip.2em

Note that $S'([\rho])\subset\widetilde{{\mathfrak A}}_n,\;H'_n\subset G^{\prime\hspace{.15ex}{\mathfrak A}}_n=U'\rtimes\widetilde{{\mathfrak A}}_n$, and spin representations $\mho^\epsilon_n\;(\epsilon=+,-)$ of $\widetilde{{\mathfrak A}}_n$ cannot be extended to $\widetilde{{\mathfrak S}}_n$ as seen in [II, \S 8.3, pp.162--163].

Inducing $\pi^0\boxdot\pi^1$ above up to $G'_n=R\big(G(m,1,n)\big)$, we get a spin IR 
$$
\Pi^{\rm VI}_{M^n}:=\mathrm{Ind}^{G'_n}_{H'_n} (\pi^0\boxdot\pi^1).
$$

\vskip.2em

{\bf Proposition 10.5.}\; 
{\it 
According to the three cases in Lemma 10.4, there hold the following \lq\lq\,tensor product\,+\,inducing up'' relations: 
\begin{eqnarray*}
\label{2017-07-17-1}
\!\!
\begin{array}{ll}
\mbox{\rm Case 1.}\quad\;&
\Pi^{\rm VI}_{M^n}\cong 
\mathrm{Ind}^{G'_n}_{G^{\prime\hspace{.15ex}{\mathfrak A}}_n}\big(\Pi^{\mho^+\!,\hspace{.15ex}{\mathfrak A}}_n\otimes \mathrm{Ind}^{G^{\prime\hspace{.15ex}{\mathfrak A}}_n}_{H'_n}(\zeta_\gamma\hspace{.2ex}{\bf 1}_{Z_3}\cdot\widetilde{\rho}_{M^n})\big),
\\[.7ex]
\mbox{\rm Case 2.}&
\Pi^{\rm VI}_{M^n}\cong 
\mathrm{Ind}^{G'_n}_{G^{\prime\hspace{.15ex}{\mathfrak A}}_n}\big(\Pi^{\mho^\epsilon\!,\hspace{.15ex}{\mathfrak A}}_n\otimes\mathrm{Ind}^{G^{\prime\hspace{.15ex}{\mathfrak A}}_n}_{H'_n}\!( \zeta_\gamma\hspace{.2ex}{\bf 1}_{Z_3}\cdot\widetilde{\rho}_{\Lambda^n})\big),
\\[.7ex]
\mbox{\rm Case 3.}&
\Pi^{\rm VI}_{M^n}\cong 
\mathrm{Ind}^{G'_n}_{G^{\prime\hspace{.15ex}{\mathfrak A}}_n}\big(\Pi^{\mho^+\!,\hspace{.15ex}{\mathfrak A}}_n\otimes\mathrm{Ind}^{G^{\prime\hspace{.15ex}{\mathfrak A}}_n}_{H'_n}\!( \zeta_\gamma\hspace{.2ex}{\bf 1}_{Z_3}\cdot\widetilde{\rho}_{\Lambda^n}\big)).
\end{array}
\end{eqnarray*}
}

{\it Proof.}\; 
Note that \;$\big((P^\epsilon_\gamma\hspace{.2ex}{\bf 1}_{Z_3})\cdot\mho^\epsilon_n|_{S'([\rho])}\big)\boxdot\widetilde{\rho}_{M^n}
\cong \big((P^\epsilon\hspace{.2ex}{\bf 1}_{Z_3})\cdot \mho^\epsilon_n\big)\big|_{H'_n}\otimes(\zeta_\gamma\hspace{.2ex}{\bf 1}_{Z_3})\cdot\widetilde{\rho}_{M^n}
$\, for $M^n=(\Lambda^n,\kappa)$ or $M^n=\Lambda^n$. 
Then we obtain by Lemma 2.3 the asserted relations, with spin IRs 
$\Pi^{\mho^\epsilon,\hspace{.15ex}{\mathfrak A}}_n\;(\epsilon =\pm)$ of $G^{\prime\hspace{.15ex}{\mathfrak A}}_n$ in (\ref{2017-07-16-21}). 
\hfill 
$\Box$\quad\;

\vskip1em

{\bf Theorem 10.6.}\; 
{\it 
Assume $n\ge 4$, $m$ even. A complete set of representatives (=\,CSR) of spin IRs of $R\big(G(m,1,n)\big)$ with spin type $\chi^{\rm VI}=(1,-1,\,1)$ 
is given by the set $\mbox{\rm {\footnotesize spin}IR}^{\rm VI}\big(G(m,1,n)\big)$ below:
\\[.5ex]
\indent
{\rm (i)}\; Case $n$ even:\;
$\mbox{\rm {\footnotesize spin}IR}^{\rm VI}\big(G(m,1,n)\big)
\\
\hspace*{21.3ex}
=\big\{\Pi^{\rm VI}_{M^n}=\mathrm{Ind}^G_H\Big(\big( (P_\gamma\hspace{.2ex}{\bf 1}_{Z_3})\cdot\nabla^{\rm II}_n\big|_{\widetilde{{\mathfrak S}}_{\boldsymbol \nu}}\big)\boxdot \tau_{M^n}\Big)\,;\, M^n \in{\mathscr Y}_n(\widehat{T}^0)\big\}.$
\\[.5ex]
\indent
{\rm (ii)}\; Case $n$ odd, and $n>m^0=m/2$:\;

\;\qquad
$\mbox{\rm {\footnotesize spin}IR}^{\rm VI}\big(G(m,1,n)\big)=
\big\{\Pi^{\rm VI}_{M^n}\,;\, 
M^n\in {\mathscr Y}^{\mathfrak A}_n(\widehat{T}^0)^{\rm ev,non}\bigsqcup{\boldsymbol Y}^{\rm sh}_{\!n}(\widehat{T}^0)^{\rm odd}\big\}$.

\;{\rm (iii)}\; Case $n$ odd, and $n\le m^0=m/2$: \;

\;\qquad 
$\mbox{\rm {\footnotesize spin}IR}^{\rm VI}\big(G(m,1,n)\big)= 
\big\{\Pi^{\rm VI}_{M^n}\,;\,M^n\in {\mathscr Y}^{\mathfrak A}_n(\widehat{T}^0)^{\rm ev,non}\bigsqcup{\boldsymbol Y}^{\rm sh}_{\!n}(\widehat{T}^0)^{\rm odd} \big\} 
\\
\hspace*{28ex}
\bigsqcup\big\{
\Pi^{\rm VI}_{M^n}\,;\,M^n=(\epsilon,\Lambda^n)\in\{+,-\}\!\times\!{\mathscr Y}^{\mathfrak A}_n(\widehat{T}^0)^{\rm ev,tri}\big\}$.

}

\vskip1.2em

{\bf 10.2. Irreducible characters of spin type $\chi^{\rm VI}$.}
\vskip.5em

{\bf 10.2.1. Case of $n\ge 4$ even.}\; 
We use spin IR $\tau_{M^n}$ of 
$\widetilde{{\mathfrak S}}_{\boldsymbol \nu}$ with  
\begin{eqnarray*}
\;\;
\begin{array}{ll}
\mbox{\rm for $d(\Lambda^n)$ even,} & 
 M^n=\Lambda^n\in {\boldsymbol Y}^{\rm sh}_{\!n}({\boldsymbol \nu})\,\subset\, {\mathscr Y}_n(\widehat{T}^0)^{\rm ev}$, and $\tau_{M^n}=\tau_{\Lambda^n},
\\[.3ex]
\mbox{\rm for $d(\Lambda^n)$ odd,} &  
M^n=(\Lambda^n,\mu)\in{\boldsymbol Y}^{\rm sh}_{\!n}({\boldsymbol \nu})^{\rm odd}\times\{\pm 1\},$ and $\tau_{M^n}=
\mathrm{sgn}^{(1-\mu)/2}\cdot \tau_{\Lambda^n}. 
\end{array}
\end{eqnarray*}

{\bf Theorem 10.7.}\; For $n\ge 4$ even, irreducible characters $\chi\big(\Pi^{\rm VI}_{M^n}\big|\,\cdot)$ are given as 
\begin{eqnarray*}
\chi\big(\Pi^{\rm VI}_{M^n}\big|\,g')=\chi\big(\Pi^{\rm II0}_n\big|\,g'\big)\cdot\chi\big(\Pi^{\rm IV}_{M^n}\big|\,g'\big)\quad(g'\in G'_n). 
\end{eqnarray*}
 
\vskip.2em

{\bf 10.2.2. Case of $n\ge 5$ odd.}\; 
The following is direct from 
 Lemma 10.5. 
\vskip1em

{\bf Lemma 10.8.}\; 
There hold the following relations: 
\begin{eqnarray*}
\;\,
\begin{array}{ll}
\mbox{\rm Case 1}.\qquad\quad &
\chi\big(\Pi^{\rm VI}_{M^n}\big)=\mathrm{Ind}^{G'_n}_{G^{\prime\hspace{.15ex}{\mathfrak A}}_n}\Big(\chi\big(\Pi^{\mho^+\!,{\mathfrak A}}_n\big)\cdot\chi\big(\mathrm{Ind}^{G^{\prime\hspace{.15ex}{\mathfrak A}}_n}_{H'_n}\!( \zeta_\gamma\hspace{.2ex}{\bf 1}_{Z_3}\cdot\widetilde{\rho}_{M^n})\big)\Big),  
\\[.5ex]
\mbox{\rm Case 2}.\qquad &
\chi\big(\Pi^{\rm VI}_{M^n}\big)= \mathrm{Ind}^{G'_n}_{G^{\prime\hspace{.15ex}{\mathfrak A}}_n}\Big(\chi\big(\Pi^{\mho^\epsilon\!,{\mathfrak A}}_n\big)\cdot\chi\big(
\mathrm{Ind}^{G^{\prime\hspace{.15ex}{\mathfrak A}}_n}_{H'_n}\!(\zeta_\gamma\hspace{.2ex}{\bf 1}_{Z_3}\cdot\widetilde{\rho}_{\Lambda^n})\big)\Big),  
\\[.5ex]
\mbox{\rm Case 3}.\qquad &
\chi\big(\Pi^{\rm VI}_{M^n}\big)= \mathrm{Ind}^{G'_n}_{G^{\prime\hspace{.15ex}{\mathfrak A}}_n}\Big(\chi\big(\Pi^{\mho^+\!,{\mathfrak A}}_n\big)\cdot\chi\big(\mathrm{Ind}^{G^{\prime\hspace{.15ex}{\mathfrak A}}_n}_{H'_n}\!(\zeta_\gamma\hspace{.2ex}{\bf 1}_{Z_3}\cdot\widetilde{\rho}_{\Lambda^n})\big)\Big).   
\end{array}
\end{eqnarray*}

\vskip.5em

According to Cases 1, 2 and 3 in Lemmas 10.4, 10.5 and 10.8, we have the following.
\vskip1em 

{\bf Theorem 10.9.}\;
{\it 
Let\, $n\ge 5$ be odd, and\, $\Pi_n$ be one of\, $\Pi^{\rm VI}_{M^n}$ in Cases 1, 2 and 3 in Lemmas 10.4 and 10.8. Then its character satisfies the following. 

{\rm (i)}\,  For $g'\not\in G^{\prime\hspace{.15ex}{\mathfrak A}}_n$,\, 
$\chi(\Pi_n|g')=0$. 

{\rm (ii)}\,  Let $g'=(d',\sigma')\in G^{\prime\hspace{.15ex}{\mathfrak A}}_n=U'\rtimes\widetilde{{\mathfrak A}}_n$. 

In case\, $\mathrm{ord}(d')=\mathrm{ord}(d)\equiv 0\;(\mathrm{mod}\;2)$,
\begin{eqnarray*}
\!\!\begin{array}{ll}
\mbox{\rm 1.}\;
M^n\!=\!(\Lambda^n,\kappa) \in {\mathscr Y}^{\mathfrak A}_n(\widehat{T}^0)^{\rm ev,non},&\!\!\!
\chi\big(\Pi^{\rm VI}_{M^n}\big|g')\!=\!\chi\big(\Pi^{\mho^+\!,{\mathfrak A}}_n\big|g'\big)\chi\big(\mathrm{Ind}^{G'_n}_{H'_n}\!( \zeta_\gamma{\bf 1}_{Z_3}\,\widetilde{\rho}_{M^n})\big|g'\big),
\\[1ex]
\mbox{\rm 2.}\;
M^n\!=\!(\epsilon,\Lambda^n)\!\in\!\{\pm\}\!\times\!{\mathscr Y}^{\mathfrak A}_n(\widehat{T}^0)^{\rm ev,tri}, \!\!\!&
\chi\big(\Pi^{\rm VI}_{M^n}\big|g')\!=\!\chi\big(\Pi^{\mho^\epsilon\!,{\mathfrak A}}_n\big|g'\big)\chi\big(\mathrm{Ind}^{G'_n}_{H'_n}\!(\zeta_\gamma{\bf 1}_{Z_3}\,\widetilde{\rho}_{\Lambda^n})\big|g'\big),
\\[1ex]
\mbox{\rm 3.}\;
M^n\!=\!\Lambda^n\in{\boldsymbol Y}^{\rm sh}_{\!n}(\widehat{T}^0)^{\rm odd},&\!\!
\chi\big(\Pi^{\rm VI}_{M^n}\big|g')\!=\!\chi\big(\Pi^{\mho^+\!,{\mathfrak A}}_n\big|g'\big)\chi\big(\mathrm{Ind}^{G'_n}_{H'_n}\!( \zeta_\gamma{\bf 1}_{Z_3}\,\widetilde{\rho}_{\Lambda^n})\big|g'\big).
\end{array}
\end{eqnarray*}

In case $\mathrm{ord}(d')\equiv 1$, with an $s'_0\in\widetilde{{\mathfrak S}}_n\setminus\widetilde{{\mathfrak A}}_n$, in Cases 1, 2 and 3 respectively, 
\begin{eqnarray*}
\begin{array}{l}
\!\chi\big(\Pi^{\rm VI}_{M^n}\big|g')\!=\!\chi\big(\Pi^{\mho^+\!,{\mathfrak A}}_n\big|g'\big)\Big(\chi\big(\mathrm{Ind}^{G'_n}_{H'_n}\!( \zeta_\gamma{\bf 1}_{Z_3}\,\widetilde{\rho}_{M^n})\big|g'\big)\!-\!\chi\big(\mathrm{Ind}^{G'_n}_{H'_n}\!( \zeta_\gamma{\bf 1}_{Z_3}\,\widetilde{\rho}_{M^n})\big|s'_0g'{s'_0}^{\,-1}\big)\Big),
\\[.9ex]
\chi(\Pi^{\rm VI}_{M^n}|g')
\!=\!\chi\big(\Pi^{\mho^\epsilon\!,{\mathfrak A}}_n\big|g'\big)\Big(\chi\big(\mathrm{Ind}^{G'_n}_{H'_n}\!( \zeta_\gamma{\bf 1}_{Z_3}\,\widetilde{\rho}_{\Lambda^n})\big|g'\big)\!-\!\chi\big(\mathrm{Ind}^{G'_n}_{H'_n}\!( \zeta_\gamma{\bf 1}_{Z_3}\,\widetilde{\rho}_{\Lambda^n})\big|s'_0g'{s'_0}^{\,-1}\big)\Big),
\\[.9ex]
\chi(\Pi^{\rm VI}_{M^n}|g')=\chi\big(\Pi^{\mho^+\!,{\mathfrak A}}_n\big|g'\big)\Big(\chi\big(\mathrm{Ind}^{G'_n}_{H'_n}\!( \zeta_\gamma{\bf 1}_{Z_3}\,\widetilde{\rho}_{\Lambda^n})\big|g'\big)\!-\!\chi\big(\mathrm{Ind}^{G'_n}_{H'_n}\!( \zeta_\gamma{\bf 1}_{Z_3}\,\widetilde{\rho}_{\Lambda^n})\big|s'_0g'{s'_0}^{\,-1}\big)\Big).
\end{array}
\end{eqnarray*}
}

\vskip.2em

{\bf Note 10.1.}\; Recall spin IRs of $\widetilde{{\mathfrak A}}_n\cap\widetilde{{\mathfrak S}}_{{\boldsymbol \nu}}$ in \S 4.4, and spin IRs 
$\Pi^{\rm IV}_{\Lambda^n}$ of $G'_n=R\big(G(m,1,n)\big)$ in \S 5. In Cases 1, 2 and 3 respectively, we have  
\begin{gather*}
\mathrm{Ind}^{G'_n}_{H'_n}\!( \zeta_\gamma\hspace{.2ex}{\bf 1}_{Z_3}\cdot\widetilde{\rho}_{M^n})=\mathrm{Ind}^{G'_n}_{H'_n}\!( \zeta_\gamma\hspace{.2ex}{\bf 1}_{Z_3}\cdot\widetilde{\rho}^{(\kappa)}_{\Lambda^n})
\cong \Pi^{\rm IV}_{\Lambda^n},
\\
 \mathrm{Ind}^{G'_n}_{H'_n}\!(\zeta_\gamma\hspace{.2ex}{\bf 1}_{Z_3}\cdot\widetilde{\rho}_{\Lambda^n})
\cong \Pi^{\rm IV}_{\Lambda^n},
\quad 
\mathrm{Ind}^{G'_n}_{H'_n}\!(\zeta_\gamma\hspace{.2ex}{\bf 1}_{Z_3}\cdot\widetilde{\rho}_{\Lambda^n})
\cong \Pi^{\rm IV}_{\Lambda^n}\oplus\big(\mathrm{sgn}_{{\mathfrak S}}\,\Pi^{\rm IV}_{\Lambda^n}\big). 
\end{gather*}

\section{\large Spin IRs of spin type $\chi^{\rm VII}=(1,1,-1)$}

We classify and construct spin IRs of 
$G(m,1,n)$ of spin type $\chi^{\rm VII}=(1,\,1,-1)$ in our standard way. Denote $G'_n=R\big(G(m,1,n)\big)=\widetilde{D}^\vee_n\rtimes \widetilde{{\mathfrak S}}_n$ and $G_n=G(m,1,n)=D_n\rtimes {\mathfrak S}_n$ as $G'=U'\rtimes S'$ and $G=U\rtimes S$ respectively. 
\vskip.6em

{\it Step 1.}\, The spin dual $\widehat{U'}^{(1,-1)}$ of $U'$, of spin type $(z_2,z_3)\to (1,-1)$, is the same as for the case of spin type $\chi^{\rm III}=(-1,\,1,-1)$ in \S 8, and the $S'$-action on $\widehat{U'}^{(1,-1)}$ is the same as in Lemma 8.1.   
So we quote Proposition 8.3 for $\widehat{U'}^{(1,-1)}/S'$.

\vskip.6em

{\it Step 2.}\; For stationary subgroups $S'([\rho])$, we quote Proposition 8.4.

\vskip.6em

 {\it Step 3.}\;  Intertwining operator $J_\rho(\sigma')\equiv {\bf 1}$ (trivial one-dimensional operator). 
\vskip.6em

{\it Step 4.} IR $\pi^1$ of $S'([\rho])$ should be of spin type $z_1\to 1$ or IR of $S([\rho])= S'([\rho])/Z_1$.
\vskip.5em

{\it Step 4-a.} {\bf Case of $\rho=Z_\gamma\;(\gamma\in\Gamma^1_n)$, $S'([\rho])=\widetilde{{\mathfrak A}}_n(\gamma)$.} 

This case corresponds to Proposition\;8.4 (i), (ii-1) and (ii-3). 
Note that $S([\rho])={\mathfrak A}_n(\gamma)={\mathfrak A}_n\bigcap{\mathfrak S}_{\boldsymbol \nu}$ with ${\boldsymbol \nu}\in P_n(\widehat{T})$, where ${\boldsymbol \nu}\leftrightarrow  \gamma$ as in \S 6.1, Step 2.  
To describe the dual of  
${\mathfrak A}_n\bigcap{\mathfrak S}_{\boldsymbol \nu}$, we utilize Notation 7.6 with ${\cal K}=\widehat{T}$, for example, ${\boldsymbol Y}_{\!n}^{\mathfrak{A}}(\widehat{T})\;:=\;{\boldsymbol Y}_{\!n}^{\mathfrak{A}}(\widehat{T})^{\rm asym}\bigsqcup {\boldsymbol Y}_{\!n}^{\mathfrak{A}}(\widehat{T})^{\rm sym}\bigsqcup {\boldsymbol Y}_{\!n}(\widehat{T})^{\rm tri}.$ 
Let $\Lambda^n=({\boldsymbol \lambda}^\zeta)_{\zeta\in\widehat{T}}\in{\boldsymbol Y}_{\!n}(\widehat{T})$ and $\Lambda^n\rightsquigarrow \gamma,{\boldsymbol \nu}$. 
 \vskip1em

\hspace*{-1ex}$\bullet$ 
In the case of distinct $\gamma_j$'s,\quad 
\vskip.5em

$S([\rho])={\mathfrak A}_n(\gamma)=\{e\}$, \;\;$\pi^1=\pi_{\Lambda^n}={\bf 1}, \;\Lambda^n=({\boldsymbol \lambda}^\zeta)_{\zeta\in\widehat{T}}\in {\boldsymbol Y}^{\mathfrak A}_{\!n}(\widehat{T})^{\rm tri}$.

\vskip.7em

\hspace*{-1ex}$\bullet$ In the contrary case, 
put 
\begin{eqnarray}
\label{2015-02-06-1-9}
&&
\,
{\boldsymbol I}_n={\bigsqcup}_{\zeta\in\widehat{T}}\,I_{n,\zeta},\;I_{n,\zeta}=\{j\in{\boldsymbol I}_n\,;\,\zeta^{(\gamma_j)}=\zeta\}, \;
{\boldsymbol \nu}=(n_\zeta)_{\zeta\in\widehat{T}},\;n_\zeta=|I_{n,\zeta}|,
\end{eqnarray}
then\,  ${\prod}_{\zeta\in\widehat{T}}{\mathfrak S}_{I_{n,\zeta}}\cong \prod_{\zeta\in\widehat{T}}{\mathfrak S}_{n_\zeta}=:{\mathfrak S}_{\boldsymbol \nu}$\, and 
\;$S([\rho])={\mathfrak A}_n(\gamma)={\mathfrak A}_n\cap{\mathfrak S}_{\boldsymbol \nu}.$ 

\vskip.3em

A CSR of the dual of ${\mathfrak S}_{\boldsymbol \nu}$ is given by means of  Lemma 7.5. The situation on restricting them to the subgroup ${\mathfrak A}_n\cap{\mathfrak S}_{\boldsymbol \nu}$ of index 2 is clarified and then a CSR for the dual of ${\mathfrak A}_n\cap{\mathfrak S}_{\boldsymbol \nu}$ is given in Lemma 7.7.

\vskip1em
 {\bf Lemma 11.1.}\; 
 {\it 
 For a ${\boldsymbol \nu}=(n_\zeta)_{\zeta\in\widehat{T}}\in P_n(\widehat{T})$, put\, $G^0={\mathfrak S}_{\boldsymbol \nu}$ and\, $H^0={\mathfrak A}_n\cap G^0={\mathfrak A}_n\cap{\mathfrak S}_{\boldsymbol \nu}$.  
 
 {\rm (i)}\, For $\Lambda^n=({\boldsymbol \lambda}^\zeta)_{\zeta\in\widehat{T}}\in {\boldsymbol Y}_{\!n}(\widehat{T})$, put ${\boldsymbol \nu}=(n_\zeta)_{\zeta\in\widehat{T}}\in P_n(\widehat{T}),\,n_\zeta=|{\boldsymbol \lambda}^\zeta|$.  
For IR $\pi_{\Lambda^n}$ of\, $G^0=\large\prod_{\zeta\in\widehat{T}}{\mathfrak S}_{n_\zeta}$, \;$\mathrm{sgn}\cdot \pi_{\Lambda^n}\cong \pi_{({}^t\Lambda^n)}$ with\, ${}^t(\Lambda^n)={}^t\Lambda^n:=({}^t{\boldsymbol \lambda}^\zeta)_{\zeta\in\widehat{T}}$. 

\vskip.3em
 
{\rm (ii)}\, Assume $n_\zeta\ge 2\;(\exists \zeta\in\widehat{T})$.  
Then $[G^0:H^0]=2$. 
 \vskip.3em 
 \noindent
{\rm (a)}\, When ${}^t\Lambda^n=\Lambda^n$, i.e., ${}^t{\boldsymbol \lambda}^\zeta ={\boldsymbol \lambda}^\zeta\;(\forall \zeta\in\widehat{T})$,\; 
$\pi_{\Lambda^n}|_{H^0}\cong \rho^{(0)}_{\Lambda^n}\oplus \rho^{(1)}_{\Lambda^n},\; 
\rho^{(0)}_{\Lambda^n}\not\cong \rho^{(1)}_{\Lambda^n}$\,, 
\vskip.4em 

\noindent
{\rm (b)} When ${}^t\Lambda^n\ne\Lambda^n$, i.e., ${}^t{\boldsymbol \lambda}^\zeta \ne{\boldsymbol \lambda}^\zeta\;(\exists \zeta\in\widehat{T})$, 
\; 
 $\rho_{\Lambda^n}:=\pi_{\Lambda^n}|_{H^0} \cong (\pi_{\,{}^t\Lambda^n})|_{H^0}$\, is irreducible. 
\vskip.4em 

{\rm (iii)}\, For a fixed ${\boldsymbol \nu}=(n_\zeta)_{\zeta\in\widehat{T}}\in P_n(\widehat{T}),\, \exists n_\zeta\ge 2,$ the set of IRs listed in (a) and (b) gives a complete set of representatives of IRs of $H^0$. 
}
\vskip1em

Put $\rho=Z_\gamma$. If\, $S([\rho])={\mathfrak A}_n(\gamma)={\mathfrak A}_n\cap {\mathfrak S}_{\boldsymbol \nu}$, then we can take the counter part $\pi^1$ of $\pi^0=Z_\gamma\cdot {\bf 1}_{S'([\rho])}$\, from a CSR of linear IRs of $S([\rho])$ given above.

\vskip1em 
 
{\it Step 4-b.}\; {\bf Case of $S'([\rho])=\widetilde{{\mathfrak A}}_n(\gamma)\bigsqcup s'_0\widetilde{{\mathfrak A}}_n(\gamma)
$}  (Proposition 8.4 (ii-2), $n=2n'$){\bf :}
\vskip.3em
 
Here, mapping down to the level of $G_n=G(m,1,n)$, we have 
$s_0=\tau_0$ (in Case (1)), $s_j\tau$ (in Case (2)) odd, and $s_0^{\,2}=e$ (unit element) (loc.cit.). 
Since 
$\Psi_0(\gamma)=\tau_0\gamma+{\boldsymbol m}^0\;(\mathrm{mod}\;m)$, there exists $\delta=(\delta_j)_{j\in{\boldsymbol I}_{n'}}\in \Gamma^1_{n'}$, of half size $n'=n/2$, such that 
\begin{eqnarray}
\label{2015-02-08-3-2}
&&
\quad
\gamma=(\delta,\,\delta\!+\!{\boldsymbol m}^0),\;\;\delta\!+\!{\boldsymbol m}^0:=(\delta_1\!+\!m^0, \delta_2\!+\!m^0, \ldots, \delta_{n'}\!+\!m^0). 
\end{eqnarray}
Then, \;$s_0\gamma=\tau_0\gamma=(\delta\!+{\boldsymbol m}^0,\delta),\;s_0^{\;{\rm VII}}\gamma=(\delta,\,\delta\!+{\boldsymbol m}^0)=\gamma $.\, 
Put $\widehat{T}^0=\{\zeta^{(a)}\;;\;0\le a<m^0\}$ and   ${\boldsymbol \nu}^0:=({\boldsymbol \nu}_\zeta)_{\zeta\in\widehat{T}^0}$. Since $n_{\zeta^{(a+m^0)}}=n_{\zeta^{(a)}}\;(0\le a<m^0)$ here, we have 
\begin{eqnarray}
\label{2015-02-08-5-2}
\quad
{\mathfrak S}_{\boldsymbol \nu}\cong{\mathfrak S}^{(l)}_{{\boldsymbol \nu}^0}\times{\mathfrak S}^{(r)}_{{\boldsymbol \nu}^0},\;\;\;
 {\mathfrak A}_n(\gamma)={\mathfrak A}_n\bigcap{\mathfrak S}_{\boldsymbol \nu}\,,
\end{eqnarray}
where the above factor ${\mathfrak S}^{(*)}_{{\boldsymbol \nu}^0}$ with $*=l$ (resp. $=r$) acts on the left (resp. right)  half of ${\boldsymbol I}_n$, denoted as ${\boldsymbol I}^{(l)}_n:=\{1,2,\ldots,n'\}$ (resp. ${\boldsymbol I}^{(r)}_n:=\{1+n',2+n',\ldots,n\}$). 

Put ${\cal S}:={\mathfrak S}_{\boldsymbol \nu},\;{\cal A}:={\mathfrak A}_n(\gamma)={\mathfrak A}_n\cap{\cal S},\;{\cal C}:=S([\rho])={\cal A}\bigsqcup s_0{\cal A}$. To get a CSR of spin IRs of ${\cal C}$, first start from a CSR for ${\cal S}$, and restrict it on the subgroup ${\cal A}$ of index 2, and then induce it up to the upper group ${\cal C}$ of index 2\,:\; ${\cal S}\searrow{\cal A}\nearrow{\cal C}$.  
\vskip1em

{\bf Lemma 11.2.} (Restriction  ${\cal S}\downarrow {\cal A}$)\;  
{\it 
For a $\Lambda^n=(\Lambda^{n',1}, \Lambda^{n',2})\in{\boldsymbol Y}_{\!n'}({\boldsymbol \nu}^0)\times{\boldsymbol Y}_{\!n'}({\boldsymbol \nu}^0)\subset{\boldsymbol Y}_{\!n}({\boldsymbol \nu})$, take an IR $\pi_{\Lambda^n}\cong \pi_{\Lambda^{n',1}}\boxtimes\pi_{\Lambda^{n',2}}$ of ${\cal S}=
{\mathfrak S}^{(l)}_{{\boldsymbol \nu}^0}\times{\mathfrak S}^{(r)}_{{\boldsymbol \nu}^0}$.

{\rm (a)}\, If\, ${}^t\Lambda^n= \Lambda^n$, then $\pi_{\Lambda^n}\cong \mathrm{sgn}\cdot \pi_{\Lambda^n}$, and \,$\pi_{\Lambda^n}|_{\cal A}\cong \rho^{(0)}_{\Lambda^n}\,\bigoplus\,\rho^{(1)}_{\Lambda^n},\;(\Lambda^n,\kappa)\in{\boldsymbol Y}_{\!n}^{\mathfrak{A}}({\boldsymbol \nu})^{\rm sym},
\\ 
\kappa=0,1.$

{\rm (b)}\, If\, ${}^t\Lambda^n\ne \Lambda^n$, then \,$\pi_{\Lambda^n}\not\cong \mathrm{sgn}\cdot \pi_{\Lambda^n}$, and \,$\pi_{\Lambda^n}|_{\cal A}$ irreducible, and denote it by $\rho_{M^n},\;M^n=\{\Lambda^n,{}^t\Lambda^n\}\in{\boldsymbol Y}_{\!n}^{\mathfrak{A}}({\boldsymbol \nu})^{\rm asym}$.

IRs coming out from (a) and (b) give a CSR of the dual of ${\cal S}$. 
}
\hfill 
$\Box$\; 

\vskip1.2em 

{\bf Lemma 11.3.} (Action of $s_0$)\;  
{\it  
In the case of Proposition 8.4\,(ii-2), put $s_0=\tau_0$\, for $n'$ odd, and $s_0=s_j\tau_0$\, for $n'$ even  and $s_j\gamma=\gamma$. 
\vskip.3em

{\rm (i)}\; 
Let $\sigma=(\sigma_1,\sigma_2)\in {\cal A}={\mathfrak A}_n(\gamma)\subset{\mathfrak S}^{(l)}_{{\boldsymbol \nu}^0}\times{\mathfrak S}^{(r)}_{{\boldsymbol \nu}^0}$. Then $s_0=\tau_0$ acts as 
\;$s_0\sigma{s_0}^{\,-1}=(\sigma_2,\sigma_1)\;(transposition)$, and  
$s_0=s_j\tau_0$ with $s_j\in {\mathfrak S}^{(l)}_{{\boldsymbol \nu}^0}$ acts as 
\,$s_0\sigma{s_0}^{\,-1}=(\sigma_2,s_j\sigma_1s_j^{\,-1}).$

{\rm (ii)}\; 
For 
\,$\Lambda^n=(\Lambda^{n',1}, \Lambda^{n',2})\in{\boldsymbol Y}_{\!n}({\boldsymbol \nu})={\boldsymbol Y}_{\!n'}({\boldsymbol \nu}^0)\times{\boldsymbol Y}_{\!n'}({\boldsymbol \nu}^0)$, put \,$\overline{\tau_0}\Lambda^n:=(\Lambda^{n',2}, \Lambda^{n',1})$ (transposition). Action of $s_0$ on spin IRs of\, ${\cal A}$ is as follows\,: 
\vskip.5em
\hspace*{10ex}
$\begin{array}{ll} 
\mbox{\it in case of\;\;} {}^t\Lambda^n=\Lambda^n\,,\qquad\quad& s_0\big(\rho^{(\kappa)}_{\Lambda^n}\big)\,\cong\,\rho^{(\kappa+1)}_{\overline{\tau_0}(\Lambda^n)}\,,
\\[.3ex] 
\mbox{\it in case of\;\;}{}^t\Lambda^n\ne\Lambda^n\,, \qquad& s_0\big(\rho_{\Lambda^n}\big)\,\cong\,\rho_{\overline{\tau_0}(\Lambda^n)}\,.
\end{array}
$
\\[1ex]
Here, in case \,$\overline{\tau_0}(\Lambda^n)\ne \Lambda^n$\,,  the superfices $\kappa'=0,1$ for $\rho^{(\kappa')}_{\overline{\tau_0}(\Lambda^n)}$ is adequately arranged.  
}
\\
~\hfill 
$\Box$\;

\vskip.5em

{\bf Lemma 11.4.} (Inducing up ${\cal A}\uparrow{\cal C}$)\; 
{\it 
Let\, $\Lambda^n=(\Lambda^{n',1}, \Lambda^{n',2})\in{\boldsymbol Y}_{\!n'}({\boldsymbol \nu}^0)\times{\boldsymbol Y}_{\!n'}({\boldsymbol \nu}^0)$. 

{\rm (i)}\; Case of\, ${}^t\Lambda^n= \Lambda^n$.\;  
For $M^n=(\Lambda^n,\kappa)\in{\boldsymbol Y}^{{\mathfrak A}}_{\!n}(\widehat{T})^{\rm sym}$, \;

\quad 
$T^{\mathfrak A}_{M^n}:=\mathrm{Ind}^{\cal C}_{\cal A}\,\rho_{M^n}$\, is irreducible and for 
$\overline{\tau_0}(M^n):=\big(\overline{\tau_0}(\Lambda^n), \kappa+1\big)\in{\boldsymbol Y}^{{\mathfrak A}}_{\!n}(\widehat{T})^{\rm sym}$,\, 

\quad $T^{\mathfrak A}_{M^n}\cong T^{\mathfrak A}_{\overline{\tau_0}(M^n)}$ (equivalent).
 \vskip.3em

{\rm (ii)}\; Case of\, ${}^t\Lambda^n\ne \Lambda^n$. 
\;
For $M^n=\{\Lambda^n,{}^t\Lambda^n\}\in{\boldsymbol Y}^{{\mathfrak A}}_{\!n}(\widehat{T})^{\rm asym}$,
\vskip.3em

\quad{\rm (ii-1)}  
In case $\overline{\tau_0}(\Lambda^n)\not\in\{\Lambda^n, {}^t\Lambda^n\}$, \,$T^{\mathfrak A}_{M^n}:=\mathrm{Ind}^{\cal C}_{\cal A}\,\rho_{M^n}$\, is irreducible and 

\hspace*{5.7ex} 
put $\overline{\tau_0}(M^n):=\{\overline{\tau_0}(\Lambda^n),{}^t\,\overline{\tau_0}(\Lambda^n)\}\in{\boldsymbol Y}^{{\mathfrak A}}_{\!n}(\widehat{T})^{\rm asym}$, then \;$T^{\mathfrak A}_{M^n}\cong T^{\mathfrak A}_{\overline{\tau_0}(M^n)}$. 
\vskip.3em

\quad{\rm (ii-2)} 
In case $\overline{\tau_0}(\Lambda^n)\in\{\Lambda^n, {}^t\Lambda^n\}$, \,$T^{\mathfrak A}_{M^n}:=\mathrm{Ind}^{\cal C}_{\cal A}\,\rho_{M^n}$\, is reducible and 

\hspace*{5.7ex} 
$\overline{\tau_0}(M^n)=M^n$.   
\,$T^{\mathfrak A}_{M^n}$ splits into two irreducible components, mutually 

\hspace*{5.7ex} 
non-equivalent\,:\;  
$T^{\mathfrak A}_{M^n}\cong T^{\mathfrak A}_{M^n;0}\bigoplus T^{\mathfrak A}_{M^n;1}.$
\vskip.3em

{\rm (iii)}\; In the case of Proposition 8.4 (ii-2) and $\rho=Z_\gamma\cdot {\bf 1}_{S'([\rho])}$, IRs listed above (modulo two  equivalence relations above) gives a CSR of the dual of ${\cal C}=S([\rho])={\mathfrak A}_n(\gamma)\bigsqcup s_0{\mathfrak A}_n(\gamma)$. 
}

~\hfill 
$\Box$\; 

\vskip.5em

Now we can give a CSR of spin IRs of $G'_n=R\big(G(m,1,n)\big)$ of spin type $\chi^{\rm VII}=(1,\,1,-1)$. 
Put $H'_n=\widetilde{D}^\vee_n(T)\rtimes S'([\rho])$. \vskip1em

{\bf List 11.5.\;  
Spin IRs of $G'_n=R\big(G(m,1,n)\big)$ of spin type $\chi^{\rm VII}$.}

\vskip.5em

{\rm (1)}\; {\bf Case of $n$ odd. }  
\,$\Gamma^{\rm III}_n=\Gamma^{\rm III,0}_n\bigsqcup\Gamma^{\rm III,1}_n.$ 
\vskip.3em 

Case of $\gamma \in \Gamma^{\rm III,0}_n=\big\{\gamma\in\Gamma^1_n;\gamma\preccurlyeq\Psi_0(\gamma),\gamma_j=\gamma_{j+1}\;(\exists j)\big\}
$:\;\; 
$S([\rho])={\mathfrak A}_n(\gamma)$. 

\quad For $M^n\!=\!\{\Lambda^n,{}^t\Lambda^n\}\!\in\!{\boldsymbol Y}^{{\mathfrak A}}_{\!n}({\boldsymbol \nu})^{\rm asym}$,\, 
$\pi^1\!=\!\rho_{M^n}$,\, 
$\Pi^{\rm VII}_{M^n}\!:=\!\mathrm{Ind}^{G'_n}_{H'_n}\big(Z_\gamma\hspace{.2ex}{\bf 1}_{S'([\rho])}\boxdot \rho_{M^n}\big)$. 

\quad For $M^n=(\Lambda^n,\kappa)\in{\boldsymbol Y}^{\mathfrak A}_{\!n}({\boldsymbol \nu})^{\rm sym}$,\; 
$\pi^1=\rho_{M^n}$,\; 
$\Pi^{\rm VII}_{M^n}:=\mathrm{Ind}^{G'_n}_{H'_n}\big(Z_\gamma\hspace{.2ex}{\bf 1}_{S'([\rho])}\boxdot \rho_{M^n}\big)$. 

\vskip.3em 

Case of $\gamma\in \Gamma^{\rm III,1}_n=\big\{\gamma\in\Gamma^1_n;\mbox{\rm $\gamma_j$ all distinct}\big\}  
$:\;\;   
$S([\rho])={\mathfrak A}_n(\gamma)=\{e\}$. 

\quad For $M^n=\Lambda^n\in{\boldsymbol Y}^{{\mathfrak A}}_{\!n}({\boldsymbol \nu})^{\rm tri}$,\,$\pi^1=\rho_{M^n}={\bf 1},$\,  
$\Pi^{\rm VII}_{M^n}:=\mathrm{Ind}^{G'_n}_{H'_n}\,\big(Z_\gamma\hspace{.2ex}{\bf 1}_{S'([\rho])}\boxdot \rho_{M^n}\big)$.

\vskip1em

{\rm (2)}\; {\bf Case of $n$ even.}\; $
\Gamma^{\rm III}_n=\Gamma^{\rm III,0}_n\bigsqcup\Gamma^{\rm III,1}_n\bigsqcup\Gamma^{\rm III,2}_n$.
\vskip.3em

Case of $\gamma\in \Gamma^{\rm III,0,\prec}_n=
\big\{\gamma\in\Gamma^1_n\,;\,\gamma\prec\Psi_0(\gamma),\;\gamma_j=\gamma_{j+1}\;(\exists j)\big\}$:\;\; 
$S([\rho])={\mathfrak A}_n(\gamma)$. 

\quad For $M^n=\{\Lambda^n,{}^t\Lambda^n\}\in{\boldsymbol Y}^{\mathfrak A}_{\!n}({\boldsymbol \nu})^{\rm asym}$,\, $\pi^1=\rho_{M^n}$,\,  
$
\Pi^{\rm VII}_{M^n}:=\mathrm{Ind}^{G'_n}_{H'_n}\big(Z_\gamma\hspace{.2ex}{\bf 1}_{S'([\rho])}\boxdot \rho_{M^n}\big). 
$ 

\quad For $M^n=(\Lambda^n,\kappa)\in{\boldsymbol Y}_{\!n}({\boldsymbol \nu})^{\rm sym}$,\; 
$\pi^1=\rho_{M^n}$,\; 
$\Pi^{\rm VII}_{M^n}:=\mathrm{Ind}^{G'_n}_{H'_n}\big(Z_\gamma\hspace{.2ex}{\bf 1}_{S'([\rho])}\boxdot \rho_{M^n}\big)$. 

\vskip.3em

Case of $\gamma\in \Gamma^{\rm III,0,=}_n=
\big\{\gamma\in\Gamma^1_n\,;\,\gamma=\Psi_0(\gamma),\;\gamma_j=\gamma_{j+1}\;(\exists j)\big\}$: 
 \vskip.3em

\qquad 
$\left\{
\begin{array}{l}
\tau_0\gamma\equiv \gamma+{\boldsymbol m}^0\;(\mathrm{mod}\;m),\quad 
S([\rho])={\mathfrak A}_n(\gamma)\bigsqcup s_0{\mathfrak A}_n(\gamma),
\\[.3ex]
{\cal A}:={\mathfrak A}_n(\gamma),\;\;{\cal C}:=S([\rho]),\;\; \Lambda^n=(\Lambda^{n',1}, \Lambda^{n',2})\in{\boldsymbol Y}_{\!n'}({\boldsymbol \nu}^0)\times{\boldsymbol Y}_{\!n'}({\boldsymbol \nu}^0).
\end{array}
\right.
$
\vskip.2em

\quad For $M^n=(\Lambda^n,{}^t\Lambda^n)\in{\boldsymbol Y}^{{\mathfrak A}}_n({\boldsymbol \nu})^{\rm asym}$,

If $\overline{\tau_0}(M^n)\ne M^n$, then \; 
$\pi^1=T^{\mathfrak A}_{M^n}=\mathrm{Ind}^{\cal C}_{\cal A}\,\rho_{M^n}, 
\;
\Pi^{\rm VII}_{M^n}:=\mathrm{Ind}^{G'_n}_{H'_n}\big(Z_\gamma\hspace{.2ex}{\bf 1}_{S'([\rho])}\boxdot T^{\mathfrak A}_{M^n}\big). 
$
\\[.3ex]
\qquad\hspace{23ex} 
$\Pi^{\rm VII}_{M^n}\cong \Pi^{\rm VII}_{\overline{\tau_0}(M^n)}$ \;(equivalent).
\vskip.3em

If $\overline{\tau_0}(M^n)= M^n$, then $\pi^1=T^{\mathfrak A}_{M^n;\,\iota}\; (\iota=0,1),\;
\Pi^{\rm VII}_{M^n;\,\iota}:=\mathrm{Ind}^{G'_n}_{H'_n}\big(Z_\gamma\hspace{.2ex}{\bf 1}_{S'([\rho])}\boxdot T^{\mathfrak A}_{M^n;\,\iota}\big). 
$
\vskip.2em

\quad For $M^n=(\Lambda^n,\kappa)\in{\boldsymbol Y}^{{\mathfrak A}}_n({\boldsymbol \nu})^{\rm sym}$, 

\hspace*{21ex}
$\pi^1=T^{\mathfrak A}_{M^n}:=\mathrm{Ind}^{\cal C}_{\cal A}\,\rho_{M^n}$,\; 
\hspace{.5ex}
$\Pi^{\rm VII}_{M^n}=\mathrm{Ind}^{G'_n}_{H'_n}(Z_\gamma\hspace{.2ex}{\bf 1}_{S'([\rho])}\boxdot T^{\mathfrak A}_{M^n}).$

\hspace*{5ex} 
$\overline{\tau_0}(M^n):=\big(\overline{\tau_0}(\Lambda^n), \kappa+1\big)\in{\boldsymbol Y}^{{\mathfrak A}}_{\!n}(\widehat{T})^{\rm sym}$, \;$\Pi^{\rm VII}_{M^n}\cong \Pi^{\rm VII}_{\overline{\tau_0}(M^n)}$ (equivalent).

\vskip.3em

Case of $\gamma\in \Gamma^{\rm III,1}_n\bigsqcup 
\Gamma^{\rm III,2}_n$\,:\quad  
$S([\rho])=\{e\}$,\; 

\quad For $M^n=\Lambda^n$ quasi-trivial,\;\; $\pi^1=\rho_{M^n}={\bf 1},$\;\;
$\Pi^{\rm VII}_{M^n}:=\mathrm{Ind}^{G'_n}_{H'_n}\,\big(Z_\gamma\hspace{.2ex}{\bf 1}_{S'([\rho])}\boxdot \rho_{M^n}\big)$. 

\vskip1.2em

{\bf Theorem 11.6.} 
$\big($Spin IRs of $R\big(G(m,1,n)\big)$ of spin type $\chi^{\rm VII}=(1,\,1,-1)\big)$ 
 
{\it 
Spin IRs of $R\big(G(m,1,n)\big),\;n\ge 4,$ in Case VII and spin IRs of the quotient group  
\;$
G'=\widetilde{G}^{\rm \,VII}_n=R\big(G(m,1,n\big)/Z_{12},\;Z_{12}=Z_1\times Z_2, 
$\; 
of spin type $z_3\to -1$ can be naturally identified. Then a complete set of representatives of the latter is given as follows\,: $M^n$ contains $\Lambda^n$, and $\Lambda^n$ determines $\gamma\in\Gamma^1_n$, i.e., $\Lambda^n\rightsquigarrow \gamma$, 
\\[.5ex]
\indent
{\rm (i)}\;\; {\bf Case of $n$ odd.} \;
\\[.5ex]
\qquad $\mbox{\rm {\footnotesize spin}IR}^{\rm VII}\big(G(m,1,n)\big)
=\big\{\Pi^{\rm VII}_{M^n}\,;\, M^n \in{\boldsymbol Y}^{{\mathfrak A}}_{\!n}(\widehat{T}),\;\gamma\in\Gamma^{\rm III}_n\big\},
$
\vskip.5em

{\rm (ii)}\; {\bf Case of $n$ even.}\; 
\\[.5ex]
\qquad $\mbox{\rm {\footnotesize spin}IR}^{\rm VII}\big(G(m,1,n)\big)=
\big\{\Pi^{\rm VII}_{M^n}\,;\,M^n\in{\boldsymbol Y}^{{\mathfrak A}}_{\!n}(\widehat{T}),\;\gamma\in \Gamma^{\rm III,0,\prec}\bigsqcup\Gamma^{\rm III,1}\bigsqcup\Gamma^{\rm III,2}\big\}
\\[.5ex]
\hspace*{15.3ex} 
\bigsqcup 
\big\{\Pi^{\rm VII}_{M^n}\cong \Pi^{\rm VII}_{\overline{\tau_0}(M^n)}\;;\,M^n\in{\boldsymbol Y}^{{\mathfrak A}}_{\!n}(\widehat{T})^{\rm asym},\;\gamma\in \Gamma^{\rm III,0,=},\,\overline{\tau_0}(M^n)\ne M^n\big\}
\\[.5ex]
\hspace*{15.3ex} 
\bigsqcup 
\big\{\Pi^{\rm VII}_{M^n;\,\iota}\;;\;\iota=0,1,\;M^n\in{\boldsymbol Y}^{{\mathfrak A}}_{\!n}(\widehat{T})^{\rm asym},\;\gamma\in \Gamma^{\rm III,0,=},\,\overline{\tau_0}(M^n)= M^n\big\}
\\[.8ex]
\hspace*{15.3ex} 
\bigsqcup 
\big\{\Pi^{\rm VII}_{M^n}\cong\Pi^{\rm VII}_{\overline{\tau_0}(M^n)} \;;\,M^n\in{\boldsymbol Y}^{{\mathfrak A}}_{\!n}(\widehat{T})^{\rm sym},\;\gamma\in \Gamma^{\rm III,0,=}\big\}.
$

}

\part{\Large 
Spin characters of infinite generalized symmetric groups $G(m,1,\infty)$}

 Applying general theory of\, \lq\lq\,Inductive limits of finite groups and their characters\,'' to the case of $G_n=G(m,1,n)\to G_\infty=G(m,1,\infty)$ and $G'_n=R\big(G(m,1,n)\big)\to G'_\infty=R\big(G(m,1,\infty)\big)$, we calculate all characters and spin characters of infinite generalized symmetric groups.

\section{\large Preliminaries to Part II}

\quad
{\bf 12.1. Supports of $f\in K(G'_\infty;\chi^{\rm Y})$ for $G'_\infty=R\big(G(m,1,\infty)\big)$.}
  
We studied in [I, \S 10] the support of $f\in K(G'_\infty;\chi)$, in particular, of $f\in E(G'_\infty;\chi)$, for each spin type $\chi\in \widehat{Z}$, from the viewpoint\; 
\lq {\it how much they are characterized by its spin type $\chi^{\rm Y}$, or by the property (\ref{2017-07-23-1})}\,''.\, 
We quote the result from [ibid., Table 10.1, p.83]. 
Take an element $g'=(d',\sigma')\in G'_\infty=\widetilde{D}^\vee_\infty\rtimes\widetilde{{\mathfrak S}}_\infty$, and put $g=(d,\sigma)=\Phi(g')\in G_\infty=D_\infty\rtimes{\mathfrak S}_\infty$ be its canonical image under $\Phi: G'_\infty\to G_\infty$. Let a decomposition into fundamental elements of $g'$ and its canonical image for $g$ be respectively 
\begin{eqnarray}
\label{2017-07-31-1}
&&\;\;
\left\{
\begin{array}{ll}
g'=z_2^{\;a}z_3^{\;b}\,\xi'_{q_1}\xi'_{q_2}\cdots\xi'_{q_r}g'_1g'_2\cdots g'_s,&
\xi'_q = \big(t'_q, (q)\big)={\widehat{\eta}_q}^{\;a_q},\;\;
g'_j=(d'_j, \sigma'_j),
\\[.5ex] 
g\,=\,\xi_{q_1}\xi_{q_2}\cdots\xi_{q_r}g_1g_2\cdots g_s,&
\xi_q = \big(t_q, (q)\big)=y_q^{\;a_q},\;\;
g_j=(d_j, \sigma_j),
\end{array}
\right.
\end{eqnarray}
where $y_q=\Phi(\widehat{\eta}_q)$ and the fundamental elements are given as 
\begin{eqnarray}
\label{2017-07-31-2}
&&
\left\{
\begin{array}{l}
\sigma_j=\Phi(\sigma'_j)\;\,\mbox{\rm a cycle},\;d_j=\Phi(d'_j),\; K_j:=\mathrm{supp}(\sigma_j) \supset \mathrm{supp}(d_j), 
\\[.7ex]
\mbox{\rm $K_j\;(j\in J:={\boldsymbol I}_s)$ are mutually disjoint,}\;Q:=\{q_1,q_2,\ldots,q_r\}.
\end{array}
\right.
\end{eqnarray}
We put \,$\ell(\sigma'_j):=\ell(\sigma_j),\,L(\sigma'_j)=L(\sigma_j):=\ell(\sigma_j)-1$. 

For each Y\,=\,odd, I, II, $\cdots$, VII, a subset ${\cal O}({\rm Y})\subset G'_\infty$ is given by Condition Y on its element $g'$ but expressed by means of $g=\Phi(g')$. Or, for $f \in K(G'_\infty;\chi^{\rm Y}),$
$$
 f(g')\ne 0\;\Longrightarrow\;g'\;\mbox{\rm satisfies Condition Y}. 
$$ 

Let ${\cal O}$ be a non-empty subset of $G'_\infty$  invariant under multiplication of the central subgroup $Z$. 
We call ${\cal O}$ {\it factorisable} if it has the property that, for any\; $g' \in {\cal O}$, 
\\[1ex]
\hspace*{8ex}
$g'=h'k',\;h',k' \in G'_\infty, \;\mathrm{supp}(h')\;{\textstyle \bigcap}\;\mathrm{supp}(k') = \varnothing \;\Longrightarrow\;h',k' \in {\cal O}\,; 
$
\\[1ex]
and {\it commutatively factorisable for type $\chi\in\widehat{Z},$} if, for any\; $g' \in {\cal O}$, 
\begin{eqnarray}
\label{2009-02-24-2}
&&
\left\{
\begin{array}{l}
g'=h'k',\;\;h',k' \in G'_\infty, \;\;\mathrm{supp}(h')\bigcap\mathrm{supp}(k') = \varnothing, 
\\[1ex]
\qquad\qquad\quad
\Longrightarrow\;\;
h',k' \in {\cal O}\;\;{\rm and}\;\;h'k'\mathrm{Ker}(\chi) =k'h'\mathrm{Ker}(\chi)\,; 
\end{array}
\right. 
\end{eqnarray}
and {\it multiplicative} if\; $ 
 h',k' \in {\cal O} \;\;\Longrightarrow\;\;h'k' \in {\cal O}.$

We call Condition\,V (=\,Condition\,VI) also as Condition\,(str) (i.e., strongest).

\vskip1.2em
\begin{center}
{\footnotesize 
 $
\begin{array}{|c||c|l|c|c|c|c|}
\hline 
\!\!
\begin{array}{c}
{\rm Subset} 
\\
{\cal O}{\rm (Y)}
\end{array}\!\!\!\!
&
\!\begin{array}{c}
\chi^{\rm Y}=
\\
(\beta_1, \beta_2, \beta_3) 
\end{array}\!
&
\!\!  
\begin{array}{c}
\mbox{\rm \hspace{5ex}(Condition Y)} 
\end{array}
\!\!\!\!
& 
\!\begin{array}{l}
 \;\;{\cal O}({\rm Y}) 
\\
\mbox{\rm Facto-}
\\
\mbox{\rm risable}
\end{array}\!\!\!\!
&
\!\!
\begin{array}{l}
\;\;{\cal O}({\rm Y}) 
\\
\mbox{\rm commutati.}
\\
\mbox{\rm factorisable}
\mbox{\rm }
\end{array}
\!\!\!\!
&
\!\!
\begin{array}{l}
\;\;{\cal O}({\rm Y}) 
\\
\mbox{\rm multipli-}
\\
\mbox{\rm cative}
\end{array}
\!\!\!\!
\\
\hline\hline
\multicolumn{2}{c}{
\begin{array}{c}
m\quad\;\;{\bf ODD}\quad 
\end{array}
}
&
\multicolumn{4}{c}{}
\\[.3ex]
\hline
\!\!\!
\begin{array}{c}
\!{\cal O}({\rm odd})\!\!
\end{array}
\!\!\!\!
&
\chi^{\rm odd}(z_1)=-1
&
\begin{array}{l}
L(\sigma_j)\equiv 0\;({\rm mod}\,2)
\\
\qquad\qquad\quad\;\; (1 \leqslant j \leqslant s)
\end{array}
\!\!
&
{\rm YES}
&
{\rm YES}
&
 {\rm NO}
 \\
 \hline 
 \hline  
 \multicolumn{2}{c}{\;\;m\quad\;\; {\bf EVEN}\quad} 
 &
 \multicolumn{4}{c}{}
 \\[.3ex]
 \hline
 {\cal O}({\rm I}) & 
 \!\!
 (-1,-1,-1)
 \!\!
  &
  \begin{array}{l} 
 \mathrm{ord} (\xi_{q_i}) \equiv 0\,({\rm mod}\,2)\,(\forall i)
 \\
 {\rm ord}(d_j) + L(\sigma_j)\equiv 0\,(\forall j)
 \end{array}
 \!\!\!\!
 &
 {\rm YES}
 &
 {\rm YES}
&
 {\rm NO}
   \\
 \hline 
 {\cal O}({\rm II}) & (-1,-1,\;1)    &
\begin{array}{l} 
\mathrm{ord} (d)  \equiv 0,\quad L(\sigma)  \equiv  0,\;
\\
\mathrm{ord} (\xi_{q_i}) \equiv 0\;(\forall i),  
\\
{\rm ord}(d_j) + L(\sigma_j)\equiv 0\,(\forall j)
\end{array}
\!\!\!\!
&{\rm NO}
&{\rm NO}
&
 {\rm NO}
\\
 \hline 
 {\cal O}({\rm III}) & (-1,\;1,-1)   &
\begin{array}{l}
 \mathrm{ord} (d) \equiv 0,
 \\
 L(\sigma_j) \equiv 0\;(\forall j) 
 \end{array}
 &{\rm NO}
 &{\rm NO}
&
 {\rm NO}
\\
 \hline 
  {\cal O}({\rm IV}) & (-1,\;1,\;1)   &
 \begin{array}{ll}
   L(\sigma_j) \equiv 0\;(\forall j),
 \\
 \mbox{\rm for }\sigma= \sigma_1\sigma_2\cdots\sigma_s 
 \end{array}
 &{\rm YES}
 &{\rm YES}
 &
 {\rm NO}
  \\
 \hline
   {\cal O}({\rm V}) & (1,-1,-1)  &
  \begin{array}{l}
  \mathrm{ord} (\xi_{q_i}) \equiv 0\;(\forall i),
  \\ 
 \mathrm{ord} (d_j) \equiv 0\;(\forall j),
 \\
 L(\sigma_j) \equiv 0\;(\forall j)
 \end{array}
 &{\rm YES}
 &{\rm YES}
 &
 {\rm NO}
 \\
  \hline  
  {\cal O}({\rm VI}) & (1,-1,\;1)      &
 \begin{array}{l}
  \mathrm{ord} (\xi_{q_i}) \equiv 0\;(\forall i),  
 \\
 {\rm ord}(d_j) \equiv 0\;(\forall j),
 \\
  L(\sigma_j)  \equiv 0\;(\forall j)
 \end{array}
 &  {\rm YES}
 &{\rm YES}
&
 {\rm NO}
\\
 \hline  
 \!
 {\cal O}({\rm VII})
 \!\!
 &  (1,\;1,-1)   &
\begin{array}{l}
  {\rm ord}(d) \equiv 0,  
  \\
  L(\sigma)\equiv 0
  \end{array}
  &{\rm NO}
  &{\rm NO}
  &
 {\rm YES}
  \\
 \hline  
  \end{array}
  $
  }
\vskip.5em

 {\bf Table 12.1.}\; Subsets ${\cal O}{\rm (Y)} \subset G'_\infty,$\, 
${\cal O}{\rm (Y)} \supset \mathrm{supp}(f),\;\forall f\in K(G'_\infty;\chi^{\rm Y})$. 
\end{center}

\vskip1.5em
We define several normal subgroups of $G_\infty=G(m,1,\infty)=D_\infty\rtimes {\mathfrak S}_\infty$. Let $S$ be a subgroup of $T={\boldsymbol Z}_m$ and put $G^S_\infty:=D^S_\infty\rtimes {\mathfrak S}_\infty$ with $D^S_\infty:=\{d\in D_\infty\,;\,P(d)\in S\}$. 
We know that any $S$ is given as $S(p):=\{t^p\,;\,t\in T\}$ for some $p|m,\,p\ge 1$. In that case, $G^S_\infty$ is denoted as $G(m,p,\infty)$. 
Together with them, we call {\it canonical  
normal subgroups} of $G_\infty$ the following:  
\begin{eqnarray}
\label{2017-11-24-11}
\left\{
\begin{array}{l}
G(m,p,\infty),\;p|m,\,p\ge 1,\;\;G^{\mathfrak A}_\infty:=D_\infty\rtimes {\mathfrak A}_\infty,\;\;
\\[.5ex]
G^{\mathfrak A}(m,p,\infty):=G^{\mathfrak A}_\infty\cap G(m,p,\infty).
\end{array}
\right.
\end{eqnarray}
For later use, we put $K_\infty:=G^{\mathfrak A}(m,2,\infty)$ for $m$ even. Then $K_\infty=D^{\rm ev}_\infty\rtimes {\mathfrak A}_\infty$ with $D^{\rm ev}_\infty:=D^{S(2)}_\infty=\{d\in D_\infty\,;\,\mathrm{ord}(d)\equiv 0\:(\mathrm{mod}\;2)\}$. 
As their full inverse images, we put 
\begin{eqnarray}
\label{2017-08-03-1}
&&
\begin{array}{l}
G^{\prime\hspace{.15ex}{\mathfrak A}}_\infty:=\Phi^{-1}(G^{\mathfrak A}_\infty)=\widetilde{D}^\vee_\infty\rtimes \widetilde{{\mathfrak A}}_\infty,
\;\;
K'_\infty:=\Phi^{-1}(K_\infty)=\widetilde{D}^{\rm ev}_\infty\rtimes \widetilde{{\mathfrak A}}_\infty,
\\[.5ex]
\hspace*{6ex}
{\rm with} 
\quad
\widetilde{D}^{\rm ev}_\infty:=\{d'\in \widetilde{D}^\vee_\infty\,;\,\mathrm{ord}(d')\equiv 0\;(\mathrm{mod}\;2)\}.
\end{array}
\end{eqnarray}
Then the subset ${\cal O}({\rm VII})$ is equal to the subgroup $K'_\infty$.
\vskip1.2em

{\bf 12.2. Criterion (EF) and Criterion (EF$^{\chi^{\rm Y}}$).}\; 
\vskip.5em 

{\bf Definition 12.2.} We call an $f \in K_1(G'_\infty)$\, {\it factorisable}\, if 
\begin{eqnarray}
\label{2017-08-01-1}
&& 
g'\in G'_\infty,\,g'=h'k',\,\mathrm{supp}(h')\cap\mathrm{supp}(k')=\emptyset\;
\Longrightarrow\;f(g')=f(h')f(k'),
\end{eqnarray}
and call $f\in K_1(G'_\infty;\chi^{\rm Y})$\, {\it $\chi^{\rm Y}$-factorisable}\, if 
\begin{eqnarray}
\label{2017-08-01-2}
&&\;\;
h'\in {\cal O}({\rm Y}),\,k'\in G'_\infty,\,\mathrm{supp}(h')\cap\mathrm{supp}(k')=\emptyset\;
\Longrightarrow\;f(h'k')=f(h')f(k').
\end{eqnarray}

The set of all factorisable $f\in K_1(G'_\infty;\chi)$ is denoted by $F(G'_\infty;\chi)$, and that of all $\chi^{\rm Y}$-factorisable $f\in K_1(G'_\infty;\chi^{\rm Y})$ by $F^{\chi^{\rm Y}}(G'_\infty;\chi^{\rm Y})$. Then, $F(G'_\infty;\chi)\subset F^{\chi^{\rm Y}}(G'_\infty;\chi^{\rm Y})$, and we ask if\; $E(G'_\infty;\chi)=F(G'_\infty;\chi)$ \,or\, $E(G'_\infty;\chi^{\rm Y})=F^{\chi^{\rm Y}}(G'_\infty;\chi^{\rm Y})$. Concerning to this, we propose two criteria to be extremal, for $f\in K_1(G'_\infty;\chi)$ and $f\in K_1(G'_\infty;\chi^{\rm Y})$ respectively, as  
\vskip.8em 

{\bf (EF)}\hspace{17.3ex} 
$f\in E(G'_\infty;\chi)\;\,\Longleftrightarrow\;f\in F(G'_\infty;\chi)\,;
$
\vskip.5em 
{\bf (EF$^{\chi^{\rm Y}}$)}\hspace{14ex} 
$f\in E(G'_\infty;\chi^{\rm Y})\;\Longleftrightarrow\;f\in F^{\chi^{\rm Y}}(G'_\infty;\chi^{\rm Y}).
$
\vskip1em 

We call spin types $\chi^{\rm Y}$,\, for Y\,=\,odd,\,I, IV, V, VI, VIII ($\chi^{\rm VIII}$ is trivial), {\it of the 1st kind}, and those for Y\,=\,II, III, VII, {\it of the 2nd kind}. 
We know from [I, Theorems 6.3] the following.  
\vskip1.1em

{\bf Theorem 12.3.}\; 
{\it 
{\rm (i)} $F(G'_\infty;\chi)\subset E(G'_\infty;\chi)$\, for any $\chi\in\widehat{Z}$, or,  
if\, $f\in K_1(G'_\infty)$ is factorisable, then it is extremal.

{\rm (ii)} For a spin type\, $\chi=\chi^{\rm Y}$ of the first kind, the criterion\, {\bf (EF)} holds. 

}
\vskip1.2em 

For a subset ${\cal O}$ of $G'_\infty$, we put 
\begin{eqnarray}
\label{2017-11-21-11}
E(G'_\infty;\chi;{\cal O}):=\{f\in E(G'_\infty)\,;\,\mathrm{supp}(f)\subset{\cal O}\}, 
\end{eqnarray}
and similarly for $F(G'_\infty;\chi;{\cal O})$. For Y\,=\,II, III and VII, an element $g'\in{\cal O}({\rm Y})$ is called $\chi^{\rm Y}$-{\it fundamental} if it cannot be decomposed into a product of non-trivial $h',k'\in{\cal O}({\rm Y})$ as $g'=h'k'$. 
From the results in [I, Theorem 11.1], [II, \S 23.3], [8, \S 10.4] 
and Table 12.1, we know

\vskip.7em

{\bf Theorem 12.4.}\;\, 
{\it 
For spin types $\chi^{\rm Y}$ of the 2nd kind (or for\, {\rm Y=II, III and VII}), there does not hold Criterion {\bf (EF)}, but 
 there holds Criterion {\bf (EF$^{\chi^{\rm Y}}$)}.
 Moreover\;  
$$
\qquad 
F\big(G'_\infty;\chi^{\rm Y}\big)=F\big(G'_\infty;\chi^{\rm Y};{\cal O}({\rm str})\big)=E\big(G'_\infty;\chi^{\rm Y};{\cal O}({\rm str})\big)\subsetneqq E\big(G'_\infty;\chi^{\rm Y}\big).
$$
}

\vskip.2em 

This result can be seen also from explicit determination of $F(G'_\infty;\chi)$ and $E(G'_\infty;\chi)$ for each $\chi\in\widehat{Z}$ (the latter is summarized in Table 22.2).
 By Theorems 12.1 and 12.4, to describe a character $f\in E(G'_\infty;\chi^{\rm Y})$ explicitly, it is enough to give its form $f(g')$ for each fundamental elements $g'\in{\cal O}({\rm Y}).$ 

\vskip1.2em

{\bf 12.3. Limiting process as $n\to\infty$.}\; 

Our fundamental ingredients here come from the study in Part I\,:   

{\bf (1)}\; Explicit form of characters of linear IRs of $G_n$, 

{\bf (2)}\; Such of spin IRs of $G'_n$ of spin type $\chi\in \widehat{Z}$, for each non-trivial $\chi^{\rm Y}$. 

Denote by ${\rm Lim}(G'_\infty)$ the set of all the pointwise limits $f_\infty=\lim_{n\to\infty}\widetilde{\chi}_{\pi_n}$ of normalized characters of IRs $\pi_n$ of $G'_n\;(n\ge 4)$. 
Applying [17, Theorem 14.3] 
 for inductive limits of compact groups to the case of $G'_n\to G'_\infty$, we have the following

\vskip1em

{\bf Theorem 12.5.}\; 
{\it 
Any character $f\in E(G'_\infty)$ is a pointwise limit of normalized characters $\widetilde{\chi}_{\pi_n}\in E(G'_n)\;(n\ge 4)$ of a series of IRs $\pi_n$ of $G'_n$, i.e., $E(G'_\infty)\subset {\rm Lim}(G'_\infty)$.  

}
 \vskip1.1em
 
In the following, we calculate all the limits $f_\infty=\lim_{n\to\infty}\widetilde{\chi}_{\pi_n}$ explicitly, and then find that,  
for Y of the 1st kind, \,${\rm Lim}(G'_\infty;\chi^{\rm Y})\subset F(G'_\infty;\chi^{\rm Y})$, and 
for Y of the 2nd kind, \,${\rm Lim}(G'_\infty;\chi^{\rm Y})\subset F^{\chi^{\rm Y}}(G'_\infty;\chi^{\rm Y})$. 
In the latter case, actually there exists a difference between $F\big(G'_\infty;\chi^{\rm Y}\big)$ and $F^{\chi^{\rm Y}}\big(G'_\infty;\chi^{\rm Y}\big)$, reflecting  the difference ${\cal O}({\rm str})\subsetneq {\cal O}({\rm Y})$. 

 \vskip1.2em

{\bf 12.4. Finite-dimensional spin IRs of $G'_\infty$.}
 
By [ibid., Theorems 12.1], the unique spin type for $G'_\infty$ which admits finite-dimensional representations is $\chi^{\rm VII}=(1,\,1,-1),$ for $m$ even. IRs with this spin type are  $\pi_{2,\zeta^{(k)}}$ given by $\zeta^{(k)}\in\widehat{T}^0\;(0\le k <m^0=m/2)$ as \;
\begin{eqnarray*}
&&
\pi_{2,\zeta^{(k)}}(z_1)=\pi_{2,\zeta^{(k)}}(z_2)=I_2,\;\pi_{2,\zeta^{(k)}}(z_3)=-I_2,
\\[.5ex]
&&
\pi_{2,\zeta^{(k)}}(r_i)=
\begin{pmatrix} 0& 1 
\\
1&0
\end{pmatrix}\;\;(i\ge 1),\quad 
\pi_{2,\zeta^{(k)}}(\widehat{\eta}_j)=
\begin{pmatrix} \zeta^{(k)}(\eta_1)\!\!& 0 
\\
0&\!\!\!-\zeta^{(k)}(\eta_1)
\end{pmatrix}\;\;(j\ge 1), 
\end{eqnarray*}
where $I_2$ denotes the $2\times 2$ identity matrix, and $\zeta^{(k)}(\eta_1)=\omega^k=e^{2\pi k/m}$. Their characters are given by 
[17, Theorem 12.2]\,: for $g'=(d',\sigma')\in G'_\infty$ with $d'=z^{\;a}_2z^{\;b}_3\,\widehat{\eta}^{\;a_1}_1\widehat{\eta}^{\;a_2}_2\cdots,$
\begin{eqnarray}
\label{2017-08-02-11}
\chi\big(\pi_{2,\zeta^{(k)}}\big|g'\big)=\left\{
\begin{array}{ll}
2\cdot (-1)^b\zeta^{(k)}(\eta_1)^{\mathrm{ord}(d')},\quad& \big(g'\in{\cal O}({\rm VII})\big),
\\[.5ex]
\quad 0\,, &{\rm \;otherwise}.
\end{array}
\right. 
\end{eqnarray}
  
The indicator function of ${\cal O}({\rm VII})=K'_\infty$ is denoted by $X_{{\cal O}({\rm VII})}=X_{K'_\infty}$. 
In the semidirect product 
$G^{\prime\hspace{.15ex}{\mathfrak A}}_\infty=\big(Z_3\times D^\wedge_\infty\big)\rtimes\widetilde{\mathfrak{A}}_\infty$, the sign character $\mathrm{sgn}_{Z_3}$ of $Z_3$ can be naturally extended to a character of $G^{\prime\hspace{.15ex}{\mathfrak A}}_\infty\supset{\cal O}({\rm VII})=K'_\infty$. Denote it by $\mathrm{sgn}^\mathfrak{A}_{Z_3}$, then the product\, $\mathrm{sgn}^\mathfrak{A}_{Z_3}\hspace*{-.12ex}\cdot\hspace*{-.1ex} X_{{\cal O}({\rm VII})}$\, can be considered naturally as a function on $G'_\infty$. 
 
\vskip1em

{\bf Proposition 12.6.}\; 
{\it 
The normalized character \,$f^{\rm VII}_k:=\chi(\pi_{2,\zeta^{(k)}})/2$ of $\pi_{2,\zeta^{(k)}}$ is in $E\big(G'_\infty;\,\chi^{\rm VII} \big)$ and not factorisable. In particular, $f^{\rm VII}_0$ satisfies 
\begin{eqnarray}
\label{2015-07-29-2}
&&
\;f^{\rm VII}_0=\mathrm{sgn}^\mathfrak{A}_{Z_3}\hspace*{-.12ex}\cdot\hspace*{-.1ex} X_{{\cal O}({\rm VII})}\,, \;\;\;\big(f^{\rm VII}_0\big)^2= X_{{\cal O}({\rm VII})}\,.
\end{eqnarray}
}
\vskip.3em

{\bf 12.5. Is the product of two characters also a character ?}\; 

Let $m$ be even. Let $f_1, f_2$ be characters of $G'_\infty$, and ask \lq\lq{\it Is the product $f_3=f_1f_2$ also a character ?}\;'' 
Let $f_i\in E(G'_\infty;\chi^{\rm Y_i})$ for $i=1,2$, then $f_3=f_1f_2\in K_1(G'_\infty;\chi^{\rm Y_3}$) with 
$\chi^{{\rm Y}_3}=\chi^{{\rm Y}_1}\chi^{{\rm Y}_2}$, and $\mathrm{supp}(f_3)\subset {\cal O}({\rm Y_1})\cap{\cal O}({\rm Y_2})\subset{\cal O}({\rm Y_3})$.  
The spin types of the 1st kind are $\big\{\chi^{\rm Y}\,;\,{\rm Y\,=\,I, IV, V, VI, VIII}\big\}$, and those of the 2nd kind are $\big\{\chi^{\rm Y}\,;\,{\rm Y\,=\,II, III, VII}\big\}$. 
\vskip1em

\begin{center}
{\small 
$
\begin{array}{|c||c|c|c|c|c||c|c|c|}
\hline
{\rm Y_1}&{\rm I}&{\rm IV}&{\rm V}&{\rm VI}&{\rm VIII}&{\rm II}&{\rm III}&{\rm VII}
\\
\hline
\hline
{\rm Y_2=II}&{\rm VII}&{\rm VI}&{\rm III}&{\rm IV}&{\rm II}&{\rm VIII}&{\rm V}&{\rm I}
\\
\hline
{\rm Y_2=III}&{\rm VI}&{\rm VII}&{\rm II}&{\rm I}&{\rm III}&*&{\rm VIII}&{\rm IV}
\\
\hline
{\rm Y_2=VII}&{\rm II}&{\rm III}&{\rm VI}&{\rm V}&{\rm VII}&*&*&{\rm VIII}
\\
\hline
\end{array}
$
\vskip.5em
{\bf Table 12.7.}\; Products of spin types\; $\chi^{{\rm Y}_3}=\chi^{{\rm Y}_1}\chi^{{\rm Y}_2}$.
}
\end{center}

\vskip.2em

If both of $\chi^{\rm Y_1}, \chi^{\rm Y_2}$ are of the 1st kind or of the 2nd kind, then the product $\chi^{\rm Y_1}\chi^{\rm Y_2}=\chi^{\rm Y_3}$ is of the 1st kind. Assume that $\chi^{\rm Y_1}$ be of the 1st kind and $\chi^{\rm Y_2}$ be of the 2nd kind as in the left half of Table 12.7, and consider the cases where the product $\chi^{\rm Y_3}$ is of the 1st kind. There exist 6 such triplets $({\rm Y_1,Y_2,Y_3})$'s. 
Containing them, we consider the cases where at least one of ${\rm Y}_i$ is V or VI. Then we see, for $f_i\in E(G'_\infty;\chi^{{\rm Y}_i}),\,i=1,2,$  
\begin{eqnarray}
\label{2017-09-11-1}
&&\quad 
%
\mathrm{supp}(f_1f_2)\subset {\cal O}({\rm Y_1})\cap{\cal O}({\rm Y_2})={\cal O}({\rm str})={\cal O}({\rm V})={\cal O}({\rm VI}). 
\end{eqnarray}

\vskip.2em

{\bf Theorem 12.8.}\; 
{\it 
{\rm (i)}\; 
Let $f_1, f_2$ be characters of $G'_\infty$ which are factorisable, then the product $f_1f_2$ is also a character. Symbolically  \;$F(G'_\infty)F(G'_\infty)\subset F(G'_\infty)$ and $F(G'_\infty;\chi^{{\rm Y}_1})\cdot 
\\ 
F(G'_\infty;\chi^{{\rm Y}_2})\subset F(G'_\infty;\chi^{{\rm Y}_3}).$ 
In case $\chi^{{\rm Y}_1}, \chi^{{\rm Y}_2}$ are both of the 1st kind, this means also that 
\begin{eqnarray}
\label{2017-09-15-1}
E(G'_\infty;\chi^{{\rm Y}_1})E(G'_\infty;\chi^{{\rm Y}_2})\subset E(G'_\infty;\chi^{{\rm Y}_3}).
\end{eqnarray}

{\rm (ii)}\; 
Assume that one of\, ${\rm Y}_i\;(i=1,2,3)$ is\, {\rm V} or\, {\rm VI}. Then there holds (\ref{2017-09-15-1}).

}

\vskip1em

{\it Proof.}\; We prove (ii). Note that $\mathrm{supp}(f_3)\subset{\cal O}({\rm str})$ and ${\cal O}({\rm str})$ is factorisable (Table 12.1). It suffices to  prove that $f_3$ is factorisable by Theorem 12.3\,(i).

If a $g'\in{\cal O}({\rm str})$ splits as $g'=h'k',\,\mathrm{supp}(h')\cap\mathrm{supp}(k')=\emptyset$, then $h',k'\in{\cal O}({\rm str})\subset{\cal O}({\rm Y}_i)$ for $i=1,2$. Hence $f_i(g')=f_i(h')f_i(k')\;(i=1,2)$ by Theorem 12.4 for $i=2$, and so  $f_3(g')=f_3(h')f_3(k').$ 
If a $g'\not\in{\cal O}({\rm str})$ splits as $g'=h'k',\,\mathrm{supp}(h')\cap\mathrm{supp}(k')=\emptyset$, then one of $h'$ and $k'$ is not in ${\cal O}({\rm str})={\cal O}({\rm Y}_1)\cap{\cal O}({\rm Y}_2)$. Suppose $h'\not\in{\cal O}({\rm str})$, then $h'\not\in {\cal O}({\rm Y}_1)$ or $h'\not\in {\cal O}({\rm Y}_2)$. Correspondingly $f_1(h')=0$ or $f_2(h')=0$, whence $f_3(h')=0$, and so $f_3(g')=0=f_3(h')f_3(k')$. Thus $f_3$ is factorisable, and so is a character. 
\hfill 
$\Box$\;
 
\vskip1em 
After calculating all the characters explicitly, we see  in \S 22.3, for which $({\rm Y_1,Y_2,Y_3})$ the assertion (\ref{2017-09-15-1}) fails, or\, 
$E(G'_\infty;\chi^{{\rm Y}_1})E(G'_\infty;\chi^{{\rm Y}_2})\not\subset E(G'_\infty;\chi^{{\rm Y}_3})$.

\section{\large Infinite symmetric group ${\mathfrak S}_\infty$ and its covering group $\widetilde{{\mathfrak S}}_\infty$}

We shortly review some results for these groups, which are  necessary for later studies.  
 \vskip1.2em

{\large\bf 13.1. Characters of infinite symmetric group ${\mathfrak S}_\infty$}
\vskip.5em
\quad 
{\bf 13.1.1. Thoma's character formula.}\;  
E.\,Thoma gave in \cite{[Tho1964]} a general character formula for the infinite symmetric group    
$\mathfrak{S}_\infty=G(1,1,\infty)$. To do so, he first proved in [ibid, Satz 1] the criterion (EF), that is,  for $f\in K_1({\mathfrak S}_\infty)$, 
\vskip.5em

(EF)\hspace{10ex} 
$f$ is extremal 
\;$\Longleftrightarrow$\; $f$ is factorisable.
\vskip.8em 

 Prepare a parameter $(\alpha,\beta)$, called Thoma parameter, given as 
\begin{eqnarray}
\label{Thoma-1} 
&&
\left\{
\begin{array}{l}
\alpha = (\alpha_i)_{i \geqslant 1}, \quad 
\alpha_1 \ge \alpha_2 \ge \alpha_3 \ge\; \ldots\; \ge 0,
\\[.3ex] 
\beta =(\beta_i)_{i \geqslant 1}, \quad \beta_1 \ge \beta_2 \ge \beta_3 \ge\; \ldots\; \ge 0; 
\\[.3ex]
\|\alpha\| + \|\beta\| \le 1,\quad 
{\rm with} \;\; 
\|\alpha\| := {\large\sum}_{i \geqslant 1} \alpha_i,\;\; 
\|\beta\| := \large\sum_{i \geqslant 1} \beta_i. 
\end{array}
\right.
\end{eqnarray}
The set of all Thoma parameters is denoted by ${\cal A}$\,:
\begin{eqnarray}
\label{2016-03-28-1}
{\cal A}:=\big\{(\alpha,\beta)\,;\,\alpha = (\alpha_i)_{i \geqslant 1},\, \beta = (\beta_i)_{i \geqslant 1},\,
\|\alpha\|+\|\beta\|\le 1\big\}. 
\end{eqnarray}

\vskip.1em

{\bf Theorem 13.1.} (Thoma)\;  
{\it  
A character of ${\mathfrak S}_\infty$ is given by a parameter $(\alpha,\beta)\in{\cal A}$ as $f_{\alpha,\beta}$ given below.  
A non-trivial $\sigma\in{\mathfrak S}_\infty$ is decomposed into a product of disjoint cycles $\sigma_j$ as $\sigma=\sigma_1\sigma_2\cdots\sigma_m$, let $\ell_j=\ell(\sigma_j)$ be the length of cycle $\sigma_j$, then 
\begin{eqnarray}
\label{Thoma-2}
&&
f_{\alpha, \beta}(\sigma)= \prod_{j\in{\boldsymbol I}_m}f_{\alpha,\beta}(\sigma_j),\quad
f_{\alpha, \beta}(\sigma_j)=\sum_{i\geqslant 1} 
\alpha_i^{\;\ell_j} + (-1)^{\ell_j - 1}\sum_{i\geqslant 1} 
\beta_i^{\;\ell_j}.
\end{eqnarray}

}

The set $E({\mathfrak S}_\infty)$ of all characters is compact with the pointwise convergence topology and is homeomorphic to ${\cal A}$ with the coordinatewise topology, under the correspondence $E({\mathfrak S}_\infty)\ni f_{\alpha,\beta}\mapsto (\alpha,\beta)\in{\cal A}$. 
 
\vskip1.2em 
{\bf 13.1.2. Limits of irreducible characters along ${\mathfrak S}_n \nearrow {\mathfrak S}_\infty$.}
 
Infinite symmetric group ${\mathfrak S}_\infty$ is an inductive limits of $n$-th symmetric groups as $n\to\infty$. The dual $\widehat{{\mathfrak S}_n}$ is parametrized by the set ${\boldsymbol Y}_{\!n}$ of Young diagrams of size $n$ 
\begin{eqnarray}
\label{Young-1}
\left\{
\begin{array}{l}
{\boldsymbol \lambda}^{(n)} = \big(\lambda^{(n)}_1, \lambda^{(n)}_2, 
\ldots, \lambda^{(n)}_n\big)\in P_n, 
\\[.5ex] 
 \lambda^{(n)}_1 \ge \lambda^{(n)}_2 \ge  
\ldots \ge \lambda^{(n)}_n\ge 0,\;\; |{\boldsymbol \lambda}^{(n)}|=n. 
\end{array}
\right.
\end{eqnarray}
Take a complete set of representatives (=\,CSR)  $\{\pi_{{\boldsymbol \lambda}^{(n)}}\,;\,{\boldsymbol \lambda}^{(n)}\in{\boldsymbol Y}_{\!n}\}$, and their normalized characters   
 $\widetilde{\chi}(\pi_{{\boldsymbol \lambda}^{(n)}}|\hspace{.2ex} \sigma)\; (\sigma \in {\mathfrak S}_n)$, For $1\le k\le n$, let the length of $k$-th row and that of $k$-th column of ${\boldsymbol \lambda}^{(n)}$ be respectively $r_k({\boldsymbol \lambda}^{(n)})$ and $c_k({\boldsymbol \lambda}^{(n)})$. Then,  
$r_k({\boldsymbol \lambda}^{(n)}) = \lambda^{(n)}_k$\, and 
\begin{eqnarray}
\label{Young-2}
\sum_{1 \leqslant k \leqslant n}\,r_k({\boldsymbol \lambda}^{(n)}) = n, 
\qquad 
\sum_{1 \leqslant k \leqslant n}\,c_k({\boldsymbol \lambda}^{(n)}) = n.
\end{eqnarray}

A.\,Vershik and S.\,Kerov gave in [38, Theorems 1, 2] 
 and in [39, Theorem 1] 
 (cf. also \cite{[Ker2003]}) very interesting result, which is the origin of the so-called asymptotic theory in group representations, as   
\vskip1em

{\bf Theorem 13.2.} (Vershik-Kerov)\;  
{\it 
Pointwise limit\,   
${\displaystyle 
\lim_{n \to \infty} \widetilde{\chi}(\pi_{{\boldsymbol \lambda}^{(n)}}; \sigma) \;\;(\sigma \in {\mathfrak S}_\infty)
}$\, exists if and only if the following limits exist  
\begin{eqnarray}
\label{Young-3}
\qquad
\lim_{n \to \infty} \frac{c_k({\boldsymbol \lambda}^{(n)})}{n} = 
\alpha_k, \quad 
\lim_{n \to \infty} \frac{r_k({\boldsymbol \lambda}^{(n)})}{n} = 
\beta_k\quad(k = 1, 2, \ldots). 
\end{eqnarray}
In that case, the limit is equal to the character  $f_{\alpha,\beta}$ of ${\mathfrak S}_\infty$ with 
$\alpha = (\alpha_k)_{k \geqslant 1}$ and $\beta = (\beta_k)_{k \geqslant 1}$. 

}

\vskip1.2em

The key of the proof of this theorem is the following approximative evaluation coming from the character formula of F.\;Murnaghan (\cite{[Mur1963]}, \cite{[VK1981]})\,:  

\vskip1em 

{\bf Lemma 13.3.}\; 
Let $\widetilde{\chi}(\pi_{{\boldsymbol \lambda}^{(n)}}|\hspace{.2ex}\cdot)$ be the normalized character of IR $\pi_{{\boldsymbol \lambda}^{(n)}}$ of ${\mathfrak S}_n$.  Let $\sigma^{(\ell)}\in{\mathfrak S}_n$ be a cycle of length $\ell\ge 2$, then there holds the following approximative evaluation along $n \to \infty$\,: with Frobenius parameter $a^{(n)}_k := r_k({\boldsymbol \lambda}^{(n)})- k\ge 0,\,
b^{(n)}_k:= c_k({\boldsymbol \lambda}^{(n)})- k\ge 0\;(1\le k\le \exists R_n)$, 
 \begin{eqnarray}
 \label{F-Murnaghan-1}
&&
 \quad
 \widetilde{\chi}\big({\boldsymbol \lambda}^{(n)}| \sigma^{(\ell)}\big) 
 = \sum_{1\leqslant k\leqslant R_n} \!\bigg(\dfrac{a^{(n)}_k}{n}\bigg)^{\!\ell} + (-1)^{\ell -1}
 \sum_{1\leqslant k\leqslant R_n} \!\bigg(\dfrac{b^{(n)}_k}{n}\bigg)^{\!\ell}  
 \; + \; O \!\bigg(\dfrac{1}{n}\bigg). 
 \end{eqnarray}
 \vskip.2em
 
In \cite{[Hir2004]}, we gave another kind of approximation of $f_{\alpha,\beta}$ as $n\to\infty$.  

\vskip1.2em

{\large\bf 13.2. Characters of universal covering group $\widetilde{{\mathfrak S}}_\infty$ of ${\mathfrak S}_\infty$}
\vskip.5em 
\quad 
{\bf 13.2.1. Products of two characters of $\widetilde{{\mathfrak S}}_\infty$.}\;  
The inductive limit of representation groups $\widetilde{{\mathfrak S}}_n$ of ${\mathfrak S}_n$ is the double covering group $\widetilde{{\mathfrak S}}_\infty:=\bigcup_{n\ge 4}\widetilde{{\mathfrak S}}_n$ of ${\mathfrak S}_\infty$ such as \;$\{e\}\to Z_1\to \widetilde{{\mathfrak S}}_\infty\stackrel{\Phi}{\to}{\mathfrak S}_\infty\to\{e\}$, where $Z_1=\langle z_1\rangle$. Characters of the central subgroup $Z_1$ of $\widetilde{{\mathfrak S}}_\infty$ are $\chi_+={\boldsymbol 1},\,\chi_-=\mathrm{sgn}_{Z_1}$ and so the set   
 $E(\widetilde{{\mathfrak S}}_\infty)$ of characters is divided into two subset according to spin types $\chi\in\widehat{Z_1}$. Put $E(\widetilde{{\mathfrak S}}_\infty;\chi):=\{f\in E(\widetilde{{\mathfrak S}}_\infty)\,;\,f(z_1\sigma')=\chi(z_1)f(\sigma')\;(\sigma'\in\widehat{{\mathfrak S}}_\infty)\}$, called {\it non-spin} or {\it spin} according as $\chi=\chi_+$ or $\chi_-$. Then they are the sets of all normalized characters of type II$_1$ non-spin or spin factor representations, 
 except two linear characters ${\bf 1}_{\mathfrak S}, \mathrm{sgn}_{\mathfrak S}$, and 
 \begin{eqnarray}
\label{2015-05-23-1}
&&
E\big(\widetilde{{\mathfrak S}}_\infty\big)={\bigsqcup}_{\chi\in\widehat{Z_1}}\, E\big(\widetilde{{\mathfrak S}}_\infty;\chi\big),\quad E\big(\widetilde{{\mathfrak S}}_\infty;\chi_+\big)\cong E({\mathfrak S}_\infty).
\end{eqnarray}

A $\sigma'\in\widetilde{{\mathfrak S}}_\infty$ is called {\it of the second kind}\, if $\sigma'$ is not conjugate to $z_1\sigma'$ (under $\widetilde{{\mathfrak S}}_\infty$), otherwise it is called {\it of the first kind}. If $\sigma'$ is odd, then it is of the first kind. In fact, take another odd $\tau'$ such that  $\mathrm{supp}(\tau')\cap\mathrm{supp}(\sigma')=\emptyset$, then $\tau'\sigma'{\tau'}^{\,-1}=z_1\sigma'$. The set of all elements of the second kind is denoted by ${\cal O}(\widetilde{{\mathfrak S}}_\infty;\chi_-)$. 
A decomposition $\sigma'=\sigma'_1\sigma'_2\cdots\sigma'_t$  of $\sigma'\in\widetilde{{\mathfrak S}}_\infty$ is called {\it cycle decomposition} if it becomes a cycle decomposition of $\sigma=\Phi(\sigma')$ on the level of ${\mathfrak S}_\infty$ through $\Phi_{\mathfrak S}:\,\widetilde{{\mathfrak S}}_\infty\to{\mathfrak S}_\infty$. 

\vskip1em

{\bf Lemma 13.4.}  
{\it 
A $\sigma'\in\widetilde{{\mathfrak S}}_\infty$ is of the second kind 
 if $\sigma'$ has a cycle decomposition with components $\sigma'_j$ all of odd lengths. 
 Let $\sigma'_{\boldsymbol \nu}\;({\boldsymbol \nu}\in P_\infty)$ (cf. \S 4.1) be a standard representative of a conjugacy class of $\widetilde{{\mathfrak S}}_\infty$ modulo $Z_1$, then $\sigma'_{\boldsymbol \nu}\in {\cal O}(\widetilde{{\mathfrak S}}_\infty;\chi_-)$ if and only if ${\boldsymbol \nu}\in O\!P_n$. 
 \hfill $\Box$
}

\vskip1em

{\bf Proposition 13.5.}\; 
{\it 
{\rm (i)}\; 
For any $f\in E(\widetilde{{\mathfrak S}}_\infty;\chi_-)$, \;$\mathrm{supp}(f)\subset {\cal O}(\widetilde{{\mathfrak S}}_\infty;\chi_-)\subset \widetilde{{\mathfrak A}}_\infty$. 

{\rm (ii)}\; For an $f\in K_1(\widetilde{{\mathfrak S}}_\infty;\chi):=\big\{f\in K_1(\widetilde{{\mathfrak S}}_\infty)\,;\,\mbox{\rm $f$ is of spin type $\chi$}\big\}$, 
\vskip.5em
\hspace*{1ex} {\bf (EF)} \hspace{10ex} 
\mbox{\it $f$ is extremal \;$\Longleftrightarrow$\; $f$ is factorisable.} 

\vskip.5em

{\rm (iii)}\; The set $E(\widetilde{{\mathfrak S}}_\infty)$ of all characters is closed under product, and for $f_1\in E(\widetilde{{\mathfrak S}}_\infty;\chi_1),
\\ 
f_2\in E(\widetilde{{\mathfrak S}}_\infty;\chi_2)$, we have\, $f_1f_2\in E(\widetilde{{\mathfrak S}}_\infty;\chi_1\chi_2)$. 

}

\vskip1em 

{\it Proof.} 
(i)\; Since \,$f(z_1\sigma')=\chi_-(z_1)f(\sigma')=-f(\sigma')$, if $\sigma'$ is of the first kind, then $f(z_1\sigma')=f(\sigma')$, and so $f(\sigma')=0$.

(ii)\; For $\chi=\chi_+$, this is already known. For 
$\chi=\chi_-$, since $\mathrm{supp}(f)\subset{\cal O}(\widetilde{{\mathfrak S}}_\infty;\chi_-)$, a cycle decomposition of $\sigma'\in\mathrm{supp}(f)$ is given as $\sigma'=\sigma'_1\sigma'_2\cdots\sigma'_t$ with $\mathrm{sgn}(\sigma'_j)=1$. Then $\sigma'_j$'s are mutually commutative and the factorisability\, $f(\sigma')=\prod_{1\leqslant j\leqslant t}f(\sigma'_j)$ can be proved just as for ${\mathfrak S}_\infty$. 

(iii)\; The product $f=f_1f_2$ is central, positive- definite and factorisable. By the criterion (EF), $f$ is extremal. 
\hfill 
$\Box$\quad\;

\vskip1.2em

{\bf 13.2.2. Spin characters of $\widetilde{{\mathfrak S}}_\infty$ (Nazarov's character formula).}

Schur's Hauptdarstellung $\Delta'_n,\,n\ge 4,$ is defined in 
\S 4.1 by the formula (\ref{2017-06-19-1})--(\ref{2017-06-19-2}), and its character is given in  
Theorem 4.1. So, its normalized character is, for a standard representative $\sigma'_{\boldsymbol \nu}\;\big( 
{\boldsymbol \nu}=(\nu_p)_{p\geqslant 1}\in O\!P_n\big)$,\; 
$
\widetilde{\chi}(\Delta'_n|\sigma'_{\boldsymbol \nu})={\prod}_{p\geqslant 1}(-2)^{-(\nu_p-1)/2}. 
$

\vskip1em 

{\bf Theorem 13.6.}\; 
{\it 
The pointwise limit function $\psi_\Delta(\cdot):=\lim_{n\to\infty} \widetilde{\chi}(\Delta'_n|\,\cdot)$ is given as follows: 
for $\sigma'_{\boldsymbol \nu}\in\widetilde{{\mathfrak S}}_n\subset\widetilde{{\mathfrak S}}_\infty$ with ${\boldsymbol \nu}=(\nu_p)_{p\geqslant 1}\in O\!P_n$\,, 
\begin{eqnarray}
\label{2017-08-05-1}
\psi_\Delta(\sigma'_{\boldsymbol \nu}) =\prod_{p\geqslant 1}(-2)^{-(\nu_p-1)/2}. 
\end{eqnarray}

}

As the support of $\psi_\Delta$, we have a subset of $\widetilde{{\mathfrak S}}_\infty$ as 
\begin{eqnarray}
\label{2017-08-31-1}
&&
{\cal O}\big(\widetilde{{\mathfrak S}}_\infty;O\!P_\infty\big):=\{z_1^{\,a}\sigma'_{\boldsymbol \nu}\,;\,a=0,1,\,{\boldsymbol \nu}\in O\!P_\infty\},\;\;O\!P_\infty:={\bigcup}_{n\geqslant 1}O\!P_n\,.
\end{eqnarray}

A multiplication map $\varphi:\,f\mapsto \psi_\Delta\!\cdot\! f$ on $E(\widetilde{{\mathfrak S}}_\infty)$ gives into-maps as 
\begin{eqnarray}
\label{2015-05-26-1}
E(\widetilde{{\mathfrak S}}_\infty;\chi_+)\;
\begin{array}{c}
\stackrel{\varphi}{\longrightarrow}
\\[-2.2ex]
\longleftarrow
\\[-2.2ex]
\mbox{\rm \footnotesize$\varphi$}
\end{array}
\;E(\widetilde{{\mathfrak S}}_\infty;\chi_-).
\end{eqnarray}

We prepare, as a parameter space for $E(\widetilde{{\mathfrak S}}_\infty;\chi_-)$,  
\begin{eqnarray}
\label{2016-03-28-11}
&&
{\cal C}:=\big\{\gamma=(\gamma_i)_{i\geqslant 1}\,;\,\gamma_1\ge\gamma_2\ge\gamma_3\ge \ldots\ge 0,\;\|\gamma\|={\sum}_{i\geqslant 1}\gamma_i\le 1\big\}.
\end{eqnarray}
For $\gamma\in{\cal C}$, put $\tfrac{1}{2}\gamma:=(\tfrac{1}{2}\gamma_i)_{i\geqslant 1}$. For $(\alpha,\beta)\in{\cal A}$,  we define $\alpha\!\vee\!\beta$ by collecting all components $\alpha_j, \beta_k$ of $\alpha,\beta$ and then reordering them according to the order $\ge$\,. 
Recall that, for an interval $K=[1,\ell]\subset {\boldsymbol N}$ with 
$\ell>1$, we put $\sigma'_K=r_1r_2\cdots r_{\ell-1}\in \widetilde{{\mathfrak S}}_\infty$. 
\vskip1em

{\bf Lemma 13.7.}\; 
{\it 
Let $f_{\alpha,\beta},f_{\gamma,{\bf 0}}\in E({\mathfrak S}_\infty)=E(\widetilde{{\mathfrak S}}_\infty;\chi_+)\,$ with ${\bf 0}=(0,0,\ldots)$. Put
\begin{eqnarray}
\label{2015-07-14-21}
\psi_\gamma:=\psi_\Delta\!\cdot\!f_{\gamma,{\boldsymbol 0}}=\varphi(f_{\gamma,{\bf 0}}),\quad \gamma\in{\cal C}.
\end{eqnarray}
Then, $\psi_\gamma\big(\sigma'_{[1,\ell]}\big)=
(-2)^{-(\ell-1)/2}\cdot{\sum}_{i\geqslant 1}\gamma^{\;\ell}_i\;\;\;(\ell>1\;{\rm odd}),
$ and 
\begin{eqnarray}
\nonumber
&&
\varphi(f_{\alpha,\beta})=\psi_\Delta\!\cdot\!f_{\alpha,\beta}=
\psi_\Delta\!\cdot\!f_{\alpha\vee\beta,{\bf 0}}=\psi_{\alpha\vee\beta}\in E(\widetilde{{\mathfrak S}}_\infty;\chi_-),
\\[1ex] 
&&
\nonumber
\label{2015-05-26-12}
\varphi^2(f_{\alpha,\beta})=\varphi(\psi_\gamma)=f_{\frac{1}{2}\gamma,\frac{1}{2}\gamma}\in E({\mathfrak S}_\infty)\;\;\;{\rm with}\;\;\gamma=\alpha\vee\beta.
\end{eqnarray}
 }
 \vskip.2em
 The following theorem is due to the detailed study of Nazarov. 
 \vskip1em 
 
{\bf Theorem 13.8.}\;([31, Theorem 3.3])  
{\it 
The set of spin characters of $\widetilde{{\mathfrak S}}_\infty$ is 
$$
\qquad
E(\widetilde{{\mathfrak S}}_\infty;\chi_-)=\big\{
\psi_\gamma=\psi_\Delta\!\cdot\!f_{\gamma,{\bf 0}}\;;\;\gamma\in{\cal C}\big\}.
\qquad 
$$ 
}
\vskip.2em

{\bf 13.2.3. Limits of spin irreducible characters along $\widetilde{{\mathfrak S}}_n\nearrow\widetilde{{\mathfrak S}}_\infty$.} 

Spin IRs $\tau_{\boldsymbol \lambda}\;({\boldsymbol \lambda}=(\lambda_1,\lambda_2,\ldots,\lambda_l)\in S\!P_n)$ of $\widetilde{{\mathfrak S}}_n$ are given in (\ref{2016-05-09-11}) and the supports of their characters $\chi(\tau_{\boldsymbol \lambda}|\,\cdot)$ are evaluated in (\ref{2016-05-10-1}). 
Put 
$l({\boldsymbol \lambda}):=l$\, and $\varepsilon({\boldsymbol \lambda})=0, 1$ according as $d({\boldsymbol \lambda}):=n - l({\boldsymbol \lambda})$ is even or odd. 

On the other hand, for ${\boldsymbol \nu}=(\nu_1,\nu_2,\ldots,\nu_t)\in P_n$, we put ${\boldsymbol \nu}'=(\nu'_1,\nu'_2,\ldots,\nu'_t, \nu'_{t+1})\in P_{n+1}$ with $\nu'_j=\nu_j\;(j\in{\boldsymbol I}_t)$ and $\nu'_{t+1}=1$, then $P_n$ is embedded into $P_{n+1}$ through ${\boldsymbol \nu}\mapsto  {\boldsymbol \nu}'$.  
Put $O\!P_\infty :=\bigcup_{n\geqslant 1}O\!P_n$. Then a CSR of conjugacy classes modulo $Z_1$ in ${\cal O}(\widetilde{{\mathfrak S}}_\infty;\chi_-)$ is given by the set $\{\sigma'_{\boldsymbol \nu}\,;\,{\boldsymbol \nu}\in O\!P_\infty\}$. 
To simplify the signs of character values, we replace $\sigma'_{\boldsymbol \nu}$ by $t_{\boldsymbol \nu}$ given as follows: put $t_j=z_1^{\;j-1}r_j\;(j\ge 1)$ and for ${\boldsymbol \nu}=(\nu_1,\nu_2,\ldots,\nu_t)\in P_n$, put  
$K_1=[1,\nu_1],\;K_j=[\nu_1+\cdots+\nu_{j-1}+1, \nu_1+\cdots+\nu_j]\;\;(j\ge 2),$ and $t_K:= t_at_{a+1}\cdots t_{b-1}$ for $K=[a,b]$, and\;  
$t_{\boldsymbol \nu}:=t_{K_1}t_{K_2}\cdots t_{K_t}\,.$  Then, 
\begin{eqnarray}
\label{2015-07-18-1}
\big(\psi_\Delta f_{\gamma,{\bf 0}}\big)(t_{\boldsymbol \nu})=
\prod_{1\leqslant j\leqslant t}\Big(2^{-(\nu_j-1)/2}
\sum_{i\geqslant 1}\gamma_i^{\;\nu_j}\Big). 
\end{eqnarray}
\vskip.2em

Put $\chi^{\boldsymbol \lambda}_{\boldsymbol \nu}:=\chi(\tau_{\boldsymbol \lambda}|\hspace{.2ex}t_{\boldsymbol \nu})$ for ${\boldsymbol \lambda}\in S\!P_n$ and ${\boldsymbol \nu}\in O\!P_n$, and\,  
$X^{\boldsymbol \lambda}_{\boldsymbol \nu}:=2^{(l({\boldsymbol \lambda})+\varepsilon({\boldsymbol \lambda})-l({\boldsymbol \nu}))/2}\chi^{\boldsymbol \lambda}_{\boldsymbol \nu}.$ Then, from the dimension formula, we have for ${\boldsymbol \nu}=(1,1,\ldots,1)=(1^n)$, 
\begin{eqnarray}
\label{2015-07-18-17}
 X^{\boldsymbol \lambda}_{(1^n)}= \dfrac{ n!  }{\lambda_1!\,\lambda_2!\,\cdots\,\lambda_l!}\prod_{1\leqslant i<j\leqslant l}\dfrac{\lambda_i-\lambda_j}{\lambda_i+\lambda_j}\,.
\end{eqnarray}

Put $\widetilde{\chi}^{\boldsymbol \lambda}_{\boldsymbol \nu}:=\chi^{\boldsymbol \lambda}_{\boldsymbol \nu}/\dim \tau_{\boldsymbol \lambda}$. Nazarov's  
calculation of its limits, as \,$n=|{\boldsymbol \lambda}|\to\infty$, is interesting and we follow its main stream. 
The next proposition plays a similar role as Murnaghan's character formula in the non-spin case (cf. Lemma 13.3). 
\vskip1em 

{\bf Proposition 13.9.} (cf. \cite{[Mor1976]} and [31, Proposition 1.7])\; 
{\it Let\, 
${\boldsymbol \lambda}=(\lambda_1,\ldots,\lambda_l)\in S\!P_n$, 
${\boldsymbol \nu}=(\nu_1,\ldots,\nu_t)\in O\!P_n$, $k>1$ odd integer. Let\, $\nu^0\in O\!P_{n-k}$ be a partition of\, $n-k$ obtained from ${\boldsymbol \nu}$ one term $\nu_{j_0}=k$ (i.e., $\Phi_{\mathfrak S}\big(t_{\nu_{j_0}}\big)$ is a cycle of length\, $k$). If $k$ satisfies the condition 
\begin{eqnarray}
\label{2015-07-23-1}
\lambda_i-\lambda_{i+1}>k\;(1\le \forall i<l),\quad \lambda_l>k,
\end{eqnarray}
a value of irreducible spin character of\, $\widetilde{{\mathfrak S}}_n$ is given from those of\, $\widetilde{{\mathfrak S}}_{n-k}$ as follows: 
\begin{eqnarray}
\label{2015-07-23-2}
X^{\boldsymbol \lambda}_{\boldsymbol \nu}=\sum_{1\leqslant i\leqslant l} 
X^{(\lambda_1,\ldots,\lambda_{i-1},\lambda_i-k,\lambda_{i+1},\ldots\lambda_l)}_{{\boldsymbol \nu}^0}.
\end{eqnarray}
}

Apply this formula for ${\boldsymbol \nu}=(k,1,\ldots,1)=(k,1^{n-k})$ and use (\ref{2015-07-18-17}), then we have  
\vskip1em 

{\bf Lemma 13.10.}\;\; 
{\it Suppose ${\boldsymbol \lambda}\in S\!P_n$ and $k$ satisfies the condition (\ref{2015-07-23-1}). Then 
\begin{eqnarray}
\label{2015-07-23-11}
&&
\;\qquad\qquad
\dfrac{X^{\boldsymbol \lambda}_{(k,1^{n-k})}}{X^{\boldsymbol \lambda}_{(1^n)}}
=2^{(k-1)/2}\,\widetilde{\chi}^{\boldsymbol \lambda}_{(k,1^{n-k})}
=
\sum_{1\leqslant i\leqslant l}
A_i(n)B_i(n),
\\
\nonumber
&&
A_i(n)
=
\frac{\lambda_i(\lambda_i-1)\cdots(\lambda_i-k+1)}{n(n-1)\cdots(n-k+1)},\;\;
B_i(n)= 
\prod_{j\ne i}\frac{(\lambda_i-k-\lambda_j)(\lambda_i+\lambda_j)}{(\lambda_i-k+\lambda_j)(\lambda_i-\lambda_j)}
\end{eqnarray}
}

{\bf Lemma 13.11.}\;  
{\it 
For a fixed odd integer $k>1$,  
assume that 
${\boldsymbol \lambda}={\boldsymbol \lambda}(n)=\big(\lambda_1(n),\lambda_2(n), \\ 
\ldots, \lambda_{l(n)}(n)\big)\in S\!P_n,\;l(n)=l\big({\boldsymbol \lambda}(n)\big),$ for $n\to\infty$ satisfy the condition (\ref{2015-07-23-1}).  Suppose 
\begin{eqnarray}
\label{2015-07-23-21-2} 
\lim_{n\to\infty} \frac{\lambda_i(n)}{n}\,=\,\gamma_i\quad(\forall i\ge 1). 
\end{eqnarray}
Then,\;\; 
${\displaystyle 
\lim_{n\to\infty}A_i(n)=\gamma_i^{\;k}, \;\; 
\lim_{n\to\infty}\sum_{1\leqslant i\leqslant l(n)}\!A_i(n)B_i(n) =\sum_{i\geqslant 1}\gamma_i^{\;k}, 
}$\, and\, $\gamma=(\gamma_j)_{j\geqslant 1}\in{\cal C}$.  
}
\vskip.8em

By Proposition 13.5\,(ii), the criterion (EF) holds for $\widetilde{{\mathfrak S}}_\infty$. Hence we obtain 

\vskip1em 

{\bf Theorem 13.12.} ([31, Theorem 3.5])\; 
{\it 
A series of normalized characters $\widetilde{\chi}(\tau_{{\boldsymbol \lambda}(n)}|\,\cdot)$ of spin IRs $\tau_{{\boldsymbol \lambda}(n)},\,{\boldsymbol \lambda}(n)\in S\!P_n$, converges pointwise if and only if the limits (\ref{2015-07-23-21-2}) exist. In that case\;
${\displaystyle  
\lim_{n\to\infty}\widetilde{\chi}(\tau_{{\boldsymbol \lambda}(n)})=\psi_{\Delta}\cdot f_{\gamma,{\bf 0}}\,. }$ 
Hence\; 
$$
{\rm Lim}(\widetilde{{\mathfrak S}}_\infty;\chi_-)=\{\psi_{\Delta}\cdot f_{\gamma,{\bf 0}}\;;\,\gamma\in{\cal C}\}=E(\widetilde{{\mathfrak S}}_\infty;\chi_-). 
$$
 
}

\section{\large Characters of infinite generalized symmetric group $G(m,1,\infty)$}

The case of non-spin character, or of spin type $\chi^{\rm VIII}=(1,\,1,\,1)$, of infinite generalized symmetric group ${\mathfrak S}_\infty({\boldsymbol Z}_m)=G(m,1,\infty)$ can be discussed similarly as in the case of non-spin character of infinite symmetric group
${\mathfrak S}_\infty=G(1,1,\infty)$. Fixing $m$, we can prove  \,${\rm Lim}(G_\infty)=F(G_\infty)=E(G_\infty)$\, for $G_\infty:=G(m,1,\infty)$.
 \vskip1em

 First we quote from \S 3 the results on irreducible characters of 
$G_n=G(m,1,n)=D_n(T)\rtimes {\mathfrak S}_n, \,T={\boldsymbol Z}_m.$ 
Let ${\boldsymbol Y}:={\bigsqcup}_{n\geqslant 0}{\boldsymbol Y}_{\!n},$ be the set of all Young diagrams, and put for ${\cal K}=\widehat{T}$ or $\widehat{T}^0$, 
\begin{gather}
\label{2015-07-24-33}
{\boldsymbol Y}(\widehat{T}):={\bigsqcup}_{n\geqslant 0}{\boldsymbol Y}_{\!n}(\widehat{T}),
\\
\label{2015-07-24-32}
{\boldsymbol Y}_{\!n}(\widehat{T}):=\Big\{\Lambda^n=({\boldsymbol \lambda}^{n,\zeta})_{\zeta\in\widehat{T}}\,;\,{\boldsymbol \lambda}^{n,\zeta}\in{\boldsymbol Y},\;{\sum}_{\zeta\in\widehat{T}}\,|{\boldsymbol \lambda}^{n,\zeta}|=n\Big\}. 
\end{gather}

The dual $\widehat{G_n}$ of $G_n$ is parametrized by ${\boldsymbol Y}_{\!n}(\widehat{T})$. For 
$\Lambda^n=({\boldsymbol \lambda}^{n,\zeta})_{\zeta\in\widehat{T}}\in {\boldsymbol Y}_{\!n}(\widehat{T})$, let $\breve{\Pi}_{\Lambda^n}$ be IR constructed in \S 3.1, then its character is given in Theorem 3.4.
With the same notation as in \S 3.2, we write down its normalized character $\widetilde{\chi}(\breve{\Pi}_{\Lambda^n}|g)$ as follows. For $g=(d,\sigma)\in G_n=D_n\rtimes{\mathfrak S}_n$, let its standard decomposition be as in (\ref{2014-03-18-21-8}) 
 \begin{eqnarray}
\label{2014-03-18-21-98}
&&
g=(d, \sigma) = \xi_{q_1}\xi_{q_2} \cdots \xi_{q_r}g_1g_2 \cdots g_s,\;\;\xi_q=\big(t_q,(q)\big), \;\;g_j=(d_j,\sigma_j), 
\end{eqnarray}
and put $J={\boldsymbol I}_s,\,Q=\{q_1,q_2,\ldots,q_r\}$. 
Here $\sigma_j$ is a cycle of length $\ell_j:=\ell(\sigma_j)$, $\mathrm{supp}(d_j)\subset \mathrm{supp}(\sigma_j)$, and $K_j:=\mathrm{supp}(\sigma_j)\;(j\in J)$ are mutually disjoint. Then 
$g \in H_n=
D_n(T)\rtimes {\displaystyle {\prod}_{\zeta\in\widehat{T}}\,{\mathfrak S}_{I_{n,\zeta}}
}$ if and only if \,$K_j\subset I_{n,\zeta}\;(\exists \zeta\in\widehat{T})$ for each $j\in J$.  
 
\vskip1.2em 

{\bf Theorem 14.1.} 
{\it 
 The normalized character is given as follows. For $\in G_n$ not conjugate to an element of $H_n$,\, $\widetilde{\chi}\big(\breve{\Pi}_{\Lambda^n}|g\big)=0$. 
For  $g\in H_n$, 
\begin{eqnarray}
\label{2015-07-25-11}
\nonumber
\qquad
&&
\qquad
\widetilde{\chi}\big(\breve{\Pi}_{\Lambda^n}|g\big)
=
\sum_{{\cal Q},{\cal J}} 
c(\Lambda^n;{\cal Q},{\cal J};g) \widetilde{X}(\Lambda^n;{\cal Q},{\cal J};g), 
\\
\nonumber
&& 
c(\Lambda^n;{\cal Q},{\cal J};g)
=
\frac{(n-|\mathrm{supp}(g)|)!}{n!}\cdot
{\displaystyle 
\prod_{\zeta\in\widehat{T}}
\frac{|{\boldsymbol \lambda}^{n,\zeta}|!}{
\Big(|{\boldsymbol \lambda}^{n,\zeta}| - |Q_\zeta| - 
\sum_{j \in J_\zeta}\ell_j\Big)!
}
}\,,
\\[.3ex]
\label{2015-07-25-12}
&&
\widetilde{X}(\Lambda^n;{\cal Q},{\cal J};g)
=
\prod_{\zeta\in\widehat{T}} 
\Big(\prod_{q\in Q_\zeta}\zeta(t_q)\cdot\prod_{j\in J_\zeta}\zeta\big(P(d_j)\big)\cdot \widetilde{\chi}\big(\pi_{{\boldsymbol \lambda}^{n,\zeta}}|(\ell_j)_{j\in J_\zeta}\big)\Big), 
\end{eqnarray}
with $\ell_j:=\ell(\sigma_j)$ and $\widetilde{\chi}\big(\pi_{{\boldsymbol \lambda}^{n,\zeta}}|(\ell_j)_{j\in J_\zeta}\big):=\chi\big(\pi_{{\boldsymbol \lambda}^{n,\zeta}}|(\ell_j)_{j\in J_\zeta}\big)/\dim\pi_{{\boldsymbol \lambda}^{n,\zeta}}$, where the pair of ${\cal Q}=(Q_\zeta)_{\zeta\in\widehat{T}}$ and ${\cal J}=(J_\zeta)_{\zeta\in\widehat{T}}$ runs over partitions of $Q$ and $J$ satisfying 
\begin{eqnarray}
\label{2014-03-17-11-3-2}
\nonumber
\mbox{\bf (Condition QJ)}\qquad\qquad\quad
|Q_\zeta|+\sum_{j\in J_\zeta}\ell_j\le |{\boldsymbol \lambda}^{n,\zeta}|\;\;\;
(\zeta\in\widehat{T}). \qquad\qquad\qquad\qquad
\end{eqnarray}

}

Now we calculate limits of non-spin irreducible characters of $G(m,1,n)$ along  
$\Lambda^n=({\boldsymbol \lambda}^{n,\zeta})_{\zeta\in\widehat{T}}\in{\boldsymbol Y}_{\!n}(\widehat{T})\;(n\to\infty)$.
\vskip1em

{\bf Lemma 14.2.}\;\;  
{\it 
A necessary and sufficient condition for that 
$
{\displaystyle 
\lim_{n \to \infty} \widetilde{\chi}\big(\breve{\Pi}_{\Lambda^n}|(d, {\bf 1})\big)
}$ 
converges for any $d\in D_\infty(T)$\, is the existence of the limit 
\\[1ex]
{\rm {\sc (Condition\;I)}}\quad 
\qquad\qquad 
${\displaystyle 
B_\zeta := 
\lim_{n\to \infty} \frac{|{\boldsymbol \lambda}^{n,\zeta}|}{n}  \quad  
(\zeta \in \widehat{T}). 
}
$
\vskip.8em
\noindent
In that case, ${\displaystyle {\sum}_{\zeta\in\widehat{T}}\,B_\zeta=1}$\; and 
\begin{eqnarray}
\label{2015-07-25-2} 
&&
\quad
\lim_{n\to\infty} 
\widetilde{\chi}\big(\breve{\Pi}_{\Lambda^n}|(d, {\bf 1})\big)
=\prod_{q\in Q} F_1(t_q),\;\;F_1(t):=\sum_{\zeta \in \widehat{T}}\; 
B_\zeta\, \zeta(t)\;\;\;(t \in T).
\end{eqnarray}

}

 Under (Condition I), we assume for a series 
$\Lambda^n=({\boldsymbol \lambda}^{n,\zeta})_{\zeta\in\widehat{T}},\,n\to\infty,$ \vskip.7em
\noindent
{\sc (Condition\;I$\Lambda$)}\quad  
for any  
$\zeta \in \widehat{T}^+ 
:= \{\zeta \in \widehat{T}\,;\, 
 B_\zeta > 0\}$, the following limits exist   
\begin{eqnarray*}
\label{lim-char-12} 
\qquad\quad\;\;
\lim_{n \to \infty}
\frac{r_k({\boldsymbol \lambda}^{n,\zeta})}{|{\boldsymbol \lambda}^{n,\zeta}|} = \alpha'_{\zeta;k}\,, 
\quad
\lim_{n \to \infty}
\frac{c_k({\boldsymbol \lambda}^{n,\zeta})}{|{\boldsymbol \lambda}^{n,\zeta}|} = \beta'_{\zeta;k}\, 
\quad(1 \le k < \infty).
\end{eqnarray*} 

In that case, put 
\begin{eqnarray}
\label{alpha-beta-myu-1}
&&
\qquad
\left\{
\begin{array}{l} 
{\displaystyle 
\alpha_{\zeta,0;k} 
:= \lim_{n \to \infty}\dfrac{r_k({\boldsymbol \lambda}^{n,\zeta})}{n}\,, 
\quad 
\alpha_{\zeta,1;k} 
:= \lim_{n \to \infty}\dfrac{c_k({\boldsymbol \lambda}^{n,\zeta})}{n}\,, 
}
\\[1.5ex]
{\displaystyle 
\mu_\zeta := B_\zeta - \!\sum_{\,\varepsilon = 0, 1} 
\|\alpha_{\zeta,\varepsilon}\|\ge 0,
\quad 
  \alpha_{\zeta,\varepsilon} 
  := (\alpha_{\zeta,\varepsilon;i})_{i\geqslant 1}\;\;(\varepsilon=0,1).
  \qquad
  }
  \end{array}
  \right.
\end{eqnarray} 
Then, for $\zeta\in\widehat{T}^+$, $\alpha_{\zeta,0;k}=B_\zeta\,\alpha'_{\zeta;k},\;
\alpha_{\zeta,1;k} = B_\zeta\,\beta'_{\zeta;k}$\,.  
Put, for $\zeta\not\in\widehat{T}^+$,\,  
$\alpha_{\zeta,0;k} := 
\alpha_{\zeta,1;k} := 0,\, \mu_\zeta := 0.$
Then we have 
 $\alpha_{\zeta,\varepsilon}\;(\varepsilon=0,1)$ and  
${\displaystyle \mu = (\mu_\zeta)_{\zeta \in \widehat{T}}\,,}$ and 
\begin{eqnarray}
\nonumber
 \label{2015-07-27-12}
{\sum}_{\zeta \in \widehat{T}}\;
{\sum}_{\varepsilon = 0, 1}\, 
\|\alpha_{\zeta,\varepsilon}\|
\;+\; \|\mu\| = 1,\quad \|\mu\| = {\sum}_{\zeta \in \widehat{T}}\,\mu_\zeta\,.
\end{eqnarray}
Collecting these data together, we take 
\,${\displaystyle 
A := \big( (
\alpha_{\zeta,\varepsilon})_{(\zeta, \varepsilon) \in 
\widehat{T} \times \{0,  1\}} 
\,;\,
\mu 
\big)
}$\,  
as a parameter for the limit function \;$f_A:={\displaystyle \lim_{n\to\infty}\widetilde{\chi}\big(\breve{\Pi}_{\Lambda^n})}$. 
For a 
$g = (d, \sigma) \in G_\infty$, take its standard decomposition (\ref{2014-03-18-21-98}), then  
\begin{eqnarray}
\label{lim-char-15}
\qquad 
&&
f_A(g) 
:= 
\prod_{q \in Q}
\Bigg\{ \sum_{\zeta \in \widehat{T}}
\Bigg(
 \sum_{\varepsilon \in \{0,  1\}}
\, \|\alpha_{\zeta,\varepsilon}\|   
\; + \;\mu_\zeta\Bigg)
\zeta(t_{q})\Bigg\} 
\\
\nonumber
\label{lim-char-16}
&& 
\qquad
\times 
\prod_{j \in J} 
\Bigg\{ 
\sum_{\zeta \in \widehat{T}} \Bigg(
\sum_{\varepsilon \in \{0,  1\}} 
\;
\sum_{i \in \boldsymbol{N}}\; 
(\alpha_{\zeta,\varepsilon;i})^{\ell_j}
\chi_\varepsilon(\sigma_j)
\cdot 
\zeta\big(P(d_j)\big)
\Bigg)
\Bigg\}, 
\end{eqnarray}
where $\chi_\varepsilon(\sigma_j) = 
{\rm sgn}_\mathfrak{S}(\sigma_j)^\varepsilon = 
(-1)^{\varepsilon(\ell_j - 1)},  
\; 
\ell_j = \ell(\sigma_j),\; 
P(d_j) = \prod_{i\in K_i}\!t_i \in T.
$
\vskip1em 

{\bf Notation 14.3.} To express sets of characters, put for ${\cal K}=\widehat{T}$ of $\widehat{T}^0$,  
\begin{gather}
\label{2016-03-22-1}
A = \big( (
\alpha_{\zeta,\varepsilon})_{(\zeta, \varepsilon) \in 
{\cal K} \times \{0,  1\}} 
\,;\,
\mu \big)  
\\[.3ex]
\nonumber
\alpha_{\zeta,\varepsilon}=(\alpha_{\zeta,\varepsilon;i})_{i\geqslant 1},\;\;
\alpha_{\zeta,\varepsilon;1}\ge \alpha_{\zeta,\varepsilon;2}\ge \ldots\ge 0,\quad 
\mu=(\mu_\zeta)_{\zeta\in{\cal K}},\;\;\mu_\zeta\ge 0,
\\[.5ex]
\label{2016-03-08-1} 
{\cal A}({\cal K}):=\big\{A=\big( (\alpha_{\zeta,\varepsilon})_{(\zeta, \varepsilon) \in {\cal K} \times \{0,  1\}} 
\,;\,\mu \big)  
\,;\;\mbox{\rm satisfies (\ref{2017-09-14-11})}\big\},
\\[.5ex]
 \label{2017-09-14-11}
\sum_{(\zeta,\varepsilon) \in {\cal K}\times\{0,1\}}\! 
\|\alpha_{\zeta,\varepsilon}\|
+ \|\mu\| = 1,\;\|\alpha_{\zeta,\varepsilon}\|=
\sum_{i\geqslant 1}\alpha_{\zeta,\varepsilon;i},\; \|\mu\| = \sum_{\zeta \in {\cal K}}\,\mu_\zeta\,.\qquad
\end{gather}

\vskip.2em

{\bf Theorem 14.4.}\; 
{\it 
 A series of normalized irreducible characters  $\widetilde{\chi}\big(\breve{\Pi}_{\Lambda^n}|g\big)$ of  
$G(m,1,n)={\mathfrak S}_n(T),\, T={\boldsymbol Z}_m,$ converges pointwise as $n\to\infty$\, if and only if\, {\rm (Condition I) + (Condition I$\Lambda$)} holds. 
In that case, determine by  
(\ref{alpha-beta-myu-1}), a parameter 
\begin{gather}
\label{2015-07-26-1}
\nonumber
A = \big( (
\alpha_{\zeta,\varepsilon})_{(\zeta, \varepsilon) \in \widehat{T} \times \{0,  1\}} \,;\,\mu \big)\in {\cal A}(\widehat{T}),
\end{gather}
then the limit function is $f_A$ in (\ref{lim-char-15}). Thus there holds 
$$
{\rm Lim}\big(G(m,1,\infty)\big)=\big\{f_A\,;\,A \in  {\cal A}(\widehat{T})\big\}=E\big(G(m,1,\infty)\big).
$$
}

\vskip1.2em

 \setcounter{section}{14}
\section{\large Heredity from $R\big(G(m,1,\infty)\big)$ and the case of $R\big(G(m,p,\infty)\big)$.}

\quad 
{\bf 15.1. Restriction of characters to normal subgroups.}\; 

Let $G$ be a Hausdorff topological group and $N$ its normal subgroup with relative topology. Denote by $K(N,G)$ the set of central positive-definite continuous functions on $N$ invariant under conjugation of $G$, $K_1(N, G)$ its convex subset of $f\in K(N,G)$ normalized as $f(e)=1$, and $E(N,G)$ the set of extremal elements of $K_1(N,G)$.  
\vskip1em

{\bf Theorem 15.1.}\;\;[14, Theorem 14.1]\; 
{\it 
The restriction map  ${\rm Res}^G_N: K_1(G)\ni F\mapsto f=F|_N\in K_1(N,G)$ maps $E(G)$ into $E(N,G)$.
}

\vskip1.2em

We proved the surjectivity of ${\rm Res}^G_N$ for $G=\mathfrak{S}_\infty(T)$ with $T$ a compact group in [14, Theorem 15.1]. Furthermore we have the following  
\vskip1em

{\bf Theorem 15.2.}\; 
{\it 
Let $G$ be a discrete group and $N$ its normal subgroup, then the restriction map ${\rm Res}^G_N$ maps $E(G)$ onto $E(N,G)$.
}

\vskip1em

{\it Proof.}\; Take an 
$f \in E(N,G)$. Extend $f$ to the whole of $G$ by putting 0 outside of $N$ and denote it by $F_f$. Then it is in $K_1(G)$. On the other hand, for discrete group $G$, $K_{\le 1}(G)$ is weakly compact and Extr$\big(K_{\le 1}(G)\big)=E(G)\bigsqcup\{0\}$. Apply integral representation theorem of Bishop\,-\,de Leeuw [2,Theorem 5.6] to the weakly compact convex set $K_{\le 1}(G)$, then there exists a probability measure $\mu_f$ on $E(G)$ such that 
$$
F_f=\int_{E(G)}F\,d\mu_f(F),\;\;\;\mbox{\rm and so }\;\;\;f=\int_{E(G)}(F|_N)\,d\mu_f(F).
$$
On the other hand, by Theorem 15.1, $F|_N \in E(N,G)$ for any 
$F\in E(G)$. Hence the above formula expresses $f\in E(N,G)$ by means of $F|_N\in E(N,G)$. Since $f$ is an extremal point, for almost all $F$ (with respect to $\mu_f$), we have $F|_N=f$. 
\hfill 
$\Box$\; 

\vskip1.2em

Now recall that a generarized symmetric group $G(m,1,n)=D_n(T)\rtimes \mathfrak{S}_n,\,T=\boldsymbol{Z}_m$, has, as a mother group, her child groups $G(m,p,n),\,p|m,\,p>1,$ inside her as  
\begin{eqnarray}
\label{2018-03-26-11}
G(m,p,n):=\{g=(d,\sigma)\in D_n(T)\rtimes \mathfrak{S}_n\,;\,{\rm ord}(d)\equiv 0\;({\rm mod}\;p)\}. 
 \end{eqnarray}

Let $y_i\;(i\in{\boldsymbol I}_n)$ be a specified generator of $i$-th component $T_i$ of $D_n(T)=\prod_{i\in{\boldsymbol I}_n}T_i,\,T_i=T$, and put 
\begin{eqnarray}
\label{2018-03-27-1}
x_1:=y_1^{\;p}\;\;\mbox{\rm ($x_1=e$ if $p=m$)},\quad x_j:=y_1^{\,-1}y_j\quad(j\in{\boldsymbol I}_n,\,>1).
\end{eqnarray}
Then $G(m,p,n)$ is generated by $\mathfrak{S}_n$ and $\{x_j\;(j\in{\boldsymbol I}_n)\}$. 
Moreover, we have 
\vskip1.2em 

{\bf Theorem 15.3.} [I, Proposition 3.4]\; 
{\it 
Let $4\le n<\infty,\,p|m,\,p>1$, and put $q:=m/p$. As an abstract group, $G(m,p,n)$ is presented as a group given by  

\; $\bullet$ set of generators:\quad $\{s_1,s_2,\ldots, s_{n-1},\,x_1,x_2,\ldots,x_n\}$, 

\; $\bullet$ set of fundamental equations:  
\begin{eqnarray*}
{\rm (ii)}&&
\quad s_i^{\;2}=e\;(i\in{\boldsymbol I}_{n-1}),\;(s_is_{i+1})^3=e\;(i\in{\boldsymbol I}_{n-2}),\;s_js_k=s_ks_j\;(|j-k|\ge 2); 
\\[1ex]
{\rm (iii)}
&&
\quad
\left\{
\begin{array}{l}
x_1^{\;q}=e,   
\;\,
x_j^{\;m} = e \quad  (j\in{\boldsymbol I}_n,>\!1),
\\
x_j x_k = x_k x_j \qquad\quad(j\ne k)\,;
\end{array}
\right. 
\\
\mbox{\rm (iv-1)}
&&
\quad
\left\{
\begin{array}{l}
s_ix_is_i^{\;-1} = x_{i+1},\;
\,
s_ix_{i+1}s_i^{\;-1} = x_i\;\;(i\in{\boldsymbol I}_{n-1},>\!1),
\\
s_ix_js_i^{\;-1} = x_j \quad  (j\ne i, i\!+\!1,\;i\in{\boldsymbol I}_{n-1},>\!1,\;j\in{\boldsymbol I}_n),
\end{array}
\right.
\\
\mbox{\rm (iv-2)}
&&
\quad
\left\{
\begin{array}{l}
s_1x_1s_1^{\;-1} = x_1x_2^{\;p},\;
\,
s_1x_2s_1^{\;-1} = x_2^{\;-1},
\\
s_1x_j s_1^{\;-1} = x_2^{\;-1}x_j \qquad (j\in{\boldsymbol I}_n,>\!2).
\end{array}
\right.
\end{eqnarray*}

}

To be simple, we consider the limiting case as $n\to\infty$, and let $G$ be $G(m,1,
\infty)$, and $N$ be one of its canonical normal subgroups in (\ref{2017-11-24-11})
\begin{eqnarray*}
\label{2017-09-21-1}
\hspace*{2.5ex} 
G(m,p,\infty),\;p|m,p>1,\;\; 
G^{\mathfrak A}_\infty:=D_\infty\rtimes{\mathfrak A}_\infty,\;\;
 G^{\mathfrak A}(m,p,\infty)=G^{\mathfrak A}_\infty\cap G(m,p,\infty). 
 \end{eqnarray*}
In this case, 
 the restriction map\, ${\rm Res}^G_N:\,F\mapsto f=F|_N$ for $F\in E(G)$ is studied in [14, \S 7 and \S\S 15--16] 
 and there proved that it is surjective onto $E(N)$. 
 We call this kind of property as {\it heredity}\, (from $G$ to $N$). 
Below in \S 15.6 we study how to express the heredity by means of a finite group action on the parameter space (see also Lemma 17.3 and Theorem 17.4). 

\vskip1.2em

{\bf 15.2. Representation groups of child groups.}

This  kind of heredity holds also in 
the spin case. Recall the exact sequence (\ref{2017-11-16-1}) defining the representation group $R(G(m,1,n)$ for $n\ge 4$: 
\begin{eqnarray}
\label{2018-03-26-1}
&&
\{e\}\longrightarrow Z\longrightarrow R\big(G(m,1,n)\big)\stackrel{\Phi\,}{\longrightarrow} G(m,1,n)\longrightarrow \{e\}\quad({\rm exact}),
\end{eqnarray}
where\;  $Z={\mathfrak M}(G_n)=\left\{
\begin{array}{ll}
Z_1,\quad & \mbox{\rm for $m$ odd}, 
\\
Z_1\times Z_2\times Z_3,\quad  & \mbox{\rm for $m$ even},
\end{array}
\right.
\quad Z_i=\langle z_i\rangle,\,z_i^{\;2}=e. $

In case $m$ is odd, 

\quad 
$R\big(G(m,1,n)\big)=D_n(T)\rtimes \widetilde{\mathfrak{S}}_n,\,D_n(T)=\langle \eta_1,\eta_2,\ldots,\eta_n\rangle,\,\eta_j\eta_k=\eta_k\eta_j,\,\eta^{\;m}_j=e$. 

In case $m$ is even, 

\quad  $R\big(G(m,1,n)\big)=\widetilde{D}^\vee_n(T)\rtimes \widetilde{\mathfrak{S}}_n,\,\widetilde{D}^\vee_n(T):=\widetilde{D}_n(T)\times Z_3\cong D^\wedge_n(T)\times Z_3,\,$

\quad $\widetilde{D}_n(T)=\langle z_2,\eta_1,\eta_2,\ldots,\eta_n\rangle,\,\eta_j\eta_k=z_2\eta_k\eta_j\;(j\ne k),\,\eta^{\;m}_j=e$,

\quad $D^\wedge_n(T)=\langle z_2,\widehat{\eta}_1,\widehat{\eta}_2,\ldots,\widehat{\eta}_n\rangle$ with \;$\widehat{\eta}_k:=z_3^{\;k-1}\eta_k,\;
r_i\widehat{\eta}_kr_i^{\,-1}=z_3\widehat{\eta}_{s_ik}$.

To give representation groups of $G(m,p,n)$, we separate cases as follows:\; 

\quad 
Case {\bf OO} for $p$ odd, $q$ odd\;\;($m$ odd); \;\; Case {\bf OE} for $p$ odd, $q$ even; \; 

\quad 
Case {\bf EO} for $p$ even, $q$ odd; \hspace{10.5ex} Case {\bf EE} for $p$ even, $q$ even.

\vskip1em

From [32], we see that the Schur multiplier $Z=\mathfrak{M}\big(G(m,p,n)\big)$ for $5\le n<\infty$ is given as follows: \,${\mathfrak M}\big(G(m,p,n)\big)={\boldsymbol Z}_2^{\;\ell(m,p,n)}$, 

in Case {\bf OO}, \,$Z={\boldsymbol Z}_2,\,\ell(m,p,n)=1,$ and we put $Z=Z_1$, with $Z_i=\langle z_i\rangle,\,z_i^{\;2}=e$, 

in Cases {\bf EO},\, $Z={\boldsymbol Z}_2^{\;2},\,\ell(m,p,n)=2,$ and we put $Z=Z_1\times Z_2$,

in Cases {\bf OE} and {\bf EE}, \,$Z={\boldsymbol Z}_2^{\;3},\,\ell(m,p,n)=3,$ and we put $Z=Z_1\times Z_2\times Z_3$.
\vskip1.2em

{\bf Theorem 15.4.} (Case OO) [I, Theorem 3.5]\, 
{\it 
 Let $5\le n<\infty$, and $m$ be odd, and $p|m,\,p>1$. Then a representation group $R\big(G(m,p,n)\big)$ is presented as

$\bullet$ set of generators: \quad $\{z_1,\,r_1,r_2,\ldots, r_{n-1},\,w_1,w_2,\ldots, w_n\}$; 

$\bullet$ canonical homomorphism $\Phi$: 
$\Phi(r_i) = s_i\;\;(i\in{\boldsymbol I}_{n-1}),\; \Phi(w_j) =x_j\;\;(j\in{\boldsymbol I}_n)$\,;

$\bullet$ set of fundamental relations: 
\begin{eqnarray*}
{\rm (i)} &&\quad z_1^{\;2}=e, \quad\mbox{\rm $z_1$ central element}, 
\\
{\rm (ii)} &&\quad r_i^{\;2} = e\;(i\in{\boldsymbol I}_{n-1}),\;\,
(r_ir_{i+1})^3 = e\;(i\in{\boldsymbol I}_{n-2}),
\;\,
r_ir_j = z_1r_jr_i \;\;(|i-j| \ge 2),
\\
{\rm (iii)} &&\quad 
\left\{
\begin{array}{l}
w_1^{\;q} = e, \;\;
w_j^{\;m} = e \;\;(j\in{\boldsymbol I}_n,>\!1),
\\
w_j w_k = w_k w_j \qquad\;(j\ne k);
\end{array}
\right. 
\\
\mbox{\rm (iv-1)}&&\quad 
{\rm for\;\,} i>1,\;\;  
\left\{
\begin{array}{l}
r_iw_ir_i^{\;-1} = w_{i+1},\;\;  
r_iw_{i+1}r_i^{\;-1} = w_i,
\\
r_iw_jr_i^{\;-1} = w_j \quad  (j\in{\boldsymbol I}_n,\,j \ne i, i\!+\!1);
\end{array}
\right.
\\
\mbox{\rm (iv-2)}&& \quad
\left\{
\begin{array}{l}
r_1w_1r_1^{\;-1} = w_1\,w_2^{\;p},\;\; 
r_1\,w_2\,r_1^{\;-1}\, =\,w_2^{\;-1},
\\ 
r_1w_j r_1^{\;-1} =w_2^{\;-1}w_j\, \qquad (j\in{\boldsymbol I}_n,\,j\!>\!2 ).
\end{array}
\right.
\end{eqnarray*}
}

\vskip.5em 
By (ii), $\langle r_i\;(i\in{\boldsymbol I}_{n-1})\rangle =\widetilde{\mathfrak{S}}_n,$\, and by 
(iii), \,$\langle w_j\;(j\in{\boldsymbol I}_n)\rangle \cong \langle x_j\;(j\in{\boldsymbol I}_n)\rangle\cong D_n({\boldsymbol Z}_m)^{S(p)}:=\{d\in D_n({\boldsymbol Z}_m)\,;\,P(d)\in S(p)\}$. 
The representation group $R\big(G(m,p,n)\big)$ can be expressed as a semidirect product as 
\vskip.3em
\hspace*{20ex} 
$
R\big(G(m,p,n)\big) \cong D_n({\boldsymbol Z}_m)^{S(p)}\rtimes \widetilde{\mathfrak{S}}_n.
$

\vskip1.2em

{\bf Theorem 15.5.} (Cases EO, OE, EE)\, [I, Theorems 3.6, 3.7, 3.8] \,
{\it 
Let $5\le n <\infty$, and $m$ be even, and $p|m,\,p>1$. Then a representation group $R\big(G(m,p,n)\big)$ is presented as

$\bullet$ set of generators: \quad $\{z_i\;(1\le i\le \ell(m,p,n)),\,r_1,r_2,\ldots, r_{n-1},\,w_1,w_2,\ldots, w_n\}$; 

$\bullet$ canonical homomorphism $\Phi$: 
$\Phi(r_i) = s_i\;\;(i\in{\boldsymbol I}_{n-1}),\; \Phi(w_j) =x_j\;\;(j\in{\boldsymbol I}_n)$\,; 

$\bullet$ set of fundamental relations: 
\begin{eqnarray*}
{\rm (i)} &&\quad z_i^{\;2}=e\;(1\le i\le \ell(m,p,n)), \;\;\mbox{\rm $z_i$ central element}, 
\\
{\rm (ii)} &&\quad \mbox{\rm the same as above,}
\\
{\rm (iii)} &&\quad 
\left\{
\begin{array}{l}
w_1^{\;q} = e, 
\;\;
w_j^{\;m} = z_2^{\;m/2} \,\;(j\in{\boldsymbol I}_n,>\!1),
\\
w_j w_k = z_2w_k w_j \;\,(j\ne k).
\end{array}
\right. 
\\
\mbox{\rm (iv-1)} &&\quad 
\mbox{\rm For $i>1$}, \;\;
\left\{
\begin{array}{ll}
r_iw_ir_i^{\;-1} = w_{i+1}\,,  
\;\,
r_iw_{i+1}r_i^{\;-1} = w_i,
\\
r_iw_1r_i^{\;-1} =z_3^{\;q-1}w_1\,,
\\
r_iw_jr_i^{\;-1} = w_j \;\;  (j \ne i, i\!+\!1,\;j\in{\boldsymbol I}_n,>\!1),
\end{array}
\right.
\\
\mbox{\rm (iv-2)} &&\quad 
\left\{
\begin{array}{l}
r_1w_1r_1^{\;-1} = z_2^{\;[p/2]}z_3^{\;q-1}w_1\,w_2^{\;p}, 
\;\;
r_1\,w_2\,r_1^{\;-1}\, =\,w_2^{\;-1},\quad
\\
r_1w_j r_1^{\;-1} = w_2^{\;-1}w_j \qquad\qquad (j\in{\boldsymbol I}_n,>\!2).
\end{array}
\right.
\end{eqnarray*}
\vskip.5em

Note that, in Case EO, since $q$ is odd, $z_3^{\;q-1}=e$ in {\rm (iv-2)} and so $z_3$ disappears. 

The representation group $R\big(G(m,p,n)\big)$ can be expressed as a semidirect products according to Case EO and Cases OE and EE respectively as 
\vskip.5em

$R\big(G(m,p,n)\big)\cong \widetilde{D}_n({\boldsymbol Z}_m)^{S(p)}\rtimes\widetilde{\mathfrak{S}}_n\,,$ \; 
$R\big(G(m,p,n)\big)\cong \big(\widetilde{D}_n({\boldsymbol Z}_m)^{S(p)}\times Z_3\big)\rtimes\widetilde{\mathfrak{S}}_n\,.$ 

}

\vskip.8em

For infinite generalized symmetric group $G(m,1,\infty)=\mathfrak{S}_\infty({\boldsymbol Z}_m),$ and its normal subgroup $G(m,p,\infty)=\mathfrak{S}_\infty({\boldsymbol Z}_m)^{S(p)},$ we take the inductive limits 
\begin{eqnarray}
\label{2014-03-07-11}
\left\{
\begin{array}{l}
{\displaystyle 
R\big(G(m,1,\infty)\big):=\lim_{n\to\infty}R(\big(G(m,1,n)\big),} \;\;
\\
{\displaystyle 
R\big(G(m,p,\infty)\big):=\lim_{n\to\infty}R(\big(G(m,p,n)\big),}
\end{array}
\right.
\end{eqnarray}

Then, $R\big(G(m,1,\infty)\big)$ is represented as in Theorem 15.2 by just putting $n=\infty$, that is, the symbols ${\boldsymbol I}_n,\,{\boldsymbol I}_{n-1},\,{\boldsymbol I}_{n-2}$ are all replaced by ${\boldsymbol I}_\infty={\boldsymbol N}$. Moreover, to represent 
$R\big(G(m,p,\infty)\big),\,p|m,\,p>1,$ by means of a set of generators and a set of fundamental relations, we just replace $n$ in Theorem 15.3 by $n=\infty$. 
We take these groups as interesting infinite discrete groups in the category of locally finite infinite groups, and wish to study their characters and spin ones, limiting processes as $n\to\infty$, convergence and divergence, calculation of limit functions and treatment from the side of probability theory. These studies were first
 opened by 
A.Vershik\,-\,S.\,Kerov ([38], [39]) and now become very fruitful. 

\vskip1.2em

{\bf 15.3. Subgroup $\Phi^{-1}\big(G(m,p,n)\big)\subset R\big(G(m,1,\infty)\big)$ and $R\big(G(m,p,n)\big)$.}\; 

 To show the heredity in the spin case, we study the relations between the full inverse image of $G(m,p,n)$ in $R\big(G(m,1,n)\big)$ and the representation group $R\big(G(m,p,n)\big)$. 

Put $G'_n:=R\big(G(m,1,n)\big),\,N'_n:=\Phi^{-1}\big(G(m,p,n)\big)$ for $5 \le n <\infty$, and we mean for $n=\infty$ their inductive limits as $n\to\infty$.

In the notation $\eta_j\;(j\in{\boldsymbol I}_n)$ in \S 15.2, we put 
\begin{eqnarray}
\label{2009-02-03-12}
\widehat{w}_1 = \eta_1^{\;p},\;
\;
\widehat{w}_j = \eta_1^{\;-1}\eta_j\quad(j\in{\boldsymbol I}_n,\,>\!1). 
\end{eqnarray}
Then $N'_n$ is generated by the union $Z\bigsqcup \{r_i\;(i\in{\boldsymbol I}_{n-1},\,\widehat{w}_i\;(i\in{\boldsymbol I}_n)\}$.

\vskip1.2em

{\bf Theorem 15.6.} (Case OO, $m$ odd)\; 
{\it 
The full inverse image  $\Phi^{-1}\big(G(m,p,n)\big)\subset R\big(G(m, 
\\ 
1,n)\big)$ is generated by the set 
$\{z_1,\;r_i\;(i\in{\boldsymbol I}_{n-1}),\, \widehat{w}_j\;(j\in{\boldsymbol I}_n)\}$. 
 Under the correspondence: $r_i\to r_i\;(i\in{\boldsymbol I}_{n-1}),\;\widehat{w}_j\to w_j\;(j\in{\boldsymbol I}_n)$, the normal subgroup 
$\Phi^{-1}\big(G(m,p,n)\big)$ is canonically isomorphic to the representation group $R\big(G(m,p,n)\big)$ (Theorem 15.4) of $G(m,p,n)$.
}

\vskip1.2em

{\bf Lemma 15.7.}\; 
{\it 
Assume $m$ be even. In the representation group $R\big(G(m,1,n)\big)$, the action of $r_i\;(i\in{\boldsymbol I}_{n-1})$ on $\widehat{w}_j\;(j\in{\boldsymbol I}_n)$ is given as 
\\[1ex]
$
\quad\mbox{\rm (iv$^0$-1)}
\quad{\rm For}\;\; i>1,\;\;
\left\{ 
\begin{array}{l}
r_i\widehat{w}_ir_i^{\;-1}=\widehat{w}_{i+1}\,,\;\;
r_i\widehat{w}_{i+1}r_i^{\;-1}=\widehat{w}_i,
\\
r_i\widehat{w}_jr_i^{\;-1}=\widehat{w}_j\qquad(j>1,\;\ne i,i\!+\!1), 
\end{array}
\right. 
\\[1ex]
\quad
\mbox{\rm (iv$^0$-2)}\qquad
\left\{ 
\begin{array}{l}
r_1\widehat{w}_1r_1^{\;-1} = z_2^{\;p(p-1)/2}z_3^{\;p}\widehat{w}_1{\widehat{w}_2}^{\;p},\;
\;
r_1\widehat{w}_2r_1^{\;-1}={\widehat{w}_2}^{\;-1},
\\
r_1\widehat{w}_jr_1^{\;-1}={\widehat{w}_2}^{\;-1}\widehat{w}_j\qquad\qquad(j>2).\end{array}
\right. 
$
}
\vskip1.2em

Note that, when $p$ is even, $z_3$ disappears from the above formulae.

\vskip1.2em
   
{\bf Theorem 15.8.} (Case OE: $p$ odd, $q$ even)\; 
{\it 
Normal subgroup $\Phi^{-1}\big(G(m,p,n)\big) \subset R\big(G(m,1,n)\big)$ is canonically isomorphic to the representation group $R\big(G(m,p,n)\big)$ (Theorem 15.5) under the correspondence $r_i \to r_i\;(i\in{\boldsymbol I}_{n-1}),\;\widehat{w}_j \to w_j\;(j\in{\boldsymbol I}_n)$.
}

\vskip1.2em

{\bf Theorem 15.9.} (Case EO, EE: $p$ even)  
{\it 
Let $p$ be even. 

{\rm (i)}\; 
 Define a subgroup $H'(m,p,n):=\langle z_1,z_2,\;r_i\;(i\in{\boldsymbol I}_{n-1}),\;\widehat{w}_j\;(j\in{\boldsymbol I}_n)\rangle \subset \Phi^{-1}\big(G(m, 
\\ 
p,n)\big)$.  Then $H'(m,p,n)$ is not normal but 
\begin{eqnarray}
\label{2018-03-29-11}
\Phi^{-1}\big(G(m,p,n)\big)=Z_3\cdot H'(m,p,n)\cong Z_3\times H'(m,p,n).
\end{eqnarray}

{\rm (ii)}\; 
Case {\rm EO:} $q=m/p$ odd.  Under the correspondence, $r_i \to r_i\;(i\in{\boldsymbol I}_{n-1}),\;\widehat{w}_j \to w_j\;(j\in{\boldsymbol I}_n)$, the subgroup $H'(m,p,n)$ is canonically isomorphic to $R\big(G(m,p,n)\big)$. 

{\rm (iii)}\; Case {\rm EE:} $q=m/p$ even. 
Under the correspondence 
 $r_i\to r_i\;(i\in{\boldsymbol I}_{n-1}),\;\widehat{w}_j\to w_j\;(j\in{\boldsymbol I}_n)$, the subgroup 
$H'(m,p,n)$ is canonically isomorphic to the quotient group $R\big(G(m,p,n)\big)\big/Z_3.$ 
}

\vskip1.2em

 {\bf Table 15.10.}\;\; For $n=\infty$, through the inductive limits, this theorem holds by replacing both ${\boldsymbol I}_{n-1},\,{\boldsymbol I}_n$ by ${\boldsymbol I}_\infty$. So we have for $5\le n\le \infty$ the following:
 \vskip1em

\begin{center}
$
\begin{array}{|c|l|}
\hline 
{\rm Case} &\quad \mbox{\rm Relation of $R\big(G(m,p,n)\big)$ and $\Phi^{-1}\big(G(m,p,n)\big)$} 
\\
\hline \hline
{\rm OO,\;OE}
&
R\big(G(m,p,n)\big)\;\cong\; \Phi^{-1}\big(G(m,p,n)\big)\vartriangleleft R\big(G(m,1,n)\big) 
\\
\hline 
{\rm EO} &
R\big(G(m,p,n)\big)\;\cong\; H'(m,p,n) \;\subset\; R\big(G(m,1,n)\big)
\\
\hline
{\rm EE} & 
R\big(G(m,p,n)\big)/Z_3\cong H'(m,p,n) \subset R\big(G(m,1,n)\big)
\\
\hline
\end{array}
$
\end{center}

\vskip1.2em

{\bf 15.4. Automorphisms of $\Phi^{-1}\big(G(m,p,\infty)\big)$ and invariant functions on it.}\; 

\vskip.5em

{\bf Lemma 15.11.}\; 
{\it 
 Let $4 \le n<\infty$, and $p|m,\,p>1$.  
 Then the normal subgroup $N_n=G(m,p,n)$ of $G_n=G(m,1,n)$ has index $[G_n:N_n]=p$ and a complete set of representatives for $G_n/N_n$ is given by $\{y_1^{\;j}\,;\,0\le j<p\}$. Two automorphism groups ${\rm Int}(N_n)$ and ${\rm Aut}_{G_n}(N_n):=\{\iota(g)|_{N_n}\,;\,g\in G_n\}$, with $\iota(g)h:=ghg^{-1}\;(h\in N_n)$, has index $(n,p):={\rm MCD}\{n,p\}$. A complete set of representatives for ${\rm Aut}_{G_n}(N_n)/{\rm Int}(N_n)$ is given by $\{\iota(y_1)^j\,;\,0\le j<(n,p)\}$.
}

\vskip1em

{\it Proof.}\; 
${\rm Int}(N_n)\cong N_n/Z_{N_n},\,{\rm Aut}_{G_n}(N_n)\cong G_n/Z_{G_n}(N_n)$, where $Z_{N_n}$ is the center of $N_n$ and $Z_{G_n}(N_n)$ is the centraliser of $N_n$ in $G_n$. Therefore,
$$
|{\rm Aut}_{G_n}(N_n)/{\rm Int}(N_n)|=|G_n/Z_{G_n}(N_n)|/|N_n/Z_{N_n}|
=|G_n/N_n|/|Z_{G_n}(N_n)/Z_{N_n}|.
$$

On the other hand, $|G_n/N_n|=p$.\, $Z_{G_n}(N_n)$ (resp. $Z_{N_n}$) consists of elements $(y_1y_2\cdots y_n)^a 
\\ 
(0\le a<m)$ $\big($resp. $(0\le a<m, na\equiv 0\;({\rm mod}\;p)\big)$. Here, for $0\le a<m,$
\vskip.5em

\qquad $ na\equiv 0\;({\rm mod}\;p)\;\Longleftrightarrow\;(n,p)\,a=kp\;(\exists k) \;\Longleftrightarrow\;a=k\cdot p/(n,p)\;(\exists k).$ 
\vskip.5em
\noindent
Hence $Z_{G_n}(N_n)/Z_{N_n}$ is represented by the set $(y_1y_2\cdots y_n)^a\;\big(0\le a<p/(n,p)\big).$ 
Therefore \,$|{\rm Aut}_{G_n}(N_n)/{\rm Int}(N_n)|=(n,p)$.  
\hfill
$\Box$

\vskip1em

{\bf Lemma 15.12.}\; 
{\it 
 Let $5 \le n<\infty$, and $p|m,\,p>1$.  

{\rm (i)}\; The normal subgroup $N'_n:=\Phi^{-1}\big(G(m,p,n)\big)$ of $G'_n=R\big(G(m,1,n)\big)$ has index $[G'_n:N'_n]=p$ and a complete set of representatives for $G'_n/N'_n$ is given by $\{\eta_1^{\;j}\,;\,0\le j<p\}$.

{\bf (ii)}\; \mbox{\rm Case $m$ odd.}\; 
$[{\rm Aut}_{G'_n}(N'_n): {\rm Int}(N'_n)]=(n,p)$, where $(n,p):={\rm MCD}\{n,p\}$. 
A complete set of representatives of 
${\rm Aut}_{G'_n}(N'_n)/{\rm Int}(N'_n)$ is given by 
\,$
\big\{\iota(\eta_1^{\;k})|_{N'_n}\,;\,0\le k<(n,p)\big\}, 
$\,  
where $\iota(h')g':=h'g'{h'}^{\,-1}$.

{\bf (iii)}\; \mbox{\rm Case $m$ even.}\;  
If $n$ is even, then the same assertion as above holds too. 

}
\vskip.8em

The proof is similar to that of Lemma 15.11.

\vskip1em 

Now we consider the inductive limit group $G(m,1,\infty)$ of mother groups $G(m,1,n)$ and the one $G(m,p,\infty)$ of child groups $G(m,p,n)$ for $p|m,\,p>1$.

\vskip1em

{\bf Lemma 15.13.}\; 
{\it 

{\rm (i)}\; 
 The normal subgroup $N=G(m,p,\infty)$ of $G=G(m,1,\infty)$ has index $[G:N]=p$ and a complete set of representatives for $G/N$ is given by $\{\eta_1^{\;j}\,;\,0\le j<p\}$. 

{\rm (ii)}\; 
Two automorphism groups ${\rm Int}(N)$ and ${\rm Aut}_G(N):=\{\iota(g)|_N\,;\,g\in G\}$ of $N$, with $\iota(g)h:=ghg^{-1}\;(h\in N)$, have also the index $p$. Nevertheless, let $g\in G$, then for any $n>1$, there exists an element $h_k\in N_k:=G(m,p,k)$ with $k>n$ such that $\iota(g)|_{N_n}=\iota(h_k)|_{N_n}$. 

{\rm (iii)}\; 
A function $f$ on $N$ is invariant under $G$ or $f(\iota(g)h)=f(h)\;(h\in N,\,g\in G)$\, if and only if it is invariant under $N$. Hence $K_1(N,G)=K_1(N)$ and $E(N,G)=E(N)$. 
}
\vskip1em 

{\it Proof.}\; 
We prove (ii). 
Let $g=(d,\sigma)\in G=D_\infty(T)\rtimes\mathfrak{S}_\infty$. Take $\eta_k^{\;a}$ so that ${\rm ord}(d)+a\equiv 0\;({\rm mod}\;p)$ and put $h_k:=g\eta_k^{\;a}$, then $h_k\in N_k$. Since $\iota(\eta_k)$ acts on $N_n$ trivially, we have $\iota(h_k)|_{N_n}=\iota(g)|_{N_n}$. 
\hfill 
$\Box$

\vskip1.2em 

The property of $\iota(g)$ in (ii) can be expressed as \lq\lq\,{\it $\iota(g)$ is locally inner on $N$\,}''.

\vskip1.2em

Now we go up to the spin case and put $G'_\infty:=R\big(G(m,1,\infty)\big),\,G'_n:=R\big(G(m,1,n)\big)$, and $N'_\infty:=\Phi^{-1}\big(G(m,p,\infty)\big),\,N'_n:=\Phi^{-1}\big(G(m,p,n)\big)$.

\vskip1em 

{\bf Lemma 15.14.}\; 
{\it 
 Let $p|m,\,p>1$.

{\rm (i)}\; 
For $5\le n\le \infty$, the normal subgroup $N'_n$ of $G'_n$ has index $[G'_n:N'_n]=p$ and a complete set of representatives for $G'_n/N'_n$ is given by $\{\eta_1^{\;j}\,;\,0\le j<p\}$. 

{\rm (ii)}\; 
If $p$ is odd, then every $\iota(g')\;(g'\in G'_\infty)$ acts on $N'_\infty$ locally inner way so that, for each $n>1$, there exists an $h'\in N'_\infty$ such that $\iota(g')|_{N'_n}=\iota(h')|_{N'_n}$. 
Accordingly $K_1(N'_\infty,G'_\infty)=K_1(N'_\infty)$ and 
$E(N'_\infty,G'_\infty)=E(N'_\infty)$ for $N'_\infty=\Phi^{-1}\big(G(m,p,\infty)\big)$. 

}
\vskip1em 

{\it Proof.} We prove (ii). 
The restriction onto $N'_\infty$ of $\iota(G'_\infty)$ is generated by $\iota(\eta_1)|_{N'_\infty}$ and $\iota(N'_\infty)$. So we check $\iota(\eta_1)$. For a fixed $n>1$, take $k>n$ and $h'_k=g'\eta_k^{\;p-1}\in N'_\infty$. Since $p-1$ is even and since $\eta_k^{\;2}$ commutes elementwise with $N'_n$, we have $\iota(g')h'=\iota(h'_k)h'\;(h'\in N'_n)$. 
\\
~\hfill 
$\Box$
\vskip1em

{\bf Lemma 15.15.}\; 
{\it 
 Let $p|m,\,p>1$ be even. Put $H'_\infty:={\displaystyle \lim_{n\to\infty}H'(m,p,n)}$, then $N'_\infty=Z_3\times H'_\infty$.  
The restriction $\iota(G'_\infty)|_{N'_\infty}$ is generated by $\iota(\eta_1)|_{N'_\infty}$ and $\iota(N'_\infty)$. 

Let $f$ be an invariant (under $N'_\infty$) function on $N'_\infty$. Then $f$ is invariant under $G'_\infty$ if and only if, for any $h'\in N'_\infty$,  
\begin{eqnarray}
\label{2018-03-31-1}
f(h')=f(z_3^{\;L(\sigma')}h'),
\end{eqnarray}
where $h'=zd'\sigma'$ with $z\in Z,\,d'\in \widetilde{D}_\infty,\,{\rm ord}(d')\equiv 0\;({\rm mod}\;p),\,\sigma'\in\widetilde{\mathfrak{S}}_\infty$. 
}
\vskip1em

{\it Proof.}\; 
We examine $\iota(\eta_1)$. For any $n>1$, take $k>n$ and put $h'_k:=\eta_1\eta_j^{\,-1}$, then $h'_k \in N'_k$ and $\iota(\eta_1)=\iota(h'_k)\iota(\eta_j)$. For $h'\in N'_n$, $\iota(\eta_1)h'=\iota(h'_k)\iota(\eta_j)h'=z_2^{\;{\rm ord}(d')}z_3^{\;L(\sigma')}\iota(h'_k)h'=\iota(h'_k)\big(z_3^{\;L(\sigma')}h'\big)$ because ${\rm ord}(d')$ is even.
Since $f$ is assumed to be invariant under $N'_\infty$, we have $f\big(\iota(\eta_1)h'\big)=f\big(z_3^{\;L(\sigma')}h'\big)$.
\hfill 
$\Box$

\vskip1.2em

{\bf 15.5. Heredity from mother group $R\big(G(m,1,\infty))$ to child groups}\; 
\vskip.5em

We summarize here the results on heredity in the case $n=\infty$. Let $p|m$ and put $q=m/p$. Put $G'=R\big(G(m,1,\infty)\big),\,$ and $\widehat{w}_1=\eta_1^{\;p},\,\widehat{w}_j=\eta_1^{\,-1}\eta_j\;(j\ge 2)$, and 
\vskip.5em 

\qquad
$N'=N'(m,p,\infty):=\Phi^{-1}\big(G(m,p,\infty)\big)=\langle Z,\, r_i\;(i\ge 1), \,\widehat{w}_j\;(j\ge 1)\rangle,\,$

\qquad
$H'=H'(m,p,\infty):=\langle z_1,z_2, \, r_i\;(i\ge 1), \widehat{w}_j\;(j\ge 1)\rangle,\,$ for $p$ even. 
\vskip.5em

 We list up the properties of the restriction map ${\rm Res}^{G'}_{N'}: E(G') \ni F\mapsto f=F|_{N'}\in E(N',G')$, according to the cases. 
Put $Z':={\mathfrak M}\big(G(m,p,\infty)\big)$.

\vskip.8em

{\bf Case OO} ($m$ odd). $N'\cong R\big(G(m,p,\infty)\big),\,Z'=Z_1,\,E(N',G')=E(N')$. 

\quad ${\rm Res}^{G'}_{N'}:E(G')\ni F\mapsto f=F|_{N'} \in E(N')$ is surjective onto $E(N')$.
\vskip.8em

{\bf Case OE} ($p$ odd, $q$ even). $N'\cong R\big(G(m,p,\infty)\big),\, Z'=Z_1\times Z_2\times Z_3$.

Take a spin type for $E(G')$ as $\chi\in\widehat{Z},\;Z=Z_1\times Z_2\times Z_3$, $\beta_i=\chi(z_i)\;(i\in{\boldsymbol I}_3)$. Then ${\rm Res}^{G'}_{N'}:E(G';\chi)\ni F\mapsto f=F|_{N'} \in E(N';\chi)$ is surjective for each spin type $\chi$. 

\vskip.8em

{\bf Case EO} ($p$ even, $q$ odd). $N'=Z_3\times H',\,H'\cong R\big(G(m,p,\infty)\big), \,Z'=Z_1\times Z_2$. 

Take a spin type $\chi'=(\beta_1,\beta_2)\in \widehat{Z'}$ of $H'$ with $\beta_i=\chi'(z_i)$ and put $\chi=(\beta_1,\beta_2,\beta_3)\in\widehat{Z}$ with $\beta_3=\chi(z_3)=1$. Then 
${\rm Res}^{G'}_{H'}:E(G';\chi)\ni F\mapsto f=F|_{H'} \in E(H';\chi')$ is surjective by Lemma 15.15.

\vskip.8em

{\bf Case EE} ($p$ even, $q$ even). $H'\cong R\big(G(m,p,\infty)\big)/Z_3$,\, $Z'=Z=Z_1\times Z_2\times Z_3$ for $R\big(G(m,p,\infty)\big)$. 

Take a spin type $\chi=(\beta_1,\beta_2,\beta_3)\in \widehat{Z'}$ of $R\big(G(m,p,\infty)\big)$. In case $\beta_3=\chi(z_3)=1$, put $\chi'=(\beta_1,\beta_2) \in (Z_1\times Z_2)^\wedge$, then 
${\rm Res}^{G'}_{H'}:E(G';\chi)\ni F\mapsto f=F|_{H'} \in E(H';\chi')$ is surjective by Lemma 15.15, and $ E(H';\chi')$ is canonically identified with $E\big(R\big(G(m,p,\infty)\big);\chi)$.

\vskip1.2em

As seen above, except a half of Case EE: the case of spin type $\chi=(\beta_1,\beta_2,\beta_3)\in\widehat{Z}$ with $\beta_3=\chi(z_3)=-1$, the restriction map gives us a sufficient information as heredity from mother to her children.
So we make some comments on the remained cases for Case EE. 

{\bf Method 15.1.}\; Note that $R\big(G(p,m,n)\big),\,5\le n\le \infty$, is expressed as a semidirect product as in Lemma 15.5. So, the first proposal is to apply the classical method for a semidirect product group to construct IRs and obtain irreducible characters, as explained in \S 2 and applied in Part I of this paper.

{\bf Method 15.2.}\; 
Take a spin type $\chi=(\beta_1,\beta_2,\beta_3)\in \widehat{Z}$ with $\beta_3=-1$. Then it is expressed as a product $\chi=\chi^0\cdot\chi^{\rm VII}$ of $\chi^0=(\beta_1,\beta_2, 1)$ and a special spin type $\chi^{\rm VII}=(1,1,-1)\in \widehat{Z}$. If two functions $f_1,f_2$ on $G'(m,p,\infty):=R\big(G(m,p,\infty)\big)$ are respectively of spin type $\chi_1,\chi_2$, then the product $f_1f_2$ is of spin type $\chi_1\chi_2$. Taking as a good model for $R\big(G(m,1,\infty)\big)$ in \S 21, with multiplication operator by $f^{\rm VII}_0$, we can try to check if any $f_1\in E\big(G'(m,p,\infty);\chi_1\big)$ and a special simple $F^{\rm VII}\in  E\big(G'(m,p,\infty);\chi^{\rm VII}\big)$ produces as products $f_1F^{\rm VII}$ all of $E\big(G'(m,p,\infty);\chi\big),\,\chi=\chi_1\chi^{\rm VII},$ or not. 
Here we can propose as a candidate of $F^{\rm VII}$ a simple character in $E\big(G'(m,p,\infty);\chi^{\rm VII}\big)$ below. By Theorems 15.3 and 15.5, we have the following

\vskip1em 
{\bf Lemma 15.16.}\; 
The quotient group $H'=G'(m,p,\infty)/(Z_1\times Z_2)$ is presented by
 
 the set of generators:\quad $\{s_i\;(i\in{\boldsymbol I}_\infty),\;x_j\;(j\in{\boldsymbol I}_\infty)\}$, 

the set of fundamental relations: 

\quad (ii)\;\;\;\, \quad relations in $\mathfrak{S}_\infty$ for $\{s_i\;(i\in{\boldsymbol I}_\infty)\}$,

\quad  (iii)\;\;\, \quad $x_1^{\;q}=e,\;x_j^{\;m}=e,\;x_jx_k=x_kx_j\;(j\ne k),$ 

\quad  (iv-1)\;\,\quad For $i>1$,\; $s_ix_js_i^{\,-1}=x_{s_i(j)}\quad(j>1),\;
s_ix_1s_i^{\,-1}=z_3x_1$, 

\quad  (iv-2)\;\quad \;$s_1x_1s_1^{\,-1}=z_3x_1x_2^{\;p},\;s_1x_2s_1^{\,-1}=x_
2^{\,-1},\;s_1x_js_1^{\,-1}=x_2^{\,-1}x_j\;(j>2). $
\hfill 
$\Box$

\vskip1.2em 

Denote by $D_\infty^{(q)}$ the subgroup generated by $\{x_j\;(j\in{\boldsymbol I}_\infty)\}$, then $H'=\big(D_\infty^{(q)}\times Z_3\big)\rtimes \mathfrak{S}_\infty$ under the action (iv-1)--(iv-2). 

Now take a character $X_{0,1}:={\bf 1}_{D^{(q)}_\infty}\otimes {\rm sgn}_{Z_3}$ of the abelian group $D_\infty^{(q)}\times Z_3$. Then its stationary subgroup in $\mathfrak{S}_\infty$ is equal to $\mathfrak{A}_\infty$. The induced representation $\pi_2$ of $X_{0,1}$ from $\big(D_\infty^{(q)}\times Z_3\big)\rtimes \mathfrak{A}_\infty$ to $H'$ is two-dimensinal.  

In fact, the representation space of $\pi_2$ has a basis $\varphi_1,\,\varphi_2$ such that $\varphi_i(k'h')=\rho'(k')\varphi_i(h') 
\\ 
 (k'\in K',\,h'\in H')$ and that $\varphi_1$ (resp. $\varphi_2$) is supported on $K'$ (resp. $K's_1$). The representation operators are, in matrix form with respect to $\{\varphi_1,\varphi_2\}$, and with $(-1)^0:=1$,  
$$
\pi(s_1)=\begin{pmatrix}0 & 1 \\ 1 & 0\end{pmatrix},\quad
\pi_2(x_j)=\begin{pmatrix}1 & 0 \\ 0 & (-1)^{\delta_{j,1}}\end{pmatrix},\quad
\pi(z_3)=\begin{pmatrix}-1 & 0 \\ 0 & -1\end{pmatrix}.
$$

Its character $\chi_{\pi_2}$ is given as 
\begin{eqnarray}
\label{2017-04-04-1}
&&
\chi_{\pi_2}(z_3^{\;a} \sigma\,x_1^{b_1}x_2^{\;b_2}x_3^{\;b_3}\cdots)=\left\{
\begin{array}{ll}
(-1)^a\,2, \qquad & \mbox{\rm if $\sigma\in \mathfrak{A}_\infty$ and $b_1$ even.}
\\
\;\;0, & \mbox{\rm otherwise}.
\end{array}
\right.
\end{eqnarray}

Taking the results about $f^{\rm VII}_0$ in \S 21 as a model, we can propose the following
\vskip1.2em

{\bf Conjecture 15.17.}\; Denote by $F^{\rm VII}$ the normalized chracter $\widetilde{\chi}_{\pi_2}$. For any spin type $\chi=(\beta_1,\beta_2,\beta_3)\in \widehat{Z}$ with $\beta_3=-1$, express it as $\chi=\chi^0\cdot\chi^{\rm VII}$ with $\chi^0=(\beta_1,\beta_2, 1)$ and $\chi^{\rm VII}=(1,1,-1)$. 
Then, any $f\in E\big(G'(m,p,\infty);\chi^0\big)$ and the $F^{\rm VII}\in  E\big(G'(m,p,\infty);\chi^{\rm VII}\big)$ produces as products $f\/F^{\rm VII}$ all of $E\big(G'(m,p,\infty);\chi\big),\,\chi=\chi^0\chi^{\rm VII}.$

\vskip1.2em

{\bf 15.6. Restriction map and symmetries on the parameter space.} 

As remarked in \S 15.1, for the normal subgroups $N$ of $G_\infty$ listed in \S 15.1, the set of characters $E(N)$ is equal to $\{f|_N\,;\,f\in E(G_\infty)\}$. Then the parameter space for $E(N)$ is determined canonically from the one ${\cal A}(T)$ for $G_\infty$.  
For this we prepare two kinds of operations (or symmetries) on $A= 
\big( (\alpha_{\zeta,\varepsilon})_{(\zeta, \varepsilon) \in \widehat{T} \times \{0,  1\}} 
\;;\;\mu\big)\in {\cal A}(T)$ as follows.
\vskip.2em

{\bf (1)}\; Involution $\tau:\,A \to {}^t\!A\,$:   
\begin{eqnarray}
\label{2015-08-09-12}
&&
\hspace*{10ex}
{}^t\!A 
:= \big( (\alpha'_{\zeta,\varepsilon})_{(\zeta, \varepsilon) \in \widehat{T} \times \{0,  1\}} ;\mu' \big),\, \;\; 
 \\[.3ex]
\nonumber
&&
\alpha'_{\zeta,\varepsilon}=\alpha_{\zeta,\varepsilon+1},\; 
 \mu'=\mu\;\;
\mbox{\rm (permutation of $\alpha_{\zeta,0}$ and $\alpha_{\zeta,1}$).}
\end{eqnarray}

{\bf (2)}\; Translation $R(\zeta_0)$ by $\zeta_0\in\widehat{T}$\,:\; 
\begin{eqnarray}
\label{2015-07-27-21}
&&
\qquad\qquad 
R(\zeta_0)A:= 
\big( (\alpha'_{\zeta,\varepsilon})_{(\zeta, \varepsilon) \in \widehat{T} \times \{0,  1\}} 
\;;\;\mu'\big),
\\
&&
\nonumber
\alpha'_{\zeta,\varepsilon} = 
 \alpha_{\zeta\zeta_0^{\;-1}\!,\,\varepsilon}
\;\big(\,(\zeta, \varepsilon) \in \widehat{T} \times \{0,1\}\big), \quad   
\mu' = (\mu'_\zeta)_{\zeta \in \widehat{T}}, 
\; \mu'_\zeta = \mu_{\zeta\zeta_0^{\;-1}}\,.
\end{eqnarray}
For $\zeta_0=\zeta^{(a)}\;(0\le a<m)$, we put $\kappa^{(a)}:=R(\zeta^{(a)})$ for brevity. 
When $m$ is even, $\kappa:=\kappa^{(m^0)}$ with $m^0=m/2$ is involutive, i.e., $\kappa^2={\rm id}$. 

These operations on ${\cal A}(T)$ are mutually commutative. 
 Denote by ${\cal G}$ the commutative group generated by them. Put 
\begin{eqnarray}
\label{2017-11-24-1}
P_\eta(f_A):=f_{\eta(A)}\quad(\eta\in{\cal G}), 
\end{eqnarray}
then this gives a linear representation of ${\cal G}$ on the linear span $\langle f_{\eta(A)}\,;\,\eta\in{\cal G}\rangle$. For a subgroup ${\cal H}$ of ${\cal G}$, the projection onto the space of trivial representation ${\bf 1}_{\cal H}$ is   
\begin{eqnarray}
\label{2017-08-10-11}
P_{\cal H} :=\dfrac{1}{|{\cal H}|}{\sum}_{\eta\in{\cal H}}P_\eta\,. 
\end{eqnarray}

Let $N$ be a one of canonical normal subgroups of $G_\infty=G(m,1,\infty)$.  
Then, by Theorem 12.6\,(i),  
the set $E(N)$ of characters of $N$ is obtained by the restriction map\;  
$
{\rm Res}^{G_\infty}_N:\,E(G_\infty)\ni f_A\mapsto f_A|_N\in E(N)$.

\vskip1.2em

{\bf Theorem 15.18.}\; 
[12, Theorem 14]
{\it \;Let $N=G^{\mathfrak A}_\infty$. 

{\rm (i)}\; 
For $A,A^1\in{\cal A}(\widehat{T})$, a necessary and sufficient condition for \,$f_A|_N=f_{A^1}|_N$\, is \,$A^1=A$\, or $A^1=\tau(A)={}^t\!A$. 
With the subgroup $\langle \tau\rangle=\{\tau, {\rm id}\}$ of ${\cal G}$,   
\begin{eqnarray}
\label{2017-08-10-14}
P_{\langle \tau\rangle}(f_A)=\left\{
\begin{array}{ll}
f_A|_N\quad & \mbox{\rm on $N$}, 
\\[.5ex]
\;\;0\quad &\mbox{\rm outside $N$.} 
\end{array}
\right.
\end{eqnarray}

{\rm (ii)}\; Let $[A]=\{\eta(A)\,;\,\eta\in\langle \tau\rangle\}=\{A,\tau(A)\}$ the conjugacy class of $A$ in ${\cal A}(\widehat{T})/\langle \tau\rangle$, then a parametrization of 
$E(N)$ is given as  
\begin{eqnarray}
\label{2016-04-23-11}
E(G^{\mathfrak A}_\infty)=\big\{f_A|_N\;;\;[A]\in {\cal A}(\widehat{T})/\langle \tau\rangle\big\}.
\end{eqnarray}

}
\vskip.2em

We refer [12, Theorem 15] and [14, Theorem 7.1 and Proposition 7.2] for $N=G(m,p,\infty),\\
 p|m,\, p>1.$ 
Put $S(p):=\{t^p\,;\,t\in T\}\subset T$, and $q=m/p$. 
\vskip1.2em 

{\bf Theorem 15.19.}\; 
{\it 
Let $N=G(m,p,\infty).$ 
For two parameters $A, A^1\in {\cal A}(\widehat{T})$, 
$$
f_A\big|_N=f_{A^1}\big|_N\;\Longleftrightarrow\; A^1=R(\zeta_0)A\quad \big(\exists \zeta_0\in\widehat{T},\,\zeta_0|_{S(p)}={\bf 1}_{S(p)}\big).
$$

Furthermore for $\zeta_0=\zeta^{(a)}$,\;\;
 
\qquad
$
\zeta^{(a)}|_{S(p)}={\bf 1}_{S(p)}\,\Longleftrightarrow\,ap\equiv 0\;(\mathrm{mod}\;m) 
\,\Longleftrightarrow\, m$  
  is a multiple of $q$. 
\vskip.5em

With the subgroup $\langle\kappa^{(q)}\rangle$ of ${\cal G}$,\, $P_{\langle\kappa^{(q)}\rangle}(f_A)$ gives $f_A|_N$ on $N$ and zero elsewhere on $G_\infty$. A parametrization of $E(N)$ is given by 
 \begin{eqnarray*}
\label{2016-04-23-1}
E\big(G(m,p,\infty)\big)=\big\{f_A\big|_N\;;\; [A]\in {\cal A}(\widehat{T})\big/\big\langle \kappa^{(q)}\big\rangle\big\}. 
\end{eqnarray*}

}
\vskip.2em

Now put 
 $N=G^{\mathfrak A}_\infty(m,p,\infty)=
G^{\mathfrak A}_\infty\cap G(m,p,\infty)$, which is equal to $K_\infty$ in case $p=2$. For this case, we refer  Proofs of Theorems 14 and 15 in \cite{[HH2005c]}.

\vskip1.2em

{\bf Theorem 15.20.}\; Put $N=G^{\mathfrak A}_\infty(m,p,\infty),\,p|m,\,p>1$, and $q=m/p$. 
Then, 
\begin{gather*}
\label{2015-07-27-11-9}
\nonumber
E\big(G^{\mathfrak A}(m,p,\infty)\big)=\big\{f_A|_N\,;\, A = \big( (
\alpha_{\zeta,\varepsilon})_{(\zeta, \varepsilon) \in \widehat{T} \times \{\,0,  1\,\}}\,;\,\mu \big)
\in{\cal A}(\widehat{T})\big\},
\\
\label{2016-04-23-21}
f_A|_N=f_{A^1}|_N\;\Longleftrightarrow\;[A]=[A^1]\in{\cal A}(\widehat{T})/\langle \tau, \kappa^{(q)}\rangle,  
\end{gather*}
where $[A]$ denotes the equivalence class of $A$ in the quotient space. Moreover 
\begin{eqnarray}
\nonumber
&&
\hspace*{6ex}
P_{\langle \tau, \kappa^{(q)}\rangle}(f_A)=\left\{
\begin{array}{ll} 
f_A|_N\quad & \mbox{\rm on $N$}, 
\\[.5ex]
\;\;0& \mbox{\rm outside $N$,}
\end{array}
\right.
\\
\nonumber
\label{2016-04-23-31}
&& 
E\big(G^{\mathfrak A}(m,p,\infty\big)=\big\{ 
f_A|_N\;;\;[A]\in {\cal A}(\widehat{T})\big/\big\langle \tau, \kappa^{(q)}\big\rangle\big\}.
\end{eqnarray}

\vskip.2em

{\bf Note 15.21.}\; 
The infinite Weyl group of type 
B$_\infty$/C$_\infty$ is $G(2,1,\infty)={\mathfrak S}_\infty({\boldsymbol Z}_2)$\, and that of type D$_\infty$ is its normal subgroup $G(2,2,\infty)$. Their characters were obtained in \cite{[HH2002]}.

\section{\large Spin characters of $R\big(G(m,1,\infty)\big)$ of spin type $z_1\to -1$}

As at the top of \S 5, put 
for $4\le n \le \infty$, 
\,$
\widetilde{G}_n:=\widetilde{G}(m,1,n):=D_n(T)\rtimes \widetilde{{\mathfrak S}}_n$, 
and 
\vskip.5em
\hspace*{2ex} 
$\widetilde{G}^{\rm \,odd}(m,1,n):=\widetilde{G}(m,1,n)\cong R\big(G(m,1,n)\big)$ if $m$ is odd, and 

\hspace*{2ex}
$\widetilde{G}^{\rm \,IV}(m,1,n):=\widetilde{G}(m,1,n)\cong R\big(G(m,1,n)\big)/Z_{23},\,Z_{23}:=\langle z_2,z_3\rangle,$ if $m$ is even.
\vskip.3em
  
Then, as is explained in \S 5, a spin IR of $\widetilde{G}^{\rm \,IV}(m,1,n)$ of spin type $\chi_-\in \widehat{Z_1}$ (non-trivial character) is essentially equal to a spin IR of $R\big(G(m,1,n)\big)$ of spin type $\chi^{\rm IV}=(-1,\,1,\,1)$, and similarly for their characters. In this section, discussing pointwise limits for $\widetilde{G}_n:=\widetilde{G}(m,1,n)$, we obtain the set ${\rm Lim}(\widetilde{G}_\infty;\chi_-)$ in two cases, Y\,=\,odd and IV, at the same time. 

\vskip.8em

{\bf 16.1. Normalized spin irreducible characters of $\widetilde{G}_n.$}

Spin IRs of $\widetilde{G}_n$ are classified in Theorem 5.1 and their characters are given in Theorem 5.2. For our present purpose, we give normalized characters of IRs $\Pi^{\rm Y}_{\Lambda^n},\,\Lambda^n\in {\boldsymbol Y}^{\rm sh}_{\!n}(\widehat{T})$, for Y\,=\,odd and IV.
An element $g'=(d,\sigma')\in \widetilde{G}_n=D_n(T)\rtimes\widetilde{{\mathfrak S}}_n,\,T={\boldsymbol Z}_m,$ with $d\in D_n(T),\,\sigma'\in\widetilde{{\mathfrak S}}_n,$ has a standard decomposition as 
\begin{eqnarray}
\label{2014-04-07-1-2}
&&\;
g'=(d, \sigma') = \xi_{q_1}\xi_{q_2} \cdots \xi_{q_r}g'_1g'_2 \cdots g'_s,\;\;\xi_q=\big(t_q,(q)\big), \;\;g'_j=(d_j,\sigma'_j)\, ;
\\
\nonumber
&&
Q:=\{q_1,q_2,\ldots,q_r\},\;J:={\boldsymbol I}_s,\;K_j:=\mathrm{supp}(\sigma'_j),\;\ell_j:=\ell(\sigma'_j)\;(j\in J).
\end{eqnarray}

\vskip.2em 

{\bf Lemma 16.1.}\; 
{\it 
The normalized character of spin IR\, 
 $\Pi^{\rm Y}_{\Lambda^n},\,\Lambda^n=\big({\boldsymbol \lambda}^{n,\zeta})_{\zeta\in\widehat{T}}\in{\boldsymbol Y}^{\rm sh}_{\!n}(\widehat{T}),$ of\, $\widetilde{G}_n$ is given as follows. Let ${\boldsymbol \nu}=(n_\zeta)_{\zeta\in\widehat{T}},\,n_\zeta=|{\boldsymbol \lambda}^{n,\zeta}|,$ and ${\boldsymbol I}_n=\bigsqcup_{\zeta\in\widehat{T}}I_{n,\zeta}$ the associated standard partition, and put\, $H'_n=D_n(T)\rtimes \widetilde{{\mathfrak S}}_{\boldsymbol \nu},\,\widetilde{{\mathfrak S}}_{\boldsymbol \nu}={\widehat{*}}_{\zeta\in\widehat{T}}\,\widetilde{{\mathfrak S}}_{I_{n,\zeta}}$. 
 
If $g'\in \widetilde{G}_n$ is not conjugate to any element in $H'_n$, then\; $\widetilde{\chi}\big(\Pi^{\rm Y}_{\Lambda^n}|g'\big)=0$. 

If $g'=(d,\sigma')\in H'_n$\,, then choose a normalized  element under conjugacy modulo $Z_1$ so that 
\begin{eqnarray}
\label{2017-08-08-11}
K_j=[a_j,b_j],\;\sigma'_j=r_{a_j}r_{a_j+1}\cdots r_{b_j-1}\;(=: \sigma'_{K_j}\;({\rm put}))\quad(\forall j\in J).
\end{eqnarray}

When $L(\sigma'):=\sum_{j\in J}\big(\ell(\sigma'_j)-1\big)\equiv 0\;(\mathrm{mod}\;2)$ or $\sigma'$ even, the character value is non-zero only when \,$L(\sigma'_j)\equiv 0\;(\mathrm{mod}\;2)\;\;(\forall j\in J)$,\footnote{See the row of Case IV (which corresponds to the case $m$ odd) in Table 9.1 in [I, p.81].} and in that case   
\begin{gather}
\label{2014-03-17-11-2-2}
\widetilde{\chi}\big(\Pi^{\rm Y}_{\Lambda^n}|g'\big)
=
\sum_{{\cal Q},{\cal J}} 
c(\Lambda^n;{\cal Q},{\cal J};g') \widetilde{X}(\Lambda^n;{\cal Q},{\cal J};g'), 
\\
\nonumber
c(\Lambda^n;{\cal Q},{\cal J};g')
=
\frac{(n-|\mathrm{supp}(g')|)!}{n!}\cdot
\prod_{\zeta\in\widehat{T}}
\frac{|{\boldsymbol \lambda}^{n,\zeta}|!}{
\Big(|{\boldsymbol \lambda}^{n,\zeta}| - |Q_\zeta| - 
\sum_{j \in J_\zeta}\ell_j\Big)!
}
\,,
\\[.3ex]
\nonumber
\label{2015-07-25-12-2}
\widetilde{X}(\Lambda^n;{\cal Q},{\cal J};g')
=
\prod_{\zeta\in\widehat{T}} 
\Big(
\prod_{q\in Q_\zeta}
\zeta(t_q)\cdot\prod_{j\in J_\zeta}\zeta\big(P(d_j)\big)\cdot 
\widetilde{\chi}\big(\tau_{{\boldsymbol \lambda}^{n,\zeta}}|(\ell_j)_{j\in J_\zeta}\big)
\Big), 
\end{gather} 
where $\tau_{{\boldsymbol \lambda}^{n,\zeta}}$ is spin IR of $\widetilde{{\mathfrak S}}_{n_\zeta}$ corresponding to ${\boldsymbol \lambda}^{n,\zeta}$, and \,${\cal Q}=(Q_\zeta)_{\zeta\in\widehat{T}},\;{\cal J}=(J_\zeta)_{\zeta\in\widehat{T}}$\, run over respectively partitions of $Q$ and $J$ satisfying 
\begin{eqnarray}
\label{2014-04-07-11-2-2}
\nonumber
\mbox{\rm (Condition QJ)}\hspace*{16ex}
|Q_\zeta|+{\sum}_{j\in J_\zeta}\ell_j\le |{\boldsymbol \lambda}^{n,\zeta}|\;\;\;
(\zeta\in\widehat{T}), 
\hspace{20ex}
\end{eqnarray}
and for the symbol $\widetilde{\chi}\big(\tau_{{\boldsymbol \lambda}^{n,\zeta}}|(\ell_j)_{j\in J_\zeta}\big)$, see Proposition 5.3. 

When $L(\sigma')\equiv 1\;(\mathrm{mod}\;2)$ or $\sigma'$ odd,\, 
$\widetilde{\chi}\big(\Pi^{\rm Y}_{\Lambda^n}|g'\big)\ne 0$\, only if\, $|\mathrm{supp}(g')|\ge n-1$.\footnote{Cf. the row of Case IV (corresponding to the case $m$ odd) in [loc.cit., Table 9.1].}

}

\vskip1em

{\bf Note 16.2.} The condition \,$|\mathrm{supp}(g')|\ge n-1$\, in the last assertion of Lemma 16.1 means that, when $L(\sigma')\equiv 1\;(\mathrm{mod}\;2)$, $\widetilde{\chi}\big(\Pi^{\rm Y}_{\Lambda^{n+k}}|g'\big)= 0$ for $k\ge 2$, and so their limit as $k\to\infty$ is zero ($n$ is fixed along with $g'$). 

\vskip1em

{\bf 16.2. Limits of normalized spin irreducible characters on $\widetilde{G}_\infty$.}

For an $f\in K(\widetilde{G}_\infty;\chi_-)$, its support is evaluated in Table 12.1 as\; 
$$
\mathrm{supp}(f)\subset{\cal O}':=\{g'=(d,\sigma')\,;\,L(\sigma'_j)=\ell(\sigma'_j)-1\equiv 0\;(\mathrm{mod}\;2,\,\forall j\in J)\}. 
$$ 
 
 On the other hand, from the results in \S 16.1, we see that any limit functions $f\in {\rm Lim}(\widetilde{G}_\infty;\chi_-)$ has its support \;$\mathrm{supp}(f)\subset {\cal O}'$\; as is expected (cf. Note 16.2).  Depending on the cases Y\,=\,odd or IV, we put for $g'=(d,\sigma')\in \widetilde{G}_\infty$,  
\begin{eqnarray}
\label{2017-09-17-1}
f^{\rm Y}_0(g'):=\psi_\Delta(\sigma'), 
\end{eqnarray}
a natural extension through\; $\widetilde{G}_\infty\to \widetilde{{\mathfrak S}}_\infty\cong \widetilde{G}_\infty/D_\infty$.  
For a 
$g'\in \widetilde{G}_\infty$, to calculate limit ${\displaystyle \lim_{n\to\infty}\widetilde{\chi}\big(\Pi^{\rm Y}_{\Lambda^n}|g'\big)}$ we follow the discussions in \S 14. 
\vskip1em 

{\bf Lemma 16.3.} 
{\it  A necessary and sufficient condition for the existence of pointwise limit\, ${\displaystyle \lim_{n\to\infty}\widetilde{\chi}\big(\Pi^{\rm Y}_{\Lambda^n}|(d,{\bf 1})\big)}$ at any point $d\in D_\infty(T)$ is  
\\[1ex]
{\rm {\sc (Condition\;I)}} 
\qquad\qquad 
${\displaystyle 
\exists\;\;
B_\zeta := 
\lim_{n\to \infty} \frac{|{\boldsymbol \lambda}^{n,\zeta}|}{n}  \quad  
(\zeta \in \widehat{T}). 
}
$
\\[1ex]
In that case,\,  
$ 
\sum_{\zeta\in\widehat{T}}\,B_\zeta=1.$
}

\vskip1em

Put\, $\widehat{T}^+:=\{\zeta\in\widehat{T}\,;\,B_\zeta>0\}$. For $\Lambda^n=({\boldsymbol \lambda}^{n,\zeta})_{\zeta\in\widehat{T}}\in{\boldsymbol Y}^{\rm sh}_{\!n}(\widehat{T}),\,$ let $n_\zeta:=|{\boldsymbol \lambda}^{n,\zeta}|$ and ${\boldsymbol \lambda}^{n,\zeta}=\big(\lambda^{n,\zeta}_1,\lambda^{n,\zeta}_2,\ldots, \lambda^{n,\zeta}_{l(n,\zeta)}\big)\in S\!P_{n_\zeta}$. 

  \vskip1em 

{\bf Theorem 16.4.}\; 
{\it 
 For a series $\Lambda^n=({\boldsymbol \lambda}^{n,\zeta})_{\zeta\in\widehat{T}}\in{\boldsymbol Y}^{\rm sh}_{\!n}(\widehat{T}),\,$ the pointwise limit ${\displaystyle \lim_{n\to\infty}\widetilde{\chi}\big(\Pi^{\rm Y}_{\Lambda^n}|g'\big)}$\, exists if and only if
\\[1ex]
{\rm ({\sc Condition} I$\Lambda$-spin)}
\hspace{4ex}
${\displaystyle 
\exists\;\;\lim_{n\to\infty} \frac{\lambda^{n,\zeta}_i}{|{\boldsymbol \lambda}^{n,\zeta}|}\,=\,\gamma'_{\zeta;i}\quad(\forall \zeta\in\widehat{T}^+,\; \forall i\ge 1)
}
$
\\[1ex]
holds, together with  \,{\rm (Condition I)}.  
For $\zeta\in\widehat{T}^+$, put 
\begin{eqnarray}
\label{2017-11-24-21}
{\boldsymbol \gamma}_\zeta:=\big(\gamma_{\zeta;i}\big)_{i\geqslant 1},\quad 
\mu:=(\mu_\zeta)_{\zeta\in\widehat{T}}
\end{eqnarray}
 with \,$\gamma_i:=B_\zeta\hspace{.15ex}\gamma'_{\zeta;i}$ and  
$\mu_\zeta:=B_\zeta-\|\gamma_\zeta\|\ge 0.$ 
For $\zeta\in \widehat{T}\setminus\widehat{T}^+$, put\,  ${\boldsymbol \gamma}_\zeta:={\bf 0},\;\mu_\zeta:=0$. Altogether 
\,
$C:=\big(({\boldsymbol \gamma}_\zeta)_{\zeta\in\widehat{T}}\,;\mu\big)
$\, 
is a parameter for the limit function $f^{\rm Y}_C(g')={\displaystyle \lim_{n\to\infty}\widetilde{\chi}\big(\Pi^{\rm Y}_{\Lambda^n}|g'\big)}$. Take $g'\in\widetilde{G}^{\rm Y}_\infty$ in (\ref{2014-04-07-1-2}). In case 
\,$L(\sigma'_j)\equiv 0\;(\forall j\in J,\;\mathrm{mod}\;2)$, it has the following form, otherwise it is zero:  
\begin{gather}
\label{2016-03-28-41}
\nonumber 
f^{\rm Y}_C(g')= f^{\rm Y}_0(g')\cdot
\prod_{q \in Q}
\Big\{ \sum_{\zeta \in \widehat{T}}
\big(
 \|\gamma_\zeta\|   
\; + \;\mu_\zeta\big)
\zeta(t_{q})\Big\} \,
\prod_{j \in J} 
\Big\{ \sum_{\zeta \in \widehat{T}}\Big( \sum_{i \geqslant 1}\;(\gamma_{\zeta;i})^{\ell_j}
\Big) 
\zeta\big(P(d_j)\big)
\Big\}. 
\end{gather}

In case $m$ is even, $\widetilde{G}_\infty=G'_\infty/Z_{23}$ with $G'_\infty:=R\big(G(m,1,\infty)\big),\,Z_{23}=\langle z_2,z_3\rangle$, and we can consider functions above naturally on $G'_\infty$.
}

\vskip1em

Put, for $C=\big(({\boldsymbol \gamma}_\zeta)_{\zeta\in\widehat{T}}\,;\mu\big)$, 
\begin{eqnarray}
\label{2017-09-17-21}
&&\;\;
A_C:=\big((\alpha_{\zeta,\varepsilon})_{(\zeta,\varepsilon)\in\widehat{T}\times\{0,1\}},\mu\big)\in{\cal A}(\widehat{T})\;\;{\rm with}\;\; \alpha_{\zeta,0}\!=\!\gamma_\zeta,\,\alpha_{\zeta,1}\!=\!{\bf 0}\!=\!(0,0,\ldots), 
\end{eqnarray}
then $f^{\rm Y}_C(g')=f^{\rm Y}_0(g') f_{A_C}(g)$, with $g=\Phi(g')$. 
For parameter spaces of limit functions (and also of characters) of $G'_\infty$ of spin type $\chi^{\rm Y}$, we prepare a notation as follows. We put ${\cal K}=\widehat{T}$\, for Y\,=\,odd, IV here and III later, and ${\cal K}=\widehat{T}^0$ for Y\,=V and VI later. 
\vskip1.2em

{\bf Notation 16.5.} A set of parameters is defined, for ${\cal K}=\widehat{T}$ or $=\widehat{T}^0$, as  
\begin{eqnarray}
\label{2016-03-28-42}
&& 
{\cal C}({\cal K})=\Big\{ C=\big(({\boldsymbol \gamma}_\zeta)_{\zeta\in{\cal K}}\,;\mu\big)\,;\,{\sum}_{\zeta\in{\cal K}}\,\|\gamma_\zeta\|+\|\mu\|=1
\Big\}, 
\qquad
\\
\nonumber
{\rm with}
&&\;\;
\gamma_\zeta=(\gamma_{\zeta,i})_{i\geqslant 1},\;\gamma_{\zeta,1}\ge \gamma_{\zeta,2}\ge \ldots \ge 0,\;\mu=(\mu_\zeta)_{\zeta\in{\cal K}}.
\end{eqnarray}

\vskip.2em 

{\bf Theorem 16.6.}\; \mbox{\rm (Characters of $G'_\infty=R\big(G(m,1,\infty)\big)$ of spin type $\chi^{\rm Y}$)}
 
{\it 
In case $m$ is odd and\, {\rm Y\,=\,odd}, or in case $m$ is even and\, {\rm Y=IV},
\begin{eqnarray*}
\;\;
{\rm Lim}(G'_\infty;\chi^{\rm Y})=\big\{f^{\rm Y}_C=f^{\rm Y}_0 f_{A_C}\;;\;C\in{\cal C}(\widehat{T})\big\}=E(G'_\infty;\chi^{\rm Y})=F(G'_\infty;\chi^{\rm Y}),
\end{eqnarray*}
where, for {\rm Y\,=\,IV}, the function $f^{\rm Y}_C$ on $\widetilde{G}_\infty=G'_\infty/Z_{23}$ is naturally lifted up to $G'_\infty$. 
}

\section{\large Spin characters of spin type $\chi^{\rm VII}=(1,\,1,-1)$}

Let $m$ be even. For $4\le n\le \infty$, put $G_n:=G(m,1,n),\,G'_n:=R\big(G(m,1,n)\big)$ and let $G^{\mathfrak A}_\infty,\, K_\infty \subset G_\infty$ and $G^{\prime\hspace{.15ex}{\mathfrak A}}_\infty,\,K'_\infty\subset G'_\infty$ be the normal subgroups introduced in \S 12.1. Spin characters of $G'_n$ of spin type $\chi^{\rm VII}$ for $n$ finite, and limiting process as $n\to\infty$, and spin characters of spin type $\chi^{\rm VII}$ of $G'_\infty$ are discussed in [I, Part II]. Here we give the set $E(G'_\infty;\chi^{\rm VII})$ by a quite different way. Note that, for $f\in K(G'_\infty;\chi^{\rm VII})$, $\mathrm{supp}(f)\subset {\cal O}({\rm VII})=K'_\infty$. 
Let the normalized character \,$f^{\rm VII}_0=\widetilde{\chi}_{\pi_{2,\zeta^{(0)}}}=\mathrm{sgn}^{\mathfrak A}_{Z_3}\cdot X_{{\cal O}({\rm VII})}$\, of 2-dimensional IR $\pi_{2,\zeta^{(0)}}$ be as in \S 12.4. We define two maps ${\cal M}$ and ${\cal N}$ by multiplication of $f^{\rm VII}_0$ as follows.

\vskip1em

{\bf Definition 17.1.}\; 
For $F\in K_1(G_\infty)$ and $f\in K_1\big(G'_\infty;\chi^{\rm VII}\big)$, and $g'\in G'_\infty,\,g=\Phi(g')\in G_\infty$,\, put   
\begin{eqnarray*}
\label{2009-03-02-11}
&&
{\cal M}(F)(g')
:= 
f^{\rm VII}_0(g') F(g),\quad 
{\cal N}(f)(g) 
:=
 f^{\rm VII}_0(g') f(g')
\\
&&
\qquad\qquad
K_1(G_\infty)
\begin{array}{c} 
\mbox{\rm\footnotesize ${\cal M}$}   
 \\[-1.4ex]
 \longrightarrow
  \\[-2.2ex]
 \longleftarrow
\\[-1.4ex]
 \mbox{\footnotesize ${\cal N}$}
 \end{array}
 K_1\big(G'_\infty;\chi^{\rm VII}\big)
\end{eqnarray*} 

\vskip.0em

{\bf Lemma 17.2.}\; 
{\it 
The map ${\cal M}$ is surjective, and the map ${\cal N}$ is injective, and both preserve linear combinations. The product ${\cal M}{\cal N}$ is the identity map on $K_1\big(G'_\infty;\chi^{\rm VII}\big)$. 
}
\vskip1em 

{\it Proof.}\; 
By Proposition 12.6, \,$\big(f^{\rm VII}_0\big)^2 =X_{{\cal O}({\rm VII})}$\,, the indicator function of ${\cal O}({\rm VII})=K'_\infty$, which contains $\mathrm{supp}(f)$ for any $f \in K_1\big(G'_\infty;\chi^{\rm VII}\big)$. This proves essentially all the assertions. 
\hfill 
$\Box$\;  

 \vskip1em 
 Recall from \S 12.1, 
 $K_\infty=G^{\mathfrak A}(m,2,\infty)=D^{\rm ev}_\infty\rtimes {\mathfrak A}_\infty$ with $D^{\rm ev}_\infty=\{d\in D_\infty\,;\,\mathrm{ord}(d)\equiv 0\:(\mathrm{mod}\;2)\}$, for $m$ even, and $K'_\infty=\Phi^{-1}(K_\infty)=\widetilde{D}^{\rm ev}_\infty\rtimes \widetilde{{\mathfrak A}}_\infty$ with $\widetilde{D}^{\rm ev}_\infty=\Phi^{-1}\big(D^{\rm ev}_\infty\big)$.
In \S 16, two involutive operations $\tau$ and $\kappa=\kappa^{(m^0)}\;(m^0=m/2)$ on $A\in {\cal A}\big(\widehat{T}\big)$ are given respectively 
 in (\ref{2015-08-09-12}) and 
(\ref{2015-07-27-21}).   
On $f_A\in E(G_\infty)$ in (\ref{lim-char-15}),  operations $P_{\langle \tau\rangle},\,P_{\langle \kappa\rangle}$ and $P_{\langle \tau,\kappa\rangle}=P_{\langle \tau\rangle}P_{\langle \kappa\rangle}$ are given also in \S 15. Let $X_{K_\infty}$ be the indicator function of $K_\infty$. 
  \vskip1.2em

{\bf Lemma 17.3.}\; 
{\it 
For $A, A^1\in{\cal A}(\widehat{T})$, 
\begin{eqnarray}
\label{2016-04-25-11}
 \nonumber
&& 
f_A|_{K_\infty}=f_{A^1}|_{K_\infty}\;\Longleftrightarrow\;A^1\in \{A,\tau A,\kappa A,\tau\kappa A\}, 
 \\
 \nonumber
&&
 P_{\langle\tau,\kappa \rangle}f_A=\tfrac{1}{4}(f_A+f_{\tau A}+f_{\kappa A} +f_{\kappa\tau A})
 =\left\{
 \begin{array}{ll}
 f_A|_{K_\infty}\quad &\mbox{\rm on $K_\infty$},
 \\[.3ex]
 \quad 0 & \mbox{\rm on $G_\infty\setminus K_\infty$},
 \end{array}
 \right.
 \\
 \nonumber
 &&
 ({\cal N}{\cal M})f_A=P_{\langle\tau,\kappa \rangle}f_A=X_{K_\infty}\cdot f_A\,.
 \end{eqnarray}
}

\vskip.2em 

{\bf Theorem 17.4.}\; 
{\it 
{\rm (i)}\; The map ${\cal M}$, when restricted on $E(G_\infty)$, gives a surjection to $E(G'_\infty;\chi^{\rm VII})$ and, for $A\in{\cal A}(\widehat{T})$,
\begin{eqnarray}
\label{2015-08-10-1}
\nonumber
E(G_\infty)\ni f_A \;\longmapsto\; f^{\rm VII}_0\hspace{.2ex}P_{\langle\tau,\kappa \rangle}f_A=f^{\rm VII}_0\,\big(f_A|_{K_\infty}\big)\in E(G'_\infty;\chi^{\rm VII}), 
\end{eqnarray}
where we understand $f^{\rm VII}_0\,\big(f_A|_{K_\infty}\big)$ is zero outside $K'_\infty$.  
Let\, $[A]=\{A, \tau A, \kappa A,\tau\kappa A\}\in{\cal A}(\widehat{T})/\langle \tau,\kappa\rangle$, then a parametrization of $
E(G'_\infty;\chi^{\rm VII})$ is given by 
$$
E(G'_\infty;\chi^{\rm VII})=\big\{P_{\langle\tau,\kappa \rangle}f_A\;;\;[A]\in{\cal A}(\widehat{T})/\langle \tau,\kappa\rangle\big\}.
$$

{\rm (ii)}\; 
The subset\, $F(G'_\infty;\chi^{\rm VII})$  consisting of factorisable elements in $E(G'_\infty;\chi^{\rm VII})$ is equal to 
\begin{eqnarray*}
\label{2015-08-10-2}
\hspace*{3ex}
E\big(G'_\infty;\chi^{\rm VII};{\cal O}({\rm str})\big)
\!\!&:=&\!\!
\big\{
f^{\rm VII}_0 f_A\,;\,A\in {\cal A}(\widehat{T}),\,\mathrm{supp}(f_A)\subset{\cal O}({\rm str})\big\},
\\[.3ex]
\!\!&=&\!\!
\big\{
f^{\rm VII}_0 f_A\,;\,A\in {\cal A}(\widehat{T}),\,A=\tau A=\kappa A=\tau\kappa A\big\},
\end{eqnarray*}
where the restrictive condition {\rm (str)} on $g=(d,\sigma)\in G_\infty=D_\infty(T)\rtimes{\mathfrak S}_\infty$  is 
\vskip.5em
\noindent
\;\;{\rm (str)}
\qquad 
$
\mathrm{ord}(\xi_{q_i})\equiv 0\;(i\in{\boldsymbol I}_r), \;\; 
L(\sigma_j)\equiv 0,\;\;\mathrm{ord}(d_j)\equiv 0\;(j\in{\boldsymbol I}_s).
$
}
\vskip1em

{\bf Note 17.5.}\; 
For 
$A=\big( (\alpha_{\zeta,\varepsilon})_{(\zeta, \varepsilon) \in \widehat{T} \times \{0,1\}};\mu \big)\in{\cal A}\big(\widehat{T}\big)$, let $\alpha_{\zeta,\varepsilon}=1$ for some $\zeta\in\widehat{T},\,\varepsilon\in\{0,1\},$ and other $\alpha_{\zeta',\varepsilon'}={\bf 0}$ and $\mu={\bf 0}$. Denote such an $A$ by $A(\zeta,\varepsilon)$. Then \,$f_{A(\zeta,\varepsilon)}(g)=X_{\zeta,\varepsilon}(g):=\zeta\big(P(d)\big)\cdot \mathrm{sgn}(\sigma)^\varepsilon\;\big({\rm for}\;g=(d,\sigma)\big)$\, is a one-dimensional character of $G_\infty$.  
The preimage in $E(G_\infty)$ under ${\cal M}$ for the normalized character of 2-dimensional spin IR 
 $\widetilde{\chi}_{\pi_{2,\zeta}}\in E(G'_\infty;\chi^{\rm VII}),\;\zeta\in\widehat{T}^0,$ is \,$\{X_{\zeta,0}\,,\,X_{\zeta\zeta^{(m^0)},0}\,,\,X_{\zeta,1}\,,\,X_{\zeta\zeta^{(m^0)},1}\}$.

\vskip1.2em

{\bf Theorem 17.6.} 
{\it The restriction map\, ${\rm Res}^{G_\infty}_{K_\infty}:\,E(G_\infty)\to E(K_\infty)$ is surjective, and the one\, ${\rm Res}^{G'_\infty}_{K'_\infty}:\,E\big(G'_\infty;\chi^{\rm VII}\big)\to E\big(K'_\infty;\chi^{\rm VII}\big)$ is bijective. Through these restriction maps we have bijective maps ${\cal M}'$ and ${\cal N}'$ between $E(K_\infty)$ and $E\big(K'_\infty;\chi^{\rm VII}\big)$ by multiplication of $f^{\rm VII}_0$ as is illustrated below:  
\vskip1em
\hspace*{10ex}
$
\begin{array}{cccll}
& 
E\big(G_\infty) 
&
\quad
&
E\big(G'_\infty;\chi^{\rm VII}\big)
\qquad\qquad\qquad
&
\\[1ex]
&
\mbox{\footnotesize ${\rm Res}^{G_\infty}_{K_\infty}$}\,\big\downarrow\,\mbox{\rm\footnotesize surj.} 
&
&
\quad \mbox{\rm\footnotesize bij.}\,\big\downarrow\,\mbox{\footnotesize ${\rm Res}^{G'_\infty}_{K'_\infty}$}
&
\\[.3ex]
&
E(K_\infty)
&
\;\;
\begin{array}{c}
\overset{{\cal M}'}{\longrightarrow}
\\[-2.2ex]
\underset{{\cal N}'}{\longleftarrow}
\end{array}
\;\;
& 
E\big(K'_\infty;\chi^{\rm VII}\big)\quad\;{\rm (bijections).}
&
\end{array}
$
}

\section{\large Spin characters of spin type $\chi^{\rm I}=(-1,-1,-1)$}

Let $G'_\infty=R\big(G(m,1,\infty)\big)=\widetilde{D}^\vee_\infty(T)\rtimes\widetilde{{\mathfrak S}}_\infty, \,\widetilde{D}^\vee_\infty(T)=Z_3\times D^\wedge_\infty(T).$ 
Take a $g'=(d',\sigma')\in G'_\infty$ (modulo $z\in Z$) and let its standard decomposition be    
\begin{eqnarray}
\label{2017-07-31-1-99}
g'=\xi'_{q_1}\xi'_{q_2}\cdots\xi'_{q_r}g'_1g'_2\cdots g'_s,&
\xi'_q = \big(t'_q, (q)\big)={\widehat{\eta}_q}^{\;a_q},\;\;
g'_j=(d'_j, \sigma'_j),
\end{eqnarray} 
and put its image $g=\Phi(g')$ as 
$g=\xi_{q_1}\xi_{q_2}\cdots\xi_{q_r}g_1g_2\cdots g_s,\,
\xi_q = \big(t_q, (q)\big)=y_q^{\;a_q},\,
g_j=(d_j, \sigma_j),
$ with $y_q=\Phi(\widehat{\eta}_q),\,\sigma_j=\Phi(\sigma'_j)\;\,\mbox{\rm a cycle},\;d_j=\Phi(d'_j),$ where    
\,
$
K_j:=\mathrm{supp}(\sigma_j) \supset \mathrm{supp}(d_j)
$\, and  
$K_j$'s are mutually disjoint. 
We put \,$J:={\boldsymbol I}_s,\,Q:=\{q_1,q_2,\ldots,q_r\}$.

For the support of $f\in K(G'_\infty;\chi^{\rm Y})$ for Y\,=\,I, we know by Table 12.1 that    
\begin{eqnarray}
\label{2015-07-29-21}
\nonumber
&&
\mathrm{supp}(f)\subset {\cal O}({\rm I}):=\big\{g'\in G'_\infty\,;\;\mbox{\rm $g'$\;satisfies (Condition I)}\big\},
\qquad
\\[.3ex]
\nonumber
\mbox{\rm (Condition I)}
&&\quad
\mathrm{ord}(\xi'_{q_i})\equiv 0\;(\forall i),\;\;\mathrm{ord}(d'_j)+L(\sigma'_j)\equiv 0\;(\mathrm{mod}\;2) \;(\forall j).
\hspace{13ex}
\end{eqnarray}
\vskip.2em

For $n\ge 4$, characters of IRs of $G'_n$ of spin type $\chi^{\rm I}$ is given in \S 6.3.2, and we apply them to calculate their pointwise limits as $n\to\infty$. The following lemma is a direct consequence of Theorem 6.10. 
Let   
$g'=(d',\sigma')\in G'_\infty=\widetilde{D}^\vee_\infty\rtimes\widetilde{{\mathfrak S}}_\infty$\, be as in (\ref{2017-07-31-1-99}) and take $z\in Z$, then \,$f^{\rm I}_0(zg')=\chi^{\rm I}(z)f^{\rm I}_0(g')$.

\bigskip

{\bf Proposition 18.1.}  
{\it 
If $g'=(d',\sigma')$ satisfies {\rm (Condition I)}, then 
\begin{eqnarray}
\label{2017-08-21-1}
f^{\rm I}_0(g'):=\lim_{n=2n'\to\infty}\widetilde{\chi}(P^0_n|g')=\prod_{j\in J}(-1)^{[(\ell_j-1)/2]}\, 2^{-(\ell_j-1)/2}.
\end{eqnarray}
If $g'$ does not satisfy {\rm (Condition I)}, then $f^{\rm I}_0(g')=0.$\, 
Moreover 
$$  
\lim_{n=2n'+1\to\infty}\widetilde{\chi}(P^\varepsilon_n|g')=f^{\rm I}_0(g'),
\quad\mbox{\rm for \,$\varepsilon=+,-$}. 
$$
}

\vskip.2em

For the general case we obtain the following results  by using Theorem 6.11.

\vskip1.2em

{\bf Lemma 18.2.}  
{\it 
Take a series 
$\Lambda^n=\big({\boldsymbol \lambda}^\zeta)_{\zeta\in{\cal K}} \in {\boldsymbol Y}_{\!n}({\cal K})$ with ${\cal K}=\widehat{T}^0$.    

If $g'\in G'_\infty$ satisfies {\rm (Condition I)}, then  
$$
\lim_{n\to\infty}\widetilde{\chi}\big(\Pi^{\rm \;I}_{\Lambda^n}\big|g'\big)
=\lim_{n\to\infty}\widetilde{\chi}\big(\Pi^{\rm I0}_n\big|g'\big)\times\lim_{n\to\infty}
\widetilde{\chi}\big(\breve{\Pi}_{\Lambda^n}|g\big).
$$ 
If $g'$ does not satisfy \,{\rm (Condition I)},  then 
\,${\displaystyle 
\lim_{n\to\infty}\widetilde{\chi}\big(\Pi^{\rm \;I}_{\Lambda^n}\big|g'\big)=0}$\, if exists. 
}

\vskip1.2em

Put $\widehat{T}^0=\{\zeta^{(a)}\in\widehat{T}\,;\,0\le a<m^0=m/2\}$, and introduce a parameter 
\begin{gather}
\label{2015-07-31-1}
A^0 = \big( (
\alpha_{\zeta,\varepsilon})_{(\zeta, \varepsilon) \in \widehat{T}^0 \times \{0,1\}} 
\,;\,\mu \big)\quad {\rm with}
\\
\nonumber
\alpha_{\zeta,\varepsilon} 
  := (\alpha_{\zeta,\varepsilon;i})_{i\geqslant 1}\,,\;\;\alpha_{\zeta,\varepsilon;1}\ge \alpha_{\zeta,\varepsilon;2}\ge \,\cdots\,\ge 0, 
\\
\nonumber
\mu:=(\mu_\zeta)_{\zeta\in\widehat{T}^0},\;\mu_\zeta\ge 0,
\end{gather}
and a parameter space \,${\cal A}(\widehat{T}^0)$ consisting of $A^0$ satisfying 
\begin{eqnarray} 
\label{2017-08-22-2}
\sum_{\zeta \in \widehat{T}^0}\,
\sum_{\varepsilon = 0, 1} 
\|\alpha_{\zeta,\varepsilon}\|
+ \|\mu\| = 1.
\end{eqnarray}
We can naturally consider ${\cal A}(\widehat{T}^0)$ as a subset of ${\cal A}(\widehat{T})$, and so we can define $f_{A^0}$ for $A^0\in{\cal A}(\widehat{T}^0)\subset{\cal A}(\widehat{T})$\, as in (\ref{lim-char-15}). 

\vskip1.2em

{\bf Theorem 18.3.}\; 
{\it 
{\rm (i)}\, 
 The set\, ${\rm Lim}(G'_\infty;\chi^{\rm I})$ of all pointwise limits 
${\displaystyle \lim_{n\to\infty}\widetilde{\chi}\big(\Pi^{\rm \;I}_{\Lambda^n}\big|\cdot\big)}$
of normalised irreducible characters of spin type $\chi^{\rm I}$ of $G'_n$ is given as 
\begin{eqnarray} 
\label{2017-08-22-3}
{\rm Lim}(G'_\infty;\chi^{\rm I})=
\big\{f^{\rm I}_0\,f_{A^0}\;;\,A^0\in{\cal A}(\widehat{T}^0)\big\}.
\end{eqnarray}

Let \,$g'=\xi'_{q_1}\xi'_{q_2} \cdots \xi'_{q_r} g'_1g'_2\cdots g'_s$\, be an element of $G'_\infty$ in 
(\ref{2017-07-31-1-99}), and let its image $g=\Phi(g')\in G_\infty$ be $g=\xi_{q_1}\xi_{q_2} \cdots \xi_{q_r} g_1g_2\cdots g_s$. 
If\, $g'\not\in{\cal O}({\rm I})$, then\, $f^{\rm I}_0(g')f_{A^0}(g)=0$. 

If\, $g'\in{\cal O}({\rm I})$, then with\, $\ell_j=\ell(\sigma'_j)$, 
\begin{alignat}{2}
\label{2015-07-31-13}
\nonumber
\quad 
&
F^{\rm \,I}_{A^0}(g'):=f^{\rm I}_0(g')f_{A^0}(g)  
=
\prod_{q \in Q}
\Bigg\{ \sum_{\zeta \in \widehat{T}^0}
\Big(
 \sum_{\varepsilon \in \{0,  1\}}
\, \|\alpha_{\zeta,\varepsilon}\|   
\; + \;\mu_\zeta\Big)
\zeta(t_{q})\Bigg\} \times
\\
\nonumber
&
\quad
\times 
\prod_{j \in J} 
\Bigg\{ (-1)^{[(\ell_j-1)/2]} 
2^{-(\ell_j-1)/2}
\sum_{\zeta \in \widehat{T}^0}\Big( \sum_{\varepsilon \in \{0,  1\}} 
\,
\sum_{i \in \boldsymbol{N}}\;(\alpha_{\zeta,\varepsilon;i})^{\ell_j}
\chi_\varepsilon(\sigma_j)
\Big) 
\zeta\big(P(d_j)\big)
\Bigg\}. 
\end{alignat}

The limit function \,$F^{\rm \,I}_{A^0}$ is factorisable, and so  
\,${\rm Lim}(G'_\infty;\chi^{\rm I})\subset F(G'_\infty;\chi^{\rm I})$. They all equal to $E(G'_\infty;\chi^{\rm I})$.  

{\rm (ii)}\; 
For $A^0,A^1\in{\cal A}(\widehat{T}^0)$,\, $F^{\rm I}_{A^0}=F^{\rm I}_{A^1}$ if and only if $A^0=A^1$. Hence the following gives a parametrization of characters of spin type $\chi^{\rm I}=(-1,-1,-1)$:   
$$
{\cal A}(\widehat{T}^0)\ni A^0\longmapsto\,F^{\rm \,I}_{A^0}=f^{\rm I}_0\,f_{A^0}\in E(G'_\infty;\chi^{\rm I}).
$$ 
}

\vskip.2em

{\it Proof.}\; A proof of the main part was given in [II, Part IV]. For (i), we remark that, by calculating ${\rm Lim}(G'_\infty;\chi^{\rm I})$, another proof for the validity of the criterion (EF) is obtained here. 

For (ii), 
assume $F^{\rm I}_{A^0}=F^{\rm I}_{A^1}$ and let $A^1= \big( (
\alpha^1_{\zeta,\varepsilon})_{(\zeta, \varepsilon) \in \widehat{T}^0 \times \{0,1\}} 
\,;\,\mu^1 \big),\,\mu^1=(\mu^1_\zeta)_{\zeta\in\widehat{T}^0}$. We may concentrate only one term corresponding to $g_j=(d_j,\sigma_j),\,\mathrm{ord}(d'_j)+L(\sigma'_j)=\mathrm{ord}(d_j)+(\ell_j-1)\equiv 0\;(\mathrm{mod}\;2)$. 
In case $\ell_j-1=2c\equiv 0\;(\mathrm{mod}\;2)$, we have 
\begin{gather*}
\;\;\;\sum_{\zeta \in \widehat{T}^0}\Big(
\sum_{i \in \boldsymbol{N}}\big((\alpha_{\zeta,0;i})^{2c+1}+(\alpha_{\zeta,1;i})^{2c+1}
\big)
\Big) 
\zeta(y^{2b})
\\
=
\sum_{\zeta \in \widehat{T}^0}\Big(
\sum_{i \in \boldsymbol{N}}\big((\alpha^1_{\zeta,0;i})^{2c+1}+(\alpha^1_{\zeta,1;i})^{2c+1}
\big)
\Big) \zeta(y^{2b}).
\end{gather*}
With $\zeta=\zeta^{(a)},\,\zeta^{(a)}(y)=\omega^a,\,\omega=e^{2\pi i/m}$, the matrix of coefficients is given as 
$$
\big(\zeta^{(a)}(y^{2b})\big)=(\omega^{2ab})\quad \mbox{\rm for $0\le a, b<m^0,$ of type $m^0\times m^0$. }
$$
Since $\det(\omega^{2ab})=\prod_{m^0>p>q\geqslant 0}(\omega^{2p}-\omega^{2q})\ne 0$, we obtain 
\begin{eqnarray}
\label{2017-08-23-1}
&&
\sum_{i \in \boldsymbol{N}}\big((\alpha_{\zeta,0;i})^{2c+1}+(\alpha_{\zeta,1;i})^{2c+1}
\big)
=
\sum_{i \in \boldsymbol{N}}\big((\alpha^1_{\zeta,0;i})^{2c+1}+(\alpha^1_{\zeta,1;i})^{2c+1}
\big).
\end{eqnarray}

In case $\ell_j-1=2c+1\equiv 1\;(\mathrm{mod}\;2)$, 
\begin{gather*}
\sum_{\zeta \in \widehat{T}^0}\Big(
\sum_{i \in \boldsymbol{N}}\big((\alpha_{\zeta,0;i})^{2c+2}-(\alpha_{\zeta,1;i})^{2c+2}
\big)
\Big) 
\zeta(y^{2b+1})=\hspace{15ex}
\\[-.5ex]
\hspace{15ex}
=\sum_{\zeta \in \widehat{T}^0}\Big(
\sum_{i \in \boldsymbol{N}}\big((\alpha^1_{\zeta,0;i})^{2c+2}-(\alpha^1_{\zeta,1;i})^{2c+2}
\big)
\Big) \zeta(y^{2b+1}).
\end{gather*} 
The matrix of coefficients is $(\omega^{a(2b+1)})$ for $0\le a, b<m^0$, and $\det(\omega^{a(2b+1)})\ne 0$. So  
\begin{gather}
\label{2017-08-23-2}
\sum_{i \in \boldsymbol{N}}\big((\alpha_{\zeta,0;i})^{2c+2}-(\alpha_{\zeta,1;i})^{2c+2}
\big)
=\sum_{i \in \boldsymbol{N}}\big((\alpha^1_{\zeta,0;i})^{2c+2}-(\alpha^1_{\zeta,1;i})^{2c+2}
\big).
\end{gather} 

Then equalities (\ref{2017-08-23-1})--(\ref{2017-08-23-2}) give us $\alpha_{\zeta,0;i}=\alpha^1_{\zeta,0;i}\;(i\in{\boldsymbol N})$ as desired. 
\hfill
$\Box$

\section{\large Spin characters of spin type $\chi^{\rm II}=(-1,-1,\,1)$}

Let the notations be as in \S 18, and recall the normal subgroups of $G_\infty$ and $G'_\infty$ as 
\begin{eqnarray*}
&&
G^{\mathfrak A}_\infty=G^{\mathfrak A}(m,1,\infty)=D_\infty\rtimes {\mathfrak A}_\infty,\quad 
K_\infty=D^{\rm ev}_\infty\rtimes {\mathfrak A}_\infty,   
\\
&& 
G^{\prime\hspace{.15ex}{\mathfrak A}}_\infty:=\Phi^{-1}(G^{\mathfrak A}_\infty)=\widetilde{D}^\vee_\infty\rtimes \widetilde{{\mathfrak A}}_\infty,
\hspace{4.7ex}
K'_\infty=\Phi^{-1}(K_\infty)=\widetilde{D}^{\rm ev}_\infty\rtimes \widetilde{{\mathfrak A}}_\infty, 
\end{eqnarray*}
where $D^{\rm ev}_\infty=\{d\in D_\infty\,;\,\mathrm{ord}(d)\equiv 0\:(\mathrm{mod}\;2)\}$,\;  
$\widetilde{D}^{\rm ev}_\infty=\Phi^{-1}(D^{\rm ev}_\infty).$  
Then ${\cal O}({\rm VII})=K'_\infty$. 

Let \,$f^{\rm VII}_0=\widetilde{\chi}_{\pi_{2,\zeta^{(0)}}}=\mathrm{sgn}^{\mathfrak A}_{Z_3}\cdot X_{{\cal O}({\rm VII})}$\, be as in Proposition 12.6. For $\chi^{\rm Y}\in\widehat{Z}$, we define a multiplication map ${\cal M}:\,K_1\big(G'_\infty;\chi^{\rm Y}\big)\ni F\mapsto f=f^{\rm VII}_0\cdot F\in K_1\big(G'_\infty;\chi^{\rm VII}\chi^{\rm Y}\big),$ and another one ${\cal N}:\,K_1\big(G'_\infty;\chi^{\rm VII}\chi^{\rm Y}\big)\ni f\mapsto f^{\rm VII}_0\cdot f\in K_1\big(G'_\infty;\chi^{\rm Y}\big)$ as 
\begin{eqnarray}
\label{2015-08-16-1}
K_1\big(G'_\infty;\chi^{\rm Y}\big)
\begin{array}{c}
\mbox{\footnotesize ${\cal M}$}
\\[-1.7ex]
\longrightarrow
\\[-2.3ex]
\longleftarrow
\\[-1.5ex]
\mbox{\footnotesize ${\cal N}$}
\end{array}
 K_1\big(G'_\infty;\chi^{\rm VII}\chi^{\rm Y}\big).
\end{eqnarray}
If 
\,$\chi^{\rm Y}=(\beta_1,\beta_2,\beta_3)$, then \,$\chi^{\rm VII}\chi^{\rm Y}=(\beta_1,\beta_2,-\beta_3)$. 
Moreover\, 
$({\cal N}{\cal M})F=X_{{\cal O}({\rm VII})}\,F,$ \,$({\cal M}{\cal N})f=X_{{\cal O}({\rm VII})}\,f$.\,  
For an  
$F\in K_1\big(G'_\infty;\chi^{\rm Y}\big),\,{\rm Y}={\rm I}\sim{\rm VII},$ there holds $\mathrm{supp}(F)\subset {\cal O}({\rm Y})$ and so \,$\mathrm{supp}\big({\cal M} F\big)\subset {\cal O}({\rm VII})\cap {\cal O}({\rm Y})=K'_\infty\cap {\cal O}({\rm Y})$. 

We pick up a special case of Theorem 15.1 for $G'_\infty$ and $K'_\infty$. 
 
\vskip1.2em 

{\bf Theorem 19.1.}\; 
{\it 
 The restriction map \,${\rm Res}^{G'_\infty}_{K'_\infty}:\,K_1(G'_\infty)\rightarrow K_1(K'_\infty)$\, is surjective and so are also for\, $E\big(G'_\infty\big)\to E\big(K'_\infty\big)$\, and for 
$E\big(G'_\infty;\chi^{\rm Y}\big)\to E\big(K'_\infty;\chi^{\rm Y}\big)$ in case of\, {\rm Y\,=\,I\,$\sim$\,VII}. 
}

\vskip1.2em

We have studied in \S 17 the case of $\chi^{\rm VIII}=(1,1,1)$ and $\chi^{\rm VII}\chi^{\rm VIII}=\chi^{\rm VII}=(1,1,-1)$ (see Theorem 17.6). In this section we study the case of $\chi^{\rm I}$ and $\chi^{\rm II}=\chi^{\rm VII}\chi^{\rm I}$. 
Let \,$g'=(d',\sigma')=\xi'_{q_1}\xi'_{q_2} \cdots \xi'_{q_r} g'_1g'_2\cdots g'_s\in G'_\infty,\,g'_j=(d'_j,\sigma'_j),$ be as in 
(\ref{2017-07-31-1-99}), with its image $g=\Phi(g')\in G_\infty$ as $g=\xi_{q_1}\xi_{q_2} \cdots \xi_{q_r} g_1g_2\cdots g_s$. 

\vskip1.2em 

{\bf Lemma 19.2.}\; 
{\it 
{\rm (i)}\; 
${\cal O}({\rm II})={\cal O}({\rm VII})\cap{\cal O}({\rm I})\subset K'_\infty$.
\vskip.3em

{\rm (ii)}\; Put \,$f^{\rm II}_0(g'):=\!\!{\displaystyle \lim_{n=2n'\to\infty}\widetilde{\chi}(P^{\rm II0}_n|g')=\!\!\lim_{n=2n'+1\to\infty}\widetilde{\chi}(\Pi^{\rm \mho^\varepsilon,{\mathfrak A}}_n|g') }$. Then,   
\begin{eqnarray}
&&
\label{2015-08-23-1}
f^{\rm\,II}_0(zg')= 
\chi^{\rm II}(z)\cdot\prod_{j\in J}(-1)^{[(\ell_j-1)/2]}\, 2^{-(\ell_j-1)/2}\;\;(z\in Z).
\end{eqnarray}
for $g'\in{\cal O}({\rm II})$, and $f^{\rm\,II}_0(g')=0$ for $g'\not\in{\cal O}({\rm II})$. Thus\, 
$f^{\rm II}_0=f^{\rm VII}_0\hspace{.15ex}f^{\rm I}_0$.
}

\vskip1em

{\it Proof.}\; (i)\; (Condition I) is given just after (\ref{2017-07-31-1-99}) and (Condition II) is given in the notation in (\ref{2017-07-31-1-99}) as 
\begin{eqnarray}
\label{2015-08-14-11}
\nonumber
\mbox{\rm (Condition II)}\qquad 
\left\{
\begin{array}{l}
\mathrm{ord}(d')\equiv 0\;(\mathrm{mod}\;2),\;\;L(\sigma')\equiv 0\;(\mathrm{mod}\;2)\,;
\\[.5ex]
\mathrm{ord}(\xi'_{q_i})\equiv 0\;(\forall i),\;\mathrm{ord}(d'_j)+L(\sigma'_j)\equiv 0\;(\mathrm{mod}\;2) \;\,(\forall j). 
\end{array}
\right.
\qquad\qquad
\end{eqnarray}

(ii)\; By Lemma 7.11 and Theorem 7.15. 
\hfill 
$\Box$
\bigskip

{\bf Theorem 19.3.} 
{\it 
Let the spin type be 
$\chi^{\rm II}=(-1,-1,\,1)$. The set of all pointwise limits as $n\to\infty$ of series of normalized spin irreducible characters in $E(G'_n;\chi^{\rm II})$ is given as 
\begin{eqnarray}
\label{2009-08-07-12-2}
\nonumber 
 {\rm Lim}\big(G'_\infty\;\chi^{\rm II}\big)
=
\big\{F^{\rm \,II}_{A^0}=f^{\rm \,II}_0 f_{A^0}\;;\; A^0\in {\cal A}\big(\widehat{T}^0\big)\big\},
\end{eqnarray}
where 
 $F^{\rm II}_{A^0}(g')=0$ outside ${\cal O}({\rm II})$, and for $g'\in{\cal O}({\rm II})$, with 
$\ell_j=\ell(\sigma'_j)$, 
\begin{alignat}{1}
\nonumber
&
F^{\rm \,II}_{A^0}(g')
=
f^{\rm II}_0(g')f_{A^0}(g)  
\\
\nonumber
&\hspace*{7ex}
=
\prod_{q \in Q}
\Bigg\{ \sum_{\zeta \in \widehat{T}^0}
\Big(
 \sum_{\varepsilon \in \{0,  1\}}
\, \|\alpha_{\zeta,\varepsilon}\|   
\; + \;\mu_\zeta\Big)
\zeta(t_{q})\Bigg\} \times
\\
\nonumber
&\hspace{3ex}
\times 
\prod_{j \in J} 
\Bigg\{ (-1)^{[(\ell_j-1)/2]} 
2^{-(\ell_j-1)/2}
\sum_{\zeta \in \widehat{T}^0}\Big( \sum_{\varepsilon \in \{0,  1\}} 
\,
\sum_{i \in \boldsymbol{N}}\;(\alpha_{\zeta,\varepsilon;i})^{\ell_j}
\chi_\varepsilon(\sigma_j)
\Big) 
\zeta\big(P(d_j)\big)
\Bigg\}. 
\end{alignat}
}
\vskip.2em

{\it Proof.} By direct calculation using Theorems 7.12,  7.17 and Lemma 19.2.
\hfill
$\Box$\;

\vskip1.2em

As we see easily, every $f\in{\rm Lim}\big(G'_\infty;\chi^{\rm II}\big)$ is not necessarily factorisable, and so the criterion (EF) does not holds. This is also a motivation to study relations between $E\big(G'_\infty\,;\,\chi^{\rm Y}\big)$ for Y\,=\,I and II. 

Recall that \,$E\big(G'_\infty\,;\,\chi^{\rm I}\big)=\big\{F^{\rm \,I}_{A^0}=f^{\rm I}_0 f_{A^0}\,;\,A^0\in {\cal A}\big(\widehat{T}^0\big)\big\}$, and $\mathrm{supp}(f^{\rm II}_0)\subset K'_\infty$.

\vskip1.2em

{\bf Lemma 19.4.}\; 
{\it 
{\rm (i)}\; In spin type $\chi^{\rm I}$, ${\rm Res}^{G'_\infty}_{K'_\infty}\big(F^{\rm\,I}_{A^0}\big)={\rm Res}^{G'_\infty}_{K'_\infty}\big(F^{\rm\,I}_{A^1}\big)$\, for\, $A^0, A^1\in{\cal A}(\widehat{T}^0)$\, if and only if\, $A^1=A^0$ or $A^1={}^t\!A^0$, where ${}^t\!A^0={}^t(A^0)=\tau A^0$ is in (\ref{2015-08-09-12}). 

{\rm (ii)}\;  The restriction map\, ${\rm Res}^{G'_\infty}_{K'_\infty}:\,E\big(G'_\infty;\chi^{\rm II}\big)\to E\big(K'_\infty;\chi^{\rm II}\big)$\,  is bijective. 
}

\bigskip

{\it Proof.}\; 
(i)\; On the set $K'_\infty={\cal O}({\rm VII})$, \;$\prod_{j\in J}\chi_1(\sigma_j)=(-1)^{L(\sigma)}=1\;(\varepsilon=1)$. So  
$$
{\displaystyle 
\;
\prod_{j\in J}\sum_{\varepsilon \in \{0,  1\}} 
\!(\alpha_{\zeta,\varepsilon+1;i})^{\ell_j}
\chi_{\varepsilon}(\sigma_j)=\prod_{j\in J}\sum_{\varepsilon \in \{0,  1\}}\!(\alpha_{\zeta,\varepsilon;i})^{\ell_j}
\chi_\varepsilon(\sigma_j).
}
$$

Using the explicit form of $F^{\rm \,I}_{A^0}=f^{\rm I}_0\cdot f_{A^0}$ in Theorem 18.3, we see easily \,
${\rm Res}^{K'_\infty}_{G'_\infty}\big(F^{\rm\,I}_{A^0}\big)={\rm Res}^{K'_\infty}_{G'_\infty}\big(F^{\rm\,I}_{{}^t\!A^0}\big)$. 

Conversely assume that 
$F^{\rm \,I}_{A^0}|_{K'_\infty}=F^{\rm \,I}_{A^1}|_{K'_\infty}$. 
We wish to calculate as in the proof of Theorem 18.3. 
However in the present case, the functions $F^{\rm \,I}_{A^k}|_{K'_\infty}\;(k =0,1)$ in the both side are not exactly factorisable, but factorisable in a restrictive sense so that the factors with $\ell_j-1=L(\sigma'_j)\equiv 0\;(\mathrm{mod}\;2)$ are counted as independent but those with $\ell_j-1\equiv 1\;(\mathrm{mod}\;2)$ should be coupled into pairs. Hence, for $\ell_j-1\equiv 0$,  
we get exactly (\ref{2017-08-23-1}). But, for any pair $(j,j')$ with $\ell_j-1\equiv\ell_{j'}-1\equiv 1$,  we get a variant of (\ref{2017-08-23-2}) as 
\begin{gather}
\label{2017-08-24-1}
\sum_{i \in \boldsymbol{N}}\big((\alpha_{\zeta,0;i})^{2c+2}-(\alpha_{\zeta,1;i})^{2c+2}
\big)
=(-1)^a\sum_{i \in \boldsymbol{N}}\big((\alpha^1_{\zeta,0;i})^{2c+2}-(\alpha^1_{\zeta,1;i})^{2c+2}
\big),
\end{gather} 
where the sign $(-1)^a=1$ or $-1$ is a fixed constant. Depending to it, we get $A^1=A^0$ or $A^1={}^t\!A^0$. 

(ii)\, Noting that, for $f\in K_1(G'_\infty;\chi^{\rm II})$, $\mathrm{supp}(f)\subset{\cal O}({\rm II})\subset K'_\infty={\cal O}({\rm VII})$, we see the maps are injective. 
\hfill 
$\Box$

\vskip1.2em

{\bf Theorem 19.5.}\; 
{\it 
{\rm (i)}\; When restricted on the sets of characters, the multiplication maps\, ${\cal M},{\cal N}$ of\, $f^{\rm\;VII}_0$ give bijective maps ${\cal M}',\;{\cal N}'$ between $E\big(K'_\infty;\chi^{\rm I}\big)$ and  
$E\big(G'_\infty;\chi^{\rm II}\big)$, $\chi^{\rm II}=(-1,-1,\,1)$, described as follows: 

\vskip1.2em
$
\begin{array}{ccccc}
 &
E\big(G'_\infty;\chi^{\rm I}\big) 
&
\quad
&
E\big(G'_\infty;\chi^{\rm II}\big)
\\[1ex]
&
\mbox{\footnotesize ${\rm Res}^{G'_\infty}_{K'_\infty}$}\,\big\downarrow\,\mbox{\rm\footnotesize surj.} 
&
&
\mbox{\rm\footnotesize bij.}\,\big\downarrow\,\mbox{\footnotesize ${\rm Res}^{G'_\infty}_{K'_\infty}$}
&
\\[.3ex]
&
E\big(K'_\infty;\chi^{\rm I}\big)
&
\begin{array}{c}
\overset{{\cal M}'}{\longrightarrow}
\\[-2.2ex]
\underset{{\cal N}'}{\longleftarrow}
\end{array}
&
E\big(K'_\infty;\chi^{\rm II}\big)
&{\rm (isom.)}
\end{array}
\begin{array}{ccccc}
 & 
  F^{\rm \,I}_{A^0}=f^{\rm I}_0\,f_{A^0} 
  &
  & 
  F^{\rm \,II}_{A^0}=f^{\rm \,II}_0\,f_{A^0} 
  & 
  \\
 &
 \downarrow
 &
 &
 \updownarrow
 &
 \\
 &
F^{\rm \,I}_{A^0}\big|_{K'_\infty}
&
  \!\!\!\rightleftarrows\!\! \!
&
F^{\rm \,II}_{A^0}\big|_{K'_\infty}
&
\end{array}
$
\\[1.5ex]

{\rm (ii)} 
The set of characters of 
 $G'_\infty$ of spin type $\chi^{\rm II}$ is parametrized as 
\begin{eqnarray}
\label{2009-08-07-12-2-9}
\nonumber 
 E\big(G'_\infty;\chi^{\rm II}\big)
=
\big\{F^{\rm \,II}_{A^0}=f^{\rm \,II}_0 f_{A^0}\;;\; [A^0]:=\{A^0,\tau A^0\}\in {\cal A}\big(\widehat{T}^0\big)/\langle \tau \rangle\big\}.
\end{eqnarray}
}
\vskip.2em

{\bf Theorem 19.6.}\; 
{\it 
A character $F^{\rm \,II}_{A^0}=f^{\rm II}_0(g')f_{A^0}\in E(G'_\infty;\,\chi^{\rm II}),\,A^0\in {\cal A}(\widehat{T}^0),$ is factorisable if and only if \;${}^t\!A^0=A^0$. 
 The set of factorisable characters of \,$G'_\infty$ of spin type $\chi^{\rm II}$ is given as 
 $
 F(G'_\infty;\,\chi^{\rm II})=\{F^{\rm \,II}_{A^0};\,A^0\in {\cal A}(\widehat{T}^0),{}^t\!A^0=A^0\}$ and also as 
 $$
\quad F(G'_\infty;\,\chi^{\rm II})=E\big(G'_\infty;\,\chi^{\rm II}; {\cal O}({\rm str})\big):=\big\{f\in E\big(G'_\infty;\,\chi^{\rm II})\,;\,\mathrm{supp}(f)\subset {\cal O}({\rm str})\big\}, 
 $$ 
where ${\cal O}({\rm str}):={\cal O}({\rm V})={\cal O}({\rm VI})$ in Table 12.1. 
}
\vskip1em

{\it Proof.}\; 
If $F^{\rm \,II}_{A^0}$ is factorisable, then $F^{\rm \,II}_{A^0}(g')=0$ if $g'\in{\cal O}({\rm II})$ has a fundamental component $g'_j=(d'_j,\sigma')$ with $\ell_j=\ell(\sigma'_j)$ even and $\mathrm{ord}(d'_j)$ odd by (Condition II). 
Put $\mathrm{ord}(d'_j)=2b+1$, then we have 
\begin{alignat*}{2}
0&=\sum_{\zeta \in \widehat{T}^0}\Bigg( \sum_{\varepsilon \in \{0,  1\}} 
\,
\sum_{i \in \boldsymbol{N}}\;(\alpha_{\zeta,\varepsilon;i})^{\ell_j}
\chi_\varepsilon(\sigma_j)
\Bigg) 
\zeta\big(P(d_j)\big)
\hspace{15ex}
\\
\hspace*{12ex} 
&=
\sum_{0\leqslant a<m^0}\Bigg( \sum_{i \in \boldsymbol{N}}\;\Big((\alpha_{\zeta,0;i})^{\ell_j}-(\alpha_{\zeta,1;i})^{\ell_j}
\Big)\Bigg)\cdot\omega^{a(2b+1)}.
\end{alignat*}
From this, we get\; 
${\displaystyle 
{\sum}_{i \in \boldsymbol{N}}\big((\alpha_{\zeta,0;i})^{\ell_j}-(\alpha_{\zeta,1;i})^{\ell_j}
\big)=0\;(\forall \zeta\in\widehat{T}^0). 
}$ 
 Hence\, $(\alpha_{\zeta,0;i})^2=(\alpha_{\zeta,1;i})^2$, and\, $\alpha_{\zeta,0;i}=\alpha_{\zeta,1;i}\;(i\in{\boldsymbol N},\,\zeta\in\widehat{T}^0)$, and finally ${}^t\!A^0=A^0$.
\hfill
$\Box$

\section{\large Spin characters of spin type $\chi^{\rm III}=(-1,\,1,-1)$ and $\chi^{\rm IV}=(-1,\,1,\,1)$}

{\bf 20.1. Method of obtaining characters.} 
\vskip.5em

The spin types still remained can be divide into two pairs as  
\vskip.4em
\quad
$\chi^{\rm IV}=(-1,\,1,\,1)$ and 
$\chi^{\rm III}=(-1,\,1,-1)=\chi^{\rm VII}\chi^{\rm VI},$ 

\quad
$\chi^{\rm VI}=(1,-1,\,1)$ and 
$\chi^{\rm V}=(1,-1,-1)=\chi^{\rm VII}\chi^{\rm VI}.$ 
\\[.6ex]
For these two pairs respectively, we can apply the method similar to those for $(\chi^{\rm VIII}, \chi^{\rm VII})$ and $ (\chi^{\rm I}, \chi^{\rm II})$, that is, using multiplication maps ${\cal M}$ and ${\cal N}$ of  \,$f^{\rm VII}_0=\widetilde{\chi}_{\pi_{2,\zeta^{(0)}}}=\mathrm{sgn}^{\mathfrak A}_{Z_3}\cdot X_{{\cal O}({\rm VII})}$\,, given as 
\begin{eqnarray}
\label{2015-08-23-11}
\nonumber 
K_1\big(G'_\infty;\chi^{\rm IV}\big)\;
\begin{array}{c}
\mbox{\rm \footnotesize ${\cal M}$}
\\[-1.2ex]
\rightleftarrows
\\[-1.2ex]
\mbox{\rm\footnotesize ${\cal N}$}
\end{array}
\;K_1\big(G'_\infty;\chi^{\rm III}\big),
\qquad 
K_1\big(G'_\infty;\chi^{\rm VI}\big)\;
\begin{array}{c}
\mbox{\rm \footnotesize ${\cal M}$}
\\[-1.2ex]
\rightleftarrows
\\[-1.2ex]
\mbox{\rm \footnotesize ${\cal N}$}
\end{array}
\;K_1\big(G'_\infty;\chi^{\rm V}\big),
\end{eqnarray}
we can determine $E\big(G'_\infty;\chi^{\rm III}\big)$ from $E\big(G'_\infty;\chi^{\rm IV}\big)$, and 
$E\big(G'_\infty;\chi^{\rm VI}\big)$ from $E\big(G'_\infty;\chi^{\rm V}\big).$
The background of this method is a relation among supports of 
$f\in K_1\big(G'_\infty;\chi^{\rm Y}\big)$, given below. Let 
\begin{eqnarray}
\label{2017-08-27-1}
&&\;\;
\left\{
\begin{array}{l}
g'=(d',\sigma')=z\,
\xi'_{q_1}\xi'_{q_2}\cdots \xi'_{q_r}g'_1g'_2 \cdots g'_s\in G'_\infty,\;z\in Z,\;g_j=(d'_j,\sigma'_j), 
\\[.5ex]
g=\Phi(g')=(d, \sigma) = \xi_{q_1}\xi_{q_2} \cdots \xi_{q_r}\,g_1g_2 \cdots g_s\in G_\infty,\;g_j=(d_j,\sigma_j),
\end{array}
\right.
\\
\nonumber
&&\quad
Q:=\{q_1,q_2,\ldots,q_r\},\;J:={\boldsymbol I}_s,\;K_j:=\mathrm{supp}(\sigma'_j),\;\ell_j:=\ell(\sigma'_j)\;(j\in J),
\end{eqnarray}
be as in (\ref{2017-07-31-1}). Then, for Y= III\,$\sim$\,VII, we see from Table 12.1, 
\vskip.5em 

(Condition III)\qquad 
$
 \mathrm{ord}(d)\equiv 0,\;\;
 L(\sigma_j):=\ell_j-1 \equiv 0\;(\forall j)\quad(\mathrm{mod}\;2)\,;$
\vskip.2em

(Condition IV)\qquad
$L(\sigma_j)=\ell_j-1 \equiv 0\;(\forall j)\,;$
\vskip.2em

(Condition V)\hspace*{5ex}  
$
\mathrm{ord}(\xi_{q_i})\equiv 0\;(\forall i),
\;\;
{\rm ord}(d_j)\!\equiv\! 0\;(\forall j),
\;\;
L(\sigma_j)\equiv 0\;(\forall j)\,;$
\vskip.2em

(Condition VI)\hspace{6ex}  
\lq\lq\,the same as above\,''\;; 

\vskip.2em    
(Condition VII)\qquad 
$
\mathrm{ord}(d)\equiv 0,\;\;
L(\sigma):=\sum_{j\in J}L(\sigma_j)\equiv 0\,.
$
\\[1ex]
Hence \hspace{9ex} 
$
{\cal O}({\rm III})={\cal O}({\rm VII})\cap{\cal O}({\rm IV}),\quad  {\cal O}({\rm VI})={\cal O}({\rm V})\subset{\cal O}({\rm VII}).
$

\vskip1.2em

{\bf 20.2. Spin characters of spin type $\chi^{\rm IV}=(-1,\,1,\,1)$.}
\vskip.5em

We use the notation ${\cal C}({\cal K})=\{ C=\big(({\boldsymbol \gamma}_\zeta)_{\zeta\in{\cal K}};\mu\big)\,;\,{\sum}_{\zeta\in{\cal K}}\,\|\gamma_\zeta\|+\|\mu\|=1\}$ in 
Notation 16.5, (\ref{2016-03-28-42}), for ${\cal K}=\widehat{T}$. 
 For spin type $\chi^{\rm IV}$ of $G'_\infty$, ${\cal O}({\rm IV})=\{g'\in G'_\infty\,;\,\ell_j\equiv 1\;(\forall j)\}$, and ${\cal K}=\widehat{T}$. 
From Theorems 16.4 and 16.6, we obtain the following. 
\vskip1.2em

{\bf Theorem 20.1.}\; 
{\it 
Let\, $G'_\infty=R\big(G(m,1,n)$ for $m$ even. Then  
\begin{gather}
\label{2017-08-28-1}
{\rm Lim}(G'_\infty;\chi^{\rm IV})=E(G'_\infty;\chi^{\rm IV})=F(G'_\infty;\chi^{\rm IV})=\big\{f^{\rm IV}_C\,;\,C\in{\cal C}(\widehat{T})\big\},
\end{gather}
Here, if $g'=z\, \xi'_{q_1}\xi'_{q_2} \cdots \xi'_{q_r}g'_1g'_2 \cdots g'_s\in G'_\infty$ in (\ref{2017-08-27-1}) is outside ${\cal O}({\rm IV})$, then\, $f^{\rm IV}_C(g')=0$, and for $g'\in {\cal O}({\rm IV})$, with $g=\Phi(g')$, 
\begin{gather}
\nonumber
\label{2016-03-28-41-2}
f^{\rm IV}_C(g')=f^{\rm IV}_0(g')\cdot f_{A_C}(g), 
\end{gather}
as in Theorem 16.6 with $A_C$ in (\ref{2017-09-17-21}), under the normalized condition
\begin{gather} 
\label{2017-08-28-2}
K_j:=\mathrm{supp}(\sigma'_j)=[a_j,b_j]\;\,\mbox{\rm an interval},\;\;\sigma'_j=\sigma'_{K_j}:=r_{a_j}r_{a_j+1}\cdots r_{b_j-1}.
\end{gather}

}

\vskip.2em

{\bf 20.3. Spin characters of spin type $\chi^{\rm III}=(-1,\,1,-1)$.}
\vskip.5em

 Recall that ${\cal O}({\rm III})={\cal O}({\rm VII})\cap {\cal O}({\rm IV})$ (cf.\;Table 12.1) and ${\cal O}({\rm VII})=K'_\infty:=\{g'=(d',\sigma')\in G'_\infty\,;\,\mathrm{ord}(d')\equiv 0,\,L(\sigma')\equiv 0\;(\mathrm{mod}\;2)\}$\, and\, $f^{\rm VII}_0=\mathrm{sgn}^{\mathfrak A}_{Z_3}\cdot X_{{\cal O}({\rm VII})}$. 
 An involutive transformation $\kappa=R(\zeta^{(m^0)})$ acts on $f^{\rm IV}_C$ as: for $g'=(d',\sigma')\in G'_\infty$,
\begin{eqnarray}
\label{2017-11-25-1}
&&\;\;
 f^{\rm IV}_{R(\zeta^{(m^0)})C}(g')=(-1)^{\mathrm{ord}(d')}f^{\rm IV}_C(g')
 =\left\{
 \begin{array}{ll}
 \;\;f^{\rm IV}_C(g'),\;\;&\mathrm{ord}(d')\equiv 0\;(\mathrm{mod}\;2),
 \\[1ex]
 -f^{\rm IV}_C(g'),\;&\mathrm{ord}(d')\equiv 1\;(\mathrm{mod}\;2).
 \end{array}
 \right.
\end{eqnarray}

Put, for an equivalence class $[C]:=\{C, \kappa C\}\in {\cal C}(\widehat{T})/\langle \kappa\rangle$, 
\begin{eqnarray}
\label{2016-04-09-11}
f^{\rm IV}_{[C]}:=P_{\langle \kappa\rangle} f^{\rm IV}_C=\tfrac{1}{2}\big(f^{\rm IV}_C+f^{\rm IV}_{\kappa C}\big).
\end{eqnarray}
Then\, 
$f^{\rm IV}_{[C]}|_{K'_\infty}=f^{\rm IV}_C|_{K'_\infty}={\rm Res}^{G'_\infty}_{K'_\infty}f^{\rm IV}_C$, and $f^{\rm IV}_{[C]}=0$  outside $K'_\infty$. Define  
 \begin{eqnarray}
\label{2016-04-09-12}
f^{\rm III}_C:=f^{\rm VII}_0f^{\rm IV}_C,\quad 
F^{\rm III}_{[C]}:=f^{\rm VII}_0f^{\rm IV}_{[C]}=f^{\rm VII}_0P_{\langle \kappa\rangle} f^{\rm IV}_C\,,
\end{eqnarray}
as functions on $G'_\infty$. Then $f^{\rm III}_C=F^{\rm III}_{[C]}$\, and is of spin type $\chi^{\rm III}=\chi^{\rm VII}\chi^{\rm IV}$.  
 
  \vskip1.2em
 
{\bf Theorem 20.2.}\; 
{\it 
{\rm (i)}\; 
Let 
$G'_\infty=R\big(G(m,1,\infty)\big),\; K'_\infty=\Phi^{-1}(K_\infty)$. Then the set of characters $E\big(G'_\infty;\chi^{\rm III}\big)$ of spin type $\chi^{\rm III}$ is obtained from  $E\big(G'_\infty;\chi^{\rm IV}\big)$ of spin type $\chi^{\rm IV}$ through multiplication operators ${\cal M}, {\cal N}$ of $f^{\rm VII}_0$ as shown below: 

\vskip1.2em
\noindent
\!\!$
\begin{array}{ccccc}
 &
E\big(G'_\infty;\chi^{\rm IV}\big) 
&
\quad
&
E\big(G'_\infty;\chi^{\rm III}\big)
\\[1ex]
&
\mbox{\footnotesize ${\rm Res}^{G'_\infty}_{K'_\infty}$}\,\big\downarrow\,\mbox{\rm\footnotesize surj.} 
&
&
\mbox{\rm\footnotesize bij.}\,\big\downarrow\,\mbox{\footnotesize ${\rm Res}^{G'_\infty}_{K'_\infty}$}
&
\\[.3ex]
&
E\big(K'_\infty;\chi^{\rm IV}\big)
&
\!\!\begin{array}{c}
\overset{{\cal M}'}{\longrightarrow}
\\[-2.2ex]
\underset{{\cal N}'}{\longleftarrow}
\end{array}
\!\!
&
E\big(K'_\infty;\chi^{\rm III}\big)
&\;
\end{array} 
\!\!\!\!
\begin{array}{ccccc}
 & 
  f^{\rm IV}_C\;\; 
  &
  & 
  F^{\rm III}_{[C]}=f^{\rm VII}_0f^{\rm IV}_{[C]}\;\;\, 
  & 
  \\
 &
 \downarrow\;\;
 &
 &
 \updownarrow\qquad\;
 &
 \\
 &
f^{\rm IV}_C\big|_{K'_\infty}\!\!=\!f^{\rm IV}_{[C]}\big|_{K'_\infty}
&
\!\!\rightleftarrows\!\! 
&
 F^{\rm III}_{[C]}\big|_{K'_\infty}\!\!=\!f^{\rm VII}_0 f^{\rm IV}_{[C]}\big|_{K'_\infty}
&
\end{array}
$
\vskip1em
\noindent 
Here ${\cal M}'$ and ${\cal N}'$ denote the restrictions of ${\cal M}$ and ${\cal N}$.

\vskip.2em

{\rm (ii)}\; The set of characters of spin type $\chi^{\rm III}$ is given as 
\begin{eqnarray*}
\label{2016-04-09-21} 
\quad
E\big(G'_\infty;\chi^{\rm III}\big)=
\big\{ 
f^{\rm III}_C=F^{\rm III}_{[C]}\;;\;[C]=\{C, \kappa C\}\in{\cal C}(\widehat{T})/\langle \kappa\rangle
\big\}, 
\end{eqnarray*}
and the set $F(G'_\infty;\chi^{\rm III})$ of factorisable characters of spin type $\chi^{\rm III}$ is equal to  
\begin{alignat*}{2}
E(G'_\infty;\chi^{\rm III};{\cal O}({\rm str}))
&:=\big\{f\in E\big(G'_\infty;\chi^{\rm III}\big)\;;\;  \mathrm{supp}(f)\subset {\cal O}({\rm str})\big\}
\\
&\,=\big\{f^{\rm III}_C=F^{\rm III}_{[C]}\;;\; 
C=\kappa C\in{\cal C}(\widehat{T})\big\}.
\end{alignat*}
}

{\it Proof.}\; (ii)\; We prove here the last equality. It is easy to see 
$F(G'_\infty;\chi^{\rm III})\supset \big\{F^{\rm III}_{[C]}\;;\; 
C=\kappa C\in{\cal C}(\widehat{T})\}.$ 
To prove the converse, take $g'=g'_1g'_2,\,g'_j=(d'_j,\sigma'_j)$ with $\ell_j\equiv 1\;(\mathrm{mod}\;2),\,d'_j=\eta^{2k_j+1}_{c_j}\;(\exists c_j\in K_j,\,0\le k_j\le m^0\;{\rm odd})$. Then, with $J=\{1,2\}$, 
\begin{alignat*}{2}
0 
&= F^{\rm III}_{[C]}(g') = 
\Bigg\{\sum_{0\leqslant a <m^0}\Big( \sum_{i \geqslant 1}\;(\gamma_{\zeta^{(a)};i})^{\ell_1}-\sum_{i \geqslant 1}\;(\gamma_{\zeta^{(a+m^0)};i})^{\ell_1}
\Big) \omega^{a(2k_1+1)}
\Bigg\}\times
\\
&\hspace{13.5ex} \times 
\Bigg\{\sum_{0\leqslant b <m^0}\Big( \sum_{i \geqslant 1}\;(\gamma_{\zeta^{(b)};i})^{\ell_2}-\sum_{i \geqslant 1}\;(\gamma_{\zeta^{(b+m^0)};i})^{\ell_2}
\Big) \omega^{b(2k_2+1)}
\Bigg\}. 
\end{alignat*}
Therefore we get, for $0\le a<m^0,\,\ell>1$ odd, 
$$
\qquad\qquad
\sum_{i \geqslant 1}\;(\gamma_{\zeta^{(a)};i})^\ell-\sum_{i \geqslant 1}\;(\gamma_{\zeta^{(a+m^0)};i})^\ell=0.\quad 
{\rm Hence}\;\;\;\gamma_{\zeta^{(a)};i}=\gamma_{\zeta^{(a+m^0)};i}.
\qquad \qquad
\Box
$$

\section{\large Spin characters of spin type $\chi^{\rm V}=(1,-1,-1)$ and $\chi^{\rm VI}=(1,-1,\,1)$}

{\bf 21.1. Spin characters of spin type $\chi^{\rm VI}=(1,-1,\,1)$.}
\vskip.5em

An element $g'\in G'_\infty=R\big(G(m,1,\infty)\big)=\widetilde{D}^\vee_\infty\rtimes\widetilde{{\mathfrak S}}_\infty$ has a standard decomposition $g'=z\,
\xi'_{q_1}\xi'_{q_2}\cdots \xi'_{q_r}g'_1g'_2 \cdots g'_s,\,z\in Z,\,g'_j=(d'_j,\sigma'_j),$  
as in (\ref{2017-08-27-1}), with $J:={\boldsymbol I}_s,\;Q:=\{q_1,\ldots,q_r\}$. As we have seen, ${\cal O}({\rm VI})={\cal O}({\rm V})={\cal O}({\rm str})\subset{\cal O}({\rm VII})=K'_\infty$, 
 where 
\begin{alignat}{2}
\label{2017-08-29-1}
\;\;
{\cal O}({\rm str})
:=\big\{g'\in G'_\infty;\,\mathrm{ord}(\xi_{q_i})\equiv 0\;(\forall i),\, \mathrm{ord}(d_j)\equiv L(\sigma_j)\equiv 0\;(\forall j)\;(\mathrm{mod}\;2)\big\}.
\end{alignat}
Moreover any element of $G'_\infty$ is conjugate modulo $Z$ to an element $g'=(d',\sigma')=\xi'_{q_1}\cdots \xi'_{q_r}\cdot 
\\ 
g'_1\cdots g'_s,$ satisfying the  
normalization (\ref{2016-05-12-11}), or, 
 \begin{eqnarray}
\label{2017-08-29-11}
&&
\left\{
\begin{array}{l}
\xi'_q = \big(t'_q, (q)\big)={\widehat{\eta}_q}^{\;a_q},\;\;
g'_j=(d'_j, \sigma'_j),\;K_j:=\mathrm{supp}(\sigma'_j) \supset \mathrm{supp}(d'_j),
\\[.7ex]
\mbox{\rm $K_j$ is a subinterval $[a_j, b_j]$ of ${\boldsymbol I}_n=[1,n]$},\;\mbox{\rm  mutually disjoint,}
\\[.7ex]  
\sigma'_j=r_{a_j}r_{a_j+1}\cdots r_{b_j-1},\;
d'_j =\widehat{\eta}_{a_j}^{\;\,\mathrm{ord}(d'_j)}\;(j\in J).
\end{array}
\right.
\end{eqnarray}

Irreducible characters of spin type $\chi^{\rm VI}$ of $G'_n=R\big(G(m,1,n)\big),\,4\le n<\infty$, are given in Theorem 10.7 for $n$ even, and in Theorem 10.9 for $n$ odd. Using them, we can calculate pointwise limits of normalized irreducible characters similarly as in \S 16.2. Moreover we can utilize a character relation in Theorem 10.7 for $n$ even such as\;  
$\chi\big(\Pi^{\rm VI}_{M^n}\big)=\chi\big(\Pi^{\rm II0}_n\big)\,\chi\big(\Pi^{\rm IV}_{M^n}\big)$, and similar relations in Theorem 10.9 for $n$ odd. Introduce a parameter space
${\cal C}(\widehat{T}^0)$ 
as in (\ref{2016-03-28-42}). 
For a $C^0\in {\cal C}(\widehat{T}^0)$, we 
define a function $f^{\rm VI}_{C^0}$ on $G'_\infty$ as follows:\, 
$f^{\rm VI}_{C^0}(zg')=\chi^{\rm VI}(z)f(g')\; 
(z\in Z,\, g'\in G'_\infty)$, and 

if $g'\not\in{\cal O}({\rm VI}) $, then $f^{\rm VI}_{C^0}(g')=0$, and 

if $g'\in{\cal O}({\rm VI})$, then for $g'=\xi'_{q_1}\xi'_{q_2} \cdots \xi'_{q_r}g'_1g'_2 \cdots g'_s,\;g'_j=(d'_j,\sigma'_j),\;\ell_j:=\ell(\sigma'_j),$ under the normalization condition (\ref{2017-08-29-11}),  
\begin{gather}
\label{2017-08-29-23} 
f^{\rm VI}_{C^0}(g')\,=\, 
\prod_{q \in Q}
\Big\{ \sum_{\zeta \in \widehat{T}^0}
\big(
 \|\gamma_\zeta\|   
\; + \;\mu_\zeta\big)
\zeta(t'_{q})\Big\} \times
\hspace{20ex}
\\
\nonumber
\hspace*{20ex}
\prod_{j \in J} 
\Big\{ \sum_{\zeta \in \widehat{T}^0}\Big((-1)^{(\ell_j-1)/2}2^{-(\ell_j-1)/2} \sum_{i \geqslant 1}\;(\gamma_{\zeta;i})^{\ell_j}
\Big) 
\zeta(d'_j)
\Big\}. 
\end{gather}

{\bf Theorem 21.1.}\; 
{\it 
There holds \;$f^{\rm VI}_{C^0}=f^{\rm II}_0f^{\rm IV}_{C^0}$\; for $C^0\in {\cal C}(\widehat{T}^0)$, and  
\begin{alignat}{1}
\nonumber
{\rm Lim}(G'_\infty;\chi^{\rm VI})
&
=\big\{f^{\rm VI}_{C^0}=f^{\rm II}_0f^{\rm IV}_{C^0}\,;\,C^0\in {\cal C}(\widehat{T}^0)\big\}
=E(G'_\infty;\chi^{\rm VI})=F(G'_\infty;\chi^{\rm VI}).
\end{alignat}
}

\vskip.2em

{\bf 21.2. Spin characters of spin type $\chi^{\rm V}=(1,-1,-1)$.}
\vskip.5em

Since \,$
{\cal O}({\rm VI})={\cal O}({\rm V})={\cal O}({\rm str})\subset{\cal O}({\rm VII}),
$ the multiplication operators ${\cal M}$ and ${\cal N}$ of $f^{\rm VII}_0=\widetilde{\chi}_{\pi_{2,\zeta^{(0)}}}=\mathrm{sgn}^{\mathfrak A}_{Z_3}\cdot X_{{\cal O}({\rm VII})}$ given as
\begin{eqnarray}
\label{2015-08-23-11-2}
K_1\big(G'_\infty;\chi^{\rm V}\big)\;
\begin{array}{c}
\mbox{\rm \footnotesize ${\cal M}$}
\\[-1.2ex]
\rightleftarrows
\\[-1.2ex]
\mbox{\rm\footnotesize ${\cal N}$}
\end{array}
\;K_1\big(G'_\infty;\chi^{\rm VI}\big),
\end{eqnarray}
are both bijective and ${\cal N}={\cal M}^{-1}$. 
Moreover\; $f^{\rm VII}_0f^{\rm I}_0=f^{\rm II}_0$\; by Lemma 19.2\,(ii), and so \,$f^{\rm I}_0=f^{\rm VII}_0f^{\rm II}_0$\, on the subset ${\cal O}({\rm II})\supset{\cal O}({\rm str})$.\vskip1.2em

{\bf Theorem 21.2.}\; 
{\it 
Let 
$G'_\infty=R\big(G(m,1,\infty)\big)$. For the set $E(G'_\infty;\chi^{\rm V})$ of characters of spin type $\chi^{\rm V}=(1,-1,-1)$, the restrictions ${\cal M}'={\cal M}|_{E(G'_\infty;\chi^{\rm VI})}$ and ${\cal N}'={\cal N}|_{E(G'_\infty;\chi^{\rm V})},$ of ${\cal M}$ and ${\cal N}$ are bijections and mutually inverse of the other: 
\begin{eqnarray*}
\label{2016-04-21-11}
E\big(G'_\infty;\chi^{\rm V}\big)\;
\begin{array}{c}
\mbox{\rm \footnotesize ${\cal M}'$}
\\[-1.2ex]
\rightleftarrows
\\[-1.2ex]
\mbox{\rm\footnotesize ${\cal N}'$}
\end{array}
\;E\big(G'_\infty;\chi^{\rm VI}\big),
\quad 
\begin{array}{ccc}
f^{\rm V}_{C^0}:=f^{\rm I}_0\,f^{\rm IV}_{C^0}
&
\!\rightleftarrows\! 
&
 f^{\rm VI}_{C^0}=f^{\rm II}_0 f^{\rm IV}_{C^0}, 
\end{array}
\end{eqnarray*}
and\; ${\rm Lim}\big(G'_\infty;\chi^{\rm V}\big)=E\big(G'_\infty;\chi^{\rm V}\big)=F\big(G'_\infty;\chi^{\rm V}\big).$ 
The following is a parametrization of characters of spin type $\chi^{\rm V}$\,: 
\begin{eqnarray}
\label{2016-04-21-21}
E\big(G'_\infty;\chi^{\rm V}\big)=\big\{f^{\rm V}_{C^0}=f^{\rm I}_0\,f^{\rm IV}_{C^0}
=\mathrm{sgn}^{\mathfrak A}_{Z_3}\!\cdot\! f^{\rm VI}_{C^0}\;;\;C^0\in{\cal C}(\widehat{T}^0)\big\}, 
\end{eqnarray}
where $\mathrm{sgn}^{\mathfrak A}_{Z_3}$ denotes the natural extension of $\mathrm{sgn}_{Z_3}$ to $G^{\prime\hspace{.15ex}{\mathfrak A}}_\infty=(Z_3\times D^\vee_\infty)\rtimes \widetilde{{\mathfrak A}}_\infty$ (cf. \S 12.4). 
}

\section{\large Parameter spaces and summary for the cases of $n=\infty$}

{\bf 22.1. Parameter spaces for characters.}\; 
The parameter space for characters of the infinite symmetric group ${\mathfrak S}_\infty=G(1,1,\infty)$ is the space of Thoma parameters: 
\begin{eqnarray}
\label{2012-08-21-1}
\nonumber
&&
{\cal A}:=\big\{(\alpha,\beta)\,;\;\alpha=(\alpha_i)_{i\geqslant 1},\;\beta=(\beta_i)_{i\geqslant 1}\;\; \mbox{\rm  satisfy (\ref{2012-08-21-2})}\big\},
\\[1ex]
\label{2012-08-21-2}
&&
\;
\left\{
\begin{array}{l}
\alpha_1\ge \alpha_2\ge \ldots \ge\alpha_i\ge \ldots \ge 0,
\;\;
\beta_1\ge \beta_2\ge\ldots \ge\beta_i\ge \ldots\ge 0,
\\[1.2ex]
\|\alpha\|+\|\beta\|\le 1,\;\;\;\|\alpha\|:={\sum}_{i\geqslant 1}\,\alpha_i,\;\;\;\|\beta\|:={\sum}_{i\geqslant 1}\,\beta_i.
\end{array}
\right.
\end{eqnarray}

The double covering group $\widetilde{{\mathfrak S}}_\infty=R\big(G(1,1,\infty)\big)$ of ${\mathfrak S}_\infty$ has its parameter space of spin characters consisting of Nazarov parameters $\gamma$ as   
\begin{eqnarray}
\label{2012-08-21-11-9}
\nonumber
&&
{\cal C}:=\big\{\gamma\,;\;\gamma=(\gamma_i)_{i\geqslant 1}\;\;\mbox{\rm satisfies (\ref{2012-08-21-12-9})}\big\},
\\[1ex]
\label{2012-08-21-12-9}
&&
\;
\gamma_1\ge \gamma_2\ge \ldots \ge\gamma_i\ge \ldots \ge 0,
\quad
\|\gamma\|\le 1.
\end{eqnarray}

For characters and spin characters of complex reflection groups $G_\infty=G(m,1,\infty)\\ 
={\mathfrak S}_\infty(T),\,T={\boldsymbol Z}_m,\,m>1$ even, we prepared parameter  
\,${\displaystyle A = \big( (\alpha_{\zeta,\varepsilon})_{(\zeta, \varepsilon) \in {\cal K} \times \{0,1\}}\,;\,
\mu \big)}$ and its set ${\cal A}({\cal K})$ in Notation 14.3, and also parameter \,$C=\big((\gamma_\zeta)_{\zeta\in{\cal K}}\,;\mu\big)$\, and its set ${\cal C}({\cal K})$ in Notation 16.5, both for ${\cal K}=\widehat{T}$ or ${\cal K}=\widehat{T}^0$. 
\vskip1em

{\bf 22.2. Summary for ${\mathfrak S}_\infty$ and $\widetilde{{\mathfrak S}}_\infty.$}\; 
We summarize parametrizations of characters and spin characters of ${\mathfrak S}_\infty$ and $\widetilde{{\mathfrak S}}_\infty$ in the following table. Recall ${\cal O}(\widetilde{{\mathfrak S}}_\infty;O\!P_\infty)$ in (\ref{2017-08-31-1}) and $\psi_\Delta$ in (\ref{2017-08-05-1}). 
\vskip1em

\begin{center}
{\small 
$
\begin{array}{|c|c|c|l|c|l|c|}
\cline{4-7}
\multicolumn{3}{l}{}
&
\multicolumn{1}{|c}{
\begin{array}{l}
\mbox{\rm parame-}
\\[-.4ex]
{\rm trization} 
\end{array}
}
&
\multicolumn{1}{|l}{
\begin{array}{l}
{\rm parameter} 
\\[-.4ex]
{\rm space} 
\end{array} 
}
&
\multicolumn{1}{|c|}{
\!\!\begin{array}{l}
\mbox{\rm subset containing} 
\\[-.1ex]
\mathrm{supp}(f) 
\end{array}\!\!\!\!
}
&
{\rm reference}
\\
\hline
\hline
\multicolumn{3}{|l}{
 \!\!\begin{array}{c}
 \\[-1.7ex]
{\mathfrak S}_\infty
\\[-2.2ex]
 ~
 \end{array}\!\!\!\!
 } 
&
\multicolumn{1}{|c}{f_{\alpha,\beta}} 
&
\multicolumn{1}{|c}{(\alpha,\beta)\in{\cal A}}
&
\multicolumn{1}{|c|}{{\mathfrak S}_\infty}
&
{\rm Th.\,13.1}
\\ 
\hline
\multicolumn{3}{|l}{
 \!\!\begin{array}{c}
 \\[-1.7ex]
 \widetilde{\mathfrak S}_\infty
\\[-2.2ex]
 ~
 \end{array}\!\!\!\!
 } 
&
\multicolumn{1}{|c}{\psi_\Delta\cdot f_{\gamma,0}} 
&
\multicolumn{1}{|c}{\gamma\in{\cal C}}
&
\multicolumn{1}{|c|}{
{\cal O}(\widetilde{{\mathfrak S}}_\infty;O\!P_\infty)\subset \widetilde{\mathfrak A}_\infty }
&
{\rm Th.\,13.8}
\\ 
\hline
\end{array}
$
}
\vskip.2em
 {\bf Table 22.1.}\; Parameter spaces of characters for ${\mathfrak S}_\infty$ and $\widetilde{{\mathfrak S}}_\infty$\; 
\end{center}

\vskip.3em

{\bf 22.3. Summary for $G_\infty=G(m,1,\infty)$ and $G'_\infty=R\big(G(m,1,\infty)\big)$.}
 
The parametrization of of characters and spin characters of infinite complex reflection groups $G_\infty=G(m,1,\infty)$, and the parameter spaces depending on spin types are summarized in the following table. Recall \,$f^{\rm II}_0=f^{\rm VII}_0f^{\rm I}_0$\, and $A_C$ in (\ref{2017-09-17-21}).

\vskip1em
\noindent
\begin{center}
{\small 
$
\begin{array}{|c|l|l|l|l|c|c|}
\hline 
\multicolumn{1}{|c}{ 
\begin{array}{c}
{\rm Case}
\\[-.3ex]
{\rm Y}
\end{array}
} 
&
\multicolumn{1}{|c}{
\begin{array}{l}
\mbox{\rm spin}
\\[-.4ex]
\mbox{\rm type} 
\end{array}
\!\!\!\!
}
&
\multicolumn{1}{|c}{
\begin{array}{l}
\mbox{\rm parame-}
\\[-.4ex]
{\rm trization} 
\end{array}
}
&
\multicolumn{1}{|c}{
\begin{array}{l}
\mbox{\rm parameter}
\\[-.4ex]
{\rm space}
\end{array}
} 
&
\multicolumn{1}{|c|}{
\!\!\!
\begin{array}{c}
 {\rm subset}\supset 
\\[-.3ex]
\mathrm{supp}(f) 
\end{array}
\!\!\!\!\!
}
&
\multicolumn{1}{|c|}{
\!\!\!
\begin{array}{c}
 \mbox{\rm criteri-\,} 
\\[-.3ex]
\mbox{\rm on\,(EF)\!} 
\end{array}
\!\!\!\!\!
}
&
\begin{array}{l}
 \mbox{\rm refer-} 
\\[-.3ex]
{\rm ence} 
\end{array}
\\
\hline
\hline 
 \multicolumn{1}{|c}{
  \!\!\!\!
\begin{array}{c}
 \mbox{\rm Y$=$\,odd}
 \\[-.4ex]
 \mbox{\rm ($m$\;odd)}
 \end{array}
\!\!\!\!\!
 }
 &
\multicolumn{1}{|c}{z_1\to -1} 
& 
 \multicolumn{1}{|c}{
  \!f^{\rm \,odd}_C\!=\!f^{\rm \,odd}_0f_{A_C}
  \!\!\!
 } 
&
\multicolumn{1}{|c}{
C\in {\cal C}(\widehat{T})
}
&
 \multicolumn{1}{|c|}{{\cal O}({\rm odd}) 
  }
  &
  {\rm YES}
&
{\rm Th.\,16.6}
\\
\hline 
\hline 
 \multicolumn{1}{|c}{
 \begin{array}{c}
 \\[-1.7ex]
 \mbox{\rm I}
 \\[-2.2ex]
 ~
 \end{array}
 }
 &
\multicolumn{1}{|c}{\!(\!-\!1,\!-\!1,\!-\!1)\!\!} 
& 
 \multicolumn{1}{|c}{
 f^{\rm I}_0\cdot f_{A^0}
 } 
&
\multicolumn{1}{|c}{
A^0\in {\cal A}\big(\widehat{T}^0\big)
}
&
 \multicolumn{1}{|c|}{{\cal O}({\rm I}) 
  }
 &
  {\rm YES}
&
{\rm Th.\,18.3}
\\
\hline
 \multicolumn{1}{|c}{
 \begin{array}{c}
 \\[-1.7ex]
 \mbox{\rm II}
 \\[-2.2ex]
 ~
 \end{array}
 }
 &
\multicolumn{1}{|c}{(\!-\!1,\!-\!1,1)} 
& 
 \multicolumn{1}{|c}{
 f^{\rm II}_0\cdot f_{A^0}
 } 
&
\multicolumn{1}{|c}{
[A^0]\in {\cal A}\big(\widehat{T}^0\big)/\langle\tau\rangle
\!}
&
 \multicolumn{1}{|c|}{
 {\cal O}({\rm II})
  }
 &
  {\rm NO}
&
{\rm Th.\,19.5}
\\
\hline
 \multicolumn{1}{|c}{
\begin{array}{c}
 \\[-1.7ex]
 \mbox{\rm III}
 \\[-2.2ex]
 ~
 \end{array}
 }
 &
\multicolumn{1}{|c}{(\!-\!1,1,\!-\!1)} 
& 
 \multicolumn{1}{|c}{
\!f^{\rm III}_C=f^{\rm VII}_0f^{\rm IV}_C\!\!\! } 
&
\multicolumn{1}{|c}{
[C]\in{\cal C}(\widehat{T})/\langle \kappa\rangle
}
&
\multicolumn{1}{|c|}{{\cal O}({\rm III})}
 &
  {\rm NO}
&
{\rm Th.\,20.2}
\\
\hline
 \multicolumn{1}{|c}{
 \begin{array}{c}
 \\[-1.7ex]
 \mbox{\rm IV}
 \\[-2.2ex]
 ~
 \end{array}
 }
 &
\multicolumn{1}{|c}{(\!-\!1,1,1)} 
& 
 \multicolumn{1}{|c}{
 f^{\rm IV}_C=f^{\rm IV}_0f_{A_C}
 } 
&
\multicolumn{1}{|c}{
C\in{\cal C}(\widehat{T})
}
&
 \multicolumn{1}{|c|}{
 {\cal O}({\rm IV})
  }
 &
  {\rm YES}
&
{\rm Th.\,20.1}
\\
\hline
 \multicolumn{1}{|c}{
\begin{array}{c}
 \\[-1.7ex]
 \mbox{\rm V}
 \\[-2.2ex]
 ~
 \end{array}
 }
 &
\multicolumn{1}{|c}{(1,\!-\!1,\!-\!1)} 
& 
 \multicolumn{1}{|c}{
 f^{\rm V}_{C^0}=f^{\rm I}_0f^{\rm IV}_{C^0}
 } 
&
\multicolumn{1}{|c}{
C^0\in{\cal C}(\widehat{T}^0)
}
&
\multicolumn{1}{|c|}{{\cal O}({\rm str})}
 &
  {\rm YES}
&
{\rm Th.\,21.2}
\\
\hline
 \multicolumn{1}{|c}{
 \begin{array}{c}
 \\[-1.7ex]
 \mbox{\rm VI}
 \\[-2.2ex]
 ~
 \end{array}
 }
 &
\multicolumn{1}{|c}{(1,\!-\!1,1)} 
& 
 \multicolumn{1}{|c}{
 f^{\rm VI}_{C^0}=f^{\rm II}_0f^{\rm IV}_{C^0}
 } 
&
\multicolumn{1}{|c}{
C^0\in{\cal C}(\widehat{T}^0)
}
&
 \multicolumn{1}{|c|}{
 {\cal O}({\rm str})
  }
 &
  {\rm YES}
&
{\rm Th.\,21.1}
\\
\hline  
 \multicolumn{1}{|c}{
\begin{array}{c}
 \\[-1.7ex]
 \mbox{\rm VII}
 \\[-2.2ex]
 ~
 \end{array}
 }
 &
\multicolumn{1}{|c}{(1,1,\!-\!1)} 
& 
 \multicolumn{1}{|c}{
 f^{\rm VII}_0\cdot f_A
 } 
&
\multicolumn{1}{|c}{
[A]\in {\cal A}\big(\widehat{T}\big)/\langle\tau,\kappa\rangle
\!\!}
&
\multicolumn{1}{|c|}{{\cal O}({\rm VII})}
 &
  {\rm NO}
&
{\rm Th.\,17.4}
\\
\hline\hline
 \multicolumn{1}{|c}{
\begin{array}{c}
 \\[-1.7ex]
 \!\!\mbox{\rm VIII}\!\!\!
 \\[-2.2ex]
 ~
 \end{array}
 }
 &
\multicolumn{1}{|c}{
\!\!\!\!
\begin{array}{c}
(1,1,1)
\\[-.2ex]
\mbox{\rm\footnotesize non spin}
\end{array}
\!\!\!\!
} 
& 
 \multicolumn{1}{|c}{
  f_A
 } 
&
\multicolumn{1}{|c}{
A\in {\cal A}\big(\widehat{T}\big)
}
&
\multicolumn{1}{|c|}{
G_\infty
} 
 &
  {\rm YES}
&
{\rm Th.\,14.3}
\\
\hline
\end{array}
$
\vskip.5em
}

\vskip.2em

\quad{\bf Table 22.2.}\; Parameter spaces for the set of spin characters $E(G'_\infty;\chi^{\rm Y})$.

{\small \hspace*{6ex} 
(Two different pairings of Cases: (I, II), (III, IV), (V, VI), (VII, VIII), 

\hspace{12ex} and (I, V), (II, VI), (III, VII), (IV, VIII), are both interesting.) 
}
\end{center}

\vskip.5em 

{\bf Lemma 22.3.} 
{\it There hold the following product formulas: 
\begin{gather}
\nonumber
(f^{\rm I}_0)^2 =f_{A^{(1)}},\;f^{\rm I}_0f^{\rm II}_0=f^{\rm VII}_0f_{A^{(1)}},\;
f^{\rm I}_0f^{\rm VII}_0=f^{\rm II}_0,\;
(f^{\rm IV})^2=f_{A^{(2)}},\;f^{\rm II}_0f^{\rm III}_0=f^{\rm I}_0f^{\rm IV}_0,
\\[.3ex]
\nonumber 
\label{2017-09-15-11}
(f^{\rm II}_0)^2 =X_{K'_\infty}f_{A^{(1)}},\;\;(f^{\rm III}_0)^2 =X_{K'_\infty}f_{A^{(2)}},\;\;(f^{\rm VII}_0)^2 =X_{K'_\infty},
\\[.4ex]
\nonumber 
f^{\rm II}_0f^{\rm VII}_0=X_{K'_\infty}f^{\rm I}_0,\;\;
f^{\rm III}_0f^{\rm VII}_0=X_{K'_\infty}f^{\rm IV}_0,
\end{gather}
where $A^{(k)}=\big((\alpha^{(k)}_{\zeta,\varepsilon})_{(\zeta,\varepsilon)\in\widehat{T}\times\{0,1\}},\mu^{(k)}\big),\,k=1,2,$ with 
\begin{gather*}
\alpha^{(1)}_{\zeta,\varepsilon;1}=1/2\;\;\mbox{\rm for $(\zeta,\varepsilon)=(\zeta^{(0)},0),\,(\zeta^{(m^0)},1),$\,  other $\alpha^{(1)}_{\zeta,\varepsilon;i}=0,\,\mu^{(1)}={\bf 0}=(0,0,\ldots)$,
}
\\
\alpha^{(2)}_{\zeta^{(0)},\varepsilon;1}=1/2\;\;\mbox{\rm for $\varepsilon=0,1,$\,  other $\alpha^{(2)}_{\zeta,\varepsilon;i}=0,\,\mu^{(2)}={\bf 0}$.
}
\end{gather*}
}
\vskip.2em

{\it Proof.} The square $(f^{\rm I}_0)^2$ and $(f^{\rm IV}_0)^2$ have trivial spin type $\chi^{\rm VIII}$ and are factorisable, and so should be written as $f_{A^{(k)}}$ with some $A^{(k)}\in{\cal A}(\widehat{T}),\,k=1,2$ respectively. Take $g'=(d',\sigma')\in{\cal O}({\rm I})$ with $\sigma'\in\widetilde{{\mathfrak S}}_\infty$ a cycle of length $\ell$. 
Then 
$(f^{\rm I}_0)^2(g')=2^{-\ell+1}=2^{-\ell}+2^{-\ell}$\; if $\mathrm{ord}(d')+(\ell-1)\equiv 0\;(\mathrm{mod}\;2)$ and $=0$ otherwise. This is realized by the above choice of $A^{(1)}$. Similarly for $(f^{\rm IV}_0)^2$. 

The rest of the proof is easy. 
\hfill 
$\Box$\quad\;

\vskip1em 

{\bf Theorem 22.4.} 
{\it 
Let $m$ be even and $G'_\infty=R\big(G(m,1,\infty)\big)$. 

{\rm (i)}\; 
For two characters $f_i\in E(G'_\infty;\chi^{{\rm Y}_i})$ for $i=1,2$, the product $f_3=f_1f_2$ belongs to $K_1(G'_\infty;\chi^{{\rm Y}_3}),\,\chi^{{\rm Y}_3}=\chi^{{\rm Y}_1}\chi^{{\rm Y}_2}$, in general. Except cases where ${\rm Y}_1, {\rm Y}_2 \in\big\{{\rm II,III, VII}\big\}$ but $\big\{{\rm Y}_1, {\rm Y}_2\big\}\ne \{{\rm II, III}\}$ (cf. Table 12.7), $f_3=f_1f_2$ is a character of $G'_\infty$, and so 
\\
$E(G'_\infty;\chi^{{\rm Y}_1})E(G'_\infty;\chi^{{\rm Y}_2})\subset E(G'_\infty;\chi^{{\rm Y}_3})$. 

{\rm (ii)}\; 
In the exceptional cases, 
\;
$E(G'_\infty;\chi^{{\rm Y}_1})E(G'_\infty;\chi^{{\rm Y}_2})\not\subset E(G'_\infty;\chi^{{\rm Y}_3})$, \;whereas\\
 $F(G'_\infty;\chi^{{\rm Y}_1})F(G'_\infty;\chi^{{\rm Y}_2})\subset F(G'_\infty;\chi^{{\rm Y}_3})$.

{\rm (iii)}\; 
The total set $E(G'_\infty)$ of characters is generated under product by the sets 
\begin{eqnarray}
\label{2017-09-16-11}
\big\{\,f_A\;\big(A\in {\cal A}(\widehat{T})\big),\;\; 
f^{\rm I}_0, \;\;f^{\rm IV}_0,\;\; f^{\rm VII}_0\,\big\}.
\end{eqnarray}

}

\vskip.2em

{\it Proof.}\; We apply Theorem 12.8, Table 22.2 and Lemma 22.3. Note that $(f^{\rm VII}_0)^2 =X_{K'_\infty}$, and $X_{K'_\infty}\in K_1(G'_\infty;\chi^{\rm VIII})$, the indicator function of the normal subgroup $K'_\infty$, is of trivial spin type and not factorisable and so not a character. Moreover $X_{K'_\infty}f=f$ if $\mathrm{supp}(f)\subset K'_\infty$. Then we can check case by case if necessary and obtain the assertion. 
\hfill
$\Box$ 

\vskip1em

{\bf 22.4. Case of canonical normal subgroups.}\; 
Let $N$ be one of normal subgroups of $G_\infty=G(m,1,\infty)=D_\infty\rtimes{\mathfrak S}_\infty,\,D_\infty=D_\infty(T),\,T={\boldsymbol Z}_m,$ given as 
$$
G(m,p,\infty),\, p|m,\,p>1,\;\; G^{\mathfrak A}_\infty=D_\infty\rtimes {\mathfrak A}_\infty,\;\; G^{\mathfrak A}(m,p,\infty)=G^{\mathfrak A}_\infty\bigcap G(m,p,\infty).
$$

We list up below the results on parametrization of characters of $N$.     
\vskip.5em

\begin{center}
{\small 
$
\begin{array}{|c|c|c|c|} 
\hline
\multicolumn{1}{|c}{
\!\!\!
\begin{array}{c}
\mbox{\rm normal subgroup} 
\\
N 
\end{array}
\!\!\!\!}
&
\multicolumn{1}{|c}{
\begin{array}{c}
{\rm parametrization} 
\end{array}
}
&
\multicolumn{1}{|c}{
\begin{array}{c}
\mbox{\rm parameter space}
\\[-.2ex]
\mbox{\rm for $E(N)$} 
\end{array} 
}
&
\multicolumn{1}{|c|}{
{\rm reference}
}
\\  
\hline 
\hline   
\multicolumn{1}{|c}{
 \begin{array}{c}
G(m,p,\infty) 
 \end{array}
 } 
&
\multicolumn{1}{|c}{
\!\!\!
\begin{array}{c}
\\[-1.7ex]
\;f_A\big|_N=\big(P_{\langle \kappa^{(q)}\rangle} f_A\big)\big|_N
\\
A=\big( (\alpha_{\zeta,\varepsilon})_{(\zeta, \varepsilon) \in \widehat{T} \times \{0,  1\}} 
;\mu \big)
\\[-1.5ex]
 ~
\end{array}
\!\!\!\!
} 
&
\multicolumn{1}{|c}{
\begin{array}{c}
[A]\in{\cal A}(\widehat{T})\big/\big\langle \kappa^{(q)}\big\rangle 
\\
q=m/p
\end{array}
}
&
\multicolumn{1}{|c|}{
\!\!\!\!
\begin{array}{l}
{\rm Th.\,15.18}
\end{array}
\!\!\!\!\!}
\\ 
\hline   
\multicolumn{1}{|c}{
\begin{array}{c}
 \\[-1.2ex] 
G^{\mathfrak A}_\infty
\\[-1.2ex]
 ~
 \end{array}
 } 
&
\multicolumn{1}{|c}{
f_A\big|_N=\big(P_{\langle \tau\rangle} f_A\big)\big|_N

} 
&
\multicolumn{1}{|c}{
[A]\in{\cal A}(\widehat{T})\big/\big\langle \tau\big\rangle
}
&
\multicolumn{1}{|c|}{
\!\!\!\!{\rm Th.\,15.17}\!\!\!\!\!
}
\\ 
\hline             
\multicolumn{1}{|c}{
\begin{array}{c}
 \\[-1.2ex] 
 G^{\mathfrak A}(m,p,\infty)
\\[-1.2ex]
 ~
 \end{array}
 } 
&
\multicolumn{1}{|c}{f_A\big|_N=\big(P_{\langle\tau, \kappa^{(q)}\rangle} f_A\big)\big|_N
} 
&
\multicolumn{1}{|c}{
[A]\in{\cal A}(\widehat{T})\big/\big\langle \tau,\,\kappa^{(q)}\big\rangle
}
&
\multicolumn{1}{|c|}{
\!\!\!\!{\rm Th.\,15.19}\!\!\!\!\!
}
\\ 
\hline
\end{array}
$
}

\vskip.5em

{\normalsize 
 {\bf Table 22.5.} Parameters of characters of normal subgroups of $G(m,1,\infty)$ 
 }
\end{center}

\vskip1em 

{\bf 22.5. Works of Stembridge [36] and of Morris-Jones \cite{[MoJo2003]}.} 

We compare here our results in Part I with the preceding results in J. Stembridge [36] and in A.O.\;Morris and H.I.\;Jones \cite{[MoJo2003]} for Part I. 
Let $m$ be even and $Z=\langle z_1,z_2,z_3\rangle={\mathfrak M}\big(G(m,1,n)\big),\,n\ge 4,$ be as in \S 1.1. We denote characters $\chi_Z\in \widehat{Z}$ as 
\\[.5ex]
\hspace*{10ex}
$\chi^{\rm Y}=(\beta_1,\beta_2,\beta_3),\;\;\beta_i:=\chi^{\rm Y}(z_i)=\pm 1$\; with Y\,=\,I, II, \ldots, VIII, 
\\[.8ex]
as in \S 1.2 (and in Table 22.2). 
On the other hand, it is denoted as $[\gamma,\lambda,\mu]$ in \cite{[MoJo2003]} and $[t_\vartheta,t_\mu,t_\nu]$ in Dudko\,-Nessonov \cite{[DuNe2008]}, with $\gamma=t_\vartheta=\beta_1,\,\lambda=t_\mu=\beta_3,\,\mu=t_\nu=\beta_2$, and 1-1 correspondence of notations in (14) of \cite{[DuNe2008]} with ours is 
{\small 
$$
\begin{array}{|c|c|c|c|c|c|c|c|}
\hline
~& & & & & & &~ 
\\[-2.2ex]
\chi^{\rm I}&\chi^{\rm II}&\chi^{\rm III}&\chi^{\rm IV}&\chi^{\rm V}&\chi^{\rm VI}&\chi^{\rm VII}&\chi^{\rm VIII}
\\
\hline 
 \theta_7 &\theta_5 &\theta_3 &\theta_1 &\theta_6 &\theta_4 &\theta_2 &\theta_0 
 \\
 \hline 
 \end{array}
$$
}

{\bf Work of J.\;Stembridge [36].} 

Let $W_n=G(2,1,n)$ be the Weyl group of type $B_n/C_n$. 
In its Introduction,  
\begin{quotation} 
\noindent
the main result of this paper is the fact that the projective representation of $W_n$ can be explicitly constructed from the linear and projective representations of symmetric groups, 
\end{quotation}
by adding two qualifications: (1) to add to this list the representation obtained by restricting the basic spin representation of $O_n$ to $W_n$, and 
(2) the irreducible decomposition of each symmetric group module to the alternating group.

In the part of Section 3 in Introduction, his approach is explained as follows. 

$\bullet$ A group $G$ is said to {\it have $K$-quotient}\, if it has a normal subgroup $N$ such that $G/N\cong K$. 
 Let $G$ be a group with ${\boldsymbol Z}_2^{\;2}$-quotient, then it has a two-dimensional projective IR $\rho$ given by anti-commuting involutions such as $\big(\begin{smallmatrix}1&0\\ 0&-1\end{smallmatrix}\big),$ $\big(\begin{smallmatrix}0&1\\ 1&0\end{smallmatrix}\big)$. Let $G_1$ and $G_2$ be two groups with ${\boldsymbol Z}_2$-quotient, then $G_1\times G_2$ has a natural  ${\boldsymbol Z}_2^{\;2}$-quotient and accordingly has a 2-dimensional projective IRs $\rho$. Start with irreducible modules $V_i$ of $G_i$, and consider $\rho\otimes(V_1\otimes V_2)$. Define the twisted product $V_1\widehat{\otimes}V_2$ to be one or more of the irreducible submodules of $\rho\otimes(V_1\otimes V_2)$. 

$\bullet$ $W_n$ has a natural ${\boldsymbol Z}^{\;2}_2$-quotient and for a projective IR $\pi$ of $W_n$, the operation $\pi\mapsto \rho\otimes \pi$\, and $\chi_\pi\mapsto \chi_\rho\cdot\chi_\pi$\, give 
transition of spin types as $\chi^\pi\mapsto \chi^\rho\cdot \chi^\pi$ (in our presentation, cf. {\bf 1.2} of this paper), where $\pi(z)=\chi^\pi(z)I\;(z\in Z)$ with $\chi^\pi\in\widehat{Z}$. 

(Note that $\rho$ is something like $\pi_{2,\zeta^{(0)}}$ in {\bf 12.4} of this paper, and $\chi^\rho$ is something like $\chi^{\rm VII}$ and the above operation works much simpler in the case of $n=\infty$, see e.g. Part II, \S 17 and \S\S 19--21 and Table 22.2 in this paper.)

$\bullet$ Let $\pi$ be a projective IR of $W_n$. 
Apply the theory of A.H.\;Clifford\;[3] to the irreducible submodules of $\rho\otimes \pi$, then 
irreducible characters of these submodules can be expressed in terms of 
difference characters of $W_n$ or one of its double covering.

\vskip1.2em

{\bf Work of A.O.\;Morris and H.I.\;Jones \cite{[MoJo2003]}.} 
 From Introduction, 
\begin{quotation}
\noindent 
$\cdots\cdots$\, all of the non-equivalent irreducible projective representations can be expressed in terms of certain \lq building blocks'. 
These are the ordinary and spin representations of the symmetric group $G(1,1,n)$, 
$\cdots\cdots\cdots.$ 
Also, required are basic spin representations $P, Q$ and $R$ of $B^m_n$ for certain 2-cocycles. 
\end{quotation}

In \S 2, the methods of the paper are prepared. In \S 2.3, there explained the case of $G$ with ${\boldsymbol Z}_2$-quotient and its difference character, and also the case of $G$ with ${\boldsymbol Z}_2^{\;2}$-quotient and its difference character. 
In \S 2.4, the Clifford algebra $C(n)$ generated by $\{1, e_j\;(j\in{\boldsymbol I}_n)\}$, subjected\; $e_j^{\;2}=1,\;e_je_k=-e_ke_j\;(j\ne k)$, is introduced, and the double coverings $Pin(n)$ and $Spin(n)$ respectively of $O(n)$ and $SO(n)$ are defined. Then basic spin representations of Clifford algebra $C(n)$ are given.

\vskip1.2em 

{\bf 22.6. \,Work of A.V.\;Dudko and N.I.\;Nessonov \cite{[DuNe2008]}.}

Concerning our Part II, we refer the preceding results of Dudko\,-\,Nessonov in \cite{[DuNe2008]} (see also [II, \S 25]). Put $B_m:={\mathfrak S}_\infty\ltimes {\boldsymbol Z}_m^{\,\infty}$. 
From Introduction, 
\begin{quotation}
\noindent
$\cdots\cdots.$\, In the present paper, a complete description of the projective factor representations of finite type is given for the group ${\mathfrak S}_\infty\ltimes {\boldsymbol Z}_m^{\,\infty}$. 
\vskip1em 
\noindent 
$\cdots$ in the proof of the classification theorem  for the projective characters on $B_m$ we show that every indecomposable characters  is explicitly determined by the character of some projective basis factor representations and an ordinary character of the group $B_m$ (see Theorem 10).  
\end{quotation}

Our method of getting all spin characters of $G(m,1,\infty)\;(\cong B_m)$ is based on the approximation theorem (Theorem 12.5) that asserts every character of $R\big(G(m,1,\infty)\big)=\lim_{n\to \infty}R\big( 
G(m,1,n)\big)$ is a pointwise limits of normalized irreducible characters $\widetilde{\chi}_{\pi_n}$ as $n\to\infty$. 

Here we make only a short remark that Lemma 22.3 and Theorem 22.4 in this paper correspond to Theorem 23 and Remark 24 in \cite{[DuNe2008]}, and give a little more detailed information, e.g., for the case $\{\chi^{\rm Y_1},\chi^{\rm Y_2}\}=\{\chi^{\rm II},\chi^{\rm III}\}$ which corresponds to $\{\theta_3,\theta_5\}$.


 \vskip2em 

{\small  
T. Hirai, 22-8 Nakazaichi-Cho, Iwakura, Sakyo-Ku, Kyoto 606-0027, Japan; 

hirai.takeshi.24e@st.kyoto-u.ac.jp
 \vskip1em
 
A. Hora, Department of Mathematics, Hokkaido University, Sapporo 060-0810, Japan; 

hora@math.sci.hokudai.ac.jp


\vskip2em 

{\small

\begin{thebibliography}{HH} 


\bibitem{[Bia1996]}
 P.\,Biane, Minimal factorization of a cycle and central 
multiplicative functions on the infinite 
symmetric groups,
J. Combin. Theory, Ser. A, {\bf 76}(1996), 197-212.

\bibitem{[Bi-Karel-de-Le]}
E.\,Bishop and K.\,de~Leeuw, 
The representation of linear functionals by 
measures on sets of extreme points, Ann. Inst. Fourier, 
{\bf 9}(1959), 305--331. 



\bibitem{[Clif1937]} 
A.H.\;Clifford, 
Representations induced in an invariant subgroup, Ann. Math., {\bf 38}(1937), 533-550. 




\bibitem{[DaMo1974]}
J.W.\;Davies and A.O.\;Morris, The Schur multiplier of the generalized symmetric group, 
 J. London Math. Soc., {\bf (2) 8}(1974), 615-620. 




\bibitem{[DuNe2008]}
A.\;Dudko and N.\;Nessonov, 
Characters of projective representations of the infinite generalized symmetric group, Sbornik: Mathematics {\bf 199}(2008), 1421-1450.



 
\bibitem{[Hir2004]}
 T.\,Hirai, Centralization of positive definite functions, weak containment of representations and
 Thoma characters for the infinite symmetric group,   
J. Math. Kyoto Univ., {\bf 44}(2004), 685-713. 


\bibitem{[Hir2013]}
T.\;Hirai, Classical method of constructing all irreducible representations of semidirect product of a compact group with a finite group, Probability and Math. Statistics, {\bf 33}(2013), 353-362.

\bibitem{[Hir2018]}
T.\;Hirai, {\it Introduction to the theory of spin representations of groups, --- from the beginning beyond the spin representations of symmetric groups  ---} (in Japanese), in S\^ugaku no Mori series, S\^ugaku Shob\^o, Tokyo, May 2018. 


\bibitem{[HH2002]}
T.\,Hirai and E.\,Hirai, 
Characters for the infinite Weyl 
groups of type B$_\infty$/C$_\infty$ and D$_\infty$, and for analogous groups, in \lq\,{\it 
 Non-Commutativity, Infinite-Dimensionality 
 and Probability at the Crossroad}\,', pp.296-317, World Scientific, 2002.

 
\bibitem{[HH2005a]}
T.\,Hirai and E.\,Hirai, 
 Character formula for wreath products of finite groups 
with the infinite symmetric group, in \lq\,{\it the Proceedings of Japanese-German Seminar on Infinite-Dimensional Harmonic Analysis III}, pp.119-139, World Scientific, 2005.  


\bibitem{[HH2005b]}
T.\,Hirai and E.\,Hirai, 
Positive definite class functions on a topological group 
and characters of factor representations, J. Math. Kyoto Univ., {\bf 45}(2005), 355-379. 




\bibitem{[HH2005c]}
T.\,Hirai and E.\,Hirai, 
Characters of wreath products of finite groups with 
the infinite symmetric group, J. Math. Kyoto Univ.,  
{\bf 45}(2005), 547-597.  
 

\bibitem{[HH2006]}
T.\,Hirai and E.\,Hirai, Character formula for wreath products of compact groups with the infinite symmetric group, in {\it Quantum Probability: the Proceedings of 25th QP Conference Quantum Probability and Related Topics  2004 in B\c{e}dlewo}, Banach Center Publications, Institute of Mathematics, Polish Academy of Sciences, {\bf 73}(2006), 207-221.  
 


\bibitem{[HH2007]}
T.\,Hirai and E.\,Hirai, Characters of wreath products of compact groups with the infinite symmetric group and characters of their canonical subgroups, J. Math. Kyoto Univ., {\bf 47}(2007), 269-320.





\bibitem{[HHH2006]}
T.\,Hirai, E.\,Hirai and A.\,Hora,, 
Realization of factor representations of finite type with emphasis on their characters for wreath products of compact groups with the infinite symmetric group, J. Math. Kyoto Univ., {\bf 46}(2006), 75-106.

 
 
 
\bibitem{[HHH2009]}
T.\;Hirai, E.\;Hirai and A.\;Hora, 
Towards projective representations and spin characters of finite and infinite complex reflection groups, {\it Proceedings of the fourth German-Japanese Symposium, Infinite Dimensional Harmonic Analysis, IV}, 
World Scientific, 2009, pp.112-128.


\bibitem{[HHH2009b]}
T.\;Hirai, E.\;Hirai and A.\;Hora, 
Limits of characters of wreath products ${\mathfrak S}_n(T)$ of a compact group $T$ with the symmetric groups and characters of ${\mathfrak S}_\infty(T)$, I, Nagoya Math. J., {\bf 193}(2009), 1-93. 


 

\bibitem{[HHH2013]} 
(quoted as [I])
T.\;Hirai, E.\;Hirai and A.\;Hora, 
Projective representations and spin characters of 
complex reflection groups $G(m,p,n)$ and $G(m,p,\infty)$, \;I, in MSJ Memoirs, Vol.\,29, Math. Soc. Japan, 2013, pp.49-122.   



\bibitem{[HHo2014]}
T.\;Hirai and A.\;Hora, 
Spin representations of twisted central products of double covering finite groups and the case of permutation groups, J. Math. Soc. Japan, {\bf 66}(2014), 1191-1226.





\bibitem{[HHoH2013a]}
T.\;Hirai, A.\;Hora and E.\;Hirai, 
Introductory expositions on projective representations of 
groups, in MSJ Memoirs, Vol.\,29, Math. Soc. Japan, 2013, pp.1-47.   


\bibitem{[HHoH2013b]}
(quoted as [II])\;  
T.\;Hirai, A.\;Hora and E.\;Hirai, 
Projective representations and spin characters of complex reflection groups $G(m,p,n)$ and $G(m,p,\infty)$, II, Case of generalized symmetric groups, ibid., pp.123-272. 


 \bibitem{[HoHu1992]}
P. Hoffman and J. Humphreys, 
{\it  Projective representations of the symmetric group, {\small $Q$-functions and shifted tableaux}}, 
Oxford Mathematical Monographs, 
 Clarendon Press, Oxford, 1992. 


\bibitem{[HoHH2008]}
A.\;Hora, T.\;Hirai and E.\;Hirai,  
Limits of characters of wreath products ${\mathfrak S}_n(T)$ of a compact group $T$ with the symmetric groups and characters of ${\mathfrak S}_\infty(T)$, II, 
From a view point of probability theory, J. Math. Soc. Japan, {\bf 60}(2008), 1187-1217. 


\bibitem{[HoH2014]}
A.\,Hora and T.\,Hirai,  
Harmonic functions on branching graph associated with the infinite wreath product of a compact group, Kyoto J. Math., {\bf 54}(2014), 775-817. 






 
 \bibitem{[Ker2003]}
 S.V.\,Kerov, 
{\it 
Asymptotic Representation Theory of the Symmetric Group and its
Applications in Analysis,} Translations of Mathematical Monographs, Vol.
219, 2003.



\bibitem{[Kle2005]}
A. Kleshchev, {\it Linear and projective representations of 
symmetric groups}, Cambridge Univ. Press, Cambridge, 2005.


\bibitem{[Mor1976]}
A.O.\;Morris, 
A survey on Hall-Littlewood functions and their applications, in 
{\it Combinatoire et Repr\'esentaion du Groupe Sym\'etrique}, 
Springer LN in Math., {\bf 579}(1976), 136-154. 


\bibitem{[Mor2004]}  
A.O.\;Morris, Projective representations of reflection groups III, Communications in Algebra, {\bf 32-7}(2004), 2679-2694.


\bibitem{[MoJo2003]}
A.\;Morris and H.\;Jones, 
Projective representations of generalized symmetric groups, 
S\'eminaire Lotharingien de Combinatoire, {\bf 50}(2003), Article B50b, 1-27. 




\bibitem{[Mur1963]}
F.D.\,Murnaghan, {\it The theory of group representations}, 
Dover Publications, Mineola, N.Y., 1963. 



\bibitem{[Naz1992]}
M.\,Nazarov, Projective representations of the infinite symmetric group, 
Advances Soviet Math., {\bf 9}(1992), 115-130. 


\bibitem{[Rea1977]}
E.W. Read, 
The projective representations of the generalized symmetric groups, 
 J. Algebras, {\bf 46}(1977), 102-133.




\bibitem{[Sch1904]}
J.\;Schur (= I.\;Schur), \"Uber die Darstellung der endlichen Gruppen durch gebrochene lineare Substitutionen, J. f\"ur die reine und angewante Mathematik, {\bf 127}(1904), 20-50.



 
\bibitem{[Sch1907]}
 J.\;Schur, Untersuchungen \"uber die Darstellung der endlichen Gruppen durch gebrochene lineare Substitutionen, 
 ibid., {\bf 132}(1907), 85-137.


\bibitem{[Sch1911]}
 J.\;Schur,  
\"Uber die Darstellung der 
symmetrischen und der alternierenden Gruppen durch 
gebrochene lineare Substitutionen, ibid., {\bf 139}(1911), 155-255. 




\bibitem{[Stem1992]}
J.\;Stembridge, Projective representations of hyperoctahedral groups, J. Algebra, {\bf 145}(1992), 396-453. 





\bibitem{[Tho1964]}
 E.\,Thoma, Die unzerlegbaren positiv-definiten 
Klassenfunktionen der abz\"{a}hlbar unendlichen,
 symmetrischen Gruppe, Math. Z., {\bf 85}(1964), 
40-61.



  

\bibitem{[VK1981]}
 A.\;Vershik and S.\;Kerov, 
 Asymptotic theory of characters of the symmetric group, Funkts. Anal. i Priloemen., {\bf 15}(1981), 15-27; English transl., Funct. Anal. Appl., {\bf 15}(1982), 246-255.


\bibitem{[VK1988]}
 A.\,Vershik and S.\,Kerov, Characters and 
 realizations of representations 
 of an infinite-dimensional Hecke algebra, 
 and knot invariants, Dokl. Acad. Nauk SSSR, {\bf 301}(1988), 777-780; 
 English transl., Soviet Math. Dokl., {\bf 15}(1989), 134-137.

 \end{thebibliography}
 
 
 }
\end{document}